\documentclass [a4paper,10.5pt,oneside]{book}

\usepackage{longtable}
\usepackage[T1]{fontenc}
\usepackage[utf8]{inputenc}
\usepackage[english]{babel}
\usepackage{graphicx}
\usepackage{mathtools}
\usepackage{indentfirst}
\usepackage{braket}
\usepackage{mathrsfs}
\usepackage{amssymb}
\usepackage{textcomp}
\usepackage{dsfont}
\usepackage{latexsym}
\usepackage{amsthm}
\usepackage{amsmath}
\DeclareMathAlphabet{\mathpzc}{OT1}{pzc}{m}{en}
\usepackage{yfonts}
\usepackage{xfrac}
\usepackage{newlfont}
\usepackage{etoolbox}
\usepackage{scalerel}
\usepackage[usestackEOL]{stackengine}
\usepackage{color}
\usepackage{fixmath}
\usepackage[margin=1in]{geometry}

\usepackage[original]{imakeidx}

\newcommand{\meg}{\leqslant}
\newcommand{\Meg}{\geqslant}
\newcommand{\eps}{\varepsilon}
\newcommand{\mi}{\mu}

\newcommand{\Lin}{\mathscr{L}}

\newcommand{\N}{\mathds{N}}
\newcommand{\Z}{\mathds{Z}}
\newcommand{\C}{\mathds{C}}
\newcommand{\R}{\mathds{R}}
\newcommand{\Hd}{\mathds{H}}
\newcommand{\Ac}{\mathcal{A}}
\newcommand{\Bc}{\mathcal{B}}
\newcommand{\Dc}{\mathcal{D}}
\newcommand{\Ec}{\mathcal{E}}
\newcommand{\Fc}{\mathcal{F}}
\newcommand{\Gc}{\mathcal{G}}
\newcommand{\Hc}{\mathcal{H}}
\newcommand{\Ic}{\mathcal{I}}
\newcommand{\Kc}{\mathcal{K}}
\newcommand{\Lc}{\mathcal{L}}
\newcommand{\Mcal}{\mathcal{M}}
\newcommand{\Nc}{\mathcal{N}}
\newcommand{\Oc}{\mathcal{O}}
\newcommand{\Pc}{\mathcal{P}}
\newcommand{\Sc}{\mathcal{S}}
\newcommand{\Ff}{\mathfrak{F}}
\newcommand{\Sf}{\mathfrak{S}}
\newcommand{\Uf}{\mathfrak{U}}

\newcommand{\Ds}{\mathscr{D}}
\newcommand{\Ms}{\mathscr{M}}
\newcommand{\Sb}{\mathbf{S}}
\newcommand{\mbeps}{\mathbold{\varepsilon}}

\DeclarePairedDelimiter{\abs}{\lvert}{\rvert}
\DeclarePairedDelimiter{\norm}{\lVert}{\rVert}
\newcommand{\Supp}[1]{\supp\left( #1\right) }
\newcommand{\Pfaff}{\mathrm{Pf}}
\newcommand{\ee}{\mathrm{e}}
\newcommand{\loc}{\mathrm{loc}}
\newcommand{\vect}[1]{\mathbf{{#1}}}
\newcommand{\dd}{\mathrm{d}}

\newcommand{\atomic}{(L)}
\newcommand{\atomics}{(L')}

\DeclareMathOperator{\card}{Card}

\DeclareMathOperator{\supp}{Supp}
\DeclareMathOperator{\tr}{Tr}
\DeclareMathOperator{\Hol}{Hol}

\newcommand{\half}{\textstyle{\frac12}}
\newcommand{\Szego}{S}
\newcommand{\CB}{}%{\color{black} }
\newcommand{\CR}{}%{\color{red} } 

\newcommand{\dashint}{\,\ThisStyle{\ensurestackMath{%
	\stackinset{c}{.2\LMpt}{c}{.5\LMpt}{\SavedStyle-}{\SavedStyle\phantom{\int}}}%
\setbox0=\hbox{$\SavedStyle\int\,$}\kern-\wd0}\int}

\newcommand{\leftexp}[2]{{\vphantom{#2}}^{#1}{#2}} 
\newcommand{\trasp}{\leftexp{t}}

\newenvironment{abstract}{%
	\clearpage\thispagestyle{empty}\vspace*{3\baselineskip}
	\begin{center}\Large\abstractname\end{center}
	\begin{quotation}
	}{\end{quotation}\clearpage}

\newtheorem{teo}{Theorem}[chapter]
\newtheorem{lem}[teo]{Lemma}
\newtheorem{prop}[teo]{Proposition}
\newtheorem{cor}[teo]{Corollary}

\theoremstyle{definition}
\newtheorem{deff}[teo]{Definition}

\theoremstyle{remark}
\newtheorem{oss}[teo]{Remark}
\newtheorem{ex}[teo]{Example}

\numberwithin{section}{chapter}
\numberwithin{equation}{chapter}

\makeindex

\date{}

\begin{document}

\frontmatter

\title{Holomorphic Function Spaces  on Homogeneous Siegel Domains}

\author{Mattia Calzi \\ Dipartimento di Matematica,\\ Universit\`a degli Studi di  Milano,\\ Via C. Saldini 50, 20133 Milano, Italy\\ {\tt mattia.calzi@unimi.it} \and Marco M.\ Peloso\footnote{Both authors were partially supported by the grant PRIN 2015 \emph{Real and Complex Manifolds: Geometry, Topology and Harmonic Analysis}, and by the grant \emph{Fractional Laplacians and subLaplacians on Lie groups and trees} of the Gruppo Nazionale per l'Analisi Matematica, la Probabilita e le loro Applicazioni (GNAMPA) of the Istituto Nazionale di Alta Matematica (INdAM).} \\ Dipartimento di Matematica,\\ Universit\`a degli Studi di  Milano,\\ Via C. Saldini 50, 20133 Milano, Italy\\ {\tt marco.peloso@unimi.it}}

\maketitle

\tableofcontents

\begin{abstract}
	We study several connected problems  of holomorphic function
	spaces  on homogeneous Siegel domains.  The main object of our study
	concerns 
	weighted mixed norm Bergman spaces on homogeneous Siegel domains of
	type II. These problems include: sampling, atomic decomposition, duality, 
	boundary values,   boundedness of the Bergman projectors.   Our
	analysis include the Hardy spaces, and suitable
	generalizations of the  classical Bloch and Dirichlet spaces.  One
	of the
	main novelties in this work is the generality of the domains under
	consideration, that is, homogeneous Siegel domains, extending many
	results from the more particular cases of the upper half-plane,
	Siegel domains of tube type over irreducible cones, or symmetric,
	irreducible Siegel domains of type II. 
	
	2010 Mathematical Subject Classification: 42B35, 32M15.
	
	Keywords: Bergman spaces, homogeneous Siegel domains, atomic decomposition, boundary values, Bergman projectors, Hardy spaces, Dirichlet space, Bloch space.
\end{abstract}

\chapter*{Introduction}

This research monograph is devoted to the study of spaces of holomorphic functions on a class of unbounded domains in several variables.  The domains we are concerned with are the homogeneous Siegel domains of type II, and the spaces we shall consider are Hardy and Bergman spaces.
On the one hand, we wish to give a uniform and systematic presentation of the main results developed in the literature on several related problems on weighted Bergman spaces on homogeneous Siegel domains of type II. 
On the other hand, we extend the aforementioned results from a variety of particular settings to the general case.

Recall that a connected open subset $D$ of $\C^n$ is said to be a homogeneous domain if the group of holomorphic automorphisms of $D$ acts transitively on $D$.\index{Homogeneous domain}
Homogeneous Siegel domains of type II have been introduced in~\cite{PiatetskiiShapiro}  as unbounded models of the bounded homogeneous domains on finite-dimensional complex vector spaces (cf.~\cite{VinbergGindikinPiatetskiiShapiro}).  

The simplest example of a Siegel domain is the upper half-plane $\C_+\coloneqq \R+i \R_+^*$ in $\C$\label{0}, which is an unbounded domain biholomorphically equivalent to the unit disc in $\C$ by means of the Cayley transform. Various spaces of holomorphic functions on $\C_+$ can be considered, such as the Bergman space ($p\in ]0,\infty]$)
\[
A^p(\C_+)=\Hol(\C_+)\cap L^p(\C_+),
\]
where $\Hol(\C_+)$ denotes the space of holomorphic functions on $\C_+$, the Hardy space 
\[
H^p(\C_+)=\Set{ f\in \Hol(\C_+)\colon\sup\limits_{y>0} \norm{f(\,\cdot\,+i y)}_{L^p(\R)}<\infty},
\]
the Dirichlet space
\[
\Dc(\C_+)=\Set{f\in \Hol(\C_+)\colon f'\in L^2(\C_+)}
\]
(modulo constant functions),
and the Bloch space
\[
\Bc(\C_+)=\Set{f\in \Hol(\C_+)\colon \sup_{z\in \C_+} \Im z \abs{f'(z)}<\infty}
\]
(modulo constant functions).

Following~\cite{RicciTaibleson}, for $p,q\in ]0,\infty]$ and $s>0$ one may also define the weighted Bergman space
\[
A^{p,q}_s(\C_+)=\Set{f\in \Hol(\C_+)\colon \int_0^\infty \left( \norm{f(\,\cdot\,+i y)}_{L^p(\R)} y^s\right) ^q \,\frac{\dd y}{y} <\infty}
\]
(resp.
\[
A^{p,q}_s(\C_+)=\Set{f\in \Hol(\C_+)\colon \sup\limits_{y>0}\left(  \norm{f(\,\cdot\,+i y)}_{L^p(\R)} y^s\right) <\infty}
\]
when $q=\infty$). 
With this definition, the mapping $f\mapsto f'$ induces an isomorphism of $A^{p,q}_s(\C_+) $ onto $A^{p,q}_{s+1}(\C_+)$,\footnote{At least when $p=q\in [1,\infty[$, this is a consequence of Cauchy estimates and a suitable version of Hardy's inequality. In the general case, this follows from the atomic decomposition described below.} so that one may define the space $ A^{p,q}_{s}(\C_+)$ for $p,q\in ]0,\infty]$ and $s\meg 0$, as the space of holomorphic functions on $\C_+$ (modulo polynomials of degree at most $[-s]$) whose derivative of order $[-s]+1$ belongs to $A^{p,q}_{s+[-s]+1}(\C_+)$. 
In this way, the space $A^{2,2}_0(\C_+)$ is canonically isomorphic to the Hardy space $H^2(\C_+)$, while $A^{2,2}_{-1/2}(\C_+)$ is  the Dirichlet space $\Dc(\C_+)$ and $A^{\infty,\infty}_0(\C_+)$ is the Bloch space.\footnote{Notice that, in Chapter~\ref{sec:7} below, the Bloch space will be denoted by $\widehat A^{\infty,\infty}_0(\C_+)$, while $A^{\infty,\infty}_0(\C_+)$ will denote the Hardy space $H^\infty(\C_+)$.}

All the above spaces have been widely studied in various directions. For example, since the space $A^{2,2}_s(\C_+)$, $s>0$, is a hilbertian space which embeds continuously into the space of holomorphic functions, it has a reproducing kernel, namely
\[
K_s(z,w)= \frac{s ( 2 i)^{2 s+1}}{2 \pi } \frac{1}{ (z-\overline w)^{2 s+1}}.
\]
One may then consider the integral operator ($s'>0$)
\[
(P_{s'} f)(z) = \int_{\C_+} K_{s'}(z,w) f(w) (\Im w)^{s'-1}\,\dd w,
\]
called the Bergman projector, and prove that it induces an endomorphism of the mixed norm space 
\[
L^{p,q}_{s}(\C_+)=\Set{f\colon \C_+\to \C  \colon f \text{ measurable,}   \int_0^\infty \left( \norm{f(\,\cdot\,+i y)}_{L^p(\R)} y^s\right) ^q \,\frac{\dd y}{y} <\infty}
\]
(modification if $q=\infty$)  if and only if $p,q\Meg 1$ and $2s'> s+2$. In addition, $P_{s'}$ reproduces the elements of $A^{p,q}_s(\C_+)$, that is
\[
P_{s'} f= f \qquad \forall f \in A^{p,q}_s(\C_+)
\] 
if and only if $2 s'> s+1+1/\min(1,p)$ ($p,q\in ]0,\infty]$). In particular, if $p,q\Meg 1$ and $2s'> s+2$, then $P_{s'}$ induces a continuous projector of $L^{p,q}_s(\C_+)$ onto $A^{p,q}_s(\C_+)$.

Furthermore, holomorphy implies that the elements of $A^{p,q}_s(\C_+)$ are `quasi-constant' (in a suitable weak sense) on the sets
\[
Q_{R,j,k}\coloneqq [2^{k R}R j,2^{k R}R(j+1)]\times [2^{k R},2^{(k+1)R}],
\]
for $j,k\in \Z$ and $R>0$, so that the sampling map
\[
\Sc\colon f \mapsto (\Im z_{R,j,k})^{s+1/p} f(z_{R,j,k})
\]
(where the $z_{R,j,k}\in Q_{R,j,k}$ are arbitrary) induces an isomorphism of $A^{p,q}_s(\C_+)$ onto a closed subspace of
\[
\ell^{p,q} =\Set{ \lambda\in \C^{\Z\times \Z}\colon \sum_{k}  \norm{(\lambda_{j,k})_j}_{\ell^p}^q<\infty }
\]
(modification if $q=\infty$). 

Observe that the only reason why the Bergman projector $P_{s'}$ fails to induce an endomorphism of $L^{p,q}_s(\C_+)$ when $p<1$ or $q<1$ is that the elements of $L^{p,q}_s(\C_+)$ are not locally integrable in general in those cases. 
In view of the preceding sampling result, this problem could be overcome `discretizing' the domain $L^{p,q}_s(\C_+)$ of $P_{s'}$, thus findining an `atomic decomposition' of $A^{p,q}_s(\C_+)$. One may then prove that the mapping
\[
\Ac\colon\ell^{p,q}\ni \lambda \mapsto \sum_{j,k} \lambda_{j,k} K_{s'}(\,\cdot\,, z_{R,j,k}) (\Im z_{R,j,k})^{2 s'-s-1-1/p}\in A^{p,q}_s(\C_+)
\]
is continuous (and has a continuous linear section if $R$ is sufficiently small) if and only if $ 2 s'>s+1+1/\min(1,p)$. 

Another interesting problem related to these spaces is that of determining the boundary values of the elements of $A^{p,q}_s(\C_+)$, that is, the limits of the functions $f_y\colon \R \ni x \mapsto f(x+iy)\in \C$ for $y\to 0^+$, for every $f\in A^{p,q}_s(\C_+)$.
It turns out the the boundary values of the elements of $A^{p,q}_s(\C_+)$ span the closed subspace of the homogeneous Besov space $\dot B^{-s}_{p,q}(\R)$ consisting of the distributions whose Fourier transform is supported in $\R_+$. 
Furthermore, the so-defined mapping 
\[
\Bc\colon A^{p,q}_s(\C_+)\ni f \mapsto \lim_{y\to 0^+} f(\,\cdot\,+i y)\in \dot B^{-s}_{p,s}(\R)
\]
is an isomorphism onto its image
For $p=q=2$, this latter result is closely related to the Paley--Wiener representation theorems, which assert that the elements of $A^{2,2}_s(\C_+)$ are the Laplace--Fourier transforms of the elements of $L^2(\abs{\,\cdot\,}^{-2 s}\cdot \Hc^1)$ supported in $\R_+$.

Finally, suitable inclusions between the spaces $A^{p,q}_s(\C_+)$, which are the counterparts of the Sobolev embeddings between the corresponding Besov spaces $B^{-s}_{p,q}$, allow to characterize the dual of $A^{p,q}_s(\C_+)$ by means of the sesquilinear form, for $p,q\in ]0,\infty[$ and $s,s'>0$,\footnote{Define $p'\coloneqq \max(1,p)'$, so that $p'=\infty$ if $p\meg 1$ and $\frac 1 p+\frac{1}{p'}=1$ if $p\Meg 1$. Define $q'$ analogously.}
\[
A^{p,q}_s(\C_+) \times A^{p',q'}_{s'}(\C_+)\ni (f,g)\mapsto \int_{\C_+} f(x+i y)\overline{g(x+i y)} y^{s+s'+(1/p-1)_+}\,\dd x \frac{\dd y}{y}.
\]
The preceding sesquilinear forms then induce continuous sesquilinear forms on $A^{p,q}_s(\C_+) \times A^{p',q'}_{s'}(\C_+)$ for all $s,s'\in \R$ which induce an antilinear isomorphism of $A^{p',q'}_{s'}(\C_+)$ onto $A^{p,q}_s(\C_+)'$. In particular, this identifies the Bloch space $A^{\infty,\infty}_0(\C_+)$ with the dual of the (unweighted) Bergman space $A^{1,1}_1(\C_+)$.

\medskip

The preceding results can be (and have been) investigated in more general contexts. 
One first observes that $\R_+^*$ is the simplest example of a symmetric cone, that is, a self-dual open cone (with respect to some chosen scalar product) not containing any affine line, on which the group of affine automorphisms acts transitively. Then, one may consider an arbitrary symmetric cone $\Omega\subseteq  \R^m $ and the associated tube domain, which is also called a Siegel domain of type I, 
\[
D=\R^m  +i \Omega.
\] 
In this case,  $D$ can be foliated as the union of the translates $
\R^m + i h$ of the `\v{S}ilov boundary' or `distinguished boundary' 
\[
b D= \R^m 
\]
of $D$, that is, the smallest closed subset of $\partial D$ on which every bounded continuous function on $\overline D$ which is holomorphic on $D$ has the same supremum as on $\overline D$.
In this case, the weight function on $\Omega$ can be replaced by a power of the characteristic function
\[
\varphi\colon D\ni z \mapsto \int_\Omega e^{i \langle h, z\rangle}\,\dd h
\] 
of $\Omega$ or, more generally, by a `generalized power function'  $\Delta^{\vect s}_\Omega: \Omega\to\R^*_+$, where $\vect s\in\R^r$ and $r$ is the \emph{rank} of the cone $\Omega$. Throughout the whole
work, we denote by $\Hc^k$  the  $k$-dimensional Hausdorff measure and by $\dd (\cdot)$ the integration in the given 
variable  with respect to  the appropriate  Hausdorff  measure.  Given the symmetric, or more generally \emph{homogeneuous}, cone $\Omega$, there exists a vector $\vect d\in \R^r$ canonically associated to $\Omega$ such that $\nu_\Omega:= \Delta_\Omega^{\vect d}\cdot\Hc^m$ is the unique (up to a multiplicative constant) Radon measure invariant under the linear transformations that preserve $\Omega$ (see Definitions~\ref{16-bis}  and~\ref{def:3.5}).
Then, we consider the spaces
\[
A^{p,q}_{\vect s}(D) = \Set{f\in \Hol(D)\colon \int_\Omega  \big( \norm{f(\,\cdot\,+ih)}_{L^p(\R^m)} \Delta^{\vect s}_\Omega(h)\big)^{q} \, \dd \nu_\Omega(h)<\infty}
\] 
(modification when $\max(p,q)=\infty$). 
In this case, the boundary values of the elements of $A^{p,q}_{\vect s}(D)$ can be naturally embedded in a space of Besov type $B^{-\vect s}_{p,q}(b D)$ under suitable assumptions, even though the corresponding (continuous and one-to-one) boundary mapping
\[
\Bc\colon A^{p,q}_{\vect s}(D)\to B^{-\vect s}_{p,q}(b D)
\]
is not an isomorphism, in general. 
Observe that, when $D=\C_+$, one has natural isomorphisms $(-\Delta)^{s'/2}\colon\dot B^s_{p,q}(\R)\to \dot B^{s-s'}_{p,q}(\R)$, where $\Delta$ denotes the second derivative on $\R$. In the general case, one may define suitable generalized Riemann--Liouville operators which induce isomorphisms $B^{\vect s}_{p,q}(b D)\to B^{\vect s-\vect{s'}}_{p,q}(b D)$. Nonetheless, such operators are in general far from being (fractional powers of) hypoelliptic operators. 
For this reason, while in the case of $\C_+$ it is relatively easy to study the spaces $A^{p,q}_s(D)$ for $s\meg 0$ by means of suitable Taylor expansions, in the general case the situation is far more involved (cf., e.g.,~\cite{VergneRossi}).

For what concerns sampling and atomic decomposition, they can be investigated replacing the sets $Q_{R,j,k}$ with the balls $B(z_{j,k},R)$ associated with some invariant distance, where  $(z_{j,k})$ is an  $R$-lattice on $D$.
In a similar way one may also study tube domains over homogeneous cones, that is, open cones not containing any affine line, on which the group of affine automorphisms acts transitively. 

\medskip

The preceding considerations may be further extended to homogeneous Siegel domains of type II, that is, homogeneous domains of the form 
\[
D=\Set{(\zeta,z)\in \C^n \times \C^m\colon \Im z-\Phi(\zeta,\zeta)\in \Omega},
\]
where $\Omega$ is an open convex cone in $\R^m$ which does not contain any affine line and
\[
\Phi\colon \C^n\times \C^n\to \C^m
\]
is a hermitian mapping such that $\Phi(\zeta,\zeta)\in \overline \Omega \setminus\Set{0} $ for every non-zero $\zeta\in \C^n$. The homogeneity of $D$ is equivalent to the homogeneity of $\Omega$ and a suitable condition on $\Phi$ (cf., e.g.,~\cite[Theorem 2.3]{Xu}). 
In this case, the \v Silov boundary 
\[
b D= \Set{(\zeta,z)\in \C^n \times \C^m\colon \Im z-\Phi(\zeta,\zeta)=0}
\]
of $D$ can be identified with $\Nc=\C^n\times \R^m$ by means of the mapping 
\[
(\zeta,x)\mapsto (\zeta,x+i\Phi(\zeta,\zeta)).
\]
As before, $D$ can be foliated as the union of the translates $b D+i h$, $h\in \Omega$, of $b D$. 
Furthermore, $\Nc$ can be endowed with the structure of a $2$-step
stratified group which acts holomorphically on $D$.    The Haar measure
on $\Nc$ is both left and right invariant and coincides with the
Lebesgue (i.e. Hausdorff) measure $\Hc^{2n+m}$.  
In this setting, we consider the spaces
\[
A^{p,q}_{\vect s}(D) = \bigg\{ \, f\in \Hol(D)\colon\\ 
\int_\Omega \bigg( \int_\Nc \abs{f(\zeta,x+i\Phi(\zeta,\zeta)+i h)}^p \, \dd (\zeta,x) \bigg)^{q/p}  \Delta^{q \vect s}_\Omega (h)  \, \dd \nu_\Omega(h)<\infty\, \bigg\}
\]
(modification when $\max(p,q)=\infty$). 

In this situation, though, a more prominent lack of symmetry can generate new phenomena. 
For example, if $\overline \Omega$ is the closed convex envelope of  $\Set{\Phi(\zeta,\zeta)\colon \zeta\in \C^n}$, then, at least for $p\in [1,\infty]$,
\[
H^p(D)= \Set{ f\in \Hol(D)\colon \sup\limits_{h\in \Omega\cap U} \int_\Nc \abs{f(\zeta,x+i\Phi(\zeta,\zeta)+i h)}^p\,\dd (\zeta,x)<\infty}
\] 
(modification if $p=\infty$) for every neighbourhood $U$ of $0$ in $\Omega$. 
On the contrary, in the classical case $\C_+=\R+i \R_+^*$, a holomorphic function $f$ on $D$ such that $\sup\limits_{0<y\meg R} \norm{f(\,\cdot\,+i y)}_{L^p(\R)}$ for every $R>0$ can grow exponentially fast at infinity. This is the case, e.g., for the function $f\colon z \mapsto \int_{-1}^1 e^{i t z}\,\dd t=2\frac{\sin z}{z}$ on $\C_+$.

In addition, in order to define the Besov type spaces  $B^{\vect s}_{p,q}(bD)$, while identifying $bD$ with $\Nc$, one needs to make use of the non-commutative Fourier transform associated with $\Nc$, which is far less manageable that in the commutative case. Even in the simple case of the Heisenberg groups, the characterization of the image of the Schwartz space under the Fourier transform provided in~\cite{Geller} is quite complicated. Consequently, working with the Fourier transform of general tempered distributions may not be convenient.

\medskip

The aforementioned problems have been studied at different levels of
generality.   We point out that 
there exists also a vector $\vect b \in \R^r$, canonically
associated with $D$, that plays an important role  in the analysis,
see Definitions~\ref{16-bis}  and~\ref{def:3.5}. We now indicate is a brief list  of papers most relevant to the
present work. 
\begin{itemize}
\item Gindikin in~\cite{Gindikin} studied the  {\it
unweighted} Bergman space $A^{2, 2}_{-\vect d/2}(D)$ and the Hardy space $H^2(D)$  for general Siegel domains of type II, established Plancherel formulae and determined the reproducing kernels;
	
\item Ogden and V\'agi in~\cite{OgdenVagi} studied the spaces $H^2(D)$  for general Siegel domains of type II in connection with the Fourier transform of the group $\Nc$; 
	
\item  Coifman and Rochberg in~\cite{CoifmanRochberg} studied the spaces $A^{p,p}_{s(\vect b+2\vect d)-\vect d/p}(D)$,
$p\in ]0,\infty[$, on symmetric \footnote{See Definition~\ref{4}. If  $D$ is symmetric, then it is also homogeneous.}  Siegel   domains of type II, and investigated atomic decomposition;
	
\item Ricci and Taibleson in~\cite{RicciTaibleson} established almost all the results considered in the discussion above for the spaces $A^{p,q}_s(\C_+)$;
			
\item Békollé and Temgoua Kagou in~\cite{BekolleKagou} studied the spaces $A^{p,p}_{\vect s}(D)$, $p\in [1,\infty[$, on homogeneous Siegel domains of type II, and investigated the properties of the Bergman projectors;
	
\item Békollé in~\cite{Bekolle} studied the spaces $A^{p,p}_{s\vect 1_r}(D)$, $p\in ]0,1]$, on irreducible symmetric Siegel domains of type II, and characterized their duals;
	
\item Bonami in~\cite{Bonami} studied the spaces $A^{p,q}_{s\vect 1_r}(D)$, $p,q\in [1,\infty[$, on irreducible symmetric tube domains, and surveyed results concerning the boundedness of the Bergman projectors and the boundary values;\footnote{Here, and in the rest of the paper, $\vect 1_r\coloneqq (1,\dots, 1)\in \N^r$.}
	
\item Békollé, Bonami, Garrig\'os, Nana, Peloso, and Ricci in~\cite{Bekolleetal} studied the spaces $A^{p,q}_{s\vect 1_r}(D)$, $p,q\in [1,\infty[$, on irreducible symmetric tube domains, and investigated their atomic decomposition and the boundedness of the Bergman projectors;
	
\item Békollé, Bonami, Garrig\'os, and Ricci in~\cite{BekolleBonamiGarrigosRicci} studied the spaces $A^{p,q}_{s\vect 1_r}(D)$ on irreducible symmetric tube domains and investigated their boundary values and the boundedness of the Bergman projectors;

	\CR 
\item Ciatti and Ricci in~\cite{CiattiRicci} studied the boundedness of the unweighted Bergman projectors on unweighted Bergman spaces on Siegel domains of type II over (not necessarily homogeneous) polyhedral cones; 

\item Debertol in~\cite{Debertol} studied the spaces $A^{p,q}_{\vect s}(D)$ on irreducible symmetric tube domains and investigated their boundary values and the boundedness of the Bergman projectors;\CB
	
\item Nana and Trojan  in~\cite{NanaTrojan} studied the spaces $A^{p,p}_{\vect s}(D)$, $p\in [1,\infty[$ on irreducible symmetric tube domains, and investigated the boundedness of the Bergman projectors; 
	
\item  Nana in~\cite{Nana} studied the spaces $A^{p,q}_{\vect s}(D)$, $p,q\in [1,\infty[$, on homogeneous Siegel domains of type II, and investigated the boundedness of the Bergman projectors;
	
\item Békollé, Ishi, and Nana   in~\cite{BekolleIshiNana} studied the spaces $A^{p,p}_{\vect s}(D)$ on homogeneous Siegel domains of type II, and investigated their atomic decomposition; 
	
\item Békollé, Gonessa, and Nana, in~\cite{BekolleGonessaNana} study the spaces $A^{p,q}_{\vect s}(D)$, $p,q\in [1,\infty]$, on irreducible symmetric tube domains, and investigate their atomic decomposition and complex interpolation spaces;
	
\item Christensen in~\cite{Christensen2} studied the spaces $A^{p,q}_{s\vect 1_r}(D)$, $p,q\in [1,\infty]$, on symmetric tube domains, and investigated a different kind of atomic decomposition thereon;
	
\item Arcozzi, Monguzzi, Peloso, and Salvatori in~\cite{ArcozziMonguzziPelosoSalvatori} studied the spaces $A^{2,2}_{s}(D)$ on the symmetric Siegel domains of type II associated with the cone $\R_+^*$, and established Paley--Wiener theorems;
\CR
\item Békollé, Gonessa, and Nana in~\cite{BekolleGonessaNana2} studied the spaces $A^{p,q}_{\vect s}(D)$, $p,q\in [1,\infty[$, on homogeneous Siegel domains of type II, and investigated the boundedness of the Bergman projectors.\footnote{\CR Notice, though, that when $D$ is not a tube domain and $p\neq q$, then the definition of $A^{p,q}_{\vect s}(D)$ adopted in~\cite{Nana,BekolleGonessaNana2} is different from ours. \CB}\CB
\end{itemize}
The preceding description of the available literature (which makes no claim to be exhaustive) should suggest how rarely the general study of the spaces $A^{p,q}_{\vect s}(D)$, $p,q\in ]0,\infty]$, $\vect s\in \R^r$, on homogeneous Siegel domains of type II has been undertaken. This is one of the main objects of the present  research  monograph.

\medskip

Before describing the structure of each chapter, we shall briefly discuss an issue concerning the notion of a \emph{weighted} Bergman space. In the literature, two different interpretations of the notion of a weight appear. 
The first one is that employed in~\cite{RicciTaibleson}, and is the one we described above and shall consistently adopt in Chapters~\ref{sec:6} and~\ref{sec:7}. In this case, the weight $\Delta^{\vect s}_\Omega$ acts as a multiplier of the function and not of the measure. With this convention, the space of boundary values of $A^{p,q}_{s}(\C_+)$ embeds canonically in the classical homogeneous Besov space $\dot B^{-s}_{p,q}(\R)$.

The second one is that employed in~\cite{BekolleKagou,Bekolle,Bonami,Bekolleetal,BekolleBonamiGarrigosRicci,NanaTrojan,Nana,BekolleIshiNana,BekolleGonessaNana,Christensen2}. 
In this case, the weight $\Delta^{\vect s}_\Omega$ acts as a multiplier of the measure, and the resulting space $\Ac^{p,q}_{\vect s}(D)$ is then equal to $A^{p,q}_{ (\vect s+\vect b/2)/q }(D)$,  where $\vect b\in \R^r$ is as above.  
Even though this choice is quite common in the literature, it still has some drawbacks. On the one hand, it fails to include the spaces $A^{p,\infty}_{\vect s}(D)$ for $\vect s\neq \vect 0$. On the other hand, it suggests considering only the Bergman projector $P_{(\vect s+\vect b/2)/2}$ on $\Ac^{p,q}_{\vect s}(D)$, thus ignoring the more general pattern that appears taking into account the full scale of the projectors $P_{\vect {s'}}$, as (more or less explicitly) done in~\cite{RicciTaibleson}.
The aforementioned issues have also undesirable (minor) consequences, such as: the inclusions between the spaces $\Ac^{p,q}_{\vect s}(D)$ are usually described in terms of the growth of the $L^p$ or $L^\infty$-norms of the function on the translates of $b D$, since most of the spaces $A^{p,\infty}_{\vect s}(D)$ are missing. 
The correspondences between the indices $p,q,\vect s$ are more cumbersome in some situations since, e.g.,  $\Ac^{p,q}_{\vect s}(D)$ has boundary values in $B^{-(\vect s+\vect b/2)/q}_{p,q}(D)$, while the Riemann--Liouville operator corresponding to the generalized power function $\Delta^{\vect {s'}}_{\Omega'}$ on $\Omega'$ maps $\Ac^{p,q}_{\vect s}(D) $ into $\Ac^{p,q}_{\vect s+q \vect{s'}}(D)$ (for $\vect{s'}$ sufficiently large).

For these reasons, we preferred to follow~\cite{RicciTaibleson} and to consider the weight as a multiplier of the function and not of the measure. We shall only depart from this choice in Chapter~\ref{sec:1}, whose results are mostly propaedeutical to the chapters which follow.

\medskip

In Chapter~\ref{sec:1}, we shall first introduce the notion of a Siegel domain $D$ of type II, and define the group structure on its \v Silov boundary $b D$. We shall then briefly indicate some basic facts about the Fourier transform on $b D$, and then apply this machinery to the study of the weighted Bergman spaces $A^{p,q}_\mi(D)$, defined as
\[
A^{p,q}_\mi(D) =  \Set{\! f\in \Hol(D)\colon \int_\Omega \left(\int_\Nc  \abs{f(\zeta,x+i \Phi(\zeta,\zeta)+i h)}^p \dd (\zeta,x) \right) ^{q/p}\dd \mi(h)<\infty}
\]
(modification when $\max(p,q)=\infty$), where $\mi$ is a positive
Radon measure on $\Omega$ satisfying some minimal assumptions. 
 
This leads to a characterization of $A^{2,q}_\mi(D)$ by means of the Fourier transform (cf.~Propositions~\ref{prop:2} and~\ref{prop:4}), and to the determination (up to a Fourier inversion) of the reproducing kernels of the Hardy space $A^{2,\infty}_\mi(D)$ and of the weighted Bergman space $A^{2,2}_\mi(D)$. 
We shall also prove some results on the inclusions between the various weighted Bergman spaces (Proposition~\ref{prop:1}) and a density result (Proposition~\ref{prop:37}). 
Because of the general nature of the `weight' $\mi$, though, these results are only `local' in a suitable sense, and will only be extended to global results by homogeneity, in Chapter~\ref{sec:6}.

The results concerning  the use of the Fourier transform of $\Nc$ are basically a simplified version of~\cite{OgdenVagi}. The main novelty here is the generality of the `weight' $\mi$, which can lead to some new phenomena in some cases: for example, it may happen that, given $f\in A^{2,2}_\mi(D)$ and $h\in \Omega$, the function $f(\,\cdot\,+i h )$ does \emph{not} belong to the Hardy space $A^{2,\infty}_\mi(D)$. This is the case, e.g., when $D=\C_+$ and $\dd \mi(y)= e^{-y}\,\dd y $.

In Chapter~\ref{sec:3}, we shall introduce the notation on homogeneous Siegel domains which we shall consistently employ in the chapters that follow, and we shall prove several technical lemmas. In this chapter we shall introduce our notation for the `generalized power functions' $\Delta^{\vect s}_\Omega$ and $\Delta^{\vect s}_{\Omega'}$ on $\Omega$ and $\Omega'$ resp.
We shall prove some results on the gamma functions on $\Omega$ and $\Omega'$, most of which are already present in the literature (cf.~\cite{Gindikin,NanaTrojan}). We shall then introduce our notation for the Riemann--Liouville operators and prove some basic facts about them.
Finally, we shall recall some basic facts on the canonical invariant K\"ahler metric and the invariant measure on $D$, and we shall prove a general version of `Kor\'anyi's lemma' (Theorem~\ref{teo:2}), extending~\cite[Theorem 1.1]{BekolleIshiNana}. 
This result is of fundamental importance to deal with atomic decomposition. We shall then define suitable lattices on $D$ and $\Omega$ and prove the quasi-constancy on invariant balls of several important functions on $D$, $\Omega$, and $\Omega'$.

In Chapter~\ref{sec:6}, we shall introduce the weighted Bergman spaces $A^{p,q}_{\vect s}(D)$ on the homogeneous Siegel domain $D$ of type II and study in full generality some of the problems we presented for $\C_+$.  This study will be completed in Chapter~\ref{sec:7}.
We shall then prove inclusion results between the Bergman spaces (Propositions~\ref{prop:5} and~\ref{prop:22}),  and density results for their intersections (Proposition~\ref{prop:8}). 
We shall then determine an explicit formula for the reproducing kernel of the Bergman space $A^{2,2}_{\vect s}(D)$ (Proposition~\ref{prop:6}) and determine its reproducing properties on other spaces (Proposition~\ref{prop:9}). We shall then briefly develop Paley--Wiener theorems for $A^{2,2}_{\vect s}(D)$, extend them below the critical index following~\cite{ArcozziMonguzziPelosoSalvatori} (Proposition~\ref{prop:54}), and determine an analogue of the classical Dirichlet space in this context (Proposition~\ref{prop:34}).

Next, we shall prove some results on sampling (Theorems~\ref{teo:3} and~\ref{teo:8}) and on atomic decomposition (Propositions~\ref{prop:16} and~\ref{prop:17}, and Theorems~\ref{teo:4} and~\ref{teo:5}), following~\cite{RicciTaibleson,Bekolleetal}. Optimality is achieved in some situations, even though the general picture is still unclear.
We then apply the preceding results to characterize the dual of $A^{p,q}_{\vect s}(D)$ under various assumptions (Propositions~\ref{prop:29} and~\ref{prop:13}), following~\cite{RicciTaibleson,Bekolleetal}.

In Chapter~\ref{sec:5}, we shall introduce the spaces of Besov type $B^{\vect s}_{p,q}(b D)$ and develop some basic facts of the associated theory, such as duality results and the compatibility with the Riemann--Liouville operators previously defined, which play the same role as the fractional Laplacian for homogeneous Besov spaces.
We shall not delve deeply into the theory of Besov spaces (thus proving results on multiplication, convolution, Fourier multipliers, or interpolation), since our main focus is in their connection with the weighted Bergman spaces.

In Chapter~\ref{sec:7}, we shall complete the study of weighted Bergman spaces started in Chapter~\ref{sec:6}.
We shall study the boundary values of the elements of $A^{p,q}_{\vect s}(D)$ (Proposition~\ref{prop:23} and Theorem~\ref{teo:6}), following~\cite{BekolleBonamiGarrigosRicci}. 
We shall then prove the equivalence of several conditions concerning the boundary values of the spaces $A^{p,q}_{\vect s}(D)$, their atomic decompositions, and the boundedness of the Riemann--Liouville operators (Corollary~\ref{cor:32}).
Finally, we shall study the boundedness of the Bergman projectors, following~\cite{BekolleBonamiGarrigosRicci} (Theorem~\ref{teo:7}).

In the Appendix~\ref{sec:app}, we shall collect a few results and notation on mixed norm spaces. After recalling the definition of the mixed norm space $L^{p,q}(\mi,\nu)$, for $p,q\in ]0,\infty]$, we shall characterize its dual in two specific situations (Propositions~\ref{prop:51}, and~\ref{prop:52}) and compare $L^{p,q}(\mi,\nu)$ with the space $L^q(\mi;L^p(\nu))$ (Propositions~\ref{prop:49} and~\ref{prop:50}). For this latter comparison to make sense, we shall briefly recall the definition and some properties of the Bochner spaces $L^p(\mi;Z)$, where $Z$ is a locally bounded $F$-space (cf.~Definition~\ref{def:7}). 

We shall adhere as far as possible to the commonly adopted notation. Possible exceptions to this rule are explicitly defined at their first occurrence and then listed in the Index of Notation, to which the reader is referred for any unspecified piece of notation.

\medskip

We conclude the introduction with some notational remarks. 
We first notice that the spaces $A^{p,p}_{\vect s}(D)$ are the (`pure norm') weighted Bergman spaces, the unweighted Bergman spaces if $\vect s=-\vect d/q$. 
In order to make the notation as uniform as possible, we shall denote by $A^{p,\infty}_{\vect 0}(D)$ the Hardy spaces $H^p(D)$. The (generalized) Dirichlet space will be denoted by $\widehat A^{2,2}_{(\vect b+\vect d)/2,\vect{s'}}(D)$, while the (generalized) Bloch space will be denoted by $\widehat A^{\infty,\infty}_{\vect 0,\vect{s'}}(D)$.\footnote{The element $\vect{s'}$ of $\N^r$ is relevant only for the definition of $\widehat A^{2,2}_{(\vect b+\vect d)/2,\vect{s'}}(D)$ and $\widehat A^{\infty,\infty}_{\vect 0,\vect{s'}}(D)$ as sets, as well as for the precise determination of their norms. As a matter of fact, changing $\vect{s'}$ gives rise to canonically isomorphic spaces.}

For a better treatment of the duality between the various weighted Bergman spaces $A^{p,q}_{\vect s}(D)$, it will be convenient to introduce some smaller spaces $A^{p,q}_{\vect s,0}(D)$, which are defined by some vanishing conditions `at infinity.'
In particular, $A^{p,q}_{\vect s}(D)=A^{p,q}_{\vect s,0}(D)$ when $p,q<\infty$, so that the difference between $A^{p,q}_{\vect s}(D)$ and $A^{p,q}_{\vect s,0}(D)$ only matters when either $p$ or $q$ equals $\infty$.

In order to deal with both $A^{p,q}_{\vect s}(D)$ and $A^{p,q}_{\vect s,0}(D)$ (and in similar circumstances), we shall make considerable use of statements of the form
\begin{center}
	`Assume that $A$ (resp.\ $B$ if $C$) holds. Then $A'$ (resp.\ $B'$ if $C'$) holds,'
\end{center}
where $A$, $B$, $C$, $A'$, $B'$, and $C'$ denote suitable properties. We believe that these statements should be clear from the context.\footnote{More precisely, the above statement is equivalent to the following lengthier one (expressed in symbols for simplicity):
	\begin{center}
		`$(A\implies A')\land [((A\land \lnot C) \lor (B \land C) )\implies ((A'\land \lnot C') \lor (B' \land C') )]$.'
\end{center} }
Sometimes, $C$ or $C'$ may be empty. We hope that this shortened notation be of help to understand the various statements.

Finally, we remark that we shall not denote explicitly the measure in integrals and Lebesgue spaces when it reduces to a suitably normalized Lebesgue measure.

\mainmatter

\chapter{Elementary Theory}\label{sec:1}

In this chapter we shall introduce our main notation on Siegel domains of type II over open convex cones not containing affine lines. 
We begin with an open convex cone $\Omega$, not containing any affine line, in a finite-dimensional real vector space $F$ of dimension $m$, and an $\Omega$-positive hermitian mapping $\Phi\colon E\times E\to F$, where $E$ is a finite-dimensional complex vector space of dimension $n$. 
One may wonder why $E$ and $F$ are not directly identified with $\C^n$ and $\R^m$, respectively, by a suitable choice of coordinates. On the one hand, this would cause some notational issues, since we shall consider on $E$ several complex structures other than the original one, giving rise to different spaces $E_\lambda$. On the other hand, we believe that this choice marks more clearly the different roles played by $E$ and $F$, which might not appear identifying $E\times F_\C$ with $\C^{n+m}$ (we denote by $F_\C$ the complexification of $F$, that is, $F\otimes_\R \C$).

Denoting by $D$ the corresponding Siegel domain, we observe that the action of the \v Silov boundary $b D$ on $D$  stems from a $2$-step nilpotent group structure on the ambient space $E\times F_\C$ which turns $D$ in a semigroup.

In Section~\ref{sec:1:2}, we recall some basic facts on the Fourier transform on $\Nc$ and establish some notation. See~\cite{OgdenVagi} for a more detailed study of the Fourier transform on $\Nc$ and for a direct comparison of the Bargmann and Schr\"odinger representations. 
We shall content ourselves with a brief exposition of the main properties of the Bargmann representations, which are particularly well adapted to the study of Siegel domains of type II.

In Section~\ref{sec:1:3}, we study a provisional version of the weighted Bergman spaces we shall study in Chapter~\ref{sec:6}. As we mentioned in the introduction, in this chapter it will be convenient to choose a positive Radon measure $\mi$ as the weight. Nonetheless, the role of $\mi$ will be of importance only when dealing with the spaces $A^{2,q}_\mi(D)$ (cf.~Section~\ref{sec:1:4}), for which the Plancherel formula allows a more detailed treatment. For other spaces, we shall only deal with `local' properties of such spaces, for which the role of $\mi$ is virtually irrelevant (provided that it satisfies some minimal requirements).

In Section~\ref{sec:1:4}, we apply the results of Section~\ref{sec:1:2} to the study of the spaces $A^{2,q}_\mi(D)$. In particular, by means of Propositions~\ref{prop:2} and~\ref{prop:4}, we shall characterize the spaces $A^{2,q}_\mi(D)$ by means of the Fourier transform, thus establishing some Paley--Wiener theorems. 
These results will then be applied to find the reproducing kernels of the Hardy space $A^{2,\infty}_\mi(D)$ (Corollary~\ref{cor:4}) and of the weighted Bergman spaces $A^{2,2}_{\mi}(D)$ (Corollary~\ref{cor:5}). 
While the procedure to find Corollary~\ref{cor:4} is essentially the same as the one of~\cite{OgdenVagi}, Corollary~\ref{cor:5} seems to be new, at least in the stated generality.

\section{General Definitions}\label{sec:1:1}

In this section we introduce the notation we shall adopt for Siegel domains of type II, the associated boundary, the group structure of the latter, and its action on the domain.

\medskip

Denote by $E$ a vector space over $\C$ of finite dimension $n$ and by $F$ a vector space over $\R$ of finite dimension $m>0$.\footnote{Notice that we do \emph{not} impose the condition $n>0$, so that our analysis applies to tube domains as well. The condition $m>0$ is nonetheless necessary to avoid trivialities.} 
Let $\Omega$ be an open convex cone with vertex $0$ in $F$ and not containing any affine lines (that is, such that its closure $\overline{\Omega}$ is proper).\label{1} 
Denote by $\Omega'$ the dual cone\index{Cone!dual} of $\Omega$, that is, the interior of the polar of $\overline \Omega$ in $F'$. More explicitly,
\[
\Omega'=\Set{\lambda\in F'\colon \langle \lambda, h \rangle>0 \quad\forall h\in \overline \Omega\setminus\Set{0}}.
\]
We shall endow $F$ with the (partial) ordering induced by $\overline \Omega$, so that `$x\meg y$' means `$y-x\in \overline \Omega$' for every $x,y\in F$.\index{Cone!induced ordering}

Throughout the whole presentation, we shall denote by 
$\langle\,\cdot\,,\,\cdot\,\rangle$ a bilinear pairing and by
$\langle\,\cdot\,\vert\,\cdot\,\rangle$  a sesquilinear pairing, without
explicitely indicating the spaces on which these pairings are
defined. Such spaces will be clear from the context and we  believe
that such abuse of notation will cause no confusion.  

Let $\Phi\colon E\times E\to F_\C$ be a positive non-degenerate hermitian mapping, that is:
\begin{enumerate}
\item[(1)] $\Phi$ is linear in the first argument; 
	
\item[(2)] $\Phi(\zeta,\zeta')=\overline{\Phi(\zeta',\zeta)}$ for every $\zeta,\zeta' \in E$; 
	
\item[(3)] $\Phi$ is non-degenerate; 
	
\item[(4)] $\Phi(\zeta)\coloneqq \Phi(\zeta,\zeta)\in \overline \Omega$ for every $\zeta\in E$.
\end{enumerate}

Observe that, assuming conditions~{(1)},~{(2)}, and~{(4)}, condition~{(3)} is equivalent to the following one:
\begin{itemize}
\item[{(3\ensuremath{'})}] $\Phi(\zeta,\zeta)\neq 0$ for every non-zero $\zeta \in E$.
\end{itemize}

\begin{deff}%%\label{def:Siegel-domain}
Define the  Siegel domain of type II  associated with the cone $\Omega$ and the mapping $\Phi$ as\label{3} 
\[
D\coloneqq \Set{(\zeta,z)\in E\times F_\C\colon \Im z-\Phi(\zeta)\in \Omega},
\]
\index{Siegel domain!of type II} and denote by \label{4}
\[
b D\coloneqq \Set{(\zeta,z)\in E\times F_\C\colon \Im z-\Phi(\zeta)=0}
\]
its \v{S}ilov boundary.

If $E=\Set{0}$, $D$ is in particular called a Siegel domain of type I or a tube
domain.\index{Siegel domain!of type I} \index{Tube domain}. 

The domain $D$ is said to be homogeneous
if  the  the group $G(D)$ of biholomorphisms of $D$ onto
itself acts transitively on $D$.\index{Siegel domain!homogeneous}

Moreover, $D$ is said to be symmetric if for every $(\zeta,z)\in D$
there exists  $\varphi\in G(D)$  such that
$\varphi\circ \varphi$ is the identity and $(\zeta,z)$ is an isolated
fixed point of $\varphi$.\index{Siegel domain!symmetric} 
Finally, $D$ is said to be irreducible if 
it is not holomorphically equivalent to a product of two 
Siegel domains. \index{Siegel domain!irreducible} 
\end{deff}

\begin{oss}
We explicitly point out that, according to our Definition~\ref{3}, the class of Siegel domains of type II includes the  class of  Siegel domains of type I, that is, we allow the space $E$ to reduce to $\Set{0}$.  There is no common agreement in the literature on this point.  For instance, in~\cite{Gindikin} and~\cite{CoifmanRochberg} the same convention is adopted, the references~\cite{Bekolle}  and~\cite{BekolleIshiNana} implicitly adopt it too. On the other hand, the works~\cite{Bonami-et-al} and~\cite{Damek-et-al} assume that $E\neq\Set{0}$ and therefore a domain of type II cannot be of type I as well.

We make the choice of including the tube type domains in the definition of Siegel domains of type II for a practical reason.  All our results are proved in both cases $E\neq\Set{0}$ and $E=\Set{0}$, unless explicitly indicated, and it becomes simpler to deal with a single class of domains.
\end{oss} \CB 

We recall that  if $D$ is symmetric, then it is also homogeneous.   In addition, recall that
  $b D$ can be characterized as the smallest closed subset $B$ of $\overline D$ such that $\sup\limits_{\overline D} \abs{f}= \sup\limits_{B} \abs{f}$ for every bounded continuous function $f\colon \overline D\to \C$ which is holomorphic on $D$  (cf.~\cite[Lemma 2.2]{Murakami}).\index{Silov boundary@\v Silov boundary}

\smallskip

The proof of the following result is straightforward and left to the reader. Notice that it both endows $b D$ with a group structure and provides a holomorphic action of $b D$ on $D$ by left translations, as we shall see below.

\begin{lem}\label{lem:3}
The mapping $\cdot\colon (E\times F_\C)\times (E\times F_\C)\to E\times F_\C$, defined by
\[
(\zeta,z)\cdot(\zeta',z')\coloneqq (\zeta+\zeta', z+z'+2 i \Phi(\zeta',\zeta))
\]
for every $\zeta,\zeta'\in E$ and for every $z,z'\in F_\C$, endows $E\times F_\C$ with a nilpotent (real) Lie group structure of step $2$ with centre $\Set{0}\times F_\C$.
The identity is $(0,0)$ and 
\[
(\zeta,z)^{-1}=(-\zeta,-z+2 i \Phi(\zeta))
\]
for every $(\zeta,z)\in E\times F_\C$.
\end{lem}

Notice that the derived subgroup $[E\times F_\C,E\times F_\C]$ of $E\times F_\C$ is the vector space generated by $\Set{0}\times\Im\Phi(E\times E)$, which is contained in, but not necessarily equal to, $\Set{0}\times F$.

\begin{deff}
	We endow $E\times F_\C$ with the product defined in Lemma~\ref{lem:3}.
\end{deff}

Observe that, if $(\zeta,z),(\zeta',z')\in E\times F_\C$, then
\[
\Im(z+z'+2 i \Phi(\zeta',\zeta))- \Phi(\zeta+\zeta')= \Im z-\Phi(\zeta)+\Im z - \Phi(\zeta'),
\]
so that 
\[
D\cdot \overline D, \overline D\cdot D\subseteq D, \qquad \overline D\cdot \overline D\subseteq \overline D, \qquad \text{and}\qquad b D\cdot b D\subseteq b D.
\] 
Thus, $D, \overline D$, and $b D$ are subsemigroups of $E\times F_\C$ and act analytically on both $D$ and $\overline D$. In addition, $b D^{-1}\subseteq bD$, so that $b D$ is a subgroup of $E\times F_\C$.

\begin{deff}
Define $\Nc$ as $E\times F$, endowed with the group structure given by the product\label{5}
\[
(\zeta,x)(\zeta',x')\coloneqq (\zeta+\zeta',x+x'+2\Im\Phi(\zeta,\zeta')).
\]
\end{deff}

Observe that the mapping
\[
\Nc \ni (\zeta,x)\mapsto (\zeta,x+i\Phi(\zeta))\in b D
\]
is an isomorphism of Lie groups.

\begin{deff}
	For every $f\colon D\to \C$ and for every $h\in \Omega$, define 
	\label{12}
	\[
	f_h\colon \Nc\ni (\zeta,x)\mapsto f(\zeta,x+i\Phi(\zeta)+i h)\in \C.
	\]
\end{deff}

\section{The Fourier Transform}\label{sec:1:2}

In this section we briefly discuss the Fourier transform on the $2$-step nilpotent group $\Nc$. Since we are mostly interested in the interaction of the Fourier transform with weighted Bergman spaces, we shall make use of the Bargmann representations, which are particularly well suited to the purpose.

\begin{deff}
Define, for every $\lambda\in F'$,
\[
B_\lambda\colon E\times E\ni (\zeta, \zeta') \mapsto \langle \lambda, \Im \Phi(\zeta, \zeta')\rangle\in \R
\]
and
\[
W\coloneqq\Set{\lambda\in F'\colon \exists \zeta\neq 0 \ \text{such
    that}\ 
  B_\lambda(\zeta,\,\cdot\,)=0 }.
\]
\end{deff}

Observe that $W$ is a proper algebraic variety of $F'$ and that $W=\emptyset$ if and only if $E=\Set{0}$.\label{6} 
In addition, $\Omega'\cap W=\emptyset$ since $B_\lambda$ is the imaginary part of the non-degenerate scalar product 
\[
E\times E\ni (\zeta, \zeta') \mapsto \langle \lambda_\C, \Phi(\zeta, \zeta')\rangle
\]
for every $\lambda\in \Omega'$ (here, $\lambda_\C=\lambda\otimes I_\C$).

We fix a scalar product on the real vector space $E\times F$ for which $E$ and $F$ are orthogonal and for which the multiplication by $i$ is a skew-adjoint endomorphism of $E$.
Observe that the Hausdorff measure $\Hc^{2n+m}=\Hc^{2n}\otimes \Hc^m$ is both left- and right-invariant on $\Nc$, so that we may choose it as the Haar measure on $\Nc$. We endow $E$ and $F$, and consequently $F'$, with the corresponding scalar products and Hausdorff measures.

\begin{deff}	
For every $\lambda\in F'\setminus W$, define $J_\lambda\in 
GL(E) $ so that 
\[
\langle \lambda_\C,\Phi(\zeta,i\zeta')\rangle=  \langle \zeta\vert J_\lambda \zeta' \rangle
\]
for every $\zeta,\zeta'\in E$.
	
In addition, define the (absolute value of the) Pfaffian\index{Pfaffian} of $J_\lambda$ (considered as an endomorphism of the real vector space underlying $E$)
\[
\abs{\Pfaff(\lambda)}\coloneqq \abs{{\det}_\C(J_\lambda)},\label{7}
\] 
where $\det_\C$ denotes the complex determinant,
and set\label{11} 
\[
\Lambda_+\coloneqq \Set{\lambda\in F'\colon \forall \zeta \neq 0\quad\langle \lambda, \Phi(\zeta)\rangle >0 }.
\]
\end{deff}

Observe that 
\[
B_\lambda(\zeta, \zeta')= \Re \langle \zeta\vert J_\lambda \zeta' \rangle
\]
for every $\zeta,\zeta'\in E$. In addition, the properties of $\Phi$ show that $J_\lambda^*=-J_\lambda$  for every $\lambda\in F'\setminus W$, and that $-i J_\lambda$ is positive if $\lambda\in \Lambda_+$, in which case $J_\lambda= i \abs{J_\lambda}$. Further,
\[
\abs{\Pfaff(\lambda)}= {\det}_\C(\abs{J_\lambda})=\abs{{\det}_\R(J_\lambda)}^{1/2}
\]
for every $\lambda\in F'\setminus W$,
 and $\abs{\Pfaff(\lambda)}=1$ for every
$\lambda\in F'$ if and only if $E=\Set{0}$. Finally,
$\Lambda_+=F'$ if and only if $E=\Set{0}$.

\CB

\begin{deff}
For $\lambda\not \in W$, define  
\[
J_\lambda'\coloneqq \abs{J_\lambda}^{-1} J_\lambda,
\]
and let $E_\lambda$ be the vector space $E$ endowed with the complex structure induced by $J_\lambda'$, i.e.,
\[(x+ i y)\cdot_\lambda \zeta\coloneqq x \zeta+  y J'_\lambda \zeta
\]
for every $x,y\in \R$ and for every $\zeta\in E_\lambda$. 
	
Define 
\[
\Phi_\lambda(\zeta,\zeta')\coloneqq \Im\Phi(J'_\lambda \zeta, \zeta')+i \Im \Phi(\zeta, \zeta') 
\] 
for every $\zeta,\zeta'\in E_\lambda$.
\end{deff}
Notice that $\Phi_\lambda$ is an $\R$-bilinear form on $E_\lambda$ which is $\C$-linear in the first argument, and that $\langle \lambda_\C, \Phi_\lambda \rangle$ is hermitian and positive. 
In addition,  $\lambda\in \Lambda_+$ if and only if $\lambda \not \in W$ and $J_\lambda'=i$ or, equivalently, $\Phi=\Phi_\lambda$. This is also equivalent to $E_\lambda=E$.

Observe, in addition, that $J'_\lambda$ and $i$ commute, so that $J'_\lambda$ is a $\C$-linear endomorphism of $E$. Since, in addition, $J'^2_\lambda=-I$, we see that $E$ is the direct sum of two orthogonal subspaces $E_{\lambda,+}$ and $E_{\lambda,-}$ such that $J'_\lambda= \pm i$ on $E_{\lambda,\pm}$. Then, $E_{\lambda}=E_{\lambda,+}$ if and only if $\lambda\in \Lambda_+$.

\medskip

Now, take $\lambda\not \in W$, and observe that $\Nc/\ker \lambda$ is isomorphic to the $(2 n+1)$-dimensional Heisenberg group (interpreting $\R$ as the $1$-dimensional Heisenberg group, by an abuse of language). Therefore, the Stone--von Neumann theorem (cf., e.g.,~\cite[Theorem 1.50]{Folland}) implies that there is (up to unitary equivalence) one and only one irreducible unitary representation $\pi_\lambda$ of $\Nc$ such that $\pi_\lambda(0,x)=\ee^{-i\langle \lambda, x\rangle} I$ for every $x\in F$. This representation is characterized by the fact that it possesses a cyclic unit vector whose corresponding diagonal coefficient is the function
\[
(\zeta,x)\mapsto \ee^{- i\langle \lambda, x\rangle - \abs{\abs{J_\lambda}^{1/2}\zeta  }^2}= \ee^{\langle \lambda_\C , -i x- \Phi_\lambda(\zeta)\rangle}
\] 
(cf.~\cite{OgdenVagi} and~\cite[Section 2]{AstengoCowlingDiBlasioSundari} for more details).
We can realize $\pi_\lambda$ as follows. 

\begin{deff}
Define $H_\lambda\coloneqq \Hol(E_\lambda)\cap
L^2(\nu_\lambda)$,\label{10} where $\nu_\lambda= \ee^{-2\langle
  \lambda, \Phi_\lambda(\,\cdot\,)\rangle}\cdot \Hc^{2 n}$, and where
$\Hol(E_\lambda)$ denotes the space of holomorphic functions on
$E_\lambda$.  Then, define \label{8}
\[
\pi_\lambda(\zeta,x)  \psi  (\omega)\coloneqq \ee^{\langle \lambda_\C, -i x +2 \Phi_\lambda(\omega, \zeta)-\Phi_\lambda(\zeta)\rangle}  \psi (\omega- \zeta)
\]
for every $\psi \in H_\lambda$, for every $\omega\in E_\lambda$, and for every $(\zeta,x)\in \Nc$.
\end{deff}

The representations $\pi_\lambda$ are also known as the Bargmann representations.\index{Bargmann representation}

\begin{deff}\label{def:9}
For $v\in E_\lambda$, define  
\[
\partial_{E_\lambda,v} \psi \coloneqq \textstyle{\frac{1}{2}} \big( 
  \partial_v-i\partial_{J'_\lambda v}  \big)    \psi  (0)
\]
for every $\R$-differentiable function $ \psi $ on $E_\lambda$, so that $\partial_{E_\lambda,v} \psi =\partial_v   \psi (0)$ for every $ \psi \in \Hol(E_\lambda)$. In particular, $\partial_{E_\lambda,v}$ is an element of the complexification of the real tangent space of $E_\lambda$ at $0$.
\end{deff}

Observe that $\partial_{E_\lambda,v}=\partial_{E,v}$ if $v\in E_{\lambda,+}$, while $\overline{\partial_{E_\lambda,v}}=\partial_{E,v}$ if $v\in E_{\lambda,-}$.

\begin{prop}\label{prop:36}
	The following hold:
\begin{itemize}
\item $\pi_\lambda$ is an  irreducible continuous unitary representation of  $\Nc$ in $H_\lambda$;
\item $e_{\lambda,0}\coloneqq \sqrt{\frac{2^n\abs{\Pfaff(\lambda)}}{\pi^n}} 
  \chi_{E_\lambda}$ is a unit vector in $H_\lambda$, and
\[
\langle \pi_\lambda(\zeta,x) e_{\lambda,0}\vert e_{\lambda,0}\rangle=\ee^{\langle \lambda_\C, -i x-\Phi_\lambda(\zeta)\rangle} 
\]
for every $(\zeta,x)\in \Nc$, 
where $\langle \cdot\,\vert\, \cdot \rangle $ denotes the inner
product in $H_\lambda$; 

\item for every $v\in E_{\lambda}$, for every $ \psi \in C^\infty(\pi_\lambda)$ and for every $\omega \in E_\lambda$.
\[
\dd\pi_\lambda(\partial_{E_\lambda,v}) \psi (\omega)=-\partial_{v}  \psi (\omega)
\]
and
\[
\dd \pi_\lambda(\overline{\partial_{E_\lambda,v}}) \psi (\omega)= 2\langle \lambda_\C, \Phi_\lambda(\omega,v)\rangle \psi (\omega).
\]	
\end{itemize}
\end{prop}

Cf.~\cite{OgdenVagi} for more details on the Fourier transform on $\Nc$, and for a connection between the representations $\pi_\lambda$ and the Schr\"odinger representations employed, e.g., in~\cite{AstengoCowlingDiBlasioSundari}.

The unit vectors $e_{\lambda,0}$ and the relative orthogonal projectors are of particular importance in the study of holomorphic functions, as we shall see in Proposition~\ref{prop:3}.

\begin{proof}
	Straightforward computations show that $\pi_\lambda$ is a unitary representation of $\Nc$ in $H_\lambda$. Continuity is established by means of~\cite[Theorem 16.28]{Yeh}. Further, arguing as in~\cite[Lemma 5.1]{Taylor}, one may prove that $\pi_\lambda$ is irreducible.
In order to prove that
\[
\langle\pi_\lambda( \zeta,x ) e_{\lambda,0}\vert e_{\lambda,0}\rangle=\ee^{\langle \lambda_\C, -ix-\Phi_\lambda(\zeta)\rangle} 
\]
 we observe  that  for all $\omega\in E$
\[
\langle \lambda, \Phi_\lambda(\omega)\rangle= \abs{\abs{J_\lambda}^{1/ 2} \omega}^2.
\]
Now, for $(\zeta,x)\in \Nc$, with some obvious change of variables we see that 
\[
\begin{split}
  \langle\pi_\lambda( \zeta,x ) e_{\lambda,0}\vert e_{\lambda,0}\rangle
  & =e^{\langle \lambda_\C, -ix-\Phi_\lambda(\zeta)\rangle} \frac{2^n\abs{\Pfaff(\lambda)}}{\pi^n}
  \int_{E} e^{-2\langle \lambda_\C, \Phi_\lambda(\omega)-\Phi_\lambda(\omega,\zeta)\rangle}\,\dd \omega\\
  &=  e^{\langle \lambda_\C, -ix-\frac{1}{2}\Phi_\lambda(\zeta)\rangle} \frac{2^n\abs{\Pfaff(\lambda)}}{\pi^n}
  \int_{E} e^{-2\langle \lambda_\C, \Phi_\lambda(\omega-\zeta/2) -i\Im \Phi_\lambda(\omega,\zeta)\rangle}\,\dd \omega  \\
  &=e^{\langle \lambda_\C, -ix-\frac{1}{2}\Phi_\lambda(\zeta)\rangle} \frac{1}{(2\pi)^n}
  \int_{E} e^{-\frac{1}{2}\abs{\omega}^2-i\Re \langle \omega\, |\, \abs{J_\lambda}^{-1/2} J_\lambda\zeta\rangle}\,\dd \omega\\
&=e^{\langle \lambda_\C, -i x-\frac{1}{2}\Phi_\lambda(\zeta)\rangle} e^{-\frac{1}{2}\abs{ \abs{J_\lambda}^{-1/2} J_\lambda\zeta  }^2  }\\
&=e^{\langle \lambda_\C, -i x-\frac{1}{2}\Phi_\lambda(\zeta)\rangle} e^{-\frac{1}{2}\langle \lambda, \Phi_\lambda(\zeta)\rangle},
\end{split}
\] 
recalling 
 that $\abs{\Pfaff(\lambda)}=\det_{\R}(\abs{J_\lambda}^{1/ 2})$, and  appling~\cite[Theorem 7.6.1]{Hormander1}). In particular, taking $(\zeta,x)=(0,0)$, we see that $\norm{e_{\lambda,0}}_{H_\lambda}=1$.
	
	The last assertions follow easily.
\end{proof}

Even though we have \emph{not} fully described the dual of $\Nc$ (when $E\neq \Set{0}$), we are now able to present the associated Plancherel formula. 
Cf.~\cite{OgdenVagi} for more details and a precise interpretation of the formalism of direct integrals.

\begin{deff}
We denote by $\Lin(H_\lambda)$ the space of endomorphisms of $H_\lambda$, and by $\Lin^2(H_\lambda)$ the space of Hilbert--Schmidt endomorphisms of $H_\lambda$. For $\lambda\in F'\setminus W$ and $f\in L^1(\Nc)$, define $\pi_\lambda(f)\in\Lin (H_\lambda)$ by setting
\[
\pi_\lambda(f)\psi = \int_\Nc  f(\zeta,x) \pi_\lambda (\zeta,x)\psi \, \dd(\zeta,x),
\]
 for $\psi\in H_\lambda$.
\end{deff}

\begin{cor}\label{cor:7}
The mapping 
\[
L^1(\Nc)\ni f \mapsto (\pi_\lambda(f))\in \prod_{\lambda \in F'\setminus W} \Lin(H_\lambda)
\]
induces an isometric isomorphism
\[
\Pc\colon L^2(\Nc)\to \frac{2^{n-m}}{ \pi^{n+m}} \int_{F'\setminus W}^\oplus \Lin^2(H_\lambda) \abs{\Pfaff(\lambda)}\dd \lambda.
\]
\end{cor}

\begin{proof}
The assertion follows from~\cite[Section 2]{AstengoCowlingDiBlasioSundari},~\cite[Theorem 5.2]{OgdenVagi}, and Proposition~\ref{prop:36} when $E\neq \Set{0}$, from standard Fourier analysis otherwise.
\end{proof}

\begin{deff}
We denote by $P_{\lambda,0}$\label{9} the self-adjoint projector of $H_\lambda$ onto the space generated by $e_{\lambda,0}$, with the notation of Proposition~\ref{prop:36}.
\end{deff}

Observe that $P_{\lambda,0}(H_\lambda)=\bigcap_{v\in E_\lambda} \ker \dd \pi_\lambda(\partial_{E_\lambda,v})$ for every $\lambda\in F'\setminus W$.

\smallskip

In the following proposition, we show that the Fourier transform of holomorphic functions is relatively simple. 
It does not say anything when $E=\Set{0}$.

\begin{prop}\label{prop:3}
	If $f\in \Hol(D)$, $h\in \Omega$ and $f_h\in L^1(\Nc)+L^2(\Nc)$, then 
	\[
	\pi_\lambda(f_h)=\chi_{\Lambda_+}(\lambda)\pi_\lambda(f_h) P_{\lambda,0}
	\]
	for almost every $\lambda\in F'\setminus W$.
\end{prop}

Cf.~\cite{Ricci,ArcozziMonguzziPelosoSalvatori} for the case $\Omega=\R_+^*$.

\begin{proof}
	Take $v\in E$ and let $\overline{Z_v}$ be the left-invariant vector field on $\Nc$ which induces $\overline{\partial_{E,v}}$ at $0$. 
	Then, 
	\[
	\overline{Z_v} f_h(\zeta,x)=\overline{\partial_{E,v}} \big[\omega\mapsto f( (\zeta,x+i\Phi(\zeta))\cdot(\omega,i h))\big]=0
	\]
	since $f$ is holomorphic and $b D$ acts holomorphically on $D$. 
	
	Fix a representative $\tau$ of the `mapping' $\lambda \mapsto \pi_\lambda(f_h)$, and observe that there is a negligible subset $N$ of $F'$, containing $W$, such that 
	\[
	\tau(\lambda) \dd \pi_\lambda(\overline{\partial_{E,v}})=0
	\]
	for every $\lambda \not \in N$ and for every $v$ in a countable dense subset of $E$, hence for every $v\in E$ by continuity.
	Then, Proposition~\ref{prop:36} implies that 
	\[
	\tau(\lambda)\langle\lambda, \Phi(\,\cdot\,,v)\rangle=0
	\]
	for every $\lambda\in \Lambda_+\setminus N$. By the arbitrariness of $v$ in $E$, it follows that $\tau(\lambda)$ vanishes at every (holomorphic) polynomial on $E_\lambda$ which vanishes at $0$, hence at every element of $H_\lambda$ which vanishes at $0$ (argue as in~\cite[Theorem 1.63]{Folland}). Therefore,	
	\[
	\tau(\lambda)=\tau(\lambda) P_{\lambda,0}
	\]
	for every $\lambda\in \Lambda_+\setminus N$.
	 
	If, otherwise, $\lambda \not \in \Lambda_+\cup N$, then we may take a non-zero $v\in E_{\lambda,-}$, so that 
	\[
	\tau(\lambda)\partial_{v}=-\tau(\lambda)\dd \pi_\lambda(\partial_{E_\lambda,v})=0
	\]
	by Proposition~\ref{prop:36}. Since the image of $\partial_{v}$ contains the set of (holomorphic) polynomials on $E_\lambda$, which is dense in $H_\lambda$ (argue as in~\cite[Theorem 1.63]{Folland}), this implies that $\tau(\lambda)=0$.
\end{proof}

\section{Bergman and Hardy Spaces}\label{sec:1:3}

In this section we shall introduce a \emph{provisional}  definition of weighted Bergman spaces. 
Most of the results proved in this section are relative to `local' properties of such spaces. For this reason, we shall also introduce the spaces $A^{p,q}_{\mi, \loc}(D)$ for a better treatment of such `local' properties. 

We observe explicitly that the spaces $A^{p,\infty}_\mi(D)$ defined below are simply the Hardy space $H^p(D)$ if $\Supp{\mi}=\Omega$.\index{Hardy space} We keep the general notation for the sake of uniformity.

\medskip

We refer the reader to the appendix for an explanation of the notation concerning Lebesgue spaces.

\begin{deff}\label{def:5}
Take a positive Radon measure $\mi$ on $\Omega$, and define, for $p,q\in ]0,+\infty]$, 
\[
A^{p,q}_\mi(D)\coloneqq \Set{f\in \Hol(D)\colon \int_\Omega \norm{f_h}_{L^p(\Nc)}^q\,\dd \mi(h)<\infty }
\]
endowed with the topology associated with the distance 
\[
(f,g)\mapsto\norm{f-g}_{A^{p,q}_\mi(D)}^{\min(1,p,q)}\coloneqq \left( \int_\Omega \norm{f_h-g_h}_{L^p(\Nc)}^q\,\dd \mi(h)\right) ^{\min(1,p,q)/q}
\]
(modifications if $q=\infty$). We define also 
\[
A^{p,q}_{\mi,\loc}(D) \coloneqq \Set{f\in \Hol(D)\colon [h\mapsto \norm{f_h}_{L^p(\Nc)}] \in L^q_\loc(\mi) },
\]
endowed with the corresponding (metrizable) topology.
We define $A^{p,q}_{\mi,0}(D)$ and $A^{p,q}_{\mi,0,\loc}(D)$ as the spaces of $f\in \Hol(D)$ such that the mapping $h\mapsto f_h$ belongs to $L^q_0(\mi;L^p_0(\Nc))$ and  $L^q_{0,\loc}(\mi;L^p_0(\Nc))$, respectively, endowed with the corresponding topology (cf.~Definition~\ref{41}).
\end{deff}

Here and throughout the paper, we denote by $L^\ell_0$ the closure of $C_c$ in $L^\ell$, so that $L^\ell_0=L^\ell$ if $\ell\in ]0,\infty[$, and $L^\ell_0=C_0$ if $\ell=\infty$.

\begin{oss}
The measure $\mi$ can be replaced by any positive Radon measure with the same support without altering the spaces $A^{p,\infty}_\mi(D)$ and $A^{p,\infty}_{\mi,\loc}(D)$. 

Indeed, the mapping 
\[
\varphi\colon h\mapsto\norm{f_h}_{L^p(\Nc)}=\sup\limits_K \norm{\chi_K f_h}_{L^p(\Nc)},
\]
where $K$ runs through the set of compact subsets of $\Nc$, is lower semi-continuous on $\Omega$. Hence, the set $\Set{h\in \Omega\colon \varphi(h)>t}$ is open in $\Omega$ for every $t\in \R$, so that $\norm{\chi_C \varphi}_{L^\infty(\mi)}= \sup\limits_{\Supp{\chi_C\cdot\mi}} \varphi$ for every closed subset $C$ of $\Omega$. In addition, 
\[
\overline{\Supp{\mi}\cap U}\subseteq\Supp{\chi_C\cdot \mi}\subseteq \Supp{\mi}\cap C
\]
if $C$ is the 
closure of an open subset $U$ of $\Omega$.
\end{oss}

The following lemma is taken from~\cite[Lemma 8.1]{OgdenVagi}. 
The functions $g^{(\eps)}$ are very useful when arguing by approximation.

\begin{lem}\label{lem:34}
	Take $\alpha\in \left]0,\frac 1 2\right[$. There is a family $(g^{(\eps)})_{\eps>0}$ of holomorphic functions on $D$ such that the following hold:
	\begin{enumerate}
		\item[\em(1)] there is a constant $C>0$ such that $\abs{g^{(\eps)}_h(\zeta,x)}\meg \ee^{-\eps C\left(\abs{\zeta}^{2\alpha}+\abs{x}^\alpha+\abs{h}^\alpha\right)}$ for every $\eps>0$, for every $h\in \Omega$, and for every $(\zeta,x)\in \Nc$;
		
		\item[\em(2)] for every $p\in ]0,\infty]$ and for every $\eps>0$ there is a constant $C_{\eps,p}>0$ such that $\norm{g^{(\eps)}_h}_{L^p(\Nc)}\meg C_{\eps,p} \ee^{-\eps C \abs{h}^\alpha  } $ for every $h\in \Omega$;
		
		\item[\em(3)] $g^{(\eps)}\to 1$ as $\eps\to 0^+$, locally uniformly on $D$. 
	\end{enumerate}
\end{lem}

In particular, $A^{p,q}_{\mi,0,\loc}(D)\neq 0$ for every $p,q\in ]0,\infty]$.

\begin{proof}
	Take a basis $\lambda_1,\dots, \lambda_m$ of $F'$ contained in $\Omega'$. Observe that, if $(\zeta,z)\in D$, then 
	\[
	\Im \langle (\lambda_j)_\C,z\rangle=\langle \lambda_j, \Im z-\Phi(\zeta)+ \Phi(\zeta)\rangle \Meg \langle \lambda_j, \Im z-\Phi(\zeta)\rangle>0
	\]
	since $\Phi(\zeta)\in \overline \Omega$. Therefore, for every $\eps>0$ we may define
	\[
	g^{(\eps)}\colon D\ni (\zeta,z)\mapsto \exp\Big( -\eps\sum_{j=1}^m \langle (\lambda_j)_\C,z\rangle^\alpha \Big) \in \C,
	\] 
	where the complex power to the exponent $\alpha$ is the unique holomorphic function on $\C \setminus (-i \R_+)$ which induces the mapping $x\mapsto x^\alpha$ on $\R_+^*$. 
	Thus, the above remarks imply that
	\[
	\abs{g^{(\eps)}(\zeta,z)}\meg \exp\Big( - \eps \cos(\alpha \pi)\sum_{j=1}^m \abs{\langle (\lambda_j)_\C, z\rangle }^{\alpha} \Big) 
	\]
	for every $(\zeta,z)\in D$, where $\cos(\alpha \pi)>0$. 
	Therefore, if $h\in \Omega$,
	\[
	\abs{g^{(\eps)}_h(\zeta,x)}\meg \exp\Big( -\eps \frac{\cos(\alpha \pi)}{ 2^{\alpha/2}} 3^{\alpha-1} \sum_{j=1}^m ( \abs{\langle \lambda_j, x\rangle}^\alpha + \langle \lambda_j, h\rangle^\alpha+\langle \lambda_j,\Phi(\zeta)\rangle ^\alpha  )  \Big) 
	\]
	for every $(\zeta,x)\in \Nc$, since
	\[
	\begin{split}
	\abs{\langle( \lambda_j)_\C, x+i \Phi(\zeta)+i h\rangle }^{\alpha}&\Meg 2^{-\alpha/2}( \abs{\langle \lambda_j, x\rangle}+\langle \lambda_j, h\rangle+\langle \lambda_j,\Phi(\zeta)\rangle   )^\alpha\\
	&\Meg 2^{-\alpha/2} 3^{\alpha-1} \left(  \abs{\langle \lambda_j, x\rangle}^\alpha + \langle \lambda_j, h\rangle^\alpha+\langle \lambda_j,\Phi(\zeta)\rangle ^\alpha \right).
	\end{split}
	\]
	Set $C'\coloneqq  \cos(\alpha \pi) 2^{-\alpha/2} 3^{\alpha-1}$. Then, for every $p\in ]0,\infty]$ there is a constant $C'_{\eps,p}>0$ such that 
	\[
	\norm{g^{(\eps)}_h}_{L^p(\Nc)}\meg C'_{\eps,p} \exp\Big(-\eps C'\sum_{j=1}^m \langle \lambda_j, h\rangle^\alpha\Big)
	\]
	for every $h\in \Omega$. The assertion follows easily.
\end{proof}

\begin{lem}\label{lem:67}
	Take $\alpha\in \left]0,\frac 1 2\right[$, $\left(g^{(\eps)}\right)_{\eps>0}$ as in Lemma~\ref{lem:34}, $p,q\in ]0,\infty]$ and $f\in A^{p,q}_{\mi,0,\loc}(D)$. Then, $g^{(\eps)} f$ converges to $f$ in $A^{p,q}_{\mi,0,\loc}(D)$ for $\eps\to 0^+$.
\end{lem}

\begin{proof}
	Assume first that $q<\infty$. If $p<\infty$, then the dominated convergence theorem shows that 
	\[
	\lim_{\eps\to 0^+}\norm{g^{(\eps)}_h f_h-f_h}_{L^p(\Nc)}=0
	\]
	for every $h\in \Omega$. 
	If, otherwise, $p=\infty$, then $f_h\in C_0(\Nc)$, so that $\norm{g^{(\eps)}_h f_h-f_h}_{L^\infty(\Nc)}\to 0 $ for $\eps\to 0^+$.
	Thus, the dominated convergence theorem shows that $g^{(\eps)} f$ converges to $f$ in $A^{p,q}_{\mi,0,\loc}(D)$ for $\eps\to 0^+$.
	
	Now, assume that $q=\infty$. Since the mapping 
	\[
	\Supp{\mi}\ni h\mapsto f_h\in L^p(\Nc)
	\]
	is then continuous, the $f_h$ stay in a compact subset of $L^p(\Nc)$ as $h$ stays in a compact subset of $\Supp{\mi}$. Now, the endomorphisms  $\varphi\mapsto g^{(\eps)}\varphi$ of $L^p(\Nc)$, as $\eps$ runs through $\R_+^*$, are equicontinuous, and converge pointwise to the identity for $\eps\to 0^+$. The assertion then follows from~\cite[Theorem 1 of Chapter X, \S 2, No.\ 4]{BourbakiGT2}.
\end{proof}

The following result, which is a simple consequence of the subharmonicity of $\abs{f}^p$, is of fundamental importance in the study of weighted Bergman spaces. 

\begin{lem}\label{lem:76}
	Take $f\in \Hol(D)$ and $p\in ]0,\infty[$. Then,
	\[
	\abs{f(\zeta,z)}^p\meg \dashint_{B_{E\times F_\C}((\zeta,z),\rho)} \abs{f(\zeta',z')}^p\,\dd (\zeta',z')
	\]
	for every $(\zeta,z)\in D$ and for every $\rho>0$ such that $\overline B_{E\times F_\C}((\zeta,z),\rho)\subseteq D$.
\end{lem}

\begin{proof}
	This follows from~\cite[Theorem 2.1.4 and Corollary 2.1.15]{Krantz}.
\end{proof}

The following result is a sort of `local' version of the inclusions between the various weighted Bergman spaces we shall prove later (cf.~Propositions~\ref{prop:5} and~\ref{prop:22}). It still has some relevant consequences.

Here and throughout the paper, if $\nu$ is a Radon measure and $f\in L^1_\loc(\nu)$, we shall denote by $f\cdot \nu$ the measure with density $f$ with respect to the measure $\nu$, so that $(f\cdot \nu)(B)=\int_B f\,\dd\nu$ for every Borel set $B$. If $\Ff$ is a subset of $L^1_\loc(\nu)$, we shall define $\Ff\cdot \nu\coloneqq \Set{f\cdot \nu\colon f \in \Ff}$.
Finally, we shall define $p'\coloneqq \max(1,p)'$ for every $p\in ]0,\infty]$, so that $p'=\infty$ if $p\in ]0,1]$ and $\frac{1}{p}+\frac{1}{p'}=1$ if $p\in [1,\infty]$.

\begin{prop}\label{prop:1}
	Take $p_1,p_2,q_1,q_2\in ]0,\infty]$ such that $p_1\meg p_2$ and $q_1\meg q_2$, and assume that $\chi_\Omega\cdot \Hc^m\in L^{q_1'}_\loc(\mi)\cdot \mi$. Then there are continuous inclusions 
	\[
	A^{p_1,q_1}_{\mi,\loc}(D)\subseteq A^{p_2,q_2}_{\mi,\loc}(D)\qquad \text{and} \qquad A^{p_1,q_1}_{\mi,0,\loc}(D)\subseteq A^{p_2,q_2}_{\mi,0,\loc}(D).
	\]
\end{prop}

The proof is based on~\cite[Proposition 2.2]{RicciTaibleson}.

\begin{proof}
	Denote by $\varrho $ the density of $\chi_\Omega \cdot \Hc^m$ with respect to $\mi$, and define $\ell\coloneqq \min (1,p_1,q_1)$ to simplify the notation.
	Define 
	\[
	\varphi\colon E\times F\times \Omega\ni (\zeta,x,h)\mapsto (\zeta, x+i\Phi(\zeta)+i h)\in D,
	\]
	and observe that $\varphi$ is a bijection of $E\times F\times \Omega$ onto $D$.
	Observe that there are $R_0>0$ and $C'>0$ such that, for every $R\in ]0,R_0]$ and for every $h\in \Omega$,
	\[
	B_{E\times F_\C}((0,i h),R/C')\subseteq \varphi(B_\Nc(0,R)\times B_F(h,R))\subseteq B_{E\times F_\C}((0,i h),C' R).
	\]
	Therefore, Lemma~\ref{lem:76} implies that there is a constant $C_{R_0}>0$ such that
	\[
	\abs{f(0,i h)}^{\ell}\meg C_{R_0} \dashint_{B_F(h,R)}\dashint_{B_\Nc(0,R)} \abs{f_{h'}(\zeta',x')}^{\ell}\,\dd (\zeta',x')\,\dd h'
	\]
	for every $h\in \Omega$ and for every $R\in ]0,R_0]$ such that $\overline B_{E\times F_\C}((0,i h),C'R)\subseteq D$.
	Hence, 
	\[
	\abs{f(\zeta,z)}^{\ell}\meg C_{R_0} \dashint_{B_F(h,R)}\dashint_{B_\Nc(0,R)} \abs{f_{h'}((\zeta,x)(\zeta',x'))}^{\ell}\,\dd (\zeta',x')\,\dd h'
	\]
	for every $(\zeta,z)\in D$ and for every $R\in ]0,R_0]$ such that $\overline B_{E\times F_\C}((0,i h),C' R)\subseteq D$, where $h\coloneqq \Im z-\Phi(\zeta)$.
	Thus, applying Minkowsky's integral inequality, the right invariance of the measure on $\Nc$ and Jensen's inequality,
	\[
	\begin{split}
	\norm{f_h}_{L^{p_1}(\Nc)}^{\min(1,q_1)}&=\norm*{\abs{f_h}^{\ell}  }_{L^{p_1/\ell}(\Nc)}^{ \min(1,q_1)/\ell}\\
		&\meg \left(  C_{R_0}\dashint_{B_F(h,R)}\norm{f_{h'}}_{L^{p_1}(\Nc)}^{\ell}\,\dd h'\right)^{ \min(1,q_1)/\ell} \\
		&\meg C_{R_0}^{ \min(1,q_1)/\ell} \dashint_{B_F(h,R)}\norm{f_{h'}}_{L^{p_1}(\Nc)}^{\min(1,q_1)}\,\dd h'.
	\end{split}
	\]
	In addition, the left invariance of the measure on $\Nc$, H\"older's inequality, and Jensen's inequality show that
	\[
	\begin{split}
	\norm{f_h}_{L^{\infty}(\Nc)}^{\min(1,q_1)}&=\norm*{\abs{f_h}^{\ell}  }_{L^{\infty}(\Nc)}^{ \min(1,q_1)/\ell }\\
	&\meg \left(\frac{C_{R_0}}{C_R'^{\ell/p_1 } }\dashint_{B_F(h,R)}\norm{f_{h'}}_{L^{p_1}(\Nc)}^{\ell}\,\dd h'\right)^{ \frac{\min(1,q_1)}{\ell} } \\
	&\meg C_{R_0}^{ \min(1,q_1)/\ell}\frac{1}{C_R'^{\frac{\min(1,q_1)}{p_1}} }\dashint_{B_F(h,R)}\norm{f_{h'}}_{L^{p_1}(\Nc)}^{\min(1,q_1)}\,\dd h',
	\end{split}
	\]
	where $C_R'\coloneqq \Hc^{2 n+m}(B_\Nc(0,R))$.
	Therefore, there is a constant $C_{p_1,q_1,p_2,R}>0$ such that
	\[
	\begin{split}
	\norm{f_h}_{L^{p_2}(\Nc)}^{\min(1,q_1)}&\meg C_{p_1,q_1,p_2,R}\dashint_{B_F(h,R)}\norm{f_{h'}}_{L^{p_1}(\Nc)}^{\min(1,q_1)}\,\dd h'\\
		&=\frac{C_{p_1,q_1,p_2,R}}{\Hc^m(B_F(0,R))}\int_{B_F(h,R)}\norm{f_{h'}}_{L^{p_1}(\Nc)}^{\min(1,q_1)} \varrho (h')\,\dd \mi(h').
	\end{split}
	\]
	Now, let $K$ be a compact subset of $\Omega$, and take $R\in ]0,R_0]$ in such a way that $\overline B_{E\times F_\C}((0,i h), C' R)\subseteq D$ for every $h\in K$. Then, setting 
	\[
	B_F(K,R)\coloneqq \bigcup_{h\in K} B_F(h,R)
	\]
	and 
	\[
	C_{p_1,q_1,p_2,q_2,R}\coloneqq \mi(K)^{1/q_2} \left( \frac{C_{p_1,q_1,p_2,R} \norm{\chi_{B_F(K,R)} \varrho  }_{L^{q_1'}(\mi)}  }{\Hc^m(B_F(0,R))} \right) ^{1/\min(1,q_1)},
	\]	
	by means of H\"older's inequality we see that
	\[
	\begin{split}
	\left( \int_K  \norm{f_h}_{L^{p_2}(\Nc)}^{q_2}\,\dd \mi(h)\right) ^{1/q_2}&\meg \mi(K)^{1/q_2} \sup\limits_{h\in K} \norm{f_h}_{L^{p_2}(\Nc)}\\
		&\meg C_{p_1,q_1,p_2,q_2,R}\left( \int_{B_F(K,R)} \norm{f}_{L^{p_1}(\Nc)}^{q_1}\,\dd \mi(h)\right)^{1/q_1},
	\end{split} 
	\]
	so that the first assertion follows.
	
	For what concerns the second assertion, take $f\in A^{p_1,q_1}_{\mi,0,\loc}(D)$ and $(g^{(\eps)})_{\eps>0}$ as in Lemma~\ref{lem:34}. Then, Lemma~\ref{lem:67} implies that $g^{(\eps)}f$ converges to $f$ in $A^{p_1,q_1}_{\mi,\loc}(D)$ as $\eps\to 0^+$, hence also in $A^{p_2,q_2}_{\mi,\loc}(D)$. 
	It will then suffice to prove that $g^{(\eps)}f\in A^{p_2,q_2}_{\mi,0,\loc}(D)$ for every $\eps>0$.
	Let us first show that $g^{(\eps)}_h f_h\in L^{p_2}_0(\Nc)$ for every $h\in \Omega$. The assertion is clear if $p_2<\infty$. 
	Otherwise, the assertion follows from the fact that $ f_h$ is continuous and bounded.
	If $q_2<\infty$, this is sufficient to conclude. 
	Otherwise, we have to show that the mapping 
	\[
	\Supp{\mi}\ni h\mapsto g^{(\eps)}_h f_h\in L^{p_2}(\Nc)
	\]
	is continuous. If $p_2<\infty$, it suffices to observe that $f$ is continuous on $D$ and that the functions $g^{(\eps)}_h f_h$, as $h$ runs through a compact subset $K$ of $D$, are uniformly bounded in absolute value by an element of $L^{p_2}(\Nc)$ (e.g., $\abs{g^{(\eps)}}\sup\limits_{h\in K} \norm{f_h}_{L^\infty(\Nc)}  $). When $p_2=\infty$, we observe that the $f_h$, as $h$ runs through a compact subset $K$ of $D$, are uniformly bounded and equicontinuous. Hence, the $ g^{(\eps)}_h f_h$ , as $h$ runs through $K$, are equicontinuous on $\Nc$ and at the point at infinity of $\Nc$. Thus, the mapping 
	\[
	\Omega\ni h\mapsto g^{(\eps)}_h f_h\in C_0(\Nc)
	\]
	is continuous by~\cite[Theorem 1 of Chapter X, \S 2, No.\ 4]{BourbakiGT2}.
\end{proof}

\begin{cor}\label{cor:6}
	Take $p,q_1,q_2\in ]0,\infty]$ and assume that $\chi_\Omega\cdot \Hc^m\in L^{q'_1}_\loc(\mi)\cdot \mi$. If $q_2\Meg q_1$, then $A^{p,q_1}_{\mi,\loc}(D)=A^{p,q_2}_{\mi,\loc}(D)$ (as topological vector spaces).
\end{cor}

In particular, if the measure $\mi$ is relatively well-behaved, then the index $q$ in the definition $A^{p,q}_{\mi,\loc}(D)$ is essentially irrelevant. 
Nonetheless, if we had simply considered the space $A^{p,\infty}_{\mi,\loc}(D)$, then the inclusion $A^{p,q}_{\mi}(D)\subseteq A^{p,q}_{\mi,\loc}(D)=A^{p,\infty}_{\mi, \loc}(D)$  would not have been clear.

\begin{cor}\label{cor:19}
	Take $p_1,p_2\in ]0,\infty]$ such that $p_1\meg p_2$.
	Assume, further, that $\Supp{\mi}=\Omega$. Then there are continuous inclusions $A^{p_1,\infty}_{\mi,\loc}(D)\subseteq A^{p_2,\infty}_{\mi,\loc}(D)$ and $A^{p_1,\infty}_{\mi,0,\loc}(D)\subseteq A^{p_2,\infty}_{\mi,0,\loc}(D)$ .
\end{cor}

\begin{proof}
	Indeed, we may assume that $\mi=\chi_\Omega \cdot \Hc^m$.
\end{proof}

\begin{cor}\label{cor:27}
	Take $p\in ]0,\infty]$, assume that $\Supp \mi=\Omega$, and take $f\in A^{p,\infty}_{\mi,\loc}(D)$. Then, $f\in A^{p,\infty}_{\mi,0,\loc}(D)$ if and only if $f_h\in L^p_0(\Nc)$ for every $h\in \Omega$.
\end{cor}

This result shows that the continuity of the mapping $h\mapsto f_h\in L^p(\Nc)$ is ensured if we know that $f_h\in L^p_0(\Nc)$ for every $h\in \Omega$.

\begin{proof}
	It will suffice to prove that, if $f\in \Hol(D)$ and $f_h\in L^p_0(\Nc)$ for every $h\in \Omega$, then $f\in A^{p,\infty}_{\mi,0,\loc}(D)$, that is, the mapping $\Omega\ni h\mapsto f_h\in L^p(\Nc)$ is continuous. Notice that we may assume that $\mi=\chi_\Omega \cdot \Hc^m$.
	Define $R_0$, $C'$, and $C_{R_0}$ as in the proof of Proposition~\ref{prop:1}, so that
	\[
	\abs{f_h(\zeta,x)}^{\min(1,p)} \meg C_{R_0} \dashint_{B_F(h,R)}\dashint_{B_\Nc(0,R)} \abs{f_{h'}((\zeta,x)(\zeta',x'))}^{\min(1,p)}\,\dd (\zeta',x')\,\dd h'
	\]
	for every $h\in \Omega$, for every $(\zeta,x)\in \Nc$, and for every $R\in ]0,R_0]$ such that 
	\[
	\overline B_{E\times F_\C}((0, i h), C' R)\subseteq D.
	\]
	Now, fix $h_0\in \Omega$ and $R_1\in ]0,R_0/2]$ in such a way that  $\overline B_{E\times F_\C}((0, i h_0), 2 C' R_1)\subseteq D$, and observe that there is a constant $C''>0$ such that
	\[
	\abs{f_h(\zeta,x)}^{\min(1,p)} \meg C'' \dashint_{B_F(h_0,2 R_1)}\dashint_{B_\Nc(0, R_1)} \abs{f_{h'}((\zeta,x)(\zeta',x'))}^{\min(1,p)}\,\dd (\zeta',x')\,\dd h'
	\]
	for every $h\in B_F(h_0, R_1)$ and for every $(\zeta,x)\in \Nc$. Define 
	\[
	\widetilde f_{h_0}\colon \Nc \ni (\zeta,x)\mapsto  \dashint_{B_F(h_0,2 R_1)}\dashint_{B_\Nc(0, R_1)} \abs{f_{h'}((\zeta,x)(\zeta',x'))}^{\min(1,p)}\,\dd (\zeta',x')\,\dd h',
	\]
	and let us prove that $\widetilde f_{h_0}^{1/\min(1,p)}\in L^p_0(\Nc)$. If $p<\infty$, it suffices to argue as in the proof of Proposition~\ref{prop:1}. 
	Then, assume that $p=\infty$, and let us prove that the mapping 
	\[
	\Omega\times \Nc\ni (h',(\zeta',x'))\mapsto R_{(\zeta',x')}f_{h'}\in C_0(\Nc)
	\]
	is measurable, where $R_{(\zeta',x')}$ denotes the right translation by $(\zeta',x')$. Since $C_0(\Nc)$ is a separable Banach space and since $(\delta_{(\zeta,x)})_{(\zeta,x)\in\Nc}$ is total (that is, generates a dense vector subspace) in the dual $\Mcal^1(\Nc)$ of $C_0(\Nc)$, by~\cite[Proposition 10 of Chapter IV, \S 5, No.\ 5]{BourbakiInt1} it will suffice to show that the mapping 
	\[
	\Omega\times \Nc\ni (h',(\zeta',x'))\mapsto \langle \delta_{(\zeta,x)},R_{(\zeta',x')}f_{h'}\rangle
	\] is measurable for every $(\zeta,x)\in \Nc$, and this is clear. 
	Therefore, the mapping 
	\[
	\Omega\times \Nc\ni (h',(\zeta',x'))\mapsto  R_{(\zeta',x')}\abs{f_{h'}}\in C_0(\Nc)
	\]
	is measurable. In addition, it is easily seen that the $R_{(\zeta',x')}\abs{f_{h'}}$, for $h'\in \overline B_F(h_0,2 R_1)$ and $(\zeta',x')\in \overline B_\Nc(0,R_1)$, stay in a bounded subset of $C_0(\Nc)$, so that $\widetilde f_{h_0}\in C_0(\Nc)$.
	
	 Then, the assertion follows from the dominated convergence theorem when $p<\infty$, and from the equicontinuity of the $f_h$, $h\in B_F(h_0,R_1)$, at the point at infinity of $\Nc$, when $p=\infty$.
\end{proof}

The following result is the extension of a well-known (and quite useful) fact concerning functions in the classical Hardy spaces $H^p(\C_+)$.  In the classical situation, it is closely related to the subharmonicity of $\abs{f}$ (for $p\Meg 1$) and to a suitable decomposition of the elements of $H^p(\C_+)$ (for $p<1$). The general case follows from the classical one by means of Fubini's theorem. 

\begin{prop}\label{prop:28}
	Take $p\in ]0,\infty]$ and assume that $\Supp{\mi}=\Omega$. Then, for every $f\in A^{p,\infty}_{\mi}(D)$, the mapping $\Omega \ni h\mapsto \norm{f_h}_{L^p(\Nc)}$ is decreasing.
\end{prop}

Recall that $F$ (hence $\Omega$) is endowed with the ordering `$x\meg y$ if and only if $y-x\in \overline \Omega$.'

\begin{proof}
	Take $h_0\in\Omega$ and a unit vector $h\in \overline\Omega$. It will suffice to prove that the mapping $\R_+^*\ni y\mapsto \norm{f_{h_0+y h}}_{L^p(\Nc)}$ is decreasing for every $f\in A^{p,\infty}_\mi(D)$.  
	Let $F_h$ be the orthogonal complement of $h$ in $F$, and fix $\zeta\in E$ and $x'\in F_h$.
	Define 
	\[
	f_{h_0,h,\zeta,x'}\colon \C_+\ni (x+i y)\mapsto f(\zeta, x'+i h_0+(x+i y)h)\in \C,
	\]
	where $\C_+= \R+i \R_+^*$. 
	Then,  $f_{h_0,h,\zeta,x'}$ is holomorphic on $\C_+$. Let us prove that $ f_{h_0,h,\zeta,x'}\in A^{p,\infty}_{\Hc^1}(\C_+)$.
	Indeed, defining $R_0$, $C'$, and $C_{R_0}$ as in the proof of Proposition~\ref{prop:1}, we see that
	\[
	\begin{split}
	&\abs{f_{h_0,h,\zeta,x'}(x+i y)}^{\min(1,p)}\\
	&\qquad\meg C_{R_0} \dashint_{B_F(h_0+y h, R)} \dashint_{B_\Nc(0,R)} \abs{f_{h'}((\zeta, x'+x h)(\zeta',x'')  )  }^{\min(1,p)}\,\dd (\zeta',x'')\,\dd h'
	\end{split}
	\]
	for every $x+i y\in \C_+$ and for every $R\in ]0, R_0]$ such that $\overline B_{E\times F_\C}((0,i  h_0), C' R)\subseteq D$. Therefore,
	\[
	\begin{split}
	\norm{f_{h_0,h,\zeta,x'}(\,\cdot\,+i y)}_{L^p(\R)}^{\min(1,p)}&\meg \frac{C_{R_0}C_R'^{1/\max(1,p)}}{C_R''^{1/\max(1,p)}} \dashint_{B_F(h_0+ y h,R)} \norm{f_{h'}}_{L^p(\Nc)}^{\min(1,p)}\,\dd h'\\
		&\meg\frac{C_{R_0}C_R'^{1/\max(1,p)}}{C_R''^{1/\max(1,p)}}\norm{f}_{A^{p,\infty}_\mi(D)}^{\min(1,p)}
	\end{split}
	\]
	for every $y>0$, where 
	\[
	C'_R\coloneqq \sup\limits_{ (\zeta',x'')\in E\times F_h} \Hc^1(B_\Nc(0,R)\cap [(\zeta',x'')+\R h])  
	\] 
	and 
	\[
	C''_R\coloneqq \Hc^{2 n+m}(B_\Nc(0,R)).
	\]
	Thus, $f_{h_0,h,\zeta,x'}\in A^{p,\infty}_{\Hc^1}(\C_+)$, so that
	\[
	\norm{f_{h_0,h,\zeta,x'}(\,\cdot\,+i y)}_{L^p(\R)}\meg \norm{f_{h_0,h,\zeta,x'}(\,\cdot\,+i y')}_{L^p(\R)}
	\]
	for every $y,y'>0$ such that $y>y'$, thanks to~\cite[Theorem 11.5]{Duren}.  The assertion follows easily.
\end{proof}

\begin{lem}\label{lem:9}
	Take $p,q\in ]0,\infty]$, and assume that $\Supp{\mi}=\Omega$. Then, $A^{p,q}_\mi(D)$ (resp.\ $A^{p,q}_{\mi,\loc}(D)$) is complete if and only if the inclusion $A^{p,q}_\mi(D)\subseteq \Hol(D)$ (resp.\ $A^{p,q}_{\mi,\loc}(D)\subseteq \Hol(D)$) is continuous.
\end{lem}

\begin{proof}
	One implication follows from the closed graph theorem (for complete metrizable topological vector spaces), the other one from the completeness of the mixed norm space $L^{p,q}(\Hc^{2 n+m},\mi)$ (cf.~Proposition~\ref{prop:62}).
\end{proof}

\begin{cor}\label{cor:2}
	Take $p,q\in ]0,\infty]$ and assume that $\chi_\Omega\cdot \Hc^m\in L^{q'}_\loc(\mi)\cdot \mi$. Then, $A^{p,q}_\mi(D)$ and $A^{p,q}_{\mi,\loc}(D)$ are complete. 
\end{cor}

\begin{proof}
	This follows from Proposition~\ref{prop:1} and Lemma~\ref{lem:9}, since $A^{\infty,\infty}_{\mi,\loc}(D)$ embeds continuously into $\Hol(D)$.
\end{proof}

The following density result is a `local' version of Proposition~\ref{prop:8}. It still has some useful consequences.

\begin{prop}\label{prop:37}
	Take $p_1,p_2,q\in ]0,\infty]$ such that $p_2\meg p_1$. Then, the space $A^{p_1,q}_{\mi,0,\loc}(D)\cap A^{p_2,q}_{\mi,0,\loc}(D)$ is dense in the space $A^{p_1,q}_{\mi,0,\loc}(D)$.
\end{prop}

The proof is based on~\cite[Lemma 8.1]{OgdenVagi}, which deals with the density of the space $A^{p_1,\infty}_{\Hc^m}(D)\cap A^{p_2,\infty}_{\Hc^m}(D)$ in the space $A^{p_1,\infty}_{\Hc^m}(D)$.

\begin{proof}
	Take $(g^{(\eps)})_{\eps>0}$ as in Lemma~\ref{lem:34}, and take $f\in A^{p_1,q}_{\mi,0,\loc}(D)$. Notice that we may assume that $p_2<p_1$. Set $s\coloneqq \frac{p_1 p_2}{p_1-p_2}$ if $p_1\neq \infty$, $s\coloneqq p_2$ if $p_1=\infty$. Then H\"older's inequality implies that
	\[
	\norm{(f g^{(\eps)})_h}_{L^{p_2}(\Nc)}\meg \norm{g^{(\eps)}_h}_{L^{s }(\Nc)}\norm{f_h}_{L^{p_1}(\Nc)}\meg C_{\eps,s} \norm{f_h}_{L^{p_1}(\Nc)} \ee^{-\eps C\abs{h}^\alpha},
	\]
	for every $h\in \Omega$. Therefore, $f g^{(\eps)}\in A^{p_1,q}_{\mi,0,\loc}(D)\cap A^{p_2,q}_{\mi,0,\loc}(D)$ for every $\eps>0$. The assertion then follows from Lemma~\ref{lem:67}.
\end{proof}

\section{Paley--Wiener Theorems}\label{sec:1:4}

In this section we shall study some global properties of the weighted Bergman spaces $A^{2,q}_\mi(D)$, making use of the Fourier transform on $\Nc$. In particular, we shall characterize the spaces $A^{2,q}_\mi(D)$ and $A^{2,q}_{\mi,\loc}(D)$ by means of the Fourier transform (thus proving analogues of the classical Paley--Wiener theorems) in Propositions~\ref{prop:2} and~\ref{prop:4}.

\begin{deff}
We denote by $\Lc \mi$ the Laplace transform of $\mi$, defined as
\[
(\Lc \mi)\CB(\lambda)= \int_\Omega \ee^{-\langle \lambda, h\rangle}\,\dd \mi(h)
\] 
for every $\lambda\in F'_\C$ such that $\Re \lambda\in  \Ds_\mi $, where\label{13}
\[
\Ds_\mi\coloneqq \Set{\lambda\in F'\colon (\Lc \mi)(\lambda)=\int_\Omega \ee^{-\langle \lambda, h\rangle}\,\dd \mi(h)<\infty}.
\]
\end{deff}

We now recall some basic facts on the Laplace transform. The proofs, which are elementary, are omitted.

\begin{lem}\label{lem:5}
	The set $ \Ds_\mi $ is a union of translates of $\overline{\Omega'}$. In addition, if $\Ds_\mi$ is not empty, then $\mi$ induces a Radon measure on $\overline \Omega$.
\end{lem}

\begin{lem}\label{lem:10}
	Assume that $\mi$ is non-zero and homogeneous of degree $d\in \R$, that is, that $(r\,\cdot\,)_*\mi=r^d \mi$ for every $r>0$. Then, the following conditions are equivalent:
	\begin{enumerate}
		\item[\em(1)] $\Ds_\mi \neq \emptyset$;
		
		\CR \item[\em(2)] $\mi$ induces a Radon measure on $\overline \Omega$;
		
		\item[\em(3)] there is an open neighbourhood $V$ of $0$ in $\overline \Omega$ such that $\mi(V\cap \Omega)<\infty$.\CB
	\end{enumerate} 
	
	Assume that $ \Ds_\mi\neq \emptyset$ and denote by $S$ be the interior of the convex hull of $\Supp{\mi}$. Then, $d<0$, $S$ is an open convex cone  which does not contain any affine line,  and $S'\subseteq  \Ds_\mi\subseteq \overline{S'}$. 
\end{lem}

Observe that, if $\chi_\Omega\cdot\Hc^m \in L^{q'}_\loc(\mi)\cdot \mi$, then $\Supp{\mi}=\Omega$.

The following result shows that the holomorphy of $f$ not only simplifies the Fourier transform of $f_h$ (cf.~Proposition~\ref{prop:3}), but also relates the Fourier transforms of $f_h$ and $f_{h'}$ for $h,h'\in \Omega$. 
If we were able to define the inverse Fourier transform of an arbitrary measurable field of operators on $\Lambda_+$, we could then define a general notion of boundary values of the elements of $A^{p,\infty}_{\mi,\loc}(D)$, at least for $p\in [1,2]$.

\begin{prop}\label{prop:2}
	Take $p\in [1,2]$, and assume that $\Supp{\mi}=\Omega$. Then, for every $f\in A^{p,\infty}_{\mi,\loc}(D)$ there is a measurable field of operators 
	\[
	\tau_f\colon F'\to \prod_{\lambda\in F'}\Lin(H_\lambda)
	\]
	such that 
	\[
	\pi_\lambda(f_h)=\ee^{-\langle \lambda,h\rangle}\tau_f(\lambda)=\chi_{\Lambda_+}(\lambda)\ee^{-\langle \lambda,h\rangle} \tau_f(\lambda) P_{\lambda,0}
	\] 
	for almost every $\lambda\not \in W$ and for every $h\in \Omega$. 
\end{prop}

The proof is based on~\cite[Theorem 19.2]{Rudin}, which deals with the case $D=\C_+$.

\begin{proof}
	Define $(g^{(\eps)})_{\eps>0}$ is as in Lemma~\ref{lem:34}.
	Then, up to replacing $f$ with $f g^{(\eps)}$ and then passing to the limit, we may assume that $p=1$.
	Take two distinct elements $h,h'$ of $ \Omega$. Let $\widetilde F$ be the orthogonal complement of $(h'-h)\R$ in $F$, so that $\widetilde F_\C$ is an algebraic complement of $(h'-h)\C$ in $F_\C$. Define $\widetilde B_F(0, \rho  )\coloneqq B_{\widetilde F}(0,\rho)+[-\rho,\rho](h'-h)$ for every $\rho>0$. 
	Then, by Cauchy's integral theorem,
	\[
	\begin{split}
	&\pi_\lambda(\chi_{B_E(0,\rho)\times \widetilde B_F(0,\rho)} f_h  )-\ee^{\langle\lambda,h'-h\rangle}\pi_\lambda(\chi_{B_E(0,\rho)\times \widetilde B_F(0,\rho)} f_{h'}  )\\
		&=-\sum_{\eps=\pm 1} \eps i  \abs{h-h'}\int_{B_E(0,\rho)\times B_{\widetilde F}(0,\rho)}  \int_0^1  f(\zeta,x+i\Phi(\zeta)+i h+ (\eps \rho+i t) (h'-h))\\
		&\qquad \times\ee^{-i (\eps \rho+t i)\langle \lambda,h'-h \rangle}\,\dd t\,\pi_\lambda(\zeta,x')  \,\dd (\zeta,x)
	\end{split}
	\]
	for every $\rho>0$ and for every $\lambda \not \in W$.
	Now, the right-hand side of the preceding formula is bounded in $\Lin(H_\lambda)$-norm by
	\[
	g(\rho)\coloneqq\sum_{\eps=\pm1} \abs{h'-h} \ee^{\langle \lambda,h'-h\rangle_+}\int_0^1 \norm{f_{h+t(h'-h)} (\,\cdot\,+\eps \rho (h'-h))}_{L^1(E\times \widetilde F)}\,\dd t.
	\]
	In addition, 
	\[
	\int_\R g(\rho)\,\dd\rho\meg 2 \ee^{\langle \lambda,h'-h\rangle_+}\sup\limits_{h''\in [h,h']}\norm{f_{h''}}_{L^1(\Nc)}<\infty,
	\]
	so that there is a sequence $\rho_k\to \infty$ such that $g(\rho_k)\to 0$ for $k\to \infty$. Hence, passing to the limit for $k\to \infty$, we see that $ \pi_\lambda( f_h  )=\ee^{\langle\lambda,h'-h\rangle}\pi_\lambda( f_{h'}  )$. 
	Therefore, it suffices to define $\tau_f(\lambda)\coloneqq \ee^{\langle\lambda,h\rangle}\pi_\lambda( f_h  )$ for some (hence every)  $h\in \Omega$ and for every $\lambda \not \in W$. The second equality follows from Proposition~\ref{prop:3}.
\end{proof}

In the following corollaries, we draw some consequences of Proposition~\ref{prop:2} in the cases of Hardy ($A^{2,\infty}_\mi(D)$) and weighted Bergman spaces ($A^{2,2}_\mi(D)$). The case of Hardy spaces has already been completely solved by at least three different methods in~\cite{Gindikin,KoranyiStein,OgdenVagi}. Our approach is basically that of~\cite{OgdenVagi}.

\begin{cor}\label{cor:1}
	Keep the hypotheses and the notation of Proposition~\ref{prop:2}, and assume that $p=2$ and that $f\in A^{2,\infty}_\mi(D)$. Then, $\tau_f(\lambda)=0$ for almost every $\lambda \not \in \Omega'$ and
	\[
	\norm{f}_{A^{2,\infty}_\mi(D)}^2= \frac{2^{n-m}}{\pi^{n+m}}\int_{\Omega'} \norm{\tau_f(\lambda)}_{\Lin^2(H_\lambda)} ^2\abs{\Pfaff(\lambda)}\,\dd \lambda .
	\] 
	In particular,  $A^{2,\infty}_\mi(D)$ is a hilbertian space. 
	In addition, there is a unique $f_0\in L^2(\Nc)$ such that $\tau_f(\lambda)= \pi_\lambda(f_0)$ for almost every $\lambda\in F'$, $\norm{f_0}_{L^2(\Nc)}=\norm{f}_{A^{2,\infty}_\mi(D)}$, and $f_h\to f_0$ in $L^2(\Nc)$ as $h\to 0$.
\end{cor}

To the best of our knowledge, the case of the weighted Bergman spaces $A^{2,2}_{\mi}(D)$ considered below has not been considered in the literature in this generality. Cf.~\cite{Gindikin} for the case $\mi=\Hc^m$, treated with a different method.

\begin{cor}\label{cor:3}
	Keep the hypotheses and the notation of Proposition~\ref{prop:2}, and assume that $p=2$ and that $f\in A^{2,2}_\mi(D)$. Then, $\tau_f(\lambda)=0$ for almost every $\lambda \not \in \frac{1}{2} \Ds_\mi \cap \Lambda_+$, and
	\[
	\norm{f}_{A^{2,2}_\mi(D)}^2= \frac{2^{n-m}}{\pi^{n+m}}\int_{\frac{1}{2} \Ds_\mi \cap \Lambda_+} \norm{\tau_f(\lambda)}_{\Lin^2(H_\lambda)} ^2 \Lc \mi(2\lambda)\abs{\Pfaff(\lambda)}\,\dd \lambda .
	\] 
\end{cor}

The following result, together with Proposition~\ref{prop:2}, completes the characterization of $A^{2,q}_\mi(D)$ and $A^{2,q}_{\mi,\loc}(D)$ by means of the Fourier transform.

\begin{prop}\label{prop:4}
	Take $q\in ]0,\infty]$ (resp.\ $q\in ]0,\infty[$), and assume that $\Supp{\mi}=\Omega$ and that $A^{2,q}_\mi(D)$ (resp.\ $A^{2,q}_{\mi,\loc}(D)$) is complete.
	Let $\tau\colon \Lambda_+\to \prod_{\lambda\in F'}\Lin(H_\lambda)$ be a measurable field of operators such that the following hold:
	\begin{enumerate}
		\item[\em(1)] $\tau(\lambda)=\tau(\lambda)P_{\lambda,0}$ for almost every $\lambda\in \Lambda_+$;
		
		\item[\em(2)] the mapping $h\mapsto \left( \int_{\Lambda_+} \norm{\tau(\lambda)}_{\Lin^2(H_\lambda)} ^2 \ee^{-2\langle \lambda,h \rangle} \abs{\Pfaff( \lambda)}\,\dd \lambda \right) ^{1/2}$ belongs to $L^q(\mi)$ (resp.\ $L^q_\loc(\mi)$).
	\end{enumerate}
	Then, there is a unique $f\in A^{2,q}_{\mi}(D)$ (resp.\  $f\in A^{2,q}_{\mi,\loc}(D)$) such that 
	\[
	\pi_\lambda(f_h)=\ee^{-\langle \lambda,h\rangle}\tau(\lambda)
	\]
	for almost every $\lambda\in \Lambda_+$ and for every $h\in \Omega$. 
	
	If, in addition, $\int_{\Lambda_+} \norm{\tau(\lambda)}_{\Lin^2(H_\lambda)} \ee^{-\langle \lambda,h \rangle} \abs{\Pfaff( \lambda)}\,\dd \lambda <\infty$ for every $h\in \Omega$, then
	\[
	f(\zeta,z)=\frac{2^{n-m}}{\pi^{n+m}} \int_{\Lambda_+} \tr(\tau(\lambda) \pi_\lambda(\zeta,\Re z)^*) \ee^{-\langle \lambda, \Im z-\Phi(\zeta)\rangle}\abs{\Pfaff(\lambda)}\,\dd \lambda
	\]
	for every $(\zeta,z)\in D$.
\end{prop}

In order to prove Proposition~\ref{prop:4}, we need two simple lemmas.

\begin{lem}\label{lem:8}
	Take $q\in ]0,\infty]$ and assume that $\Supp{\mi}=\Omega$. Let $\tau\colon  \Omega'\to \prod_{\lambda\in \Omega'}\Lin(H_\lambda)$ be a measurable field of operators such that the mapping 
	\[
	h\mapsto\left( \int_{\Omega'} \norm{\tau(\lambda)}_{\Lin^2(H_\lambda)} ^2 \ee^{-2\langle\lambda,h\rangle}\abs{\Pfaff(\lambda)}\,\dd \lambda\right) ^{1/2}
	\]
	belongs to $L^q_\loc( \mi)$.
	Then, the mapping 
	\[
	g_p\colon h\mapsto\int_{\Omega'} \norm{\tau(\lambda)}_{\Lin^2(H_\lambda)} ^p \ee^{-p\langle \lambda,h \rangle} \abs{\Pfaff( \lambda)}\,\dd \lambda  
	\]
	is  finite, continuous, and decreasing on $\Omega$ for every $p\in [1,2]$.
\end{lem}

\begin{proof}
	Observe that there is a $\mi$-negligible subset $N$ of $\Omega$ such that $ g_2(h)$ is finite for every $h\in \Omega \setminus N$, thanks to~\cite[Proposition 14 of Chapter IV, \S 5, No.\ 9]{BourbakiInt1}. 
	Then, take $h\in \Omega$, and observe that there is $h'\in \Omega\cap(h-\Omega )\setminus N$ since the support of $\mi$ is the whole of $\Omega$. 
	It then follows that the mapping $\Omega'\ni\lambda \mapsto \ee^{-\langle\lambda, h-h'\rangle }$ stays in $L^1(\nu)\cap L^\infty(\nu)$, where $\nu=\abs{\Pfaff}\cdot \Hc^m$. 
	Hence, $g_p(h)$ is finite for every $p\in [1,2]$. In addition, $g_p(h')\meg g_p(h)$ for every $h,h'\in \Omega$ such that $h'-h\in \Omega$, since $\ee^{-\langle \lambda,h'\rangle} \meg \ee^{-\langle \lambda, h\rangle}$ for every $\lambda\in \Omega'$. Hence,  $g_p$ is decreasing. Continuity follows from the dominated convergence theorem.
\end{proof}

\begin{lem}\label{lem:7}
	Let $U$ be an open subset of a hilbertian space $X$ over $\C$ of finite dimension $k$,  let $\nu$ be a Radon measure on a Hausdorff space $Y$, and let $Z$ be a Banach space over $\C$. Take $p\in [1,\infty]$ and take $f\in \Hol(U;L^p(\nu;Z))$. Then, there is a $(\Hc^{2k}\otimes \nu)$-measurable mapping $g\colon U\times Y\to Z$ such that $g(\,\cdot\,,y)$ is holomorphic for every $y\in Y$ and $g(x,\,\cdot\,)$ is a representative of $f(x)$ for every $x\in U$.
\end{lem}

\begin{proof}
	Identify $X$ with $\C^k$.
	Take $x_0\in X$ and take $r_{x_0}>0$ so that the Taylor series of $f$ at $x_0$ converges absolutely to $f$ on $\overline B(x_0,r_{x_0})$. For every $\alpha\in \N^k$, choose a representative $f_{x_0,\alpha}$ of $\frac{1}{\alpha!}\partial^\alpha f(x_0)$, and let $N_{x_0}$ be the $\nu$-negligible set of $y\in Y$ such that $\sum_{\alpha\in \N^k} \norm{f_{x_0,\alpha}(y)} r_{x_0}^{\abs{\alpha}} =\infty$. Then, define	 $g_{x_0}\colon B(x_0,r_{x_0})\times Y\to Z$ so that 
	\[
	g_{x_0}(x,y)=\sum_{\alpha\in \N^k} f_{x_0,\alpha}(y) (x-x_0)^\alpha
	\]
	for every $(x,y)\in B(x_0,r_{x_0})\times (Y\setminus N_{x_0})$, while $g_{x_0}(x,y)=0$ for every $(x,y)\in B(x_0,r_{x_0})\times N_{x_0}$. 
	Then, $g_{x_0}$ is $(\Hc^{2 k}\otimes \nu)$-measurable, $g_{x_0}(\,\cdot\,,y)$ is holomorphic for every $y\in Y$ and $g_{x_0}(x,\,\cdot\,)$ is a representative of $f(x)$ for every $x\in B(x_0,r_{x_0})$. 
	
	Now, take $x_0,x_1\in U$, and observe that, since $B(x_0,r_{x_0})\cap B(x_1,r_{x_1})$ is separable, there is a $\nu$-negligible subset $N_{x_0,x_1}$ of $Y$ such that $g_{x_0}(x,y)=g_{x_1}(x,y)$ for every $x\in B(x_0,r_{x_0})\cap B(x_1,r_{x_1})$ and for every $y\in Y\setminus N_{x_0,x_1}$.
	Next, observe that there is a countable subset $C$ of $U$ such that 
	\[
	U=\bigcup_{x\in C} B(x,r_x).
	\]
	Then, set $N\coloneqq \bigcup_{x_0,x_1\in C} N_{x_0,x_1}$, and define $g(x,y)=g_{x_0}(x,y)$ for every $x\in B(x_0,r_{x_0})$, for every $x_0\in C$, and for every $y\in Y\setminus N$, while $g(x,y)=0$ for every $(x,y)\in U\times N$. The preceding remarks show that $g$ is well-defined, that $N$ is $\nu$-negligible, and that $g$ satisfies the conditions of the statement.
\end{proof}

\begin{proof}[Proof of Proposition~\ref{prop:4}.]
	\textsc{Step I.} Assume first that either $q=\infty$ or $\tau$ is compactly supported in $\Lambda_+$. Observe that, if $q=\infty$, then our assumptions imply that $\tau$ is supported on $\overline{\Omega'}$.
	Define $\widetilde f_0\colon \Omega+i F \to L^2(\Nc)$ so that, for every $h\in \Omega+i F$, 
	\[
	\pi_\lambda(\widetilde f_0(h))=\chi_{\Lambda_+}(\lambda)\ee^{-\langle \lambda_\C,h\rangle}\tau(\lambda)
	\]
	for almost every $\lambda\in F'\setminus W$, thanks to Corollary~\ref{cor:7}.
	Let us prove that $\widetilde f_0$ is holomorphic. With the notation of Corollary~\ref{cor:7}, this is equivalent to showing that the mapping 
	\[
	h\mapsto \ee^{-\langle (\,\cdot\,)_\C, h \rangle}\tau\in \frac{2^{n-m}}{\pi^{n+m}}\int_{\Lambda_+}^\oplus \Lin^2(H_\lambda) \abs{\Pfaff(\lambda)}\,\dd \lambda
	\] 
	is holomorphic, and this latter fact follows easily from Morera's theorem, thanks our assumptions.
	Next, take $v\in E$ and let $\overline{Z_v}$ be the left-invariant vector field on $\Nc$ which induces $\overline{\partial_{E,v}}$ at the origin (cf.~Definition~\ref{def:9}). Then, 
	\[
	\pi_\lambda(\overline{Z_v} \widetilde f_0(h))= -\ee^{-\langle \lambda,h\rangle}\tau(\lambda) \dd\pi_\lambda(\overline{Z_v})=0
	\]
	for almost every $\lambda\in F'$,\footnote{To see this, convolve $\widetilde f_0(h)$ on the right with a smooth approximate identity, and pass to the limit.} so that $\overline{Z_v}\widetilde f_0(h)=0$ (in the sense of distributions) for every $h\in \Omega+i F$. 
	Next, take $w\in F$, and observe that
\[ 
\pi_\lambda(\partial_w (\widetilde f_0(h)))
=i\langle \lambda, w\rangle \ee^{-\langle \lambda_{\C},h\rangle}\tau(\lambda)=- i\partial_w\big( \ee^{-\langle \lambda_{\C},\,\cdot\,\rangle}\tau(\lambda)\big)(h)=-i\pi_\lambda\big( (\partial_{w}\widetilde f_0)(h)\big)
\] 
for almost every $\lambda\in \Lambda_+$ and for every $h\in \Omega+i F$. Now, observe that Lemma~\ref{lem:7} shows that there is a measurable mapping $\widetilde f_1\colon (\Omega+i F)\times \Nc\to \C$ such that $\widetilde f_1(\,\cdot\,,(\zeta,x))$ is holomorphic for every $(\zeta,x)\in \Nc$, while $\widetilde f_1(h,\,\cdot\,)$ is a representative of $\widetilde f_0(h)$ for every $h\in \Omega+i F$. 
Define $\widetilde f_2\colon D\to \C$ so that 
\[
\widetilde f_2(\zeta,z)\coloneqq \widetilde f_1(\Im z- \Phi(\zeta), (\zeta, \Re z))
\]
for every $(\zeta,z)$, so that $\widetilde f_2$ is measurable and locally square-integrable. In addition, $\overline{Z_v} \widetilde f_2=0$ for every $v\in E$ (interpreting $\overline{Z_v}$ as the left-invariant vector field on $E\times F_\C$ which induces $\overline{\partial_{E,v}}$ at $(0,0)$), while
\[
\partial_{i w} \widetilde f_2= i \partial_w \widetilde f_2
\]
	for every $w\in F$. It then follows that $\widetilde f_2$ satisfies Cauchy--Riemann equations on $D$ (in the sense of distributions), so that it has a holomorphic representative $f\in A^{2,q}_{\mi}(D)$ (resp.\ $f\in A^{2,q}_{\mi,\loc}(D)$) which necessarily satisfies the conditions of the statement.
	
\textsc{Step II.} Now, consider the general case. We may assume that $q<\infty$.
Let $( \Sigma_j )$ be an increasing sequence of compact subsets of $\Lambda_+$ whose union is $\Lambda_+$. Then,~{step I} implies that for every $j$ there is $\widetilde f_j\in A^{2,q}_{\mi}(D)$ (resp.\ $\widetilde f_j\in A^{2,q}_{\mi,\loc}(D)$) such that 
\[
\pi_\lambda((\widetilde f_j)_h)=\chi_{\Sigma_j }(\lambda)\ee^{-\langle \lambda,h\rangle}  \tau(\lambda)
\]
for almost every $\lambda\in \Lambda_+$ and for every $h\in \Omega$. In addition, by dominated convergence we see that $(f_j)$ is a Cauchy sequence in $A^{2,q}_{\mi}(D)$ (resp.\ $A^{2,q}_{\mi,\loc}(D)$), so that it converges to a limit $f$ in $A^{2,q}_{\mi}(D)$ (resp.\ $A^{2,q}_{\mi,\loc}(D)$). The function $f$ then satisfies the conditions of the statement.
	
	\textsc{Step III.} Observe that $\tr(\tau(\lambda)\pi_\lambda(\zeta,\Re z)^* )= \tr(\tau(\lambda) (\pi_\lambda(\zeta,\Re z) P_{\lambda,0})^*)$, and that 
	\[
	\norm{\pi_\lambda(\zeta,\Re z) P_{\lambda,0}}_{\Lin^2(H_\lambda)}^2=\tr(P_{\lambda,0})=1
	\]
	for every $\lambda\in \Lambda_+$ and for every $(\zeta,z)\in D$.
	Consequently, if 
	\[
	\int_{\Lambda_+} \norm{\tau(\lambda)}_{\Lin^2(H_\lambda)} \ee^{-\langle \lambda,h \rangle} \abs{\Pfaff( \lambda)}\,\dd \lambda <\infty
	\]
	for every $h\in \Omega$, then
	\[
	f(\zeta,z)\coloneqq\frac{2^{n-m}}{\pi^{n+m}} \int_{\Lambda_+} \tr(\tau(\lambda) \pi_\lambda(\zeta,\Re z)^*) \ee^{-\langle \lambda, \Im z-\Phi(\zeta)\rangle}\abs{\Pfaff(\lambda)}\,\dd \lambda
	\]
	for every $(\zeta,z)\in D$, thanks to Corollary~\ref{cor:7}. 	
\end{proof}

\begin{cor}\label{cor:4}
	Assume that $\Supp{\mi}=\Omega$ and define
	\[
	S((\zeta,z),(\zeta',z'))=\frac{2^{n-m}}{\pi^{n+m}} \int_{\Omega'}  \ee^{\langle \lambda_\C, i(z-\overline{z'})+2 \Phi(\zeta,\zeta')\rangle} \abs{\Pfaff(\lambda)}\,\dd \lambda
	\]
	for every $(\zeta,z),(\zeta',z')\in D$. Then, $S$ is the reproducing kernel of $A^{2,\infty}_\mi(D)$, that is, 
	\[
	f(\zeta,z)=\langle f\vert S(\,\cdot\,,(\zeta,z))\rangle_{A^{2,\infty}_\mi(D)}
	\]
	for every $f\in A^{2,\infty}_\mi(D)$ and for every $(\zeta,z)\in D$.
\end{cor}

This gives an alternative proof of~\cite[Theorem 5.3]{Gindikin}.
The function $S$  is also known as the Cauchy--Szeg\H o kernel.\index{Cauchy--Szeg\H o kernel}

\begin{proof}
	For every $(\zeta',z')\in D$, denote by $ S_{(\zeta',z')}  $ the unique element of the hilbertian space $A^{2,\infty}_\mi(D)$ such that $\langle f \vert  S_{(\zeta',z')}  \rangle=f(\zeta',z')$ for every $f\in A^{2,\infty}_\mi(D)$ (cf.~Proposition~\ref{prop:1} and Corollaries~\ref{cor:2} and~\ref{cor:1}). Set $x'\coloneqq \Re z'$ and ${h'}\coloneqq \Im z'-\Phi(\zeta')$, and observe that
	\[
	\begin{split}
	f(\zeta',z')=f_{h'}(\zeta',x')&= \frac{2^{n-m}}{\pi^{n+m}} \int_{\Omega'} \tr(\pi_\lambda(f_{h'}) \pi_{\lambda}(\zeta',x')^*) \abs{\Pfaff(\lambda)}\,\dd \lambda\\
		&=\frac{2^{n-m}}{\pi^{n+m}} \int_{\Omega'} \tr(\pi_\lambda(f_0) (\pi_{\lambda}(\zeta',x') P_{\lambda,0})^*) \ee^{-\langle \lambda,{h'}\rangle}\abs{\Pfaff(\lambda)}\,\dd \lambda,
	\end{split}
	\]
	where $f_0$ is the limit of $f_h$ in $L^2(\Nc)$ for $h\to 0$ (cf.~Proposition~\ref{prop:3} and Corollary~\ref{cor:1}). 
	Conversely, Proposition~\ref{prop:4} and Lemma~\ref{lem:8} show that, for every $f_0\in L^2(\Nc)$ such that $\pi_\lambda(f_0)=\chi_{\Omega'}(\lambda)\pi_\lambda(f_0)P_{\lambda,0}$ for almost every $\lambda\in F'\setminus W$, there is a unique $f\in A^{2,\infty}_\mi(D)$ such that $f_0=\lim\limits_{h\to 0} f_h$ in $L^2(\Nc)$. By the arbitrariness of $f_0$, it then follows that 
	\[
	\pi_\lambda(( S_{(\zeta',z')}  )_h)=\ee^{-\langle \lambda,h+{h'}\rangle}\pi_{\lambda}(\zeta',x') P_{\lambda,0}
	\]
	for every $h\in \Omega$ and for almost every $\lambda\in \Omega'$. Hence, Proposition~\ref{prop:4} and Lemma~\ref{lem:8} show that $(S_{(\zeta',z')} )_h(\zeta,x)$ equals
	\[
	\frac{2^{n-m}}{\pi^{n+m}} \int_{\Omega'} \tr(P_{\lambda,0} \pi_\lambda( \zeta'-\zeta, x'-x-2\Im \Phi(\zeta, \zeta') ) P_{\lambda,0}) \ee^{-\langle \lambda , h+{h'}\rangle  }\abs{\Pfaff(\lambda)}\,\dd \lambda
	\]
	by the Fourier inversion formula, for (almost) every $(\zeta,x)\in \Nc$ and for every $h\in \Omega$. The assertion follows from Proposition~\ref{prop:36}.
\end{proof}

The following result, which is based on the reproducing properties of the Poisson--Szeg\H o kernel, has some interesting consequences which do not arise in the study of tube domains. 
In particular, if $\Lambda_+\setminus \Omega'$ is negligible, $p\in [1,\infty]$,  and $f\in A^{p,\infty}_{\mi,\loc}(D)$, then $f(\,\cdot\,+i h)\in A^{p,\infty}_\mi(D)$ for every $h\in \Omega$, a fact which is clearly false in tube domains.

\begin{cor}\label{cor:8}
	Take $p\in [1,\infty]$ and assume that $\Supp{\mi}=\Omega$. Assume that $\Lambda_+\setminus\Omega'$ is negligible and take $f\in A^{p,\infty}_{\mi,\loc}(D)$.
	Then, the mapping $h\mapsto \norm{f_h}_{L^p(\Nc)}$ is decreasing on $\Omega$. 
\end{cor}

Observe that the assumption that $\Lambda_+\setminus \Omega'$ is negligible (for $\Hc^m$) means that the closed convex envelope of $\Phi(E)$ is $\overline \Omega$. Equivalently, $\overline{\Lambda_+}=\overline{\Omega'}$.

\begin{proof}
With the notation of Corollary~\ref{cor:4}, define (the Poisson--Szeg\H o kernel\index{Poisson--Szeg\H o kernel}), for every $(\zeta,z)\in D$, 
\[
P_{(\zeta,z)} \coloneqq S((\zeta,z),(\zeta,z))^{-1} \abs{S(\,\cdot\,,(\zeta,z))_0}^2,
\]
where $S(\,\cdot\,,(\zeta,z))_0$ denotes the limit of $S(\,\cdot\,,(\zeta,z))_h$ in $L^2(\Nc)$ as $h\to 0$. Hence, $P_{(\zeta,z)}$ is positive, has integral $1$, and 
\[
\lim_{(\zeta,z)\to(\zeta_0,x_0+i\Phi(\zeta_0))} \langle P_{(\zeta,z)}\cdot \Hc^{2 n+m}, \varphi\rangle= \varphi(\zeta_0,x_0)
\]
for every bounded continuous function $\varphi$ on $\Nc$ and for every $(\zeta_0,x_0)\in \Nc$ (cf.~\cite[Theorem 2.4]{Koranyi}). 	In addition, 
\[
P_{(\zeta',x'+i\Phi(\zeta'))\cdot(\zeta,z)}=P_{(\zeta,z)}((-\zeta',-x')\,\cdot\,)
\]
for every $(\zeta,z)\in D$ and for every $(\zeta',x')\in \Nc$. Further, if $p>1$ and $f\in A^{p,\infty}_\mi(D)$, then 
\[
f(\zeta,z+i h)=\langle f_h\vert P_{(\zeta,z)}\rangle
\]
for every $(\zeta,z)\in D$ and for every $h\in \Omega$ (cf.~\cite[Proposition 2.7]{Koranyi}).
	
	Assume first that $p\meg 2$, so that $f\in A^{2,\infty}_{\mi,\loc}(D)$ by Proposition~\ref{prop:1}. By means of Proposition~\ref{prop:2}, we then see that the mapping $h\mapsto \norm{f_h}_{L^2(\Nc)}$ is decreasing on $\Omega$, so that $f(\,\cdot\,+i h)\in A^{2,\infty}_\mi(D)$ for every $h\in \Omega$. Then, $f(\zeta,z+i h)=\langle f_h\vert P_{(\zeta,z)}\rangle$ for every $(\zeta,z)\in D$ and for every $h\in \Omega$ by the preceding remarks. Further,
\[
P_{(\zeta,x+i\Phi(\zeta)+i h)}(0)=P_{(0,i h)}(-\zeta,-x)
\]
for every $(\zeta,x)\in \Nc$ and for every $h\in \Omega$, so that
\[
f_{h+h'}= f_h* P_{(0,i h')}
\]
for every $h,h'\in \Omega$ (cf.~the proof of~\cite[Proposition 2.7]{Koranyi}). Using Young's inequality, we then see that $\norm{f_{h+h'}}_{L^p(\Nc)}\meg \norm{f_h}_{L^p(\Nc)}$. 
	
	The case $p>2$ is treated assuming first that $f\in A^{2,\infty}_{\mi,\loc}(D)\subseteq A^{p,\infty}_{\mi,\loc}(D)$ and then using Proposition~\ref{prop:37} to get the conclusion, when $p<\infty$. If $p=\infty$, then take $(g^{(\eps)})_{\eps>0}$ as in Lemma~\ref{lem:34}. Then, the above arguments show that $(g^{(\eps)}f)_{h+h'}= (g^{\eps} f)_h* P_{(0,i h')}$ for every $h,h'\in \Omega$ and for every $\eps>0$, whence $f_{h+h'}= f_h* P_{(0,i h')} $ passing to the limit. Thus, the mapping $h\mapsto \norm{f_h}_{L^\infty(\Nc)}$ is decreasing. 
\end{proof}

\begin{cor}\label{cor:5}
Assume that $A^{2,2}_\mi(D)$ is complete and that $\Supp{\mi}=\Omega$. Define $K_\mi \colon D\times D\to \C$ so that
	\[
	f(\zeta',z')=\langle f\vert K_\mi  (\,\cdot\,,(\zeta',z'))\rangle_{A^{2,2}_\mi(D)}
	\]
	for every $f\in A^{2,2}_\mi(D)$ and for every $(\zeta',z')\in D$.
	Then, 
	\[
	\pi_\lambda(K_\mi (\,\cdot\,,(\zeta',z'))_h)= \chi_{\frac 1 2 \Ds_\mi\cap \Lambda_+}(\lambda) \frac{1}{\Lc \mi(2\lambda)} \ee^{-\langle \lambda,h+\Im z'-\Phi(\zeta')\rangle} \pi_\lambda(\zeta',\Re z') P_{\lambda,0}
	\]
	for almost every $\lambda\not \in W$, for every $(\zeta',z')\in D$, and for every $h\in \Omega$.
	
	If, in addition, $\frac 1 2  \Ds_\mi \cap \Lambda_+\subseteq \overline{\Omega'}$, then
	\[
	K_\mi  ((\zeta,z),(\zeta',z'))=\frac{2^{n-m}}{\pi^{n+m}} \int_{ \frac 1 2  \Ds_\mi \cap \Omega'}  \ee^{\langle \lambda_\C, i(z-\overline{z'})+2  \Phi(\zeta,\zeta')\rangle} \frac{\abs{\Pfaff(\lambda)}}{\Lc \mi(2 \lambda)}\,\dd \lambda
	\]
	for every $(\zeta,z),(\zeta',z')\in D$.
\end{cor}

The function $ K_\mi  $, which is the reproducing kernel of the weighted Bergman space $A^{2,2}_{\mi}(D)$, is also known as the (weighted) Bergman kernel.\index{Bergman kernel}

Recall that $ \Ds_\mi $ is the set of $\lambda\in F'$ for which $\Lc \mi(\lambda)$ is finite.

\begin{proof}
	Notice that $ K_\mi $ is well defined thanks to Lemma~\ref{lem:9}.
	Define $K_{(\zeta',z')}\coloneqq   K_\mi  (\,\cdot\,,(\zeta',z'))$ for every $(\zeta',z')\in D$.  Set $x'\coloneqq \Re z'$ and ${h'}\coloneqq \Im z'-\Phi(\zeta')$, and observe that, with the notation of Corollary~\ref{cor:3},
	\[
	\begin{split}
	f(\zeta',z')&= \frac{2^{n-m}}{\pi^{n+m}} \int_{\frac1 2 \Ds_\mi \cap\Lambda_+} \tr(\pi_\lambda(f_{h'}) \pi_{\lambda}(\zeta',x')^*) \abs{\Pfaff(\lambda)}\,\dd \lambda\\
	&=\frac{2^{n-m}}{\pi^{n+m}} \int_{\frac 1 2  \Ds_\mi \cap \Lambda_+} \tr(\tau_f(\lambda) (\pi_{\lambda}(\zeta',x') P_{\lambda,0})^*) \ee^{-\langle \lambda,{h'}\rangle}\abs{\Pfaff(\lambda)}\,\dd \lambda,
	\end{split}
	\] 
	(cf.~Proposition~\ref{prop:3} and Corollary~\ref{cor:3}), but also
	\[
	f(\zeta',z')= \frac{2^{n-m}}{\pi^{n+m}} \int_{\frac1 2  \Ds_\mi \cap \Lambda_+} \tr( \tau_f(\lambda) \tau_{ K_{(\zeta',z') }}(\lambda)^*) \Lc\mi(2\lambda) \abs{\Pfaff( \lambda)}\,\dd \lambda.
	\]	
	In addition, take a compact subset $ \Sigma $ of the interior of $\frac 1 2  \Ds_\mi \cap\Lambda_+$ and observe that Proposition~\ref{prop:3} implies that there is a unique $f_{\Sigma } \in A^{2,2}_{\mi}(D)$ such that $\pi_\lambda((f_{ \Sigma } )_h)=\chi_{ \Sigma } (\lambda) \ee^{-\langle \lambda,h\rangle}  P_{\lambda,0}$ for almost every $\lambda\not \in W$. 
	Taking into account the arbitrariness of $\Sigma $ and observing that $\Lc \mi(2\lambda)$ is finite and non-zero for every $\lambda \in \frac 1 2  \Ds_\mi $ since the support of $\mi$ is the whole of $\Omega$,  it follows that 
	\[
	\pi_\lambda(( K_{(\zeta',z')})_h)  =\frac{1}{\Lc\mi(2\lambda)}\ee^{-\langle \lambda,h+{h'}\rangle}\pi_{\lambda}(\zeta',x') P_{\lambda,0}
	\]
	for every $h\in \Omega$ and for almost every $\lambda\in \frac1 2  \Ds_\mi \cap \Lambda_+$.
	
	Now, if $\frac 1 2  \Ds_\mi \cap\Lambda_+ \subseteq \overline{\Omega'}$, then Proposition~\ref{prop:4} and Lemma~\ref{lem:8}  imply that  $(B_{(\zeta',z')})_h(\zeta,x)$ equals
	\[
	\frac{2^{n-m}}{\pi^{n+m}} \int_{\frac 1 2  \Ds_\mi \cap \Omega'} \tr(P_{\lambda,0} \pi_\lambda( \zeta'-\zeta, x'-x -2 \Im \Phi(\zeta, \zeta')) P_{\lambda,0}) \ee^{-\langle \lambda , h+{h'}\rangle  }\frac{\abs{\Pfaff(\lambda)}}{\Lc \mi(2 \lambda)}\,\dd \lambda
	\]
	for every $(\zeta,x)\in \Nc$ and for every $h\in \Omega$. Thus, the assertion follows from Proposition~\ref{prop:36}.
\end{proof}

\section{The Kohn Laplacian on $\Nc$}

We now illustrate some connections  between the Hardy spaces $A^{2,\infty}_{\Hc^m}(D)$ and the Kohn Laplacian on $\Nc$ (associated with the chosen scalar product on $E$). We refer the reader to~\cite{PelosoRicciJFA,PelosoRicci2} for further information on this topic,  and references therein. 

Take $q\in \Set{0,\dots, n}$, and denote by $\Lambda_q(\Nc)$ the pull-back under the canonical projection $\Nc\to E$ of the $(0,q)$-exterior power of the cotangent bundle of $E$. If $\zeta_1,\dots, \zeta_n$ are coordinate functions on $E$ associated with a fixed orthonormal basis $(v_1,\dots, v_n)$ of $E$, a section of $\Lambda_q(\Nc)$ is therefore a mapping of the form
\[
\sum_{\alpha} a_\alpha \dd \overline\zeta^\alpha,
\]
where $a_\alpha$ is a (real-valued) function on $\Nc$ and $\dd\overline \zeta^\alpha=\dd\overline{\zeta_{\alpha_1}}\wedge\cdots\wedge \dd\overline{\zeta_{\alpha_q}}$ for every strictly increasing multi-index $\alpha$ in $\N^q$.
The scalar product on $E$ then induces a scalar product on the space of compactly supported continuous sections of $\Lambda_q(\Nc)$
\[
(\omega_1,\omega_2)\mapsto \int_\Nc \langle \omega_1(\zeta,x),\omega_2(\zeta,x)\rangle\,\dd (\zeta,x).
\]
Now, for every smooth section $\omega=\sum_{\alpha}\omega_\alpha \dd \overline\zeta^\alpha$ of $\Lambda_q(\Nc)$, define
\[
\bar \partial_b \omega\coloneqq \sum_\alpha \sum_{k=1}^n \overline{ Z_k}(\omega_\alpha) \dd \overline{\zeta_k}\wedge\dd \overline\zeta^\alpha,
\]
where $\overline{Z_k}$ is the left-invariant vector field on $E$ which induces $\partial_{E,v_k}$ at the origin.
Denote by $\bar \partial_b^*$ the formal adjoint of $\bar \partial_b$ with respect to the scalar product previously defined. The left-invariant differential operator  acting on smooth sections of $\Lambda_q(\Nc)$
\[
\Box_b^{(q)}\coloneqq \bar \partial_b \bar \partial_b^*+ \bar \partial_b^* \bar \partial_b
\]
is called the Kohn Laplacian. In~\cite[Proposition 2.1]{PelosoRicciJFA}, an explicit expression of $\Box_b^{(q)}$ has been provided. In particular,
\[
\Box_b^{(0)}=\Lc+i\sum_{k=1}^n \partial_{\Phi(v_k)},
\]
where $\Lc\coloneqq -\frac 1 2\sum_{k=1}^n (\overline{Z_k} Z_k+Z_k \overline{Z_k})$. 

\begin{prop}
The kernel of $\Box_b^{(0)}$ in $L^2(\Nc)$ is 
\[
\widetilde H^2(\Nc)\coloneqq \Set{f\in L^2(\Nc)\colon \pi_\lambda(f)=\chi_{\Lambda_+}(\lambda) \pi_\lambda(f)P_{\lambda,0} \text{ for almost every $\lambda\in F'\setminus W$}}.
\]
If $\Lambda_+\setminus\Omega'$ is negligible, then $\widetilde H^2(\Nc)$ is the set of boundary values of the Hardy space $A^{2,\infty}_{\Hc^m}(D)$.
\end{prop}

\begin{proof}
We observe that \CB
\[
\pi_\lambda(\Box_b^{(0)} f  )=\pi_\lambda(f)\dd \pi_\lambda(\Box_b^{(0),\dag})
\]
for almost every $\lambda\in F'\setminus W$ and for every $f\in L^2(\Nc)$ such that $\Box_b^{(0)} f\in L^2(\Nc)$, where $\Box_b^{(0),\dag}$ denotes the formal transpose of $\Box_b^{(0)}$, that is
\[
\Box_b^{(0),\dag}=\Lc-i\sum_{k=1}^n \partial_{\Phi(v_k)}.
\]
Now,
\[
\dd\pi_\lambda(\Box_b^{(0),\dag})=\dd \pi_\lambda(\Lc)-\sum_{k=1}^n \langle \lambda, \Phi(v_k)\rangle I_{H_\lambda}
\]
for every $\lambda\in F'\setminus W$. Take $\lambda\in F'\setminus W$, and let $v_{\lambda,1},\dots,v_{\lambda,n}$ be an orthonormal basis of $E$ which diagonalizes $- i J_\lambda$. Let $\theta_{\lambda,1},\dots, \theta_{\lambda,n}$ be the corresponding eigenvalues, and let $Z_{\lambda,1},\dots, Z_{\lambda,n}$ be the left-invariant vector fields on $\Nc$ which induce $\partial_{E,v_1},\dots,\partial_{E,v_n}$ at the origin. 
%Notice that we may assume that $v_{\lambda,1},\dots, v_{\lambda,n_\lambda}\in E_{\lambda,+}$ while $v_{\lambda,n_\lambda+1},\dots, v_{\lambda,n}\in E_{\lambda,+}$ for some $n_\lambda\in \Set{0,\dots,n}$, so that $- i J_\lambda v_{\lambda,k}=\abs{J_\lambda} v_{\lambda,k}$ if $k\meg n_\lambda$, while  $- i J_\lambda v_{\lambda,k}=-\abs{J_\lambda} v_{\lambda,k}$ if $k> n_\lambda$. 
Then, it is readily verified that
\[
\Lc=-\frac 1 2 \sum_{k=1}^n (\overline{Z_{\lambda,k}}Z_{\lambda,k}+Z_{\lambda,k}\overline{Z_{\lambda,k}}),
\]
so that Proposition~\ref{prop:36} implies that, if we define $e_{\lambda,\alpha}\coloneqq \prod_{k=1}^n \langle \lambda, \Phi_\lambda(\,\cdot\,,v_{\lambda,k})\rangle^{\alpha_k}$ for every $\alpha\in \N^n$, then
\[
\dd \pi_\lambda(\Lc) e_{\lambda,\alpha}= \sum_{k=1}^{n} ( (1+2 \alpha_k) \abs{\theta_{\lambda,k}} )  e_{\lambda, \alpha}. 
\]
In addition, it is not hard to see that $(e_{\lambda,\alpha} )_{\alpha\in \N^n}$ is a complete orthogonal system in $H_\lambda$ (cf.~\cite[Theorem 1.63]{Folland}).
Therefore, $\dd \pi_\lambda(\Box_b^{(0),\dag})$ is a positive self-adjoint operator with purely discrete spectrum, its least eigenvalue is 
\[
\tr_\C\abs{J_\lambda}-\sum_{k=1}^n \langle \lambda, \Phi(v_k)\rangle= \tr_\C(\abs{J_\lambda}+i J_\lambda),
\]
and the corresponding eigenprojector is $P_{\lambda,0}$.
It then follows that $\ker \dd\pi_\lambda(\Box_b^{(0),\dag})=\Set{0}$ for every $\lambda\in F'\setminus (W\cup \Lambda_+)$, while $\ker \dd \pi_\lambda(\Box_b^{(0),\dag})=\C e_{\lambda,0}$ if $\lambda\in \Lambda_+$. 
Therefore, the assertion follows.

In particular, if $\Lambda_+\setminus\Omega'$ is negligible, then $\widetilde H^2(\Nc)$ is the set of boundary values of the Hardy space $A^{2,\infty}_{\Hc^m}(D)$, thanks to Corollary~\ref{cor:1} and Proposition~\ref{prop:4}.
\end{proof}

\section{Notes and Further Results}

\paragraph{1.6.1} The material collected in this chapter is far from being complete or comprehensive. Several relevant topics, such as non-tangential or admissible limits at the \v Silov boundary of the elements of the Hardy spaces, or the properties of the Poisson kernel, are not treated in detail.
Only Paley--Wiener theorems are treated in depth, even though the description of the `boundary values' of the elements of the weighted Bergman spaces $A^{2,q}_{\mi}(D)$ is only given up to a Fourier transform (cf.~Propositions~\ref{prop:2}). For the weights considered in Chapter~\ref{sec:6}, we shall characterize these spaces as suitable Besov spaces (at least for $q\meg 2$). One may wonder if an analogous assertion holds for more general weights. Cf.~\cite{Rychkov} and the references therein for a treatment of Besov and Triebel--Lizorkin spaces with weights in suitable local Muckenhoupt classes. 

\paragraph{1.6.2} The measures $\Hc^{2n+m}\otimes\mu$ are precisely
the Radon measures on $\Nc\times\Omega$ that are invariant under the 'horizontal translations' by $\Nc$.
Thus, under the identification of $D$ with $\Nc\times\Omega$, such measures constitute the analogue of the translation invariant measures on $\C_+$, or, in its bounded realization, the rotational invariant measures on the unit disk.

\paragraph{1.6.3}
Most of the results in this chapter are proved under the assumption that $\supp(\mi)=\Omega$.  The case when $\supp(\mi)\subsetneqq\Omega$ is certainly of
interest, but the situation is by far more complicated.  An interesting space of holomorphic functions on $\C_+$ was studied in~\cite{PS1} in connection with the so-called
\emph{M\"untz--Sz\'asz problem for the Bergman space}, and was denoted
by $\Ms^2$.  In this case $\mu$ is an atomic measure on $\R^*_+$, with
$\supp(\mu)= \frac12 \N$  and, in our notation, 
$\Ms^2= A^{2,2}_\mu(\C_+)\cap A^{2,\infty}_{\Hc^2,\loc}(\C_+)$, endowed with the $A^{2,2}_\mu(\C_+)$-norm.  It was shown in~\cite{PS1} that the space
$A^{2,2}_\mu(\C_+)$ contains wildly behaved functions, such as $\exp(i\ee^{2\pi i z})$.  See also~\cite{PS2} and references therein for a further discussion on functions in $A^{2,2}_\mu(\C_+)$ for a general $\mi$.

\chapter{Homogeneous Siegel Domains of Type II}\label{sec:3}

In this chapter, we shall introduce our notation for various objects related to homogeneous Siegel domains of type II and recall some basic facts.
 
In Section~\ref{sec:2} we shall recall some basic facts on $T$-algebras and the associated homogeneous cones. We shall generally avoid to make use of the associated formalism as long as possible, and reformulate several definitions in a more conceptual way in Sections~\ref{sec:4:1} and~\ref{sec:3:1}. In particular, we shall define the generalized power functions $\Delta_\Omega^{\vect s}$ and $\Delta_{\Omega'}^{\vect s}$ on $\Omega$ and $\Omega'$, respectively. 
In the remainder of Section~\ref{sec:3:1} and in Section~\ref{sec:3:2}, we shall report several results on gamma and beta functions related to these generalized power functions. These results will be of importance in the development of the theory of weighted Bergman spaces (cf.~Section~\ref{sec:6}).

Finally, in Section~\ref{sec:3:3}, we shall recall some basic properties of the Bergman metric on $D$ and define some related metrics on $\Omega$ and $\Omega'$. Notice that these latter metrics are different from the canonical invariant Riemannian metric on general homogeneous cones defined in~\cite[I.4]{FarautKoranyi}, unless $E=\Set{0}$. Nonetheless, since both metrics are invariant, the difference is of minor importance. 
We shall then introduce a suitable notion of lattices on $D$ which is adapted to the decomposition of $D$ as a union of translates of $b D$. This restriction, which is analogous to that employed in~\cite{RicciTaibleson}, is necessary to deal with mixed norm Bergman spaces. For less general Bergman spaces, more general lattices may be employed. We shall then present several quasi-constancy results for various relevant functions on $D$, $\Omega$, and $\Omega'$.

\section{$T$-Algebras}\label{sec:2}

In this brief section we recall some notions from~\cite{Vinberg}. Cf.~also~\cite{NanaTrojan}.
Notice that we shall no longer make use in the sequel of the notation of this chapter, unless explicitly stated, in favour of a more conceptual one.

\medskip

For every $r\in \N^*$, we denote by $\Delta_r$ the monoid $\Set{0}\cup\Set{1,\dots,r}^2$ endowed with the product defined by $(j,k)(j',k')=(j,k')$ if $k=j'$ and $0$ otherwise, while $(j,k)0=0(j,k)=0$ for every $j,j',k,k'=1,\dots, r$. 

\begin{deff}
A $T$-algebra\index{T-algebra@$T$-algebra} of rank $r$ is a (not necessarily associative) finite-dimensional graded algebra $A$ over $\R$ of type $\Delta_r$ (cf.~\cite[Definition 1 of Chapter III, \S 3, No.\ 1]{BourbakiA1}), endowed with a linear involution $^*$ such that the following hold:
\begin{enumerate}
\item[(1)] $A^*_{j,k}=A_{k,j}$ for every $j,k=1,\dots,r$ and $(ab)^*=b^* a^*$ for every $a,b\in A$;

\item[(2)] $A_0=\Set{0}$ and $\dim A_{j,j}=1$ for every $j=1,\dots,r$;

\item[(3)] for every $j=1,\dots, r$ there is $e_j\in A_{j,j}$ such that left and right multiplication by $e_j$ induce the identity on $A_{j,k}$ and $A_{k,j}$, respectively, for every $k=1,\dots, r$;

\item[(4)] if we define $\tr a=\sum_{j=1}^r \langle e_j', a\rangle$ for every $a\in A$, where $e'_j$ is the unique graded linear functional on $A$ which takes the value $1$ at $e_j$ ($j=1,\dots,r$), then the mapping $(a,b)\mapsto \tr(a b)$ is symmetric;

\item[(5)] for every $a,b,c\in A$, $\tr(a(bc))=\tr((ab)c)$;

\item[(6)] the symmetric\footnote{Observe that $\tr(a)=\tr(a^*)$ for every $a\in A$, since the only automorphism of the field $\R$ is the identity, and since the subalgebras $A_{j,j}$ of $A$ are isomorphic to $\R$. Then,~{(4)} shows that $\tr(a^* b)=\tr(b^* a)$ for every $a,b\in A$.} bilinear mapping $(a,b)\mapsto \tr(a^*b)$ is positive and non-degenerate;

\item[(7)] setting $T\coloneqq\bigoplus_{j\meg k} A_{j,k}$, one has $t(u w)=(t u)w$ and $t( u u^*)=(t u) u^*$ for every $t,u,w\in T$.
\end{enumerate}

The vector space subjacent to $A$ may then be endowed with the scalar product $\langle a,b\rangle\coloneqq \tr(a^*b)$, $a,b\in A$. 
\end{deff}

Notice that $A_{j,k}$ is orthogonal to $A_{j',k'}$ if $(j,k)\neq (j',k')$. Indeed, if $a\in A_{j,k}$ and $b\in A_{j',k'}$, then $a^*b=0$ if $j\neq j'$, and $b a^*=0$ if $k'\neq k$, whence $\langle a,b\rangle=\tr(a^*b)=\tr(b a^*)=0$.

\begin{deff}
Define 
\[
\begin{aligned}
T_+&\coloneqq \Set{a\in T\colon \langle e_j',a\rangle>0 \quad \forall j=1,\dots, r }, &
H&\coloneqq\Set{a\in A\colon a=a^*},\\
C(A)&\coloneqq \Set{t t^*\colon t\in T_+}, &
C'(A)&\coloneqq \Set{t^* t\colon t\in T_+}.
\end{aligned}
\]
The cones $C(A)$ and $C'(A)$ are said to have rank $r$.\index{Rank!of a homogeneous cone}
\end{deff}

The following result is proved in~\cite{Vinberg}. Recall that an open convex cone $C$ not containing any affine line is said to be homogeneous if the group of affine automorphisms of $C$ acts transitively on $C$.\index{Cone!homogeneous}

\begin{teo}\label{teo:1}
The following hold:
\begin{enumerate}
\item[\em(1)] $T_+$, endowed with the product induced by $A$, is a Lie group;

\item[\em(2)] the algebra $A'$ with the same operations as $A$ and the graduation $A'_{j,k}\coloneqq A_{r+1-j,r+1-k}$, $j,k=1,\dots,r$, is a $T$-algebra;

\item[\em(3)] $C(A)$ and $C'(A)$ are open convex cones, and are dual to one another with respect to the scalar product on $H$;

\item[\em(4)] $(t x) t^*= t(x t^*)$ and $(t^* x) t = t^* (x t)$ for every $t\in T_+$ and for every $x\in H$;

\item[\em(5)] the mappings 
\[
(t,x)\mapsto t x t^* \qquad \text{and}\qquad (t,x)\mapsto t^* x t
\]
are simply transitive left and right actions of $T_+$ on $C(A)$ and $C'(A)$, respectively, which are dual to one another with respect to the scalar product on $H$.
\end{enumerate}

Finally, every homogeneous open convex cone not containing any affine line is isomorphic to $C(A)$ for some $T$-algebra $A$. 
\end{teo}

The following lemmas are the analogues of well-known facts concerning (classical) triangular groups.

\begin{lem}\label{lem:68}
Define $T_1$ as the set of $t\in T_+$ such that $\langle e'_j, t\rangle=1$ for every $j=1,\dots, r$. Then, $[T_+,T_+]=T_1$ and $T_1$ is a nilpotent subgroup of $T$.
\end{lem}

\begin{proof}
Recall that, by Theorem~\ref{teo:1}, $T_+$ is a Lie group. In addition, define a mapping $X$ from $T$ into the set of left-invariant vector fields on $T_+$ so that $(X(t) f)(t')=f'(t')\cdot (t' t)$ for every $t\in T$, for every $t'\in T_+$, and for every $f\in C^1(T_+)$. Then, it is easily verified that $X$ is an isomorphism of $T$ onto the Lie algebra of $T_+$, and that $[X(t),X(t')]=X(t t'-t' t)$ for every $t,t'\in T$. 
In addition, it is also clear that the subspace $T_0\coloneqq \bigoplus_{j<k}A_{j,k}$ of $T$ identifies with the Lie algebra of $T_1$, which is clearly nilpotent. 
Since $T_+$ and $T_1$ are connected, by~\cite[Proposition 4 of Chapter III, \S 9, No.\ 2]{BourbakiLie2} it will suffice to show that $[T,T]=T_0$. Now, take $\alpha\in \R^r$ and $\alpha_j\neq \alpha_k$ whenever $j\neq k$, and define $t_\alpha\coloneqq \sum_{j=1}^r\alpha_j e_j$. Then, for every $t\in T_0$, 
\[
t t_\alpha-t_\alpha t= \sum_{j<k} (\alpha_k-\alpha_j) t_{j,k},
\]
so that the mapping $t\mapsto t t_\alpha-t_\alpha t$ is an automorphism of $T_0$. Hence, $[T,T]=T_0$ and the assertion follows.
\end{proof}

\begin{lem}\label{lem:69}
For every $\vect{s}  \in \C^r$ the mapping 
\[
\Delta^{\vect{s}} \colon T_+\ni  t\mapsto \prod_{j=1}^r \langle e_j', t\rangle^{2 s_j} \in \C^*
\]
 is a homomorphism of Lie groups. Conversely, every homomorphism of Lie groups from $T_+$ into $\C^*$ is of this form. 
\end{lem}

The definition of $\Delta^{\vect s}$ may seem peculiar. Nonetheless, since we are more interested in the cone $C(A)$ than in the group $T_+$, we define the homomorphism $\Delta^{\vect s}$  in such a way that it induces the mapping $\sum_{j=1}^r t_j e_j\mapsto \prod_{j=1}^r t_j^{s_j}$ when transferred on $C(A)$ by means of the bijection $t\mapsto t t^*$.  Observe also that 
$\Delta^{\vect{s}} \colon T_+ \to \R_+^*$ if  and only if $ \vect{s}  \in \R^r$.  

\begin{proof}
It is clear that $\Delta^{\vect{s}}$ is a homomorphism, since $\langle e_j', t t'\rangle= \langle e_j' ,t \rangle \langle e_j', t'\rangle$ for every $t,t'\in T$ and for every $j=1,\dots, r$. 
Conversely, let $\Delta$ be a homomorphism of $T_+$ into $\C^*$. Since $B \coloneqq \bigoplus_{ j=1 }^r \R_+^* e_j$ is a subgroup of $T_+$ isomorphic to $(\R_+^*)^r$, there is a unique ${\vect{s}}\in \C^r$ such that $\Delta(t)=\prod_{j=1}^r \langle e_j', t\rangle^{2 s_j}$ for every $t\in B$. Denote by $T_1$ the subgroup of $t\in T$ such that $\langle e_j',t\rangle=1$ for every $j=1,\dots, r$, and observe that $B T_1=T_+$ and that $[T_+,T_+]=T_1$ by Lemma~\ref{lem:68}, so that $\Delta(T_1)=\Set{1}$. Then, $\Delta=\Delta^{\vect{s}}$.
\end{proof}

\section{Notation}\label{sec:4:1}

In this brief section we establish the main notation for homogeneous cones and homogeneous Siegel domains of type II. 

As mentioned in  Section~\ref{sec:2}, we shall interpret the formalism of $T$-algebras in a more conceptual way, fixing a group $T_+$ which acts simply transitively on $\Omega$ and $\Omega'$ and `base points' $e_\Omega$ and $e_{\Omega'}$ in $\Omega$ and $\Omega'$.   
We shall also define two vectors $\vect m$ and $\vect{m'}$ which are closely related to the geometry of $\Omega$ and $\Omega'$, and will play an important role in the sequel.

We shall then impose a condition on $D$ which is equivalent to (affine) homogeneity. This condition implies that the Pfaffian $\abs{\Pfaff}$ induces a character on $T_+$, which will play an important role in the sequel. We shall therefore introduce a proper piece of notation.
\medskip

Define $E$, $F$, $\Omega$, $\Phi$, $D$, $bD$, and $\Nc$ as in Section~\ref{sec:1:1}, and assume further that $\Omega$  is a homogeneous cone. Further, define $\pi_\lambda$, $H_\lambda$, and $\abs{\Pfaff}$ as in Section~\ref{sec:1:2}.
By Theorem~\ref{teo:1}, there are a $T$-algebra $A$ and an isomorphism $\Psi$ of the corresponding space $H$ onto $F$ such that $\Psi(C(A))=\Omega$. Hence, $\trasp \Psi(\Omega')$ is the image of $C'(A)$ under the isomorphism of $H$ onto its dual $H'$ induced by the scalar product of $H$. 
Define 
\[
e_\Omega\coloneqq \Psi\left( \sum_{j=1}^r e_j\right) \qquad \text{and}\qquad e_{\Omega'}=\trasp \Psi^{-1}\left(\sum_{j=1}^r e_j'\right).
\]
We shall assume that $F$ carries the scalar product induced by that of $H$, and we shall say that $r$ is the rank of $\Omega$.\label{16}
Then, $T_+$\label{14} acts simply transitively on the left on $\Omega$ and on the right on $\Omega'$, and the actions of $T_+$ on $\Omega$ and $\Omega'$ are transpose of one another. 

 \begin{deff} \label{16-bis} 
  Define $m_{j,k}\coloneqq \dim A_{j,k}$ for every $j,k=1,\dots, r$, and set\label{15}
\[
\vect m\coloneqq \left(\sum_{k>j} m_{j,k}  \right)_{j=1,\dots, r} \qquad \text{and} \qquad \vect{m'}\coloneqq \left(\sum_{k<j} m_{j,k}\right)_{j=1,\dots, r}.
\]
We then define $\Delta^{\vect s}$, for every $\vect s\in \C^r$, as in Lemma~\ref{lem:69}. 
\end{deff} 

We shall assume that for every $t\in T_+$ there is $g\in GL(E)$ such that\label{2}
\[
t\cdot \Phi=\Phi\circ (g\times g).
\]
This condition is equivalent to assuming that $D$ is a homogeneous Siegel domain (cf.~\cite[Theorem 2.3]{Xu} and~\cite[Proposition 2.2]{Murakami}).

\emph{We shall no longer make use of the symbols $A$, $H$, $e_j$, $e'_j$ ($j=1,\dots, r$), and $\Psi$ with the above meaning, unless explicitly stated.}

\begin{lem}\label{lem:70}
There is $\vect b\in \R^r$ such that $\Delta^{-\vect b}(t)=\abs{\det_{\R}(g)}$ for every $t\in T_+$ and for every $g\in GL(E)$ such that $t\cdot \Phi=\Phi\circ (g\times g)$.
\end{lem}

\begin{proof}
 Observe first that the group $U(\Phi)$ of elements of $GL(E)$ which preserve $\Phi$ is closed in $GL(E)$ and contained in the unitary group associated with the scalar product $\langle (e_{\Omega'})_\C, \Phi\rangle$ on $E$. 
 Hence, $U(\Phi)$ is compact, so that $\abs{\det g}=1$ for every $g\in U(\Phi)$. 
 Therefore, we may define a function $\det_\Phi\colon T_+\to \R_+^*$ so that $\det_\Phi(t)=\abs{\det_\R (g)}$ for every $t\in T_+$ and for every  $g\in GL(E)$ such that $t\cdot \Phi=\Phi\circ (g\times g)$. It is easily seen that $\det_\Phi$ is a homomorphism of $T_+$ into $\R_+^*$, so that $\det_\Phi=\Delta^{-\vect{b}}$ for some $\vect{b}\in \R^r$. 
\end{proof}

\begin{deff}\label{17}
We define $\vect{b}$ as in Lemma~\ref{lem:70}.
\end{deff}

We conclude this section with some examples and remarks.

\begin{ex}\label{ex:1}
	If $D$ is an irreducible symmetric Siegel domain of type II, then there are $a,b\in \N$ such that $a=m_{j,k}$ for every $j,k\in \Set{1,\dots, r}$ such that $j\neq k$, and such that $\vect b=-b\vect{1}_r$. In addition, $D$ is uniquely determined, up to a biholomorphism, by the three parameters $a$, $b$, and $r$. 
	
	Then, 
	\[
	\vect m= ( a(r-j) )_{j=1,\dots, r}, \qquad\vect{m'}= (a (j-1))_{j=1,\dots, r},
	\]
	and 
	\[
	\vect d=-( 1+a(r-1)/2 ) \vect 1_r.
	\]
	In particular, $m=-\sum_{j=1}^r d_j= r(1+a(r-1)/2)$.
	Observe that, in this situation, $\Delta^{\vect 1_r}$ and its powers play a central role in the theory of such domains.   For this reason, it is quite common, in the literature, to consider only irreducible symmetric Siegel domains of type II (or simply tube domains), and to consider only the spaces $A^{p,q}_{s\vect 1_r}(D)$, with the notation of Chapter~\ref{sec:6}.
	
	In addition to that, irreducible symmetric Siegel domains (or, more precisely, their bounded counterparts) have been completely classified by \'E.\ Cartan in~\cite{Cartan}. See also~\cite{Arazy} for a brief exposition of the classification.
\end{ex}

\begin{ex}\label{ex:2}
	If $r=1$, then $D$ is necessarily irreducible and symmetric, and is biholomorphically equivalent to the unit ball in $E\times F_\C$. In this case, $\vect d=-1$ and $\vect b=-n$. In addition, $\Nc$ is isomorphic to $\R$ if $n=0$, and to the $(2n+1)$-dimensional Heisenberg group $\Hd_n$  if $n>0$.
\end{ex}

\begin{oss}\label{oss:tube-domain}
We observe that $D$ is a tube domain, that is, $E=\Set{0}$, if and only if $\vect b=\vect 0$.
\end{oss}

\section{Gamma Functions}\label{sec:3:1}

In this section we introduce our notation for generalized power functions.
Observe that, while the definition of generalized power functions on $\Omega$ is quite uniform (up to the notation),  at least two ways to define generalized power functions on $\Omega'$ exist in the literature.
On the one hand, in e.g.~\cite{Gindikin,BekolleKagou}, the generalized power functions on $\Omega'$ are defined identifying $\Omega'$ with the cone $C(A')$, with the notation of Section~\ref{sec:2}, thus reversing the order of the variables. 
On the other hand, in e.g.~\cite{NanaTrojan}, the generalized power
functions on $\Omega'$ are defined by means of the canonical right
action of $T_+$ on $\Omega'$ in complete analogy with  the case of
$\Omega$.  
We shall follow this latter convention. Thus, when comparing our results with those presented in the literature, it will sometimes be necessary to reverse the order of the variables. 

We then identify the invariant measures on $T_+$, $\Omega$, and $\Omega'$, as well as the modular function on $T_+$.
Then, we report some basic results concerning the gamma  and beta functions on the cones $\Omega$ and $\Omega'$ introduced in~\cite{Gindikin}. 

We conclude this section with some basic results on the generalized Riemann--Liouville operators on $\Omega$ and $\Omega'$, following~\cite{Gindikin}. 

\medskip

\begin{deff}\label{def:3.5}
For every $\vect{s}\in \C^r$, define the generalized power functions $\Delta_\Omega^{\vect{s}}$ and $\Delta_{\Omega'}^{\vect{s}}$ on $\Omega$ and $\Omega'$, respectively, so that\index{Generalized power functions}
\[
\Delta_\Omega^{\vect{s}}(e_\Omega)=1, \qquad \Delta_{\Omega'}^{\vect{s}}(e_{\Omega'})=1,
\]
\[
\Delta_\Omega^{\vect{s}}(t\cdot h)=\Delta^{\vect{s}}(t)\Delta_\Omega^{\vect{s}}(h), \qquad \text{and} \qquad \Delta_{\Omega'}^{\vect{s}}( \lambda\cdot t)=\Delta^{\vect{s}}(t)\Delta_{\Omega'}^{\vect{s}}(\lambda)
\]
for every $t\in T_+$, for every $h\in \Omega$, and for every $\lambda\in \Omega'$. Recall that $\Delta^{\vect s}$ is the character of $T_+$ defined as in Lemma~\ref{lem:69}.\label{19} 

 With $\vect{m}$ and $\vect{m}'$ as in Definition~\ref{16-bis}, 
we set  \label{20} 
\[
\vect{d}\coloneqq -\left(\vect{1}_r+\half \vect{m}+\half \vect{m}'\right),
\]
 and 
\[
\nu_\Omega\coloneqq  \Delta_\Omega^{\vect{d}}\cdot \Hc^m \qquad \text{and} \qquad \nu_{\Omega'} \coloneqq \Delta^{\vect{d}}_{\Omega'}\cdot \Hc^m.
\] 
\end{deff}

In the following lemma, we define the Hausdorff measure $\Hc^m $ on $T_+$ with reference to the formalism of  Section~\ref{sec:2},  using the identifications of Section~\ref{sec:4:1}.

\begin{lem}\label{lem:19}
The Hausdorff measure $\Hc^m$ on $T_+$ is relatively invariant. Its left and right multipliers are 
\[
\Delta^{(\vect 1_r+\vect m)/2} \qquad \text{and} \qquad\Delta^{(\vect 1_r+\vect{m'})/2},
\]
respectively. In particular, the modular function on $T_+$ is 
\[
\Delta^{(\vect{m'}-\vect{m})/2}.
\]
In addition, 
\[
(\iota_\Omega)_*( \Delta^{-(\vect 1_r+\vect m)/2}\cdot \Hc^m )= 2^{-(m+r)/2}\Delta_\Omega^{\vect d}\cdot \Hc^m
\]
and
\[
(\iota_{\Omega'})_*(\Delta^{-(\vect 1_r+\vect{m'})/2}\cdot \Hc^m)=2^{-(m+r)/2}\Delta_{\Omega'}^{\vect d}\cdot \Hc^m,
\]
where
\[
\iota_\Omega\colon T_+ \ni t \mapsto t\cdot e_\Omega\in \Omega \qquad \text{and} \qquad \iota_{\Omega'}\colon T_+ \ni t \mapsto  e_{\Omega'}\cdot t\in \Omega'.
\]
\end{lem}

In particular, $\nu_\Omega$ and $\nu_{\Omega'}$ are two $T_+$-invariant measures on $\Omega$ and $\Omega'$, respectively.

\begin{proof}
Observe  first that the product on $T_+$ is bilinear, so that $\Hc^m$ is relatively invariant.
Then, take $\varphi \in C_c(T_+)$ and $t=\sum_{j=1}^r \alpha_j e_j$ for some $\alpha\in (\R_+^*)^r$, and observe that
\[
\begin{split}
\int_{T_+} \varphi(t t')\,\dd t'= \int_{T_+} \varphi\left(\sum_{j\meg k} \alpha_j t'_{j,k}\right)\,\dd t'= \alpha^{-(\vect 1_r+ \vect m)}\int_{T_+}\varphi(t')\,\dd t',
\end{split}
\]
so that $\Delta^{(\vect 1_r+\vect m)/2}$ is the left multiplier of $\Hc^m$. Analogously, $\Delta^{(\vect 1_r+\vect{m'})/2}$ is the right multiplier of $\Hc^m$, so that the modular function of $T_+$ is $\Delta^{(\vect{m'}-\vect{m})/2}$. 

Then, denote by $J$ the Jacobian of the mapping $\iota_\Omega\colon T_+\ni t\mapsto t t^*\in \Omega$. 
Observe that 
\[
\iota_\Omega'(t) v= v t^*+t v^* 
\]
for every $t\in T_+$ and for every $v\in T$ (identifying $T$ with the tangent space of $T_+$ at $t$). Define $P\colon A\to T$ by $P(a)_{j,k}\coloneqq a_{j,k}$ for $j\meg k$ and $P(a)_{j,k}=0$ for $j>k$, and observe that, for every $v\in T$ and for every $w\in H$,
\[
\langle \iota_\Omega'(t) v  ,  w\rangle= 2\langle v ,   w t\rangle=\langle v  , 2P(w t )\rangle,
\] 
so that $\iota_\Omega'(t)^* w=2 P( w t)$. 
Therefore, 
\[
\iota_\Omega'(t)^* \iota_\Omega'(t) v=2P((v t^*) t+(t v^*) t  )
\]
and, in particular, $\iota_\Omega'(e)^* \iota_\Omega'(e)=2I+2\sum_{j=1}^r e'_j$, whence 
\[
J(e)=2^{(m+r)/2}.
\] 
Now, the area formula shows that 
\[
(\iota_\Omega)_*(J\cdot \Hc^m)=\chi_\Omega\cdot\Hc^m.
\]
In addition, since $T_+$ acts on $F$ by linear mappings, the measure $\Hc^m$ on $\Omega$ is relatively invariant, so that 
\[
J(t)=J(t e)=\Delta^{\vect s}(t) J(e)
\]
for some $\vect s\in \R^r$ and for every $t\in T_+$. 
Now, to determine $\vect s$ it suffice to compute $J(t)$ for $t=\sum_{j=1}^r \alpha_j e_j$, $\alpha\in (\R_+^*)^r$. 
However, in this case 
\[
2P((v t^*) t+(t v^*) t)_{j,k}= 2(1+\delta_{j,k}) v_{j,k} \alpha_k^2
\]
for every $j\meg k$, so that 
\[
J(t)=\Delta^{(\vect 1_r+\vect {m'})/2}(\alpha) J(e).
\] 
Thus, 
\[
(\iota_\Omega)_*(\Delta^{-(\vect 1_r+\vect m)/2} \cdot \Hc^m)= 2^{-(m+r)/2} \Delta_\Omega^{\vect d  }\cdot \Hc^m.
\]
The other assertion is proved similarly
\end{proof}

Recall that $\Lc$ denotes the Laplace transform.

\begin{prop}\label{prop:58}
The measure $\Delta_\Omega^{\vect{s}}\cdot \nu_\Omega$ induces a Radon measure on $\overline \Omega$ if and only if $\Re\vect{s}\in \frac{1}{2}{\vect{m}}+(\R_+^*)^r$. In this case,
\[
\Lc(\Delta^{\vect s}_\Omega\cdot \nu_\Omega)= (2\pi)^{\frac{m-r}{2}} \prod_{j=1}^r \Gamma\left( s_j-\frac{m_j}{2} \right) \Delta^{-\vect s}_{\Omega'}.
\]

Analogously, $\Delta_{\Omega'}^{\vect{s}}\cdot \nu_{\Omega'}$ induces a Radon measure on $\overline{\Omega'}$ if and only if $\Re\vect{s}\in \frac{1}{2}\vect{m}'+(\R_+^*)^r$. In this case,
\[
\Lc(\Delta^{\vect s}_{\Omega'}\cdot \nu_{\Omega'})=(2\pi)^{\frac{m-r}{2}} \prod_{j=1}^r \Gamma\left( s_j-\frac{m'_j}{2} \right)\Delta^{-\vect s}_\Omega.
\]
\end{prop}

Cf.~\cite[Theorem 2.1]{Gindikin} for the original results, where a different normalization of the scalar product is used.

\begin{proof}
Observe first that, by Lemma~\ref{lem:19},
\[
\Lc(\Delta^{\vect s}_\Omega\cdot \nu_\Omega)(e_{\Omega'})= 2^{(m+r)/2}\int_{T_+} \ee^{-\tr(t t^*)} \Delta^{\vect s-(\vect{1}_r+\vect m)/2}(t)\,\dd t, 
\]
in the sense that one side of the equality is defined if and only if the other one is, and then equality holds. Now, $\tr(t t^*)=\abs{t}^2=\sum_{j\meg k} \abs{t_{j,k}}^2$, while  $ \Delta_\Omega^{\vect s-(\vect{1}_r+\vect m)/2}(t)=\prod_{j=1}^r t_{j,j}^{2 \vect s-\vect 1 _r-\vect m} $, so that the preceding integral is defined if and only if $\Re\vect s-\frac 1 2\vect m\in (\R_+^*)^r $ and, in this case,
\[
2^r\int_{T_+} \ee^{-\tr(t t^*)} \Delta_\Omega^{\vect s-(\vect{1}_r+\vect m)/2}(t)\,\dd t= \pi^{(m-r)/2} \prod_{j=1}^r \Gamma(s_j-m_j/2)  ,
\]
since $2\int_0^\infty \ee^{-x^2} x^{2 s_j-m_j}\,\frac{\dd x}{x}= \Gamma(s_j-m_j/2)$ ($j=1,\dots, r$) and $\int_\R \ee^{-x^2}\,\dd x=\sqrt \pi$. 
The fact that, when $\Re\vect{s}\in \frac{1}{2}{\vect{m}}+(\R_+^*)^r$,
\[
\Lc(\Delta^{\vect s}_\Omega\cdot \nu_\Omega)= (2 \pi)^{\frac{m-r}{2}} \prod_{j=1}^r \Gamma\left( s_j-\frac{m_j}{2} \right) \Delta^{-\vect s}_{\Omega'}
\]
follows by homogeneity. The second assertion is proved similarly.
\end{proof}

\begin{deff}\label{21}\index{Gamma function!of a homogeneous cone}
We define, for $\Re\vect{s}\in \frac{1}{2}{\vect{m}}+(\R_+^*)^r$
\[
\Gamma_\Omega(\vect{s})\coloneqq \int_{\Omega} \ee^{-\langle e_{\Omega'},h\rangle}\Delta^{\vect{s}}_\Omega(h)\,\dd \nu_\Omega(h)=(2\pi)^{\frac{m-r}{2}} \prod_{j=1}^r \Gamma\left( s_j-\frac{m_j}{2} \right) 
\]
and, for $\Re\vect{s}\in \frac{1}{2}\vect{m}'+(\R_+^*)^r$,
\[
\Gamma_{\Omega'}(\vect{s})\coloneqq \int_{\Omega'} \ee^{-\langle \lambda, e_{\Omega}\rangle}\Delta^{\vect{s}}_{\Omega'}(\lambda)\,\dd \nu_{\Omega'}(\lambda)=(2\pi)^{\frac{m-r}{2}} \prod_{j=1}^r \Gamma\left( s_j-\frac{m'_j}{2} \right).
\]
\end{deff}

As a corollary to Proposition~\ref{prop:58} we now present some elementary results on beta functions.

\begin{cor}\label{cor:9}
Take $\vect s,\vect{s'}\in \C^r$ and $h\in \Omega$. Then, the function $h'\mapsto \Delta_\Omega^{\vect s+\vect{d}}(h-h')\Delta_\Omega^{\vect{s'}}(h')$ is $\nu_\Omega$-integrable on $\Omega \cap (h-\Omega)$ if and only if $\Re\vect s, \Re\vect{s'}\in \frac 1 2 \vect m+(\R_+^*)^r$, in which case
\[
\int_{\Omega\cap (h-\Omega)} \Delta_\Omega^{\vect s+\vect{d}}(h-h')\Delta_\Omega^{\vect{s'}}(h')\,\dd \nu_\Omega(h')= \frac{\Gamma_\Omega(\vect s)\Gamma_\Omega(\vect{s'})}{\Gamma_\Omega(\vect s+ \vect{s'})} \Delta_\Omega^{\vect s+\vect{s'}+\vect d}(h).
\] 
\end{cor}

An analogous statement holds for $\Omega'$. 
Cf.~\cite[Theorem 2.2]{Gindikin} for the original result.

\begin{proof}
The first assertion follows easily from Proposition~\ref{prop:58}. For what concerns the second assertion,  assume that $\vect s, \vect{s'}\in \frac 1 2 \vect m+(\R_+^*)^r$ and observe that the above equality is equivalent to  $(\Delta_\Omega^{\vect s}\cdot \nu_\Omega)*(\Delta_\Omega^{\vect{s'}}\cdot \nu_\Omega) =\frac{\Gamma_\Omega(\vect s)\Gamma_\Omega(\vect{s'})}{\Gamma_\Omega(\vect s+ \vect{s'})} \Delta_\Omega^{\vect s+\vect{s'}}\cdot \nu_\Omega$. Since the Laplace transforms of both sides of the asserted equality are equal on $\Omega'+i F$ by Proposition~\ref{prop:58}, the assertion follows.
\end{proof}

\begin{cor}\label{cor:10}
Take $\vect s,\vect{s'}\in \C^r$ and $h\in \Omega$. Then, the function $h'\mapsto \Delta_\Omega^{\vect s}(h+h')\Delta_\Omega^{\vect{s'}}(h')$ is $\nu_\Omega$-integrable on $\Omega$ if and only if $\Re\vect{s'}\in \frac 1 2 \vect m+(\R_+^*)^r$ and $\Re \vect s+\Re \vect{s'}\in - \frac 1 2 \vect{m'}-(\R_+^*)^r$, in which case
\[
\int_{\Omega} \Delta_\Omega^{\vect s}(h+h')\Delta_\Omega^{\vect{s'}}(h')\,\dd \nu_\Omega(h')= \frac{\Gamma_\Omega(\vect s')\Gamma_{\Omega'}(-\vect s-\vect{s'})}{\Gamma_{\Omega'}(-\vect s)} \Delta_\Omega^{\vect s+\vect{s'}}(h).
\] 
\end{cor}

An analogous statement holds for $\Omega'$.
Cf.~\cite[Proposition 2.6]{Gindikin} for the original result (which is flawed by an incorrect computation of the modular function on $T_+$). See also~\cite[Lemma 4.19]{NanaTrojan} for an alternative proof, under the additional assumption $\vect{s}\in -\frac 1 2 \vect{m'}-(\R_+^*)^r$.

\begin{proof}
Define a mapping $I\colon h\mapsto h^{-1}$ on $\Omega$ so that $(t\cdot e_\Omega)^{-1}\coloneqq t^{-1}\cdot e_\Omega$ for every $t\in T_+$. Then, it is easily verified that $I(e_\Omega+\Omega)=\Omega \cap (e_\Omega-\Omega)$. In addition, Lemma~\ref{lem:19} shows that $I_*(\nu_\Omega)= \Delta_\Omega^{(\vect m-\vect{m'})/2}\cdot \nu_\Omega$. Therefore,
\[
\begin{split}
&\int_{\Omega} \Delta_\Omega^{\vect s}(e_\Omega+h')\Delta_\Omega^{\vect{s'}}(h')\,\dd \nu_\Omega(h')= \int_{e_\Omega+\Omega} \Delta_\Omega^{\vect s-\vect d}(h')\Delta_\Omega^{\vect{s'}+\vect{d}}(h'-e_\Omega)\,\dd \nu_\Omega(h')\\
&\qquad=\int_{\Omega\cap (e_\Omega-\Omega)} \Delta_\Omega^{\vect s-\vect d-(\vect m-\vect{m'})/2}(h'^{-1})\Delta_\Omega^{\vect{s'}+\vect{d}}(h'^{-1}-e_\Omega)\,\dd \nu_\Omega(h')\\
&\qquad=\int_{\Omega\cap (e_\Omega-\Omega)} \Delta_\Omega^{\vect d-\vect s+(\vect m-\vect{m'})/2-\vect{s'}-\vect d  }(h')\Delta_\Omega^{\vect{s'}+\vect d}(e_\Omega-h')\,\dd \nu_\Omega(h')\\
&\qquad=\int_{\Omega\cap (e_\Omega-\Omega)} \Delta_\Omega^{\vect {s'}+\vect d}(e_\Omega-h')\Delta_\Omega^{-\vect s-\vect{s'}+(\vect m-\vect{m'})/2}(h')\,\dd \nu_\Omega(h'),
\end{split}
\]
where each integral is defined if and only if the other ones are, and then all the equalities hold. Therefore, the assertion follows from Proposition~\ref{prop:58}, Corollary~\ref{cor:9}, and homogeneity.
\end{proof}

It will be relevant, e.g.~in the study of Riemann--Liouville operators, to know which generalized power functions on $\Omega$ and $\Omega'$ are actually polynomials. Cf.~\cite{Ishi} for a more detailed study of this subject.

\begin{deff}\label{22}
We denote by $\N_\Omega$ and $\N_{\Omega'}$ the sets of $\vect s\in \R^r$ such that $\Delta_\Omega^{\vect s}$ and $\Delta_{\Omega'}^{\vect s}$ are polynomials, respectively.

To simplify the notation, if $\vect s\in \C^r$ and $a\in \R_+^*$, we define
\[
a^{\vect s}\coloneqq a^{\sum_{j=1}^r s_j}=\ee^{(\sum_{j=1}^r s_j)\log a}.
\]
By an abuse of language, we also define 
\[
i^z\coloneqq \ee^{\pi i z/2}
\]
for every $z\in \C$. We define $i^{\vect s}$ accordingly.
\end{deff}

In particular, $a^{-\vect d}=a^m$ and $a^{-\vect b}=a^n$.

\begin{lem}\label{lem:71}
$\N_\Omega$ and $\N_{\Omega'}$ are sub-monoids of $\N^r$ and are isomorphic to $\N^r$.
\end{lem}

In particular, $\N_\Omega$ and $\N_{\Omega'}$ are cofinal in $\R^r$ for $\meg$.

\begin{proof}
Clearly, $\vect 0\in \N_\Omega\cap \N_{\Omega'}$ and both $\N_\Omega$ and $\N_{\Omega'}$ are stable under the sum. Thus, $\N_\Omega$ and $\N_{\Omega'}$ are sub-monoids of $\R_+^r$. To see that they are contained in $\N^r$, keep the notation of Chapter~\ref{sec:2} (identifying $F$ with $H$), and observe that 
\[
\Delta_\Omega^{\vect s}\left(\sum_{j=1}^r \alpha_j e_j\right)=\Delta_\Omega^{\vect s}\left(\left(\sum_{j=1}^r \sqrt{\alpha_j} e_j\right)\cdot e_\Omega\right)=\prod_{j=1}^r  \alpha_j^{s_j}
\]
for every $\alpha\in (\R_+^*)^r$. 
Finally, $\N_\Omega$ is isomorphic to $\N^r$ thanks to~\cite[Theorem 2.2]{Ishi}. The assertions concerning $\N_{\Omega'}$ are proved analogously.
\end{proof}

Combining the fact that polynomials are entire functions with Proposition~\ref{prop:58}, we are then able to show that the generalized power functions on $\Omega$ and $\Omega'$ always extend to holomorphic functions on $\Omega+ i F$ and $\Omega'+ i F'$, respectively.
Notice that, if $\vect s\in \Z^r$, then $\Delta_\Omega^{\vect s}$ and $\Delta_{\Omega'}^{\vect s}$ are rational functions, but may have singularities on $\partial \Omega+ i F$ and $\partial \Omega'+ i F'$, respectively, in general.
We denote by $\Sc'(F)$ and $\Sc'(F')$ the spaces of tempered distributions on $F$ and $F'$, respectively.

\begin{cor}\label{cor:20}
For every $\vect s\in \C^r$, the functions $\Delta^{\vect s}_\Omega$ and $\Delta^{\vect s}_{\Omega'}$ extend to holomorphic functions on $\Omega+i F$ and $\Omega'+i F'$. In addition, the mappings
\[
\C^r\ni\vect s\mapsto \Delta^{\vect s}_\Omega(h+i\,\cdot\,)\in \Sc'(F)
\]
and
\[
\C^r\ni\vect s\mapsto \Delta^{\vect s}_{\Omega'}(\lambda+i\,\cdot\,)\in \Sc'(F')
\]
are holomorphic for every $h\in \Omega$ and for every $\lambda\in \Omega'$.
\end{cor}

\begin{proof}
Take $\vect s\in \C^r$. By Lemma~\ref{lem:71}, we may take $\vect{s'}\in \N_\Omega$ such that $\vect{s'}-\Re \vect{s}\in \frac1 2 \vect{m'}+(\R_+^*)^r$, so that Proposition~\ref{prop:58} implies that
\[
\Delta_\Omega^{\vect s}=\Delta_\Omega^{\vect{s'}} \frac{1}{\Gamma_{\Omega'}(\vect {s'}-\vect s)} \Lc(\Delta_{\Omega'}^{\vect{s'}-\vect s}\cdot \nu_{\Omega'}),
\]
so that $\Delta_\Omega^{\vect s}$ extends to a holomorphic function on $\Omega+ i F$. In addition, it is clear that the mapping
\[
\vect s\mapsto \Delta_\Omega^{\vect s}(h+i\,\cdot\,)= \Delta_\Omega^{\vect{s'}}(h+i\,\cdot\,) \frac{1}{\Gamma_{\Omega'}(\vect {s'}-\vect s)} \Fc_{F'}(\ee^{-\langle\,\cdot\,,h\rangle}\Delta_{\Omega'}^{\vect{s'}-\vect s}\cdot \nu_{\Omega'})\in \Sc'(F)
\]
is holomorphic on the set of $\vect s\in \C^r$ such that $\Re \vect{s}-\vect{s'}\in -\frac1 2 \vect{m'}-(\R_+^*)^r$. By the arbitrariness of $\vect{s'}$, the assertions concerning $\Omega$ follow. The remaining assertions are proved similarly.
\end{proof}

Recall that, when $\Omega=\R_+^*$, the Riemann--Liouville potentials $I^s$ on $\R$ are defined so that
\[
I^s=\begin{cases}
\frac{1}{\Gamma(s)} (\,\cdot\,)^{s-1}\cdot\Hc^1 & \text{if $\Re s>0$}\\
\Fc^{-1}( (i \,\cdot\,)^{-s}  ) & \text{if $\Re s <1$}.
\end{cases}
\]
In the general case, an analogous definition by cases is no longer
possible, in full generality, since stronger conditions on $\vect s$
are needed  to ensure local integrability of $\Delta^{\vect
  s+\vect d}_\Omega$ and $\Delta^{-\vect s}_{\Omega'}$. For this
reason, it will be convenient to define the operators  $I^{\vect
  s}_\Omega$ by means of their Fourier transform. To do this, though,
we need to give a (distributional) meaning to $\lim\limits_{\lambda
  \to 0} \Delta^{-\vect s}_{\Omega'}(\lambda+ i \,\cdot\,)$. This is
the purpose of the following result.

\begin{lem}\label{lem:13}
Take $\vect s\in \C^r$. Then, 
\[
J^{\vect s}_\Omega \coloneqq\lim_{\substack{h\in \Omega\\ h\to 0}}
\Delta_{\Omega}^{\vect s}(h+i\,\cdot\,) \qquad \text{and} \qquad
J^{\vect s}_{\Omega'} \coloneqq\lim_{\substack{\lambda\in \Omega'\\ \lambda\to 0}}
\Delta_{\Omega}^{\vect s}(h+i\,\cdot\,)
\]
are well defined elements of $\Sc'(F)$ and $\Sc'(F')$, respectively. In addition, 
\[
J^{\vect s}_\Omega= i^{\vect s}\Delta_\Omega^{\vect s}  \qquad \text{and}\qquad J^{\vect s}_{\Omega'}= i^{\vect s} \Delta^{\vect s}_{\Omega'}
\]
on $\Omega$ and $\Omega'$, respectively, and the mappings
\[
\C^r\ni\vect s\mapsto J^{\vect s}_\Omega\in \Sc'(F) \qquad \text{and} \qquad \C^r\ni\vect s\mapsto J^{\vect s}_{\Omega'}\in \Sc'(F')
\]
are holomorphic.
\end{lem}

\begin{proof}
We prove only the assertions concerning $J^{\vect s}_\Omega$. If $-\Re \vect s\in \frac{1}{2}\vect{m'}+(\R_+^*)^r $, then Proposition~\ref{prop:58} shows that
\[
\begin{split}
\Delta_{\Omega}^{\vect s}(h+i\,\cdot\,)&=\frac{1}{\Gamma_{\Omega'}(-\vect s)} \Lc(\Delta_{\Omega'}^{-\vect s} \cdot \nu_{\Omega'})(h+i\,\cdot\,)=\frac{1}{\Gamma_{\Omega'}(-\vect s)} \Fc_{F'}(\ee^{-\langle\,\cdot\,,h\rangle}\Delta_{\Omega'}^{-\vect s} \cdot \nu_{\Omega'}),
\end{split}
\]
so that $J_\Omega^{\vect s}=\frac{1}{\Gamma_{\Omega'}(-\vect s)} \Fc_{F'}(\Delta_{\Omega'}^{-\vect s} \cdot \nu_{\Omega'})$. 

Now, let $U$ be the set of $(z,h)\in \C\times (\Omega+i F)$ such that $z\not \in \R_-$ and $z h\in \Omega+i F$, and let us prove that $U$ is a connected open subset of $\C\times (\Omega +i F)$. Indeed, $U$ is clearly open. Further, observe that for every $z\in \C\setminus \R_-$ the set $U_z$ of $h\in F_\C$ such that $(z,h)\in U$ is the convex set $(\Omega+i F)\cap z^{-1}(\Omega+i F) $, and is therefore connected. 
In addition, take $z_0, z_1\in \C\setminus \R_-$ such  that $\Im z_0 \Im z_1>0$ if either $\Re z_0\meg0$ or $\Re z_1\meg0$. 
Let us prove that there are $h\in \Omega$ and $h'\in F$ such that 
\[
[(1-t)\Re z_0+ t\Re z_1] h-[(1-t)\Im z_0+ t \Im z_1] h'\in\Omega 
\]
for every $t\in [0,1]$. 
Indeed, if $\Re z_0,\Re z_1>0$, it suffices to take $h\in \Omega$ and $h'=0$. If, otherwise, $\Re z_0\meg 0$ or $\Re z_1\meg 0$, it suffices to take $h\in \Omega$ and $h'= R h$ for some $R$ such that $\Re z_0-R \Im z_0, \Re z_1-R \Im z_1>0$, which is possible since $\Im z_0 \Im z_1>0$.
Thus, 
\[
[z_0(h+i h'), z_1(h+i h')]\subseteq\Omega+i F. 
\]
By the arbitrariness of $z_0$ and $z_1$, it follows that $U$ is connected.
Now, observe that the mapping
\[
\Psi\colon U\ni (z,h)\mapsto \Delta^{\vect s}_\Omega(z h)-z^{\vect s}\Delta^{\vect s}_\Omega(h)\in \C
\]
is holomorphic, where $\C\setminus \R_-\ni z \mapsto z^{\vect s}\in \C$ is the unique holomorphic function which equals $x \mapsto \ee^{\vect s\log x}$ on $\R_+^*$. Observe that $\Psi$ vanishes on $\R_+^*\times \Omega$ by the homogeneity of $\Delta^{-\vect s}_{\Omega'}$, so that it vanishes on the whole of $U$ by holomorphy. Thus,
\[
\Delta^{\vect s}_{\Omega}(h+i h')= i^{\vect s} \Delta^{\vect s}_{\Omega}(h'-i h)
\]
for every $h,h'\in \Omega$, so that $J^{\vect s}_\Omega= i^{\vect s} \Delta^{\vect s}_{\Omega}$ on $\Omega$.

In the general case, take $\vect{s'}\in \N_\Omega$ such that $\vect{s'}-\Re\vect{s}\in \frac 1 2\vect{m'}+(\R_+^*)^r $. 
Then, clearly
\[
J^{\vect s}_\Omega= \Delta_\Omega^{\vect s'}(i\,\cdot\,) J^{\vect s-\vect{s'}}_\Omega,
\] 
so that $J^{\vect s}_\Omega= i^{\vect s} \Delta^{\vect s}_{\Omega}$ on $\Omega$.

For what concerns holomorphy, fix some $\vect{s'}\in \N_\Omega$. Then, the mapping $\vect{s}\mapsto \Delta_\Omega^{\vect{s'}}(i\,\cdot\,)J^{\vect s-\vect{s'}}_\Omega\in \Sc'(F)$ is holomorphic on the set of $\vect s$ such that $\Re\vect{s}\in \vect{s'}-\frac 1 2\vect{m'}-(\R_+^*)^r$, so that the assertion follows from Lemma~\ref{lem:71}.
\end{proof}

\begin{deff}\label{23}
For every $\vect s\in \C^r$, we define Riemann--Liouville potentials\index{Riemann--Liouville potentials} on $F$ and $F'$ as follows:
\[
\begin{aligned}
I_\Omega^{\vect{s}}&=\lim_{\substack{\lambda\in \Omega'\\\lambda\to 0}}\Fc_F^{-1}( \Delta_{\Omega'}^{-\vect{s}}(\lambda+i\,\cdot\,)) \\ 
I_{\Omega'}^{\vect{s}}&= 
\lim_{\substack{h\in \Omega\\h\to 0}}\Fc^{-1}_{F'}(\Delta_{\Omega}^{-\vect{s}}(h+i \,\cdot\,)). 
\end{aligned}
\]
\end{deff}

We collect in the following result some elementary properties of the Riemann--Liouville potentials previously defined.
\begin{prop}\label{prop:38}
Take $\vect s,\vect{s'}\in \C^r$. Then, the following         properties hold:
\begin{enumerate}
\item[\em(1)] $I_\Omega^{\vect{s}}=\frac{1}{\Gamma_\Omega(\vect{s})} \Delta_\Omega^{\vect{s}}\cdot \nu_\Omega$ if $\Re \vect s \in\frac{1}{2} \vect{m}+(\R_+^*)^r$, while $I_{\Omega'}^{\vect{s}}=\frac{1}{\Gamma_{\Omega'}(\vect{s})}\Delta^{\vect s}_{\Omega'}\cdot \nu_{\Omega'}$ if $\Re \vect s \in \frac{1}{2}\vect{m}'+(\R_+^*)^r$;
\item[\em(2)] $\Lc I_\Omega^{\vect{s}}= \Delta_{\Omega'}^{-\vect s}$ on $\Omega'+i F'$ and $\Lc I_{\Omega'}^{\vect{s}}= \Delta_{\Omega}^{-\vect s}$ on $\Omega+i F$;
\item[\em(3)] $I_{\Omega}^{\vect s}$ is supported in $\overline \Omega$ and   $I^{\vect s}_{\Omega'}$ is supported in $\overline {\Omega'}$;
\item[\em(4)] $I_\Omega^{\vect{s}}$ is supported at $0$ if and only if $\vect{s}\in -\N_{\Omega'}$, while  $I_{\Omega'}^{\vect s}$ is supported at $0$ if and only if $\vect{s}\in -\N_\Omega$;
\item[\em(5)] $I_\Omega^{\vect{s}}*I_\Omega^{\vect{s'}}= I_\Omega^{\vect s+ \vect{s'}}$ and  $I_{\Omega'}^{\vect{s}}*I_{\Omega'}^{\vect{s'}}= I_{\Omega'}^{\vect s+ \vect{s'}}$.
\end{enumerate}
\end{prop}

Since clearly $I^{\vect 0}_\Omega=\delta_0$ and $I^{\vect{0}}_{\Omega'}=\delta_0$, it then follows that $I^{\vect{s}}_\Omega$ and $I^{\vect{s}}_{\Omega'}$ are fundamental solutions of the differential operators $I^{-\vect s}_\Omega$ and $I^{-\vect s}_{\Omega'}$ for every $\vect s\in \N_{\Omega'}$ and for every $\vect s\in \N_\Omega$, respectively.

\begin{proof}
(1) This follows from Proposition~\ref{prop:58}.

(2) This follows from~(1) and Proposition~\ref{prop:58} when $\Re\vect s$ is sufficiently large. The general case then follows by holomorphy, thanks to Lemma~\ref{lem:13}.

(3) This follows from~(1) when $\Re \vect s \in\frac{1}{2} \vect{m}+(\R_+^*)^r$ and $\Re \vect s \in \frac{1}{2}\vect{m}'+(\R_+^*)^r$, respectively, and by holomorphy in the general case.

(4) This follows from the definition of $\N_\Omega$ and $\N_{\Omega'}$. 

(5) This follows e.g.\ from Corollary~\ref{cor:9} when $\Re \vect s,\Re \vect{s'}\in \frac{1}{2}\vect{m}+(\R_+^*)^r$ and $\Re \vect s,\Re \vect{s'}\in \frac{1}{2}\vect{m}'+(\R_+^*)^r$, respectively. The general case follows by holomorphy, since the space 
\[
\Dc'_\Omega(F)\coloneqq \Set{T\in \Dc'(F)\colon \Supp{T}\subseteq \overline\Omega}
\]
(where $\Dc'(F)$ denotes the space of distributions on $F$)is a commutative and associative convolution algebra (cf.~\cite[Theorems XIII, XIV, and XIV bis]{Schwartz}).
\end{proof}

When convolving smooth functions with $I^{\vect s}_{\Omega}$, we shall often restrict ourselves to the case $\vect s\in -\N_{\Omega'}$, so that $I^{\vect s}_\Omega$ is a differential operator and no issues occur. 
With a deeper analysis, one may sometimes extend such results to more general values of $\vect s$.

\begin{prop}\label{prop:39}
For every $\vect s\in \C^r$, for every $\vect{s'}\in \N_{\Omega'}$, and for every $h\in \Omega+i F$,
\[
(\Delta_\Omega^{\vect s}*I_\Omega^{-\vect{s'}})(h)=\left(\vect s +
  \half   \vect{m'} \right)_{\vect s'} \Delta_\Omega^{\vect s-\vect{s'}}(h),
\]
where $\left(\vect s +\frac 1 2 \vect{m'} \right)_{\vect s'}=\prod_{j=1}^r (s_j+\frac 1 2 m'_j)\cdots(s_j-s_j'+\frac 1 2 m'_j+1)$.
\end{prop}

\begin{proof}
By holomorphy, it will suffice to prove the assertion for $h\in \Omega$ and $\Re\vect s\in \vect{s'}+\vect d +\frac 1 2 \vect m+(\R_+^*)^r$. Then, Proposition~\ref{prop:38} shows that
\[
\begin{split}
\Delta_\Omega^{\vect s}*I_\Omega^{-\vect{s'}}&=\Gamma_\Omega(\vect s-\vect d) I_\Omega^{\vect s-\vect d} *I_\Omega^{-\vect{s'}}\\
&=\Gamma_\Omega(\vect s-\vect d) I_\Omega^{\vect s-\vect{s'}-\vect d}\\
&=\frac{\Gamma_\Omega(\vect s-\vect d)}{\Gamma_\Omega(\vect s-\vect{s'}-\vect d)}\Delta_\Omega^{\vect s-\vect{s'}},
\end{split}
\] 
whence the result, since $\frac{\Gamma_\Omega(\vect s-\vect d)}{\Gamma_\Omega(\vect s-\vect{s'}-\vect d)}=\left(\vect s +\frac 1 2 \vect{m'} \right)_{\vect s'}$.
\end{proof}

We now describe further properties of the generalized power functions $\Delta_\Omega^{-\vect b}$ and $\Delta_{\Omega'}^{-\vect b}$ in terms of the preceding results.

 As we observed in Example~\ref{ex:1}, $\vect b, \vect d\in \R_- \vect 1_r$ when $D$ is irreducible and symmetric. As a consequence,  some objects in the theory of Bergman spaces simplify in that case (cf., e.g.,~\cite{Arazy}). This need no longer be  the case if $D$ is either not irreducible or not symmetric.

\begin{prop}\label{prop:59}
The following properties hold:
\begin{itemize}
\item $-\vect b\in \N_{\Omega'}$;

\item $\abs{\Pfaff(\lambda)}=\abs{\Pfaff(e_{\Omega'})}\Delta_{\Omega'}^{-\vect b}(\lambda)$ for every $\lambda\in \Omega'$;

\item $\Phi_*(\Hc^{2 n})= \frac{\pi^n}{\abs{\Pfaff(e_{\Omega'})}} I^{-\vect b}_\Omega$.
\end{itemize}
\end{prop}

\begin{proof}
Observe first that, with the notation of Section~\ref{sec:1:2},
\[
J_{\lambda\cdot t}=g^* J_\lambda g
\]
for every $\lambda\in \Omega'$, for every $t\in T_+$, and for every $g\in GL(E)$ such that $t\cdot \Phi=\Phi\circ (g\times g)$. Therefore, under the same assumptions,
\[
\abs{\Pfaff(\lambda\cdot t)}={\det}_\C(\abs{J_{\lambda\cdot t}})=\abs{{\det}_\C(g)}^2 {\det}_\C(\abs{J_{\lambda}})= \Delta^{-\vect b}(t) \abs{\Pfaff(\lambda)}
\]
by the definition of $\vect b$. Hence, $\abs{\Pfaff(\lambda)}= \abs{\Pfaff(e_{\Omega'})} \Delta_{\Omega'}^{-\vect b}(\lambda)$ for every $\lambda\in \Omega'$.

Next, observe that
\[
\abs{\Pfaff(e_{\Omega'})}\Delta_{\Omega'}^{-\vect b}(\lambda) =\abs{\Pfaff(\lambda)}=\abs{{\det}_\C(J_\lambda)} ={\det}_\C(-i J_\lambda)
\]
for every $\lambda\in \Omega'$. Since the mapping
\[
F'\setminus W\ni\lambda \mapsto {\det}_\C(-i J_\lambda)\in \C
\]
is clearly polynomial, it then follows that $-\vect b\in \N_{\Omega'}$.

Finally, observe that $\Phi_*(\Hc^{2 n})$ is a positive Radon measure on $F$ which is concentrated on $\Phi(E)\subseteq\overline\Omega$, since the positive quadratic form $\Phi$ is a proper mapping. In addition, a simple change of variables shows that
\[
\Lc(\Phi_*(\Hc^{2 n}))(\lambda)= \int_{E} \ee^{-\langle \lambda, \Phi(\zeta)\rangle}\,\dd \zeta= \frac{1}{\abs{\Pfaff(\lambda)}} \int_E \ee^{-\abs{\zeta}^2}\,\dd \zeta= \frac{\pi^n}{\abs{\Pfaff(\lambda)}}
\]
for every $\lambda\in \Omega'$, since $\langle \lambda, \Phi(\zeta)\rangle= \abs*{\abs{J_\lambda}^{1/2}\zeta}^2$  for every $\lambda\in \Omega'$ and for every $\zeta\in E$.
Now, observe that Proposition~\ref{prop:38} implies that
\[
(\Lc I^{-\vect b}_\Omega)(\lambda)=\Delta_{\Omega'}^{\vect b}(\lambda)
\]
for every $\lambda\in \Omega'$, so that
\[
\Phi_*(\Hc^{2 n})=  \frac{\pi^n}{\abs{\Pfaff(e_{\Omega'})}} I^{-\vect b}_\Omega
\]
by the uniqueness of the Laplace transform.
\end{proof}

\section{Evaluation of Some Norms}\label{sec:3:2}

This section contains several technical results which will be very useful in the following chapters. 
We mainly report or extend the results present in the literature to compute the $L^p$ norms of several functions associated with the generalized powers functions on $\Omega$ and $\Omega'$.

While the results concerning the study of $L^p$ norms for $p<\infty$ are generally reduced to the computations of $L^1$ norms, the study of $L^\infty$ norms is slightly subtler and will be addressed by some passages to the limit.

\begin{lem}\label{lem:15}
Take $\vect s\in \C^r$ and $h\in \Omega$. Then, $\Delta_\Omega^{\vect s}(h+i\,\cdot\,)\in L^1(F)$ if and only if $ \Re\vect s\in \vect d-\frac 1 2\vect{m'}-(\R_+^*)^r$. In this case, 
\[
\int_F \abs{\Delta_\Omega^{\vect s}(h+i x)}\,\dd x= (4 \pi)^m 2^{\Re\vect s} \frac{\Gamma_{\Omega'}(\vect d-\Re\vect s)}{\abs{\Gamma_{\Omega'}(-\vect s/2)}^2} \Delta_\Omega^{\Re\vect s-\vect d}(h) .
\]
\end{lem}

Again, an analogous statement holds for $\Omega'$. 
The proof is based on~\cite[Lemma 4.20]{NanaTrojan}, which deals with the case $\vect s \in -\vect{m'}-(\R_+^*)^r$.

\begin{proof}
Observe first that $\Delta_\Omega^{\vect s}(h+i x)=\Delta_\Omega^{\vect s/2}(h+i x)^2$ for every $x\in F$, so that $\Delta_\Omega^{\vect s}(h+i\,\cdot\,)\in L^1(F) $ if and only if $\Delta_\Omega^{\vect s/2}(h+i\,\cdot\,)\in L^2(F)$. Now, 
\[
\Delta_\Omega^{\vect s/2}(h+i\,\cdot\,)=\frac{1}{\Gamma_{\Omega'}(-\vect s/2)} \Fc_{F'}( \ee^{-\langle \,\cdot\,,h\rangle} \Delta_{\Omega'}^{-\vect s/2}\cdot \nu_{\Omega'} )
\]
when $\Re\vect s\in - \vect{m'}-(\R_+^*)^r$, thanks to Proposition~\ref{prop:58}. By Proposition~\ref{prop:38} and analyticity,
\[
\Delta_\Omega^{\vect s/2}(h+i\,\cdot\,)=  \Fc_{F'}( \ee^{-\langle \,\cdot\,,h\rangle} I^{-\vect s/2}_{\Omega'}  ) 
\]
on $\Omega'$ for every $\vect s\in \C^r$. 
Observe, in addition, that Proposition~\ref{prop:38} implies that 
\[
I^{-\vect s/2}_{\Omega'} = \frac{1}{\Gamma_{\Omega'}(-\vect s/2)} \Delta_{\Omega'}^{-\vect s/2}\cdot \nu_{\Omega'}
\]
on $\Omega'$, and that $I^{-\vect s/2}_{\Omega'} \neq 0$.
Therefore, Proposition~\ref{prop:58} shows that 
\[
\ee^{-\langle \,\cdot\,,h\rangle} I^{-\vect s/2}_{\Omega'}  \in L^2(\Omega')\cdot \nu_{\Omega'}
\]
if and only if $-\Re\vect s+\vect d\in \frac 1 2 \vect{m'}+(\R_+^*)^r$, in which case
\[
\begin{split}
\int_F \abs{\Delta_\Omega^{\vect s}(h+i x)}\,\dd x&=\frac{(2\pi)^m}{\abs{\Gamma_{\Omega'}(-\vect s/2)}^2} \int_{\Omega'} \ee^{-\langle  \lambda,2 h\rangle} \Delta_{\Omega'}^{\vect d-\Re \vect s}(\lambda)\,\dd \nu_{\Omega'}(\lambda)  \\
& = (2 \pi)^m 2^{\Re\vect s-\vect d} \frac{\Gamma_{\Omega'}(\vect d-\Re\vect s)}{\abs{\Gamma_{\Omega'}(-\vect s/2)}^2}\Delta_\Omega^{\Re\vect s-\vect d}(h).
\end{split}
\]
The assertion follows since $2^{-\vect d}=2^m$.
\end{proof}

We now present an analogue of Corollary~\ref{cor:10}. The proof is inspired by that of~\cite[Lemma 4.19]{NanaTrojan}.

\begin{lem}\label{lem:74}
Take $\vect s\in \C^r$ such that $\Re \vect s\in - \frac 1 2 \vect{m'}-(\R_+^*)^r$,  and take $h\in \Omega$. Then, the function $h'\mapsto \Delta_\Omega^{\vect s}(h+h')$ is $\Phi_*(\Hc^{2 n})$-integrable on $\Omega$ if and only if $\Re \vect s\in \vect b- \frac 1 2 \vect{m'}-(\R_+^*)^r$. In this case,
\[
\int_{\Omega} \Delta_\Omega^{\vect s}(h+h')\,\dd \Phi_*(\Hc^{2 n})(h')= \frac{\pi^n\Gamma_{\Omega'}(\vect b-\vect s)}{\abs{\Pfaff(e_{\Omega'})}\Gamma_{\Omega'}(-\vect{s})} \Delta_\Omega^{\vect s-\vect b}(h).
\] 
\end{lem}

For the sake of simplicity, we shall not characterize the set  $\Sb$ of $\vect s\in \R^r$ such that  the function $h'\mapsto \Delta_\Omega^{\vect s}(h+h')$ is $\Phi_*(\Hc^{2 n})$-integrable on $\Omega$ for every $h\in \Omega$. We simply observe that, if $-\vect b \in \frac 1 2 \vect{m}+(\R_+^*)^r$, then Propositions~\ref{prop:58} and~\ref{prop:38} show that $\Sb= \vect b- \frac 1 2 \vect{m'}-(\R_+^*)^r$, while $S=\R^r$ if $E=\Set{0}$ since, in this case, $\Phi_*(\Hc^{2 n})=\delta_0$.
The interested reader may recover the general case observing that  $\Sb$  is an open subset of $\R^r$, that the mapping
\[
 \Sb +i \R^r \ni \vect s \mapsto \int_{\Omega} \Delta_\Omega^{\vect s}(h+h')\,\dd \Phi_*(\Hc^{2 n})(h')\in \C
\]
is holomorphic, and that the mapping
\[
\R^r\ni \vect s\mapsto  \int_{\Omega} \Delta_\Omega^{\vect s}(h+h')\,\dd \Phi_*(\Hc^{2 n})(h')\in [0,\infty]
\]
is continuous by Fatou's lemma. 

\begin{proof}
Assume first that $\vect s\in - \frac 1 2 \vect{m'}-(\R_+^*)^r$, and observe that Propositions~\ref{prop:58} and~\ref{prop:59} and Tonelli's theorem imply that
\[
\begin{split}
\int_{\Omega} \Delta_\Omega^{\vect s}(h+h')\,\dd \Phi_*(\Hc^{2 n})(h')&=\frac{1}{\Gamma_{\Omega'}(-\vect{s})}\int_\Omega \int_{\Omega'} \ee^{-\langle\,\cdot\,, h+h'\rangle} \Delta_{\Omega'}^{-\vect s}\,\dd \nu_{\Omega'}\,\dd \Phi_*(\Hc^{2 n})(h')\\
&=\frac{\pi^n}{\abs{\Pfaff(e_{\Omega'})}\Gamma_{\Omega'}(-\vect{s})} \int_{\Omega'}   \ee^{-\langle \lambda, h\rangle} \Delta_{\Omega'}^{\vect b-\vect s}(\lambda) \,\dd \nu_{\Omega'}(\lambda)\\
&=\frac{\pi^n\Gamma_{\Omega'}(\vect b-\vect s)}{\abs{\Pfaff(e_{\Omega'})}\Gamma_{\Omega'}(-\vect{s})} \Delta_\Omega^{\vect s-\vect b}(h)
\end{split}
\]
if and only if  $ \vect s\in - \frac 1 2 \vect{m'}-(\R_+^*)^r$, while the first integral is $\infty$ otherwise. The assertion for general $\vect s$ follows, since $\abs{\Delta_\Omega^{\vect s}}=\Delta_\Omega^{\Re \vect s}$ on $\Omega$.
\end{proof}

The preceding results allow to compute $L^p$ `norms' of several functions associated with the generalized power functions $\Delta_\Omega^{\vect s}$ and $\Delta_{\Omega'}^{\vect s}$, $p\in ]0,\infty[$. We now deal with $L^\infty$ norms.

\begin{deff}\label{24}
For every $\vect s,\vect{s'}\in \R_+^r$, we define
\[
\vect s^{\vect{s'}}\coloneqq \prod_{j=1}^r s_j^{s'_j},
\]
with the convention $0^0=1$.
\end{deff}

\begin{lem}\label{lem:18}
Take $\vect s\in \C^r$. Then, $\Delta^{\vect s}_\Omega$ is bounded on the bounded subsets of $\Omega$ if and only if $\Re\vect s\in \R_+^r$. In this case,
\[
\norm{\ee^{-\langle \lambda,\,\cdot\,\rangle} \Delta_\Omega^{\vect s} }_{L^\infty(\nu_\Omega)}=(\Re\vect s/\ee)^{\Re\vect s}  \Delta_{\Omega'}^{-\Re\vect s}(\lambda)
\]
for every  $h\in \Omega$.
\end{lem}

An analogous result holds for $\Omega'$.

\begin{proof}
Keep the notation of Chapter~\ref{sec:2}, and let $B$ be a bounded subset of $\Omega$. Then, $\sup\limits_{t \cdot e_{\Omega}\in B} \langle e_j', t\rangle$ is finite for every $j=1,\dots, r$, so that $\Delta^{\vect s}_\Omega$ is bounded on $B$ whenever $\Re \vect s\in \R_+^r$. Conversely, assume that $\Re\vect s \not \in \R_+^r$, and take $j\in \Set{1,\dots,r}$ so that $\Re s_j<0$. Then, the sequence $(t_j/k+\sum_{j'\neq j} e_{j'})\cdot e_{\Omega}$ is bounded in $\Omega$, but 
\[
\lim_{k\to \infty} \Delta^{\vect s}_\Omega\left(\left(t_j/k+\sum_{j'\neq j} e_{j'}\right)\cdot e_{\Omega}  \right)=\infty.
\]

Now, assume that $\Re \vect s\in (\R_+^*)^r$, and observe that, if $p$ is (finite and) large enough, Proposition~\ref{prop:58} implies that
\[
\norm{\ee^{-\langle \lambda,\,\cdot\,\rangle} \Delta^{\Re\vect s}_\Omega}_{L^p(\nu_\Omega)}=\Gamma_\Omega(p\Re\vect s)^{1/p} p^{-\Re\vect s} \Delta_{\Omega'}^{-\Re\vect s}(\lambda)
\]
for every $\lambda\in \Omega'$. Therefore, passing to the limit for $p\to \infty$ and using Stirling's formula,
\[
\norm{\ee^{-\langle \lambda,\,\cdot\,\rangle} \Delta_\Omega^{\vect s} }_{L^\infty(\nu_\Omega)}=(\Re\vect s/\ee)^{\Re\vect s}  \Delta_{\Omega'}^{-\Re\vect s}(\lambda)
\]
for every $\lambda\in \Omega'$. The assertion for $\Re\vect s\in \R_+^r$ follows by continuity. 
\end{proof}

\begin{lem}\label{lem:23}
Take $\vect{s},\vect{s'}\in \C^r$ and $h\in\Omega$. Then, the function $h'\mapsto \Delta_\Omega^{\vect s}(h+h') \Delta_\Omega^{\vect{s'}}(h')$ is bounded on $\Omega$ if and only if $\Re\vect{s'}\in \R_+^r$ and $ \Re\vect s+\Re\vect{s'}\in -\R_+^r$, in which case
\[
\sup\limits_{h'\in \Omega} \abs{\Delta_\Omega^{\vect s}(h+h') \Delta_\Omega^{\vect{s'}}(h')}=\frac{(\Re\vect{s'})^{\Re\vect{s'}}(-\Re\vect s-\Re\vect{s'})^{-\Re\vect s-\Re\vect{s'}} }{(-\Re\vect s)^{-\Re\vect s}}\Delta_\Omega^{\Re \vect s+\Re \vect{s'}}(h).
\]
In addition, $\Delta_\Omega^{\vect s}(h+\,\cdot\,) \Delta_\Omega^{\vect{s'}}\in C_0(\Omega) $ if and only if $\Re \vect{s'}\in (\R_+^*)^r$ and $\Re \vect s+\Re \vect{s'}\in -(\R_+^*)^r$.
\end{lem}

An analogous result holds for $\Omega'$.

\begin{proof}
By homogeneity, we may assume that $h=e_\Omega$. Then, define a mapping $I\colon \Omega\ni h\mapsto h^{-1}\in \Omega$ such that $I(t\cdot e_\Omega)= t^{-1}\cdot e_\Omega$ for every $t\in T_+$, so that $I(e_\Omega+\Omega)= \Omega \cap (e_\Omega-\Omega)$. Thus,
\[
\begin{split}
\sup\limits_{h'\in \Omega }\abs{\Delta_\Omega^{\vect s}(e_\Omega+h') \Delta_\Omega^{\vect{s'}}(h')}&= \sup\limits_{h'\in e_\Omega+\Omega} \Delta_\Omega^{\Re \vect s} (h') \Delta_\Omega^{\Re \vect{s'}}(h'-e_\Omega)\\
&= \sup\limits_{h'\in \Omega \cap (e_\Omega-\Omega)} \Delta_\Omega^{-\Re \vect s}(h') \Delta_\Omega^{\Re \vect{s'}}(h'^{-1}-e_\Omega)\\
&=\sup\limits_{h'\in \Omega \cap (e_\Omega-\Omega)} \Delta_\Omega^{-\Re \vect s-\Re \vect{s'}}(h') \Delta_\Omega^{\Re \vect{s'}}(e_\Omega-h'),
\end{split}
\]
which is finite if and only if $\Re \vect{s'}\in \R_+^r$ and $\Re \vect s+\Re \vect{s'}\in -\R_+^r$, thanks to Lemma~\ref{lem:18}. 
Now, if $ \vect s\in (\R_+^*)^r$ and $\vect{s}+ \vect{s'}\in -(\R_+^*)^r$, then for $p$ (finite and) sufficiently large we have, by Corollary~\ref{cor:10},
\[
\norm{ h'\mapsto \Delta_\Omega^{\vect s}(h+h') \Delta_\Omega^{\vect{s'}}(h')  }_{L^p(\nu_\Omega)}= \left(\frac{\Gamma_\Omega(p\vect{s'})\Gamma_{\Omega'}(-p(\vect s+\vect{s'})) }{\Gamma_{\Omega'}(-p\vect s)}  \right)^{1/p} \Delta_\Omega^{\vect s+\vect{s'}}(h),
\]
so that, passing to the limit for $p\to \infty$ and using Stirling's formula,
\[
\sup\limits_{h'\in \Omega} \abs{\Delta_\Omega^{\vect s}(h+h') \Delta_\Omega^{\vect{s'}}(h')}= \frac{\vect{s'}^{\vect{s'}}(-\vect s-\vect{s'})^{-\vect s-\vect{s'}} }{(-\vect s)^{-\vect s}} \Delta^{\vect s+\vect{s'}}_\Omega(h).
\]
The assertion for general $\vect s$ and $\vect{s'}$ follows by continuity.

To conclude, observe that 
\[
\lim_{h'\to \partial(\Omega\cap (e_\Omega-\Omega))}
\Delta_\Omega^{-\Re \vect s-\Re \vect{s'}}(h') \Delta_\Omega^{\Re \vect{s'}}(e_\Omega-h')= 0
\]
if and only if $\Re \vect{s'}\in (\R_+^*)^r$ and $\Re \vect s+\Re
\vect{s'}\in -(\R_+^*)^r$, so that  $\Delta_\Omega^{\vect
  s}(h+\,\cdot\,) \Delta_\Omega^{\vect{s'}}\in C_0(\Omega)$   by the preceding remarks.
\end{proof}

\begin{cor}\label{cor:12}
Take $\vect s\in \R^r$. Then, $\Delta^{\vect s}_\Omega$  is decreasing on $\Omega$ if and only if $\vect s\in - \R_+^r$.
\end{cor}

An analogous result holds for $\Omega'$. Recall that $F$ is endowed with the ordering induced by $\overline \Omega$, so that $x\meg y$ if and only if $ y-x\in \overline \Omega$, for every $x,y\in F$.

\begin{proof}
Apply Lemma~\ref{lem:23} with $\vect{s'}=0$.
\end{proof}

\begin{lem}\label{lem:24}
Take $\vect s\in \C^r$ and $h\in \Omega$. Then, the function $ x \mapsto \abs{\Delta^{\vect s}_\Omega(h+i x)}$ is bounded on $F$ if and only if $\Re\vect s\in -\R_+^r$, in which case
\[
\sup\limits_{x\in F} \abs{\Delta^{\vect s}_\Omega(h+i x)}=\frac{(-\Re s)^{-\Re s}}{\abs{(-\vect s)^{-\vect s}}} \Delta_\Omega^{\Re\vect s}(h).
\]
In addition, $\abs{\Delta^{\vect s}_\Omega(h+i \,\cdot\,)}\in C_0(F)$ if and only if $\Re\vect s\in - (\R_+^*)^r$.
\end{lem}

An analogous result holds for $\Omega'$.

\begin{proof}
Assume first that $\Re\vect s\in -(\R_+^*)^r$. Then Lemma~\ref{lem:15} shows that, if $p$ is large enough,
\[
\norm{\Delta^{\vect s}_\Omega(h+i\,\cdot\,)}_{L^p(\Nc)}=2^{\Re\vect s} \left((4\pi)^m \frac{ \Gamma_{\Omega'}(\vect d-p\Re \vect s) }{\abs{\Gamma_{\Omega'}(-(p/2)\vect s)}^2}  \right)^{1/p} \Delta^{\Re\vect s-\vect d/p}_\Omega(h),
\]
so that, passing to the limit for $p\to \infty$, and using Stirling's formula,
\[
\sup\limits_{x\in F} \abs{\Delta^{\vect s}_\Omega(h+i x)}=\frac{(-\Re s)^{-\Re s}}{\abs{(-\vect s)^{-\vect s}} } \Delta_\Omega^{\Re\vect s}(h).
\]
Then, by continuity we see that the same happens also for every $\vect s\in- \R_+^r$. 

Now, assume that $\Delta^{\vect s}_\Omega(h+i\,\cdot\,)$ is bounded on $F$. Then, 
\[
 \Delta^{\vect s+\vect{s'}}_{\Omega } (h+i \,\cdot\,)=\Delta^{\vect s}_\Omega(h+i\,\cdot\,)\Delta^{\vect {s'}}_\Omega(h+i\,\cdot\,)
\]
is integrable on $F$ for every $\vect{s'}\in \vect d-\frac 1 2 \vect{m'}-(\R_+^*)^r$, thanks to Lemma~\ref{lem:15}, so that Lemma~\ref{lem:15} implies that $\Re\vect{s}+\vect{s'}\in  \vect d-\frac 1 2 \vect{m'}-(\R_+^*)^r$. By the arbitrariness of $\vect{s'}$, it then follows that $\Re\vect{s}\in -\R_+^r$.

Next, assume that $\vect s\in -(\R_+^*)^r$, and define 
\[
f\colon F+i\Omega \ni (x+i h')\mapsto \Delta^{\vect s}_\Omega(h'-i x)\in\C,
\]
so that, with the notation of Section~\ref{sec:1:3}, $f\in A^{p,1}_{\nu_\Omega,\loc}(F+i \Omega)$ for $p$ (finite and) sufficiently large. In addition, $\Delta^{\vect s}_\Omega(h+i\,\cdot\,)=f_h(-\,\cdot\,)  $. Since Proposition~\ref{prop:1} shows that $f\in A^{\infty,\infty}_{\nu_\Omega,0,\loc}(F+i \Omega)$, this implies that $\Delta^{\vect s}_\Omega(h+i\,\cdot\,)\in C_0(F)$.

Finally, assume that $\vect s\in -\R_+^r\setminus (-\R_+^*)^r$, and let us prove that $\Delta^{\vect s}(h+i\,\cdot\,)\not\in C_0(F)$. Indeed, keep the notation of Section~\ref{sec:2}, and define 
\[
t_\rho\coloneqq r e_j+\sum_{k\neq j} e_k
\]
for every $\rho>0$, where $j\in\Set{1,\dots, r}$ is chosen so that $s_j=0$. Then,  
\[
\Delta^{\vect s}_\Omega(h+i t_\rho \cdot x)=\Delta^{\vect s}_\Omega(t_\rho^{-1}\cdot h+i x)
\]
for every $x\in F$ and for every $r>0$. Arguing as in the proof of Lemma~\ref{lem:13}, we then see that, for every $x\in \Omega$,
\[
\Delta^{\vect s}_\Omega(h+i t_\rho \cdot x)= i^{\vect s}\Delta^{\vect s}_\Omega(x-i t_\rho^{-1}\cdot h),
\]
which converges to $ i^{\vect s}\Delta^{\vect s}_\Omega(x-i P(h))$ as $r\to +\infty$, where $P$ is the projector of $H$ onto $H\cap \left(\sum_{k,\ell\neq j} A_{k,\ell}\right)$ with kernel $H\cap \left(\sum_{k=1}^r (A_{j,k}+A_{k,j})\right) $. Since $\Delta^{\vect s}_\Omega$ is holomorphic on $\Omega+i F$ and does not vanish identically, there is $x\in \Omega$ such that $ \Delta^{\vect s}_\Omega(x-i P(h))\neq 0$, so that $\lim\limits_{r\to +\infty}\Delta^{\vect s}_\Omega(h+i t_\rho \cdot x) \neq 0$ even though $\lim\limits_{r\to +\infty } t_\rho \cdot x=\infty$. Thus, $\Delta^{\vect s}_\Omega(h+i \,\cdot\,)\not \in C_0(F) $.
\end{proof}

The following functions are closely related to the reproducing kernels of the spaces $A^{2,2}_{\vect{s'}}(D)$ that we shall define in Section~\ref{sec:6} (cf.~Proposition~\ref{prop:6}).

\begin{deff}\label{25}
Take $\vect s\in \C^r$, and define 
\[
B^{\vect s}_{(\zeta',z')}(\zeta,z)\coloneqq \Delta_\Omega^{\vect s}\left( \frac{z-\overline{z'}}{2i} -\Phi(\zeta,\zeta') \right)
\]
for every $((\zeta,z),(\zeta',z'))\in (D\times \overline D)\cup (\overline D \times D)$. 
\end{deff}

Observe that $B^{\vect s}_{(\zeta',z')}(\zeta,z) B^{-\vect s}_{(\zeta',z')}(\zeta,z)=1$, so that $B^{\vect s}_{(\zeta',z')}(\zeta,z)\neq 0$ for every $((\zeta,z),(\zeta',z'))\in (D\times \overline D)\cup (\overline D \times D)$.

We conclude this section with the computation of the $L^{p,q}_{\vect s}(D)$ norms of the functions  $B^{\vect{s'}}_{(\zeta',z')}$  (cf.~Definition~\ref{26}).

\begin{lem}\label{lem:21}
Take $p\in ]0,\infty]$, $\vect s\in \C^r$, and $((\zeta',z'),h)\in (D\times \overline \Omega)\cup(\overline D\cup \Omega)$. Then, $\left(B^{\vect s}_{(\zeta',z')}\right)_h\in L^p_0(\Nc)$ (resp.\ $\left(B^{\vect s}_{(\zeta',z')}\right)_h\in L^p(\Nc)$)  if and only if $\Re\vect s\in \frac 1 p (\vect b+\vect d)-\frac{1}{2p} \vect{m'}-(\R_+^*)^r$ (resp.\ $\Re\vect s\in -\R_+^r$ if $p=\infty$). In this case,
\[
\norm*{ \left(B^{\vect s}_{(\zeta',z')}\right)_h }_{L^p(\Nc)}=C_{\vect s,p}\Delta_\Omega^{\Re\vect s-(\vect b+\vect d)/p}(h+\Im z'-\Phi(\zeta')),
\]
where 
\[
C_{\vect s,p}\coloneqq \left(\frac{(4 \pi)^m \pi^n \Gamma_{\Omega'}(\vect b+\vect d-p\Re\vect{s})}{\abs{\Pfaff(e_{\Omega'})} \abs{\Gamma_{\Omega'}(-(p/2)\vect{s})}^2}   \right)^{1/p} 
\]
if $p<\infty$, while $C_{\vect s,\infty}\coloneqq \frac{(-\Re s)^{-\Re s}}{\abs{(-\vect s)^{-\vect s}} }$.
\end{lem}

\begin{proof}
Assume first that $p=1$, and define $h'\coloneqq \Im z'-\Phi(\zeta')$. Then,
\[
\begin{split}
\int_{\Nc} &\abs*{\Delta_\Omega^{\vect s}\left( \frac{(x+i\Phi(\zeta)+i h)-\overline{z'}}{2i}-\Phi(\zeta,\zeta') \right)} \,\dd (\zeta,x)\\
&=2^{-\Re\vect s}\int_{E\times F}  \abs{\Delta_\Omega^{\vect s}\left( h+h'+\Phi(\zeta-\zeta')-(x-x')i-\Im\Phi(\zeta,\zeta')i \right)}\,\dd (\zeta,x)\\
&=2^{-\Re\vect s}\int_{\Omega} \int_{F}  \abs{\Delta_\Omega^{\vect s}\left( h+h'+h''- ix \right)}\,\dd x \,\dd \Phi_*(\Hc^{2 n})( h'')\\
&=\frac{(4 \pi)^m\Gamma_{\Omega'}(\vect d-\Re\vect s)  }{\abs{\Gamma_{\Omega'}(-\vect s/2)}^2}\int_{ \Omega}\Delta_\Omega^{\Re\vect s-\vect d}(h+h'+h'')\,\dd \Phi_*(\Hc^{2 n}) (h'')\\
&=  \frac{(4\pi)^m \pi^n\Gamma_{\Omega'}(\vect b+\vect d-\Re\vect s)}{\abs{\Pfaff(e_{\Omega'})}\abs{\Gamma_{\Omega'}(-\vect s/2)}^2}\Delta_\Omega^{\Re\vect s-\vect b-\vect d}(h+h')
\end{split}
\]
by Lemmas~\ref{lem:15} and~\ref{lem:74},  provided that $\Re\vect s \in \vect b+\vect d-\frac 1 2\vect{m'}-(\R_+^*)^r$. The first integral is $\infty$ otherwise. 
The  assertion for $p<\infty$ follows since $\norm*{ \left(B^{\vect s}_{(\zeta',z')}\right)_h }_{L^p(\Nc)}=\norm*{ \left(B^{p\vect s}_{(\zeta',z')}\right)_h }_{L^1(\Nc)}^{1/p}$.

Next, assume that $p=\infty$. Then,  Corollary~\ref{cor:12} and Lemma~\ref{lem:24} imply that $B^{\vect{s}}_{(\zeta',z')}\in L^\infty(\Nc)$ if and only if $\Re\vect{s}\in -\R_+^r$, and that, in this case,
\[
\norm*{\left(B^{\vect{s}}_{(\zeta',z')}\right)_h}_{L^\infty(\Nc)}=C_{\vect{s},\infty} \Delta_\Omega^{\vect{s}}(h+\Im z'-\Phi(\zeta')).
\]
In addition, Lemma~\ref{lem:24} again shows that  $\left(B^{\vect s}_{(\zeta',z')}\right)_h\in C_0(\Nc)$ if and only if $\vect s\in - (\R_+^*)^r$.
\end{proof}

\begin{deff}\label{26}
Take $\vect s\in \R^r$ and $p,q\in ]0,\infty]$. 
We define $L^{p,q}_{\vect s}(D)$ as the Hausdorff space associated with the space
\[
\Set{f\colon D\to \C\colon f \text{ is measurable}, \int_\Omega  \left( \Delta_\Omega^{\vect s}(h) \norm{f_h}_{L^p(\Nc)} \right) ^q\,\dd \nu_\Omega(h)<\infty  }
\]
(modification when $q=\infty$), endowed with the corresponding  structure of a locally bounded $F$-space. 
We define $L^{p,q}_{\vect s,0}(D)$ as the closure of $C_c( D)$ in $L^{p,q}_{\vect s}(D)$.
\end{deff}

Notice that $L^{p,q}_{\vect s,0}(D)=L^{p,q}_{\vect s}(D)$ if (and only
if) $p,q<\infty$.

\begin{prop}\label{prop:60}
Take $\vect s\in \R^r$, $\vect{s'}\in \C^r$, $p,q\in ]0,\infty ]$, and $(\zeta',z')\in D$. Then, $B^{\vect{s'}}_{(\zeta',z')}\in L^{p,q}_{\vect{s},0}(D)$ (resp.\ $B^{\vect {s'}}_{(\zeta',z')}\in L^{p,q}_{\vect{s}}(D)$) if and only if the following conditions are satisfied:
\begin{itemize}
\item $\vect s\in \frac{1}{2 q} \vect{m}+(\R_+^*)^r$ (resp.\ $\vect s\in \R_+^r$ if $q=\infty$);

\item $\Re\vect{s'}\in\frac 1 p (\vect b+\vect d)-\frac{1}{2p} \vect{m'}-(\R_+^*)^r $ (resp.\ $\Re\vect{s'}\in- \R_+^r$ if $p=\infty$);

\item $\vect s+\Re\vect{s'}\in \frac 1 p (\vect b+\vect d)-\frac{1}{2 q} \vect{m'}-(\R_+^*)^r$ (resp.\ $\vect s+\Re\vect{s'}\in \frac{1}{p}(\vect b+\vect d)- \R_+^r$ if $q=\infty$).
\end{itemize}
In this case,
\[
\norm*{B^{\vect{s'}}_{(\zeta',z')}}_{L^{p,q}_{\vect{s}}(D)}=C_{\vect{s'},p}C_{\vect s,\vect{s'},p,q} \Delta_\Omega^{\vect s+\vect{s'}-(\vect b+\vect d)/p}(\Im z'-\Phi(\zeta')),
\]
where $C_{\vect{s'},p}$ is that of Lemma~\ref{lem:21}, and
\[
C_{\vect s,\vect{s'},p,q} \coloneqq \left( \frac{ \Gamma_\Omega(q\vect s) \Gamma_{\Omega'}((q/p)(\vect b+\vect d)-q(\vect s+\Re\vect{s'}) ) }{\Gamma_{\Omega'}((q/p)(\vect b+\vect d)-q\Re\vect{s'})}\right) ^{1/q},
\]
if $q<\infty$, while 
\[
C_{\vect s,\vect{s'},p,\infty} \coloneqq  \frac{\vect s^{\vect s} ((\vect b+\vect d)/p-\vect s-\Re\vect{s'})^{(\vect b+\vect d)/p-\vect s-\Re\vect{s'}} }{((\vect b+\vect d)/p-\Re\vect{s'})^{(\vect b+\vect d)/p-\Re\vect{s'}} }.
\]
\end{prop}

\begin{proof}
Take $h\in \Omega$. Then Lemma~\ref{lem:21} shows that $\left(B^{\vect{s'}}_{(\zeta',z')}\right)_h\in L^p_0(\Nc)$ (resp.\ $\left(B^{\vect{s'}}_{(\zeta',z')}\right)_h\in L^p(\Nc)$) if and only if $\Re\vect{s'}\in \frac 1 p (\vect b+\vect d)-\frac{1}{2p} \vect{m'}-(\R_+^*)^r$ (resp.\ $\Re\vect{s'}\in -\R_+^r$ if $p=\infty$), and that, in this case,
\[
\norm*{(B^{\vect{s'}}_{(\zeta',z')})_h}_{L^p(\Nc)}= C_{\vect{s'},p} \Delta_\Omega^{\Re\vect{s'}-(\vect b+\vect d)/p}(h+\Im z'-\Phi(\zeta')).
\]
If $q<\infty$, then Corollary~\ref{cor:10} shows that $B^{\vect{s'}}_{(\zeta',z')}\in L^{p,q}_{\vect s,0}(D) $ (resp.\ $B^{\vect{s'}}_{(\zeta',z')}\in L^{p,q}_{\vect s}(D) $) if and only if the conditions of the statement hold and that, in this case,
\[
\norm*{B^{\vect{s'}}_{(\zeta',z')}}_{A^{p,q}_{\vect s}(D) }= C_{\vect{s'},p}C_{\vect s,\vect{s'},p,q} \Delta_\Omega^{\vect s+\Re\vect{s'}-(\vect b+\vect d)/p}(\Im z'-\Phi(\zeta')).
\]
Finally, if $q=\infty$, then Lemma~\ref{lem:23} shows that $B^{\vect{s'}}_{(\zeta',z')}\in L^{p,\infty}_{\vect s,0}(D) $ (resp.\ $B^{\vect{s'}}_{(\zeta',z')}\in L^{p,\infty}_{\vect s}(D) $) if and only if  the conditions of the statement hold and that, in this case,
\[
\norm*{B^{\vect{s'}}_{(\zeta',z')}}_{L^{p,\infty}_{\vect s}(D) }= C_{\vect{s'},p}C_{\vect s,\vect{s'},p,\infty} \Delta_\Omega^{\vect s+\Re\vect{s'}-(\vect b+\vect d)/p}(\Im z'-\Phi(\zeta')).
\]
The assertion follows.
\end{proof}

\section{Lattices and Quasi-Constancy}\label{sec:3:3}

This section deals with the various tools needed for the `discretization' of several problems related to weighted Bergman spaces (and Besov spaces).

We first recall the definition of the Bergman metric and use this metric to induce an invariant metric on $\Omega$ (which is not, in general, the `canonical' invariant metric on a homogeneous cone). We also fix an invariant measure on $D$.  
We then deal with the quasi-constancy of several functions on $D$, $\Omega$, $\Omega'$ on invariant balls. 

We conclude this section with a suitable definition of lattices on $D$, $\Omega$, and $\Omega'$, and with a proof of their main properties.

\medskip

Denote temporarily with $B$ the composite of the unweighted Bergman kernel (that is, the reproducing kernel associated with the unweighted Bergman space $A^{2,2}_{\Hc^m}(D)$, with the notation of Section~\ref{sec:1:3})  and the diagonal mapping of $D$. Then, Corollary~\ref{cor:5} and Propositions~\ref{prop:58} and~\ref{prop:59} show that
\[
B(\zeta,z)= \frac{\abs{\Pfaff(e_{\Omega'})} \Gamma_{\Omega'}(-\vect b-2 \vect d)  }{4^m \pi^{n+m}\Gamma_\Omega(-\vect d) } \Delta_\Omega^{\vect b+2\vect d}\left(\Im z-\Phi(\zeta) \right)
\]
for every $(\zeta,z)\in D$.  The Bergman metric $\vect{k}$ is then defined by\index{Bergman metric}
\[
\vect{k}_{(\zeta,z)} v v'\coloneqq \partial_v\overline{\partial_{v'}}(\log B)(\zeta,z)
\]
for every $(\zeta,z)\in D$ and for every $v,v'\in E\times F_\C$. Notice that $\vect{k}$ is a K\"ahler metric. We denote by $\vect{g}$ its real part, that is, the associated Riemannian metric.

Simple computations lead to the following result.

\begin{lem}\label{lem:73}
For every $(\zeta,z)\in D$ and for every $(v,w),(v',w')\in E\times F_\C$,
\[
\begin{split}
&\vect{k}_{(\zeta,z)}((v,w)(v',w'))= \left(\frac{(\Delta^{\vect b+2\vect  d}_\Omega)''(h)}{\Delta_\Omega^{\vect b+2\vect d}(h)}-\left(\frac{(\Delta^{\vect b+2\vect  d}_\Omega)'(h)  }{\Delta_\Omega^{\vect b+2 \vect d}(h) }\right)^{\otimes 2}   \right)\cdot \\
&\qquad \cdot\left(-\frac{i}{2}w-\Phi(v,\zeta)  \right)\cdot\left( \frac{i}{2}\overline{w'}-\Phi(\zeta,v') \right)-\frac{(\Delta^{\vect b+2\vect  d}_\Omega)'(h)}{\Delta_\Omega^{\vect b+2\vect d}(h)}\cdot \Phi(v,v'),
\end{split}
\]
where $h\coloneqq \Im z-\Phi(\zeta)$.
\end{lem}

\begin{deff}\label{27}
We endow $D$ with the Riemannian distance $d$ associated with $\vect g$, and we define 
\[
\nu_D\coloneqq \frac{4^m \pi^{n+m}\Gamma_\Omega(-\vect d) }{\abs{\Pfaff(e_{\Omega'})} \Gamma_{\Omega'}(-\vect b-2 \vect d)  } B\cdot \Hc^{2 n+2 m}.
\]
We denote with $B((\zeta,z),\rho)$ the open ball with centre $(\zeta,z)\in D$ and radius $\rho>0$ relative to the distance $d$.
\end{deff}

The constant in the definition of $\nu_D$ has been added for computational reasons. Indeed, with these conventions,
\[
\int_D f\,\dd \nu_D= \int_\Omega \int_\Nc f_h(\zeta,x)\,\dd (\zeta,x)  \Delta_\Omega^{\vect b+\vect d}(h)\,\dd\nu_\Omega(h)
\]
for every $f\in L^1(\nu_D)$.

\begin{prop}\label{prop:26}
Denote by $G(D)$ the group of biholomorphisms of $D$ onto itself. Then, the following hold:
\begin{enumerate}
\item[\em(1)] $G(D)$ acts transitively on $D$;

\item[\em(2)] $\vect k$, $\vect g$, and the associated distance $d$ are $G(D)$-invariant;

\item[\em(3)] the measure $\nu_D$ is $G(D)$-invariant;

\item[\em(4)] the metric space $(D,d)$ is complete. In particular, the closed balls relative to $d$ are compact.
\end{enumerate}
\end{prop}

\begin{proof}
{(1)} This follows from the fact that $G(D)$ contains the left translations by elements of $b D$ and the automorphisms of the form $(\zeta,z)\mapsto (g\zeta, t\cdot z)$ for $t\in T_+$ and $g\in GL(E)$ such that $t\cdot \Phi=\Phi\circ (g\times g)$. 

{ (2)} This follows from~\cite[Proposition 1.4.15]{Krantz}.

{ (3)} This follows from~\cite[Proposition 1.4.12]{Krantz}.

{ (4)} The first assertion follows from~\cite[Theorem 1.3]{Nakajima}, while the second assertion follows from the Hopf--Rinow theorem (cf.~\cite[Corollary 6.7]{Lang}). 
\end{proof}

\begin{deff}\label{29}
We endow $\Omega$ with the Riemannian metric $\vect g_\Omega$ induced by $\vect g$ through the embedding $h\mapsto (0,i h)$ of $\Omega$ in $D$.  We denote with $d_\Omega$ the associated distance.
We denote with $B_\Omega(h,\rho)$ the open ball with centre $h\in \Omega$ and radius $\rho>0$ relative to the distance $d_\Omega$.

We endow $\Omega'$ with the Riemannian metric $\vect g_{\Omega'}$ and the associated distance $d_{\Omega'}$ obtained by means of the identification $\Omega \ni t\cdot e_\Omega\mapsto e_{\Omega'}\cdot t^{-1}\in \Omega'$, where $t\in T_+$. 
We denote with $B_{\Omega'}(\lambda,\rho)$ the open ball with centre $\lambda\in \Omega'$ and radius $\rho>0$ relative to the distance $d_{\Omega'}$.
\end{deff}

In this way, the distances on $\Omega$ and $\Omega'$ are $T_+$-invariant.

The following result allows us to compare the invariant balls on $D$ and $\Omega$. In particular, it shows that, even though $\vect{g}_\Omega$ was defined by means of an embedding, the distance on $\Omega$ is actually a `quotient' of the distance on $D$.

\begin{lem}\label{lem:45} 
The mapping $\varrho \colon D\ni (\zeta,z)\mapsto \Im z-\Phi(\zeta)\in \Omega$ is a Riemannian submersion. In addition, for every $(\zeta,z),(\zeta',z')\in D$, 
\[
d_\Omega(\varrho (\zeta,z),\varrho (\zeta',z'))\meg d((\zeta,z),(\zeta',z')),
\]
with equality if $(\zeta,z),(\zeta',z')\in (\zeta'',z'')\cdot (\Set{0}\times i\Omega) $ for some $(\zeta'',z'')\in b D$.
\end{lem}

Recall that, if $(M_1,\vect g_1)$ and $(M_2, \vect g_2)$ are Riemannian manifolds and $f\colon M_1\to M_2$ is a submersion, then $f$ is a Riemannian submersion if $T_{x_1}(f)$ induces an isometry of $\ker T_{x_1}(f)^\perp$ onto $T_{f(x_1)}(M_2)$ for every $x_1\in M_1$.

Observe that, in particular, the closed balls of $\Omega$ are compact, so that the distance $d_\Omega$ induces the uniformity of a complete space on $\Omega$.

\begin{proof}
Take $(v,w)\in E\times F$ and $w'\in i F$, and observe that 
\[
\vect k_{(0,i h)}(v,w)(0,w')\in i \C
\]
by Lemma~\ref{lem:73}, since $\Delta^{\vect b+2\vect d}_\Omega$ is real on $\Omega$. Thus, $(v,w)$ and $(0,w')$ are orthogonal for $\vect g_{(0,i h)}$.
Since $E\times F=\ker \varrho '(0,i h)$ for every $h\in \Omega$, this implies that $i F$ is the $\vect{g}_{(0, i h)}$-orthogonal complement of $\ker \varrho '(0, i h)$. In addition, the restriction of $\varrho '(0, i h)$ to $i F$ is an isometry between $\vect{g}_{(0, i h)}$ and $(\vect{g}_{\Omega})_h$ by definition, so that $\varrho $ is a Riemannian submersion at $(0, i h)$ for every $h\in \Omega$. By the invariance of $\vect k$ and $\varrho $ under the action of $b D$ (cf.~Proposition~\ref{prop:26}), it follows that $\varrho $ is a Riemannian submersion on the whole of $D$.

Thus, if $\gamma\colon [0,1]\to D$ is a piecewise smooth curve joining two points $(\zeta,z)$ and $(\zeta',z')$ of $D$, then $\varrho \circ \gamma\colon [0,1]\to \Omega$ is a piecewise smooth curve joining $h\coloneqq \varrho (\zeta,z)$ and $h'\coloneqq\varrho (\zeta',z')$, and
\[
\norm{(\varrho  \circ \gamma)'(t)}_{\vect{g_\Omega}}\meg \norm{\gamma'(t)}_{\vect g}
\]
for every $t$ where $\gamma$ is differentiable. 
Consequently,
\[
d_\Omega(\varrho (\zeta,z),\varrho (\zeta',z'))\meg d((\zeta,z),(\zeta',z'))
\]
for every $(\zeta,z),(\zeta',z')\in D$, by the arbitrariness of $\gamma$. 
To prove the last assertion, we may assume that $(\zeta,z)=(0,i h)$ and $(\zeta',z')=(0, i h')$, by homogeneity.
If $\gamma$ is a piecewise smooth curve joining $h$ and $h'$ in $\Omega$, then $(0,i \gamma)$ is a piecewise smooth curve joining $(0,i h)$ and $(0,i h')$, $\gamma=\varrho \circ (0, i \gamma)$, and  
\[
\norm{\gamma'(t)}_{\vect{g_\Omega}}= \norm{(0, i \gamma)'(t)}_{\vect g},
\]
whence  
\[
d_\Omega(h,h')\Meg d((0,i h),(0,i h')).
\]
Equality holds by the preceding remarks.
\end{proof}

We now prove some quasi-constancy lemmas for several functions on $D$, $\Omega$, and $\Omega'$.
The following result, also known as Kor\'anyi's lemma, is of fundamental importance to deal with atomic decomposition of weighted Bergman spaces. 
It first appeared in~\cite[Lemma 2.3]{CoifmanRochberg} for symmetric Siegel domains of type II and $\vect s=\vect b+2\vect d$. It has then been extended in~\cite[Theorem 1.1]{BekolleIshiNana} to general homogeneous Siegel domain of type II, for $\vect s=\vect b+2\vect d$. 
By an inductive procedure, we extend this result to the case of general $\vect s$.

\begin{teo}\label{teo:2}
 Take $\vect s\in \C^r$ and $R>0$. Then, there is a constant $C>0$ such that
\[
\abs*{ \frac{B^{\vect s}_{(\zeta,z)}(\zeta',z')}{B^{\vect s}_{(\zeta,z)}(\zeta'',z'')} -1 }\meg C d((\zeta',z'),(\zeta'',z''))
\]
for every $(\zeta,z),(\zeta',z'),(\zeta'',z'')\in D$ such that $d((\zeta',z'),(\zeta'',z''))\meg R$.
\end{teo}

\begin{oss}
This result may be extended in two ways. On the one hand, one may prove the same assertion under the weaker assumption $(\zeta,z)\in \overline D$, by continuity. On the other hand, arguing by symmetry, one may prove that
\[
\abs*{ \frac{B^{\vect s}_{(\zeta_1,z_1)}(\zeta'_1,z'_1)}{B^{\vect s}_{(\zeta_2,z_2)}(\zeta_2',z_2')} -1 }\meg C (d((\zeta_1,z_1),(\zeta_2, z_2))+d((\zeta_1',z'_1),(\zeta_2',z_2')))
\]
for every $(\zeta_1,z_1),(\zeta'_1,z'_1),(\zeta_2,z_2),(\zeta_2',z_2')\in D$ such that 
\[
d((\zeta_1,z_1),(\zeta_2,z_2)),d((\zeta'_1,z'_1),(\zeta_2',z_2'))\meg R.
\]
\end{oss}

\begin{proof}
We prove the assertion by induction on $r$. Even though the case $r=1$ can be treated aside, since the following construction is meaningful (though trivial) for $r=1$, we may artificially introduce the rather trivial case $r=0$  and proceed with the inductive step. Observe, though, that the we do \emph{not} consider the trivial domain $D=\Set{0}$ elsewhere.
Then, assume that $r\Meg 1$ and that the assertion holds for every homogeneous Siegel domain of type II over a homogeneous cone of rank $<r$. Keep the notation of Section~\ref{sec:2}, so that $F=H$ and $\Omega=C(A)$. Define 
\[
A^{(r-1)}\coloneqq\bigoplus_{j,k\Meg 2} A_{j,k},
\]
and endow $A^{(r-1)}$  with the structure of a $T$-algebra of rank $r-1$ with the operations of $A$ and the graduation $A^{(r-1)}_{j,k}\coloneqq A_{j+1,k+1}$, for $j,k=1,\dots, r-1$. 
Let $T_+^{(r-1)}$, $H^{(r-1)}$ and $\Omega^{(r-1)}=C(A^{(r-1)})$ be the associated spaces.

Let $\pi^{(r-1)}\colon A\to A^{(r-1)}$ be the linear mapping which
induces the identity on $A^{(r-1)}$ and vanishes on $A_{j,k}$ for
every $j,k=1,\dots,r$ such that  either $j=1$ or $k=1$. Thus, $\pi^{(r-1)}$ may be considered as the self-adjoint projector of $A$ onto $A^{(r-1)}$. 
Notice that, even though $\pi^{(r-1)}$ need \emph{not} be a morphism of algebras, it is easily verified that
\[
\pi^{(r-1)}(a^*)=\pi^{(r-1)}(a)^*, \qquad \pi^{(r-1)}(t t')=\pi^{(r-1)}(t)\pi^{(r-1)}(t'),
\]
and
\[
\pi^{(r-1)}(t x t^*)=\pi^{(r-1)}(t)\pi^{(r-1)}(x)\pi^{(r-1)}(t)^*
\]
for every $a\in A$, for every $x\in H$, and for every $t,t'\in T$. In particular, $\pi^{(r-1)}(H)=H^{(r-1)}$, $\pi^{(r-1)}$ induces a homomorphism of $T_+$  onto $T^{(r-1)}_+$, and $\pi^{(r-1)}(C)=\Omega^{(r-1)}$.
Therefore, by an abuse of notation, $\Delta^{\vect{s'}}_{\Omega^{(r-1)}}\circ \pi^{(r-1)}=\Delta^{(0,\vect{s'})}_C$ on $C$ for every $\vect{s'}\in \C^{r-1}$, whence 
\[
\Delta^{\vect{s'}}_{\Omega^{(r-1)}}\circ \pi^{(r-1)}_\C=\Delta^{(0,\vect{s'})}_C
\]
on $C+ i H$ by holomorphy. 

Next, observe that $\pi^{(r-1)}_\C\circ \Phi$ is a $\Omega^{(r-1)}$-positive hermitian mapping of $E$ into $H^{(r-1)}$.
Let $E^{(r-1)}$ be the quotient of $E$ by the radical of $\pi^{(r-1)}_\C\circ \Phi$, and let $\pi_E\colon E\to E^{(r-1)}$ be the canonical projection. We denote by $\Phi^{(r-1)}$ the non-degenerate $\Omega^{(r-1)}$-positive hermitian mapping of $E^{(r-1)}$ into $H^{(r-1)}$ induced by $\pi^{(r-1)}_\C\circ \Phi$. 
Observe that for every $t\in T_+^{(r-1)}$ there is $g\in GL(E)$ such that $t\cdot \Phi=\Phi\circ (g\times g)$, so that 
\[
t\cdot (\pi^{(r-1)}_\C\circ \Phi)=(\pi^{(r-1)}_\C\circ\Phi)\circ (g\times g).
\]
In particular, $g$ preserves the radical of $\pi^{(r-1)}_\C\circ \Phi$, so that it induces an element $g'$ of $ GL(E^{(r-1)})$ such that $t\cdot \Phi^{(r-1)}=\Phi^{(r-1)}\circ (g'\times g')$.  
Therefore, 
\[
D^{(r-1)}\coloneqq \Set{(\zeta',z')\in E^{(r-1)}\times H^{(r-1)}_\C\colon \Im z'-\Phi^{(r-1)}(\zeta')\in \Omega^{(r-1)}}
\]
is a homogeneous Siegel domain of type II, and $D^{(r-1)}=(\pi_E\times \pi^{(r-1)}_\C )(D)$. 
Define
\[
B_{(\zeta'_2,z'_2)}^{\vect{s'},(r-1)}(\zeta'_1,z_1')\coloneqq \Delta_{\Omega^{(r-1)}}^{\vect{s'}}\left( \frac{z_1'-\overline{z'_2}}{2 i}-\Phi^{(r-1)}(\zeta'_1,\zeta'_2)  \right)
\]
for every $(\zeta'_1,z'_1),(\zeta'_2,z'_2)\in D^{(r-1)}$ and for every $\vect{s'}\in \C^{r-1}$. In addition, let $d^{(r-1)}$ be the Bergman distance on $D^{(r-1)}$. Then, the inductive assumption implies that for every $\vect{s'}\in \C^{r-1}$ there is a constant $C_{\vect{s'}}>0$ such that
\[
\abs*{\frac{B_{(\zeta'_0,z'_0)}^{\vect{s'},(r-1)}(\zeta'_1,z_1')}{B_{(\zeta'_0,z'_0)}^{\vect{s'},(r-1)}(\zeta'_2,z_2')}-1}\meg C_{\vect{s'}} d^{(r-1)}((\zeta'_1,z'_1),(\zeta'_2, z'_2))
\]
for every $(\zeta'_0,z'_0), (\zeta'_1,z'_1),(\zeta'_2, z'_2)\in D^{(r-1)}$ such that $d^{(r-1)}((\zeta'_1,z'_1),(\zeta'_2, z'_2))\meg R$. Further, the preceding remarks show that, for every $(\zeta_0,z_0),(\zeta_1,z_1)\in D$, 
\[
B^{(0,\vect{s'})}_{(\zeta_0,z_0)}(\zeta_1,z_1)=B^{\vect{s'},(r-1)}_{(\pi_E(\zeta_0),\pi^{(r-1)}_\C(z_0))}(\pi_E(\zeta_1),\pi^{(r-1)}_\C(z_1)).
\]
Then, take $\vect s\in \C^r$, and observe that there are $\vect{s'}\in \C^{r-1}$ and $\tau\in \C$ such that $\vect s= (0,\vect{s'})+\tau(\vect b+2 \vect d)$, so that
\[
B^{\vect s}_{(\zeta_0,z_0)}(\zeta_1,z_1)=B^{\tau(\vect b+2 \vect d)}_{(\zeta_0,z_0)}(\zeta_1,z_1) B^{\vect{s'},(r-1)}_{(\pi_E(\zeta_0),\pi^{(r-1)}_\C(z_0))}(\pi_E(\zeta_1),\pi^{(r-1)}_\C(z_1))
\]
for every $(\zeta_0,z_0),(\zeta_1,z_1)\in D$.
Now, select an $n$-dimensional real subspace $\widetilde E$ of $E$, and interpret $E$ as the complexification of $\widetilde E$. Endow $E$ with the corresponding conjugation.
Then, observe that, since $\overline D\times D$ is convex, hence simply connected, for every $\tau\in \C$ there is a holomorphic function
\[
b_\tau \colon D\times \overline D\to \C
\] 
such that
\[
B_{\overline{(\zeta,z)}}^{\tau(\vect b+2\vect d)}(\zeta',z')= \ee^{b_\tau((\zeta,z),(\zeta',z'))}
\]
for every $\overline{(\zeta,z)},(\zeta',z')\in D$.
Then, it is readily verified that for every $\tau\in \C$ there is constant $c_\tau\in \C$ such that
\[
b_\tau=\tau b_1+c_\tau
\]
Now, by~\cite[the first proof of Theorem 1.1]{BekolleIshiNana} there is a constant $\widetilde C>0$ such that
\[
\abs{ b_1((\zeta,z),(\zeta',z'))-b_1((\zeta,z),(\zeta'',z''))}\meg \widetilde C d((\zeta',z'),(\zeta'',z''))
\]
for every $(\zeta,z)\in \overline D$ and for every $(\zeta',z'),(\zeta'',z'')\in D$ such that $ d((\zeta',z'),(\zeta'',z''))\meg R$.
Therefore, for every $\tau\in \C$ there is a constant $\widetilde C_\tau>0$ such that
\[
\abs*{\frac{B^{\tau(\vect b+2\vect d)}_{(\zeta,z)}(\zeta',z')  }{B^{\tau(\vect b+2\vect d)}_{(\zeta,z)}(\zeta'',z'') }-1  }\meg \widetilde C_\tau d((\zeta',z'),(\zeta'',z''))
\]
for every $(\zeta,z)\in \overline D$ and for every $(\zeta',z'),(\zeta'',z'')\in D$ such that $ d((\zeta',z'),(\zeta'',z''))\meg R$.

Next, observe that the mapping $\pi_E\times \pi^{(r-1)}_\C\colon D\to D^{(r-1)}$ is analytic, hence locally Lipschitz.
Then, Proposition~\ref{prop:26} and the preceding remarks show that there is a constant  $C>0$ such that
\[
\abs*{ \frac{B^{\vect s}_{(\zeta,z)}(\zeta',z')}{B^{\vect s}_{(\zeta,z)}(0,i e_\Omega)} -1 }\meg C d((\zeta',z'),(0, i e_\Omega))
\]
for every $(\zeta,z),(\zeta',z')\in D$ such that $d((\zeta',z'),(0, i e_\Omega))\meg R$. The assertion then follows by homogeneity.
\end{proof}

\begin{cor}\label{cor:34}
Take $\vect s\in \C^r$. Then, there are two constants $R,C>0$ such that
\[
\abs*{\frac{\Delta_\Omega^{\vect s}(h+h')}{\Delta^{\vect s}_\Omega(h+h'')}-1}\meg C d_\Omega(h',h'')
\]
for every $h\in \overline \Omega$ and for every $h',h''\in \Omega $ such that $d_\Omega(h',h'')\meg R$.

If, in addition, $\vect s\in \R^r$, then for every $R'>0$ there is a constant $C'>0$ such that
\[
\frac{1}{C'}\Delta^{\vect s}_\Omega(h+h'')\meg\Delta_\Omega^{\vect s}(h+h')\meg C' \Delta^{\vect s}_\Omega(h+h'')
\]
for every $h\in \overline \Omega$ and for every $h',h''\in \Omega $ such that $d_\Omega(h',h'')\meg R'$.
\end{cor}

An analogous result holds for $\Omega'$.

\begin{proof}
When $h\in \Omega$, it suffices to apply Theorem~\ref{teo:2} with $(\zeta,z)=(0, i h)$, $(\zeta',z')=(0, i h')$, and $(\zeta'',z'')=(0, i h'')$, thanks to Lemma~\ref{lem:45}. The general case follows by continuity.
\end{proof}

\begin{lem}\label{lem:48}
Take $R>0$. Then, there is a constant $C>0$ such that, for every $h,h'\in \Omega$ and for every $\lambda,\lambda'\in \Omega'$ such that $d_\Omega(h,h')\meg R$ and $d_{\Omega'}(\lambda,\lambda')\meg R$,
\[
\frac{1}{C}\meg \frac{\langle \lambda, h\rangle}{\langle \lambda',h'\rangle}\meg C.
\]
\end{lem}

The proof is based on~\cite[Theorem 2.45]{Bekolleetal}, which deals with symmetric cones.

\begin{proof}
Observe that $\Omega=(\Omega')'$, so that $\langle \lambda,h\rangle>0$ for every non-zero $\lambda\in \overline{\Omega'}$ and for every $h\in \Omega$.
By compactness (cf.~Proposition~\ref{prop:26}), there is a constant $C_1>0$ such that
\[
\frac{1}{C_1}\meg \frac{\langle \lambda, h\rangle}{\langle \lambda,e_\Omega\rangle}\meg C_1
\]
for every $\lambda\in \overline{\Omega'}$ such that $\abs{\lambda}=1$ and for every $h\in \overline B_\Omega(e_\Omega,R)$. Then, the same holds for every $\lambda\in \Omega'$. Next, take $h,h'\in \Omega$ such that $d_\Omega(h,h')\meg R$, and choose $t\in T_+$ so that $h'=t\cdot e_\Omega$. Then, $d_\Omega(t^{-1}\cdot h, e_\Omega )\meg R$ and $ \frac{\langle \lambda, h\rangle}{\langle \lambda,h'\rangle}= \frac{\langle \lambda\cdot t, t^{-1}\cdot h\rangle}{\langle \lambda\cdot t,e_\Omega\rangle}$, so that
\[
\frac{1}{C_1}\meg \frac{\langle \lambda, h\rangle}{\langle \lambda,h'\rangle}\meg C_1
\]
for every $\lambda\in \Omega'$. Applying the above arguments to $\Omega'$, we see that there is a constant $C_2>0$ such that
\[
\frac{1}{C_2}\meg \frac{\langle \lambda, h\rangle}{\langle \lambda',h\rangle}\meg C_2
\]
for every $h\in \Omega$ and for every $\lambda,\lambda'\in \Omega'$ such that $d_{\Omega'}(\lambda, \lambda')\meg R$. Then, the assertion follows with $C=C_1 C_2$.
\end{proof}

\begin{cor}\label{cor:18}
Take $R>0$. Then, there is a constant $C>0$ such that, for every $h,h'\in \Omega$ and for every $\lambda,\lambda'\in \Omega'$ such that $d_\Omega(h,h')\meg R$ and $d_{\Omega'}(\lambda,\lambda')\meg R$,
\[
\frac{1}{C}\meg \frac{\abs{h}}{\abs{h'}}\meg C \qquad \text{and} \qquad \frac{1}{C}\meg \frac{\abs{\lambda}}{\abs{\lambda'}}\meg C.
\]
\end{cor}

\begin{proof}
We prove only the second assertion. The first one is proved analogously. Observe that the mapping $\overline{\Omega'}\ni \lambda \mapsto \langle \lambda, e_\Omega \rangle \in \R_+$ is (positively) homogeneous of degree $1$ and vanishes only at $0$ (since $e_\Omega\in\Omega=(\Omega')'$). Therefore, there is a constant $C_1>0$ such that
\[
\frac{1}{C_1}\abs{\lambda}\meg \langle \lambda, e_\Omega \rangle\meg C_1 \abs{\lambda}
\]
for every $\lambda\in \overline{\Omega'}$. The assertion follows from Lemma~\ref{lem:48}.
\end{proof}

\begin{lem}\label{lem:50}
There is a constant $C>0$ such that
\[
\frac{1}{\sqrt r}\abs{t\cdot e_\Omega}\meg \norm{t\,\cdot\,}\meg C \abs{t\cdot e_\Omega}
\]
and such that
\[
\norm{t\,\cdot\,}^{-1}\meg \norm{t^{-1}\,\cdot\,}\meg C\frac{\norm{t\,\cdot\,}^{r-1}  }{\Delta^{\vect 1_r}(t)}
\]
for every $t\in T_+$.
\end{lem}

An analogous result hold for $\Omega'$.

The first assertion has been proved, with $C=1$, in~\cite[Lemma 2.12]{BekolleBonamiGarrigosRicci} for the case of symmetric cones.

\begin{proof}
Since $\abs{e_{\Omega}}=\sqrt r$, it is clear that $\frac{1}{\sqrt r}\abs{t\cdot e_\Omega}\meg \norm{t\,\cdot\,}$. Analogously, it is clear that $\norm{t\,\cdot\,}^{-1}\meg \norm{t^{-1}\,\cdot\,}$.
Then, keep the notation of Section~\ref{sec:2}, and observe that there is a constant $C_1>0$ such that
\[
\abs{a a'}\meg C_1\abs{a}\abs{a'}
\]
for every $a,a'\in A$. Therefore, for every $x\in V$ and for every $t\in T_+$,
\[
\abs{t\cdot x}=\abs{t x t^*}\meg C_1^2 \abs{t}\abs{x}\abs{t^*}=C_1^2\abs{x} \abs{t}^2
\]
since clearly $\abs{t}=\abs{t^*}$. Therefore, $\norm{t\,\cdot\,}\meg C_1\abs{t}^2$. Now, 
\[
\abs{t}^2=\tr(t t^*)=\tr(t\cdot e_\Omega)=\langle t\cdot e_\Omega , e_\Omega\rangle\meg \sqrt r \abs{t\cdot e_\Omega},
\] 
whence the first assertion.

Finally, observe that 
\[
\begin{split}
t^{-1}&=\sum_{k=1}^r (-1)^{k-1} \sum_{j_1<\cdots<j_k} \frac{ t_{j_1,j_2}\cdots t_{j_{k-1}j_k}    }{\prod_{\ell=1}^k \langle e'_{j_\ell},t\rangle  } \\
&=\frac{1}{\Delta^{\vect 1_r/2}(t)}\sum_{k=1}^r (-1)^{k-1} \sum_{j_1<\cdots<j_k}  t_{j_1,j_2}\cdots t_{j_{k-1}j_k}    \prod_{\ell\not \in \Set{j_1,\dots,j_k}} \langle e'_{j_\ell},t\rangle  
\end{split}
\]
for every $t\in T_+$. Observe that $\abs{e_j}=1$, so that $0< \langle e'_j,t\rangle=\tr(e_j t)\meg \abs{t}$ for every $t\in T_+$ and for every $j=1,\dots, r$. Therefore, 
\[
\abs{t^{-1}}\meg\frac{(1+C_1)^r-1 }{C_1^2\Delta^{\vect 1_r/2}(t)  }\abs{t}^{r-1},
\] 
so that
\[
\norm{t^{-1}\,\cdot\,}\meg  \frac{((1+C_1)^r-1)^2 }{C_1^2 \Delta^{\vect 1_r}(t)} \abs{t}^{2(r-1)}
\]
and the result since $\abs{t}^2\meg \sqrt r \abs{t\cdot e_\Omega}\meg r\norm{t\,\cdot\,}$.
\end{proof}

\begin{deff}
Take $\delta>0$, $R>1$, and let $(h_k)_{k\in K} $ be a family of elements of $\Omega$. Then, we say that $(h_{k})$ is a $(\delta,R)$-lattice\index{Lattice!on a homogeneous cone} if the following hold:
\begin{itemize}
\item the balls $B_\Omega(h_k,\delta)$ are pairwise disjoint;

\item the balls  $B_\Omega(h_k,R\delta)$ cover $D$.
\end{itemize}
We define $(\delta,R)$-lattices on $\Omega'$ in an analogous fashion. 
\end{deff}

Since we want to deal with mixed norm spaces on $D$, we shall define lattices on $D$ more carefully. 
Heuristically, the index set must have a product structure in order to keep track of the mixed norm structure of the functions spaces on $D$ we shall consider.
In order to simplify the computations, we shall also require that each lattice on $D$ induce, in a rather straightforward way, a lattice on $\Omega$. More general choices are possible.

\begin{deff}
Take $\delta>0$, $R>1$, and let  $(\zeta_{j,k},z_{j,k})_{j\in J,k\in K} $ be a family of elements of $D$. Then, we say that $(\zeta_{j,k},z_{j,k})$ is a $(\delta,R)$-lattice if the following hold:\index{Lattice!on a homogeneous Siegel domain}
\begin{itemize}
\item  there is a $(\delta,R)$-lattice $(h_k)_{k\in K}$   on $\Omega$ such that
$h_k := \Im  z_{j,k}-\Phi(\zeta_{j,k})$  for every  $j\in J$ and every $k\in K$; 
\item the balls $B((\zeta_{j,k},z_{j,k}),\delta)$ are pairwise disjoint;
\item the balls  $B((\zeta_{j,k},z_{j,k}),R\delta)$ cover $D$.
\end{itemize}
\end{deff}

\begin{lem}\label{lem:32}
Take $\delta>0$. Then, there is a $(\delta,4)$-lattice on $D$. 
\end{lem}

Consequently, there are $(\delta,4)$-lattices on $\Omega$ and $\Omega'$ (actually, one may prove that there are $(\delta,2)$-lattices on $\Omega$ and $\Omega'$, as in~\cite[Lemma 2.6]{BekolleBonamiGarrigosRicci}).

The argument is classical. See, for example,~\cite[Lemma 2.6]{BekolleBonamiGarrigosRicci}.

\begin{proof}
Let $J$ be a maximal subset of $b D+i e_\Omega$  such that the balls $B((\zeta,z),\delta)$, for $(\zeta,z)\in J$, are pairwise disjoint. Observe that, by maximality, 
\[
b D+i e_\Omega\subseteq \bigcup_{(\zeta,z)\in J} B((\zeta,z),2 \delta).
\]
Next, choose a maximal subset $K$ of $T_+$ such that the sets 
\[ 
U  (b D+i t \cdot e_\Omega,\delta)\coloneqq \bigcup_{(\zeta,z)\in b D+i t \cdot e_\Omega} B((\zeta,z),\delta),
\] 
for $t\in K$, are pairwise disjoint, so that, by maximality (and the $b D$-invariance of $d$), 
\[
D=\bigcup_{t\in K} 
U  (b D+i t \cdot e_\Omega,2\delta).
\]
Then, for every $k\in K$ choose $g_k\in GL(E)$ such that $k\cdot \Phi=\Phi\circ (g_k\times g_k)$, and define $(\zeta_{j,k},z_{j,k})\coloneqq (g_k \times k) j$ for every $(j,k)\in J\times K$. 
Let us prove that the balls $B((\zeta_{j,k},z_{j,k}),\delta) $ are pairwise disjoint. 
Indeed, take $(j,k),(j',k')\in J\times K$ such that $(j,k)\neq (j',k')$. If $k\neq k'$, then 
\[
B((\zeta_{j,k},z_{j,k}),\delta)\subseteq  U  (b D+ i k\cdot
e_\Omega,\delta) \; \text{and} \;
B((\zeta_{j',k'},z_{j',k'}),\delta)\subseteq U  (b D+ i k'\cdot e_\Omega,\delta),
\]
so that $B((\zeta_{j,k},z_{j,k}),\delta)\cap B((\zeta_{j',k'},z_{j',k'}),\delta)=\emptyset$.
If, otherwise, $k=k'$, then $j\neq j'$, so that
\[
B((\zeta_{j,k},z_{j,k}),\delta)=(g_k,k)(B(j,\delta))\quad \text{and}\quad B((\zeta_{j',k},z_{j',k}),\delta)=(g_{k},k)(B(j',\delta))
\]
are disjoint. 
Finally, observe that $U (b D+i e_\Omega,2\delta)\subseteq \bigcup_{j\in J} B(j,4\delta) $, so that the balls $B((\zeta_{j,k},z_{j,k}),4 \delta) $ cover $D$. Then, $(\zeta_{j,k},z_{j,k})_{j\in J, k\in K}$ is a $(\delta,4)$-lattice.
\end{proof}

\begin{prop}\label{prop:56}
Take $\delta_0>0$ and $R_0>1$. Then, there is $N\in \N$ such that, for every  $(\delta,R)$-lattice  $(\zeta_{j,k},z_{j,k})_{j\in J,k\in K} $ on $D$ with $\delta\in ]0,\delta_0]$ and $R\in ]1, R_0]$, every ball $B((\zeta_{j,k},z_{j,k}),R\delta)$ intersects at most $N$ balls $B((\zeta_{j',k'},z_{j',k'}),R\delta)$, $(j',k')\in J\times K$.
\end{prop}

The argument is classical. See, for example,~\cite[Lemma 2.6]{BekolleBonamiGarrigosRicci}.
\CR Notice that a similar result holds for lattices on $\Omega$ and $\Omega'$.\CB

\begin{proof}
Indeed, fix $(j,k)\in J\times K$, and let $L$ be the set of $(j',k')\in J\times K$ such that $B((\zeta_{j,k},z_{j,k}),R\delta)\cap B((\zeta_{j',k'},z_{j',k'}),R\delta)\neq \emptyset$. Then, 
\[
\bigcup_{(j',k')\in L} B((\zeta_{j',k'},z_{j',k'}),R\delta)\subseteq B((\zeta_{j,k},z_{j,k}), 2 R \delta),
\]
so that, by Proposition~\ref{prop:26},
\[
\card(L) \nu_D(B((0,i e_\Omega),\delta))\meg \nu_D(B((0, i e_\Omega),2 R \delta)).
\]
Thus, it suffices to take $N$ greater than $ \sup\limits_{0<\delta\meg \delta_0  } \frac{ \nu_D(B((0, i e_\Omega),2 R_0 \delta))  }{ \nu_D(B((0, i e_\Omega),\delta))}$, which is finite since $\nu_D$ has a smooth density with respect to $\Hc^{2 n+2 m}$.
\end{proof}

\section{Notes and Further Results}\label{sec:4:5}

\paragraph{2.6.1}  The study of generalized power functions on homogeneous cones and of the associated gamma functions was first developed in~\cite{Gindikin}, to which the reader is referred for some further developments which are not treated here. 

Of particular interest is the study of the generalized Riemann--Liouville potentials $I^{\vect s}_\Omega$. Here we discuss only the case $\vect s\in - \N_{\Omega'}$ for simplicity, so that $I^{\vect s}_\Omega$ is supported on $\Set{0}$. 
As we noted in Proposition~\ref{prop:38},  the differential operator $\Ic^{\vect s}_\Omega$ with convolution kernel $I^{\vect s}_\Omega$ has an explicit fundamental solution, namely $I^{-\vect s}_\Omega$. 
In addition, since $I^{-\vect s}_\Omega$ is supported in $\overline \Omega$, convolution with $I^{-\vect s}_\Omega$ with distributions supported on a translate of $\overline \Omega$ is always possible, thus giving an inverse to the operator $\Ic^{\vect s}_\Omega$ on the space of such distributions. 
Nonetheless, it is sometimes relevant to solve the equation $\Ic^{\vect s}_\Omega f=g$ for more general $g$, especially in a constructive way (e.g., for $g$ analytic). 

If $\Omega=\R_+^*$, then subtracting a suitable Taylor polynomial to $g$ is often sufficient to define a formal perturbation of $g* I^{-\vect s}_\Omega$ (which need  not be defined) which gives a solution of the equation $\Ic^{\vect s}_\Omega f=g$. 
For more general $\Omega$, the situation is more complicated, since partial `indefinite integrations' over suitable subcones of $\overline \Omega$ are needed. Cf.~\cite{Gindikin} for further details.

\paragraph{2.6.2} \label{sec:2:6:2} Another interesting problem concerning the Riemann--Liouville potentials is the determination of the set $\Gc(\Omega)$ of $\vect s\in\R^r$ for which $I^{\vect s}_\Omega$ is a positive measure (this set is called the `Gindikin--Wallach set'\index{Gindikin--Wallach set} in~\cite{Ishi2}). From Proposition~\ref{prop:38}  it follows that  $\Gc(\Omega)$   contains $\frac 1 2 \vect m+(\R_+^*)^r$, but  in fact  it is strictly larger. The set  $\Gc(\Omega)$   has been first described by Gindikin~\cite{Gindikin2} for general homogeneous cones. See also~\cite{VergneRossi,Ishi2,ArazyUpmeier} for more details.

Here we content ourselves with a brief description of $\Gc(\Omega)$, based on~\cite{Ishi2}. For every $ \mbeps \in \Set{0,1}^r$, define
\[
\vect m^{( \mbeps)}\coloneqq \left( \sum_{k>j} \eps_k m_{j,k} \right)_{j=1,\dots, r},
\]
so that $\vect m= \vect m^{(\vect1_r)}$.
Then, endowing $\R^r$ with the componentwise product,
\[
\Gc(\Omega)=\bigcup_{\mbeps\in \Set{0,1}^r} \left( \frac 1 2 \vect m^{(\mbeps)}+\mbeps(\R_+^*)^r \right),
\]
where the union is disjoint. In addition, $\overline \Omega$ splits into the disjoint union of its $T_+$-orbits $\Oc_\mbeps$, for $\mbeps\in \Set{0,1}^r$, defined so that
\[
\Oc_\mbeps= T_+\cdot \sum_{j=1}^r \eps_j e_j,
\]
with the notation of Section~\ref{sec:2}. Then, $\Oc_{\vect 1_r}=\Omega$. Finally, for every $\vect s\in  \frac 1 2 \vect m^{(\mbeps)}+\mbeps(\R_+^*)^r$ the positive measure $I^{\vect s}_\Omega$ is concentrated on $\Oc_\mbeps$ and can be explicitly described.

Similar results hold for $\Omega'$.

\paragraph{2.6.3}  Using the results of~2.6.2, we are now able to give a necessary and sufficient condition for $\Lambda_+\setminus \Omega'$ to be negligible.

\begin{prop}\label{prop:80}
	Take $\vect s\in \Gc(\Omega)$ (cf.~2.6.2). Then, the following conditions are equivalent:
	\begin{enumerate}
		\item[\em(1)] $\vect s\in (\R_+^*)^r$;
		
		\item[\em(2)] the closed convex hull of $\Supp{I^{\vect s}_\Omega}$ is $\overline \Omega$.
	\end{enumerate}
\end{prop}

The proof is based on~\cite[Theorem 1.6]{Garrigos}, which deals with the case in which $\Omega$ is an irreducible symmetric cone. 

\begin{proof}
	(1) $\implies$ (2). Keep the notation of Section~\ref{sec:2}. 
 Take $\mbeps\in \Set{0,1}^r$ such that $I^{\vect s}_\Omega$ is concentrated on $\Oc_\mbeps =T_+\cdot e_{\mbeps}$, where $e_\mbeps\coloneqq\sum_{j'=1}^r \eps_{j'} e_{j'}$. In other words, $\vect s\in \frac 1 2 \vect{m}^{(\mbeps)}+\mbeps (\R_+^*)^r$.
 Take $j\in \Set{1,\dots,r}$ and let us prove that $e_j\in\overline{\Oc_\mbeps }$. 
	If $\eps_j=1$, then it is clear that $( e_{\vect 1_r}/k+e_j)\cdot e_\mbeps\to e_j$ for $k\to \infty$, so that $e_j\in\overline{\Oc_\mbeps }$.
	If, otherwise, $\eps_j=0$, then $m^{(\mbeps)}_j>0$, so that there is $k\in \Set{j+1,\dots, r}$ such that $m_{j,k}>0$ and $\eps_k=1$. Then, take  $t\in A_{j,k}$ so that $\abs{t}=1$, and observe that $( e_{\vect 1_r}/k+t)\cdot e_\mbeps\to t t^*=e_j$. 
	
	Therefore, $T_+\cdot e_j\subseteq \overline{\Oc_\eps}$ for every $j=1,\dots, r$. Since
	\[
	\Omega=T_+\cdot e_{\vect 1_r}\subseteq \sum_{j=1}^r T_+\cdot e_j,
	\]
	we then see that $\Omega$ is contained in the closed convex hull of $\Oc_\eps$, whence~(2).
	
	(2) $\implies$ (1). Assume by contradiction that $s_j=0$ for some $j\in \Set{1,\dots,r}$. Then, $\eps_j=0$ and $0=m^{(\mbeps)}_j=\sum_{k>j} \eps_k m_{j,k}$, so that for every $k\in \Set{j,\dots,r}$ either $\eps_k=0$ or $A_{j,k}=\Set{0}$. Consequently, 
	\[
	\langle e'_j, t\cdot e_{\mbeps}\rangle=\sum_{k\Meg j}\eps_k \langle e'_j, t_{j,k}(t_{j,k})^*\rangle=0
	\]
	for every $t\in T_+$. It then follows that the closed convex hull of $\Oc_\mbeps$ is contained in the kernel of $e'_j$: contradiction.
\end{proof}

\begin{cor}
	The following conditions are equivalent:
	\begin{enumerate}
		\item[\em(1)] $\Lambda_+\setminus \Omega'$ is negligible;
		
		\item[\em(2)] $\vect b\in (\R_-^*)^r$.
	\end{enumerate}
\end{cor}

\begin{proof}
	This follows from the fact that $\overline {\Lambda_+}$ is the polar of  $\Supp{\Phi_*(\Hc^{2 n})}=\Supp{I^{-\vect b}_\Omega}=\Oc_{-\vect b}$ (cf.~Proposition~\ref{prop:59}) and from Proposition~\ref{prop:80}.
\end{proof}

\chapter[Weighted Bergman Spaces I]{Weighted Bergman Spaces: Sampling, Atomic Decomposition, and Duality}\label{sec:6}

In this chapter we develop at length the theory of weighted Bergman spaces on homogeneous Siegel domains of type II. We recall that the provisional notation of Sections~\ref{sec:1:3} and~\ref{sec:1:4} will no longer be used.

In Section~\ref{sec:6:1}, we shall introduce our notation for the weighted Bergman spaces $A^{p,q}_{\vect s}(D)$, and prove some inclusions (Propositions~\ref{prop:5} and~\ref{prop:22}) between them. We shall also characterize the values of $p,q,\vect s$ for which the space $A^{p,q}_{\vect s}(D)$ is  non-trivial  (Proposition~\ref{prop:11}) and prove some (weak) density results between the various weighted Bergman spaces (Proposition~\ref{prop:8}).

In Section~\ref{sec:6:2}, we shall translate in the new notation the Paley--Wiener theorems developed in Section~\ref{sec:1:4}, determine an explicit expression of the reproducing kernel of $A^{2,2}_{\vect s}(D)$ (Proposition~\ref{prop:6}), and determine its reproducing properties on the various weighted Bergman spaces (Proposition~\ref{prop:9}). 
By means of these results, we shall then extend the definition of the Bergman spaces $A^{2,2}_{\vect s}(D)$ beyond the critical index $\frac 1 4 \vect m$ and prove Paley--Wiener theorems for such spaces (Proposition~\ref{prop:54}).

In Section~\ref{sec:6:3}, we shall deal with sampling results (Theorems~\ref{teo:3} and~\ref{teo:8}), with the aid of the results of Section~\ref{sec:3:3} 

In Section~\ref{sec:6:4}, we shall deal with atomic decomposition for the spaces $A^{p,q}_{\vect s}(D)$. As mentioned in the introduction, by atomic decomposition we mean the possibility of expressing each element $f$ of $A^{p,q}_{\vect s}(D)$ as a sum of the form 
\[
\sum_{j,k} \lambda_{j,k} B^{\vect{s'}}_{(\zeta_{j,k},z_{j,k})} \Delta_\Omega^{(\vect b+\vect d)/p-\vect s-\vect{s'}}(h_k),
\]
where $(\zeta_{j,k},z_{j,k})$ is a lattice on $D$, $(h_k)$ is the associated lattice on $\Omega$, and 
\[
\norm{f}_{A^{p,q}_{\vect s}(D)}\sim \norm{\lambda}_{\ell^{p,q}}
\]
(cf.~Definition~\ref{36}).
Thus, the validity of atomic decomposition consists (virtually) of two halves: the continuity of the map 
\[
\Psi\colon\lambda \mapsto \sum_{j,k} \lambda_{j,k} B^{\vect{s'}}_{(\zeta_{j,k},z_{j,k})} \Delta_\Omega^{(\vect b+\vect d)/p-\vect s-\vect{s'}}(h_k)
\]
from $\ell^{p,q}$ into $A^{p,q}_{\vect s}(D)$, and the surjectivity of $\Psi$ (which implies that $\Psi$ induces an isomorphism of a quotient of $\ell^{p,q}$ onto $A^{p,q}_{\vect s}(D)$). 
We then say that property $\atomic^{p,q}_{\vect s}$ holds if $\Psi$ is continuous, and that property $\atomics^{p,q}_{\vect s}$ holds if $\Psi$ is continuous and onto (hence a strict morphism, by the open mapping theorem). 
Recall that, if $X$ and $Y$ are topological vector spaces, then a strict morphism $T\colon X\to Y$ is a continuous linear mapping which induces an isomorphism $X/\ker T\to T(X)$.\index{Strict morphism}
In order to give sufficient conditions for the validity of properties $\atomic^{p,q}_{\vect s}$ and $\atomics^{p,q}_{\vect s}$, it is convenient to prove the validity of the stronger property $\atomic^{p,q}_{\vect s,+}$, that is, the continuity of the mapping 
\[
\lambda \mapsto \sum_{j,k} \lambda_{j,k} \abs{B^{\vect{s'}}_{(\zeta_{j,k},z_{j,k})}} \Delta_\Omega^{(\vect b+\vect d)/p-\vect s-\vect{s'}}(h_k)
\]
from  $\ell^{p,q}$ into $A^{p,q}_{\vect s}(D)$ (Theorem~\ref{teo:5}). It turns out that property $\atomic^{p,q}_{\vect s,+}$ implies not only property $\atomic^{p,q}_{\vect s}$, as one would expect, but also property $\atomics^{p,q}_{\vect s}$ (Theorem~\ref{teo:4}).

In Section~\ref{sec:6:5}, we show some interesting connections between property $\atomics^{p,q}_{\vect s}$ and the characterization of the dual of $A^{p,q}_{\vect s,0}(D)$ (cf.~Definition~\ref{34}). On the one hand, if property $\atomics^{p,q}_{\vect s}$ holds, then the sesquilinear form on $A^{p,q}_{\vect s,0}(D)\times A^{p',q'}_{(\vect b+\vect d)/\min(1,p)-\vect s-\vect{s'}}$
\[
(f,g)\mapsto \int_D f(\zeta,z) \overline{g(\zeta,z)} \Delta^{-\vect{s'}}_\Omega(\Im z-\Phi(\zeta))\,\dd \nu_D(\zeta,z)
\]
induces an isomorphism of $A^{p',q'}_{(\vect b+\vect d)/\min(1,p)-\vect s-\vect{s'}}$ onto $A^{p,q}_{\vect s,0}(D)'$ (Proposition~\ref{prop:13}). On the other hand, Theorem~\ref{teo:3} implies that the space $A^{p,q}_{\vect s,0}(D)'$ satisfies an analogue of property $\atomics^{p',q'}_{(\vect b+\vect d)/\min(1,p)-\vect s-\vect{s'},\vect{s'}}$ (Proposition~\ref{prop:13}). In particular, property $\atomics^{p',q'}_{(\vect b+\vect d)/\min(1,p)-\vect s-\vect{s'},\vect{s'}}$ holds if and only if the space $A^{p',q'}_{(\vect b+\vect d)/\min(1,p)-\vect s-\vect{s'}}(D)$ is canonically isomorphic to $A^{p,q}_{\vect s,0}(D)'$.

\medskip

We keep the hypotheses and the notation of Chapter~\ref{sec:3}.

\section{Weighted Bergman Spaces}\label{sec:6:1}

In this section we define (mixed norm) weighted Bergman spaces, and present their basic properties. 
In particular, we characterize the values of $\vect s$ for which $A_{\vect s}^{p,q}(D)$ is not trivial (Proposition~\ref{prop:11}) and we prove some (weak) density results (Proposition~\ref{prop:8}).

Throughout the chapter, we shall often indicate some minor modifications needed to treat the case $\max(p,q)=\infty$. As in Chapter~\ref{sec:3}, it turns out that the conditions on $\vect s$ are more uniform for the spaces $A^{p,q}_{\vect s,0}(D)$ (whose elements vanish at $\infty$ in a suitable sense), even though the larger spaces $A^{p,q}_{\vect s}(D)$ are somewhat more natural. 
In order to give a general treatment of duality between weighted Bergman spaces, though, it is worthwhile dealing thoroughly with both spaces.

\begin{deff}
Take $\vect s\in \R^r$ and $p,q\in ]0,\infty]$. We define the weighted Bergman spaces (cf.~Definition~\ref{26})\index{Bergman space!weighted}
\[
A^{p,q}_{\vect s}(D)\coloneqq \Hol(D)\cap L^{p,q}_{\vect s}(D)\qquad \text{and} \qquad A^{p,q}_{\vect s,0}(D)\coloneqq \Hol(D)\cap L^{p,q}_{\vect s,0}(D),
\]
with the topology induced by $L^{p,q}_{\vect s}(D)$.\label{34}
\end{deff}

Recall that $L^{p,q}_{\vect s,0}(D)$ is the closure of $C_c(D)$ in
$L^{p,q}_{\vect s}(D)$, and differs from $L^{p,q}_{\vect s}(D)$ (if
and) only if $\max(p,q)=\infty$.

Observe that, with the notation
of Chapter~\ref{sec:1}, $A^{p,q}_{\vect s}(D) = A^{p,q}_\mi(D)$, where
\[
\mi = \Delta_\Omega^{q\vect s+\vect d}\cdot \Hc^m .
\]
Moreover, the \emph{unweighted} Bergman spaces are the spaces $A^{p,p}_{-\vect d/p}(D)$.

In the following result we provide simple inclusions between the spaces $A^{p,q}_{\vect s}(D)$. As we shall see in Section~\ref{sec:6:6}, these inclusions are closely related to the Sobolev embeddings of the Besov spaces considered in Chapter~\ref{sec:5}.
We shall deal with the spaces  $A^{p,q}_{\vect s,0}(D)$ in Proposition~\ref{prop:22}.

\begin{prop}\label{prop:5}
Take $\vect{s_1},\vect{s_2}\in \R^r$ and $p_1,q_1,p_2,q_2\in ]0,\infty]$ such that 
\[
p_1\meg p_2, \qquad q_1\meg q_2, \qquad \text{and} \qquad \vect{s_2}=\vect{s_1}+\left(\frac{1}{p_2}-\frac{1}{p_1}\right)(\vect{b}+\vect{d}) .
\]
Then, there is a continuous inclusion 
\[
A^{p_1,q_1}_{\vect{s_1}}(D)\subseteq A^{p_2,q_2}_{\vect{s_2} }(D).
\]
\end{prop}

\begin{proof}
We assume that $q_1<\infty$ and leave the (inessential) modifications for the case $q_1=\infty$ to the reader.
By Proposition~\ref{prop:1}, there are a compact subset $K$ of $\Omega$ and a constant $C>0$ such that
\[
\norm{f_{e_\Omega}}_{L^{p_2}(\Nc)} \meg C \left(\int_K \left(\Delta^{\vect{s_1}}_\Omega(h) \norm{f_{ h}}_{L^{p_1}(\Nc)}\right)^{q_1}\,\dd \nu_\Omega(h)\right) ^{1/q_1}
\]
for every $f\in A^{p_1,q_1}_{\nu_\Omega,\loc}(D)$ (with the notation of Chapter~\ref{sec:1:3}, cf.~Definition~\ref{def:5}).
Then, take $t\in T_+$, choose $g\in GL(E)$ so that $t\cdot \Phi =\Phi\circ (g\times g)$, and observe that the $\C$-linear mapping $(\zeta,z)\mapsto (g\zeta,t\cdot z)$ preserves $D$. 
Applying the preceding estimate to the mapping $(\zeta,z)\mapsto f(g\zeta,t\cdot z)$, which clearly belongs to $A^{p_1,q_1}_{\nu_\Omega,\loc}(D)$, we obtain
\[
\begin{split}
&\Delta^{(\vect b+\vect d)/p_2}(t)\norm{f_{t\cdot e_\Omega}}_{L^{p_2}(\Nc)}\\
&\qquad  \qquad\qquad\meg C\Delta^{(\vect b+\vect d)/p_1}(t) \left(\int_K \left(\Delta^{\vect{s_1}}_\Omega(h) \norm{f_{t\cdot h}}_{L^{p_1}(\Nc)}\right)^{q_1} \,\dd \nu_\Omega(h)\right) ^{1/q_1}\\
&\qquad  \qquad\qquad= C\Delta^{(\vect b+\vect d)/p_1-\vect {s_1}}(t) \left(\int_{ t\cdot K}\left(  \Delta^{\vect {s_1}}_\Omega(h) \norm{f_{ h}}_{L^{p_1}(\Nc)}\right) ^{q_1} \,\dd \nu_\Omega(h)\right) ^{1/q_1}.
\end{split}
\]
Thus,
\[
\sup\limits_{h\in \Omega} \abs*{\Delta_{\Omega}^{\vect {s_2}}(h) \norm{f_h}_{L^{p_2}(\Nc)}} \meg C \norm*{f}_{A^{p_1,q_1}_{\vect {s_1}}(D)}, 
\]
and
\[
\left(\int_\Omega \left( \Delta_{\Omega}^{\vect {s_2}}(h) \norm{f_h}_{L^{p_2}(\Nc)}\right) ^{q_1} \,\dd \nu_\Omega(h)\right)^{1/q_1}  \meg C\nu_\Omega(K^{-1})^{1/q_1} \norm*{f}_{A^{p_1,q_1}_{\vect {s_1}}(D)},
\]
where $K^{-1}\coloneqq \Set{t\cdot e_\Omega\colon t^{-1}\cdot e_\Omega\in K}$.
The assertion follows by means of H\"older's inequality. 
\end{proof}

\begin{cor}\label{cor:26}
Take $\vect s\in \R^r$, $p,q\in ]0,\infty]$, and $f\in A^{p,q}_{\vect s}(D) $. Then, the mapping $h\mapsto \norm{f_h}_{L^p(\Nc)}$ is decreasing on $\Omega$. 
\end{cor}

\begin{proof}
Take $(g^{(\eps)})_{\eps>0}$ as in Lemma~\ref{lem:34}, and observe that $g^{(\eps)} f(\,\cdot\,+i h)\in A^{p,\infty}_{\vect 0}(D)$ for every $h\in \Omega$ and for every $\eps>0$, since $f\in A^{p,\infty}_{\vect s}(D)$ by Proposition~\ref{prop:5}. Then, Proposition~\ref{prop:28} shows that
\[
\norm{g^{(\eps)} f_{h+h'}}_{L^p(\Nc)}\meg \norm{g^{(\eps)}f_{h+h''}}_{L^p(\Nc)}
\]
for every $\eps>0$ and for every $h,h',h''\in \Omega$ such that $h'-h''\in \overline \Omega $. Passing to the limit for $\eps>0$, we then obtain
\[
\norm{f_{h+h'}}_{L^p(\Nc)}\meg \norm{f_{h+h''}}_{L^p(\Nc)}
\]
for every $h,h',h''\in \Omega$ such that $h'-h''\in \overline \Omega $.\footnote{Indeed, since $\norm{g^{(\eps)}_h}_\infty =1$ for every $h\in \Omega$, it is readily seen that $\norm{f_h}_{L^p(\Nc)}\meg \liminf\limits_{\eps\to 0^+}\norm{f_h g^{(\eps)}_h}_{L^p(\Nc)}\meg \norm{f_h}_{L^p(\Nc)}$ for every $h\in \Omega$.} Then assertion follows by the arbitrariness of $h$.
\end{proof}

In Section~\ref{sec:6:5}, we shall investigate the conditions on $\vect s$ and $\vect{s'}$ which ensure that the sesquilinear form defined in the following result induces an antilinear isomorphism of $ A^{p',q'}_{\vect{s'}}(D)$ onto $A^{p,q}_{\vect s,0}(D)'$.

\begin{cor}\label{cor:13}
Take $\vect s,\vect{s'}\in \R^r$ and $p,q\in ]0,\infty]$. Then, the sesquilinear form on $A^{p,q}_{\vect s}(D)\times A^{p',q'}_{\vect{s'}}(D)$
\[
(f,g)\mapsto \int_D f(\zeta,z) \overline{g(\zeta,z)} \Delta_\Omega^{\vect s+\vect{s'}-(\vect b+\vect d)/\min(1,p) }(\Im z-\Phi(\zeta))\,\dd \nu_D (\zeta,z)
\]
is well defined and continuous.
\end{cor}

\begin{proof}
It suffices to observe that $A^{p,q}_{\vect s}(D)\subseteq A^{\max(1,p),\max(1,q)}_{\vect s-(1/p-1)_+(\vect b+\vect d)}(D)$ continuously by Proposition~\ref{prop:5}.
\end{proof}

We now characterize the values of $p,q, \vect s$ for which
$A^{p,q}_{\vect s}(D)$ is non-trivial. 

\begin{prop}\label{prop:11}
 Take $\vect s\in \R^r$ and $p,q\in ]0,\infty]$. Then,   $A^{p,q}_{\vect s,0}(D)\neq \Set{0}$ (resp.\   $A^{p,q}_{\vect s}(D)\neq \Set{0}$) 
if and only if  $\vect s\in \frac{1}{2 q}\vect m+(\R_+^*)^r$ (resp.  $\vect s\in \R_+^r$ if $q=\infty$). 
\end{prop}

The proof is based on~\cite[Proposition 3.8]{Bekolleetal}, which deals with the case in which $p,q\Meg 1$, $\vect s\in \R\vect 1_r$ and $D$ is an irreducible symmetric tube domain.

\begin{proof}
\textsc{Step I.} Assume that $\vect s\in \frac{1}{2 q}\vect m+(\R_+^*)^r$ (resp.\ $\vect s\in \R_+^r$ if $q=\infty$), and observe that the functions $g^{(\eps)}$, $\eps>0$, of Lemma~\ref{lem:34} belong to $A^{p,q}_{\vect s,0}(D)$ (resp.\ $A^{p,q}_{\vect s}(D)$).  Thus, $A^{p,q}_{\vect s,0}(D)\neq \Set{0}$ (resp.\ $A^{p,q}_{\vect s}(D)\neq \Set{0}$). 

\textsc{Step II.} Assume first that $q<\infty$ and that $\vect s\not\in \frac{1}{2 q}\vect m+(\R_+^*)^r$, and take $f\in A^{p,q}_{\vect s}(D)$. Then, Corollary~\ref{cor:26} implies that the mapping $h\mapsto \norm{f_h}_{L^p(\Nc)}$ is decreasing on $\Omega$, hence everywhere $0$ since $\Delta_\Omega^{q\vect s } \cdot \nu_\Omega$ does not induce a Radon measure on $\overline \Omega$ (by Proposition~\ref{prop:58}). 

Then, assume that $q=\infty$. Keep the notation of Chapter~\ref{sec:2}.
If $\vect s\in \R_+^r$ but $s_j=0$ for some $j\in \Set{1,\dots,r}$, then define $t_k\coloneqq e_j/k+\sum_{j'\neq j} e_{j'}$, so that $\Delta^{\vect s}(t_k\cdot e_\Omega)=1$ for every $k\in \N$, and $ t_k\cdot e_\Omega$ converges to the point at infinity of $\Omega$ as $k\to \infty$. 
If $f\in A^{p,\infty}_{\vect s,0}(D)$, then 
\[
\lim_{k\to \infty} \norm{f_{t_k\cdot e_\Omega }}_{L^p(\Nc)}=0,
\]
which implies that $f=0$ since $f$ is holomorphic and the mapping $h\mapsto \norm{f_h}_{L^p(\Nc)}$ is decreasing on $\Omega$.
Finally, assume that $\vect s\not \in \R_+^r$ and take $f\in A^{p,\infty}_{\vect s}(D)$. Take $j\in \Set{1,\dots,r}$ such that $s_j<0$, and define $t_k$ as above, so that $\Delta^{\vect s}(t_k)\to +\infty$ as $k\to \infty$. Arguing as above, we then see that $\norm{f_{t_k\cdot e_\Omega}}_{L^p(\Nc)}\to 0 $ as $k\to \infty$, so that $f=0$. 
\end{proof}

\begin{lem}\label{lem:12}
Take $\vect s\in \R^r$, $p,q\in ]0,\infty]$, and $f\in A^{p,q}_{\vect s,0}(D)$ (resp.\ $f\in A^{p,q}_{\vect s}(D)$). Then, the mapping 
\[
\overline{\Omega} \ni h\mapsto  f(\,\cdot\,+ i h)\in A^{p,q}_{\vect s}(D)
\]
is continuous (resp.\ for the weak topology $\sigma(A^{p,q}_{\vect s}(D), L^{p',q'}_{-\vect s+(1/p-1)_+(\vect b+\vect d)}(D))$).
\end{lem}

\begin{proof}
Indeed, Proposition~\ref{prop:1} and Corollary~\ref{cor:26} show that $\norm{f_h}_{L^p(\Nc)}$ is a decreasing function of $h\in \Omega$ which is also continuous if $f\in A^{p,q}_{\vect s,0}(D)$, so that the assertion is easily established.
\end{proof}

\begin{prop}\label{prop:22}
Take $\vect{s_1},\vect{s_2}\in (\R_+^*)^r$ and $p_1,q_1,p_2,q_2\in ]0,\infty]$ such that
\[
p_1\meg p_2, \qquad q_1\meg q_2, \qquad \text{and} \qquad \vect{s_2}=\vect{s_1}+\left(\frac{1}{p_2}-\frac{1}{p_1}\right)(\vect{b}+\vect{d}) .
\]
Then, there is a continuous inclusion 
\[
A^{p_1,q_1}_{\vect{s_1},0}(D)\subseteq A^{p_2,q_2}_{\vect{s_2},0 }(D).
\]
\end{prop}

\begin{proof}
Take $f\in A^{p_1,q_1}_{\vect{s_1},0}(D)$, and observe that Proposition~\ref{prop:1} implies that $f_h\in L^{p_2}_0(\Nc)$ for every $h\in \Omega$ and that the mapping $\Omega\ni h\mapsto f_h\in L^{p_2}_0(\Nc)$ is continuous. 
Then, take $(g^{(\eps)})_{\eps>0}$ as in Lemma~\ref{lem:34}, and observe that $f(\,\cdot\,+i h) g^{(\eps)}\in A^{p_2,q_2}_{\vect{s_2},0}(D) $ for every $h\in \Omega$ and for every $\eps>0$ (cf.~Proposition~\ref{prop:11}). In addition, using Lemmas~\ref{lem:67} and~\ref{lem:12} we see that $f(\,\cdot\,+i h) g^{(\eps)}$ converges to $f$ in $A^{p_1,q_1}_{\vect{s_1}}(D)$, hence in $A^{p_2,q_2}_{\vect{s_2}}(D)$ by Proposition~\ref{prop:5}, as $h\to 0$ and $\eps\to 0^+$. The assertion follows thanks to Proposition~\ref{prop:5}.
\end{proof}

\begin{cor}\label{cor:14}
Take $\vect s\in \R^r$, $p,q\in ]0,\infty]$, and $f\in A^{p,q}_{\vect s}(D) $. Then, the following conditions are equivalent:
\begin{enumerate}
\item[\em(1)] $f\in A^{p,q}_{\vect s,0}(D)$;

\item[\em(2)] $f_h\in L^p_0(\Nc)$ for every $h\in \Omega$ (and $\Delta_\Omega^{\vect s}(h)\norm{f_h}_{L^p(\Nc)}\to 0$ as $h\to \infty$ if $q=\infty$);

\item[\em(3)] $f_h\in L^p_0(\Nc)$ and the mapping $h\mapsto \Delta_\Omega^{\vect s}(h) f_h$ belongs to $C_0(\Omega;L^p_0(\Nc))$.
\end{enumerate}
\end{cor}

\begin{proof}
{ (1) $\implies$  (3).} This follows from Proposition~\ref{prop:22}.

{ (3) $\implies$ (2).} Obvious.

{ (2) $\implies$ (1).} This follows from Proposition~\ref{prop:50} if $q<\infty$. If, otherwise, $q=\infty$, then it suffices to show that the mapping $\Omega\ni h\mapsto f_h\in L^p_0(\Nc)$ is continuous, but this follows from Corollary~\ref{cor:27}.
\end{proof}

\begin{prop}\label{prop:8}
Take $\vect{s_1}, \vect{s_2}\in \R^r$ and $p_1,p_2,q_1,q_2\in ]0,\infty]$, and assume that $\vect{s_2}\in \frac{1}{2 q_2}\vect m +(\R_+^*)^r  $ (resp.\ $\vect{s_2}\in \R_+^r$ if $q_2=\infty$). Then, 
\[
A^{p_1,q_1}_{\vect{s_1},0}(D)\cap A^{p_2,q_2}_{\vect{s_2},0}(D) \quad \text{(resp.\ $A^{p_1,q_1}_{\vect{s_1},0}(D)\cap A^{p_2,q_2}_{\vect{s_2}}(D)$)} \quad\text{is dense in} \quad A^{p_1,q_1}_{\vect{s_1},0}(D).
\]
In addition,
\[
A^{p_1,q_1}_{\vect{s_1}}(D)\cap A^{p_2,q_2}_{\vect{s_2},0}(D) \quad \text{(resp.\ $A^{p_1,q_1}_{\vect{s_1}}(D)\cap A^{p_2,q_2}_{\vect{s_2}}(D)$)} \quad\text{is dense in} \quad A^{p_1,q_1}_{\vect{s_1}}(D)
\]
for the weak topology $\sigma(A^{p_1,q_1}_{\vect{s_1}}(D), L^{p_1',q_1'}_{-\vect{s_1}+(1/p_1-1)_+(\vect b+\vect d)}(D))$.
\end{prop}

\begin{proof}
Take $f\in A^{p_1,q_1}_{\vect{s_1} ,0}(D)$ and take $g^{(\eps)}$ as in Lemma~\ref{lem:34} for some choice of $\alpha\in ]0,1/2[$. Then, Lemma~\ref{lem:12} implies that $f(\,\cdot\,+i h)g^{(\eps)}\in A^{p_1,q_1}_{\vect{s_1},0}(D)$ for every $h\in \Omega$ and for every $\eps>0$. In addition, set $p_3\coloneqq\frac{p_1 p_2}{p_1-p_2}$ if $p_2<p_1$, and $p_3\coloneqq \infty $  otherwise, so that H\"older's inequality implies that
\[
\norm{(f(\,\cdot \,+i h) g^{(\eps)})_{h'}}_{L^{p_2}(\Nc)}\meg C_{p_3,\eps} \norm{f_{h+h'}}_{L^{\max(p_1,p_2)}(\Nc)} \ee^{-\eps C \abs{h'}^\alpha},
\]
for every $h,h'\in \Omega$ and for every $\eps>0$, with the notation of Lemma~\ref{lem:34}. Now,  Corollary~\ref{cor:26} implies that the mapping 
\[
h'\mapsto  \norm{f_{h'}}_{L^{\max(p_1,p_2)}(\Nc)}
\]
is decreasing on $\Omega$, so that $f(\,\cdot \,+i h) g^{(\eps)}\in A^{p_2,q_2}_{\vect{s_2},0}(D)$ (resp.\ $f(\,\cdot \,+i h) g^{(\eps)}\in A^{p_2,q_2}_{\vect{s_2}}(D)$) for every $h\in \Omega$ and for every $\eps>0$.
Next, it is easily verified that $f(\,\cdot\,+i h) g^{(\eps)}$ converges to $f(\,\cdot\,+i h)$ in $A^{p_1,q_1}_{\vect{s_1},0}(D)$ as $\eps\to 0^+$, while $f(\,\cdot\,+ i h)$ converges to $f$ in $A^{p_1,q_1}_{\vect{s_1},0}(D)$ as $h\to 0$ thanks to Lemma~\ref{lem:12}.
The second assertion is proved similarly.
\end{proof}

\section{Paley--Wiener Theorems}\label{sec:6:2}

In this section we translate some results of Section~\ref{sec:1:4} into the new formalism of this chapter (Proposition~\ref{prop:6}). 
We then show that the weighted Bergman kernels associated with the spaces $A^{2,2}_{\vect s}(D)$ have nice reproducing properties for more general weighted Bergman spaces (Proposition~\ref{prop:9}), and that Riemann--Liouville operators induce isometries between the spaces $A^{2,2}_{\vect s}(D)$, up to a constant (Proposition~\ref{prop:7}). Even though this latter result will be considerably extended in Section~\ref{sec:6:6}, we need Proposition~\ref{prop:7} to give a (somewhat formal but) reasonable extension of Proposition~\ref{prop:6} for $\vect s \not \in \frac 1 2 \vect{m}+ (\R_+^*)^r$ (Proposition~\ref{prop:54}).
We conclude this section presenting some invariance properties of the so-defined space $\widehat A^{2,2}_{(\vect b+\vect d)/2}(D)$ (Proposition~\ref{prop:34}), which may therefore be interpreted as a generalization of the classical Dirichlet space.

\begin{deff}\label{def:10}
For every $\vect s\in \R^r$, define\label{35}
\[
\Lc^2_{\vect s}(\Omega')\coloneqq \Set{ \tau\in\int_{\Omega'}^\oplus \Lin^2(H_\lambda)\Delta_{\Omega'}^{-2 \vect s-\vect b}(\lambda)\,\dd \lambda \colon  \tau=\tau\, P_{\,\cdot\,,0} }.
\]
\end{deff}

\begin{prop}\label{prop:6}
Take $\vect{s}\in \frac{1}{4}\vect{m}+ (\R_+^*)^r$.  Then, there is a unique isometric isomorphism $\Pc\colon A^{2,2}_{\vect s}(D)\to\Lc^2_{\vect s}(\Omega')$ such that 
\[
\Pc f(\lambda) 
  =\sqrt{\frac{ \abs{\Pfaff(e_{\Omega'})} \Gamma_{\Omega}(2 \vect s)}{ 4^{ \vect s}2^{m-n}\pi^{n+m}} } \ee^{\langle \lambda, h\rangle} \pi_\lambda(f_h)
\]
for every $f\in A^{2,2}_{\vect s}(D)$, for every $h\in \Omega$, and for almost every $\lambda\in \Omega'$.\footnote{Recall that $\ee^{\langle \lambda, h\rangle} \pi_\lambda(f_h)$ does not depend on $h$ for almost every $\lambda$, thanks to Proposition~\ref{prop:2}.}
In addition, the reproducing kernel of $A^{2,2}_{\vect s}(D)$ is the
mapping 
\[
K_{\vect s}((\zeta,z),(\zeta',z')) := \CB \frac{  \abs{\Pfaff(e_{\Omega'})}\Gamma_{\Omega'}(2 \vect s-\vect b-\vect{d})  }{4^{m}\pi^{n+m} \Gamma_\Omega(2\vect s)}  B^{\vect{b}+\vect d-2 \vect s}_{(\zeta',z')}(\zeta,z).
\]
\end{prop}

\begin{oss} 
Thus, $K_{\vect s}((\zeta,z),(\zeta',z'))$ is the
weighted Bergman kernel for the weighted Bergman space $A^{2,2}_{\vect
  s}(D) = \Hol(D)\cap L^{2,2}_{\vect s}(D)$ and the inner product is given by
\[
  \begin{split}
\langle f\,\vert\, g\rangle_{A^{2,2}_{\vect s}(D)}
  & = \int_\Omega \int_\Nc f_h(\zeta,x)\overline{g_h(\zeta,x)} \dd
  (\zeta,x)\, \Delta_\Omega^{2{\vect s}+{\vect d}}(h) \dd h\\
  & =
\int_D f(\zeta,z) \overline{g(\zeta,z)}
  \Delta_\Omega^{2\vect s-\vect  b-\vect d}  (\Im z -\Phi(\zeta)) \dd \nu_D (\zeta,z),
\end{split}
\]  
see Definition~\ref{26}.
\end{oss}

\begin{proof}
Observe first that Corollary~\ref{cor:3} shows that
\[
\norm{f}_{A^{2,2}_{\vect s}(D)}^2= \frac{2^{n-m}}{\pi^{n+m}} \int_{\Omega'} \norm*{e^{\langle \lambda, e_\Omega\rangle} \pi_\lambda(f_{e_\Omega})}_{\Lin^2(H_\lambda)}^2 \Lc(\Delta_\Omega^{2 \vect s}\cdot \nu_{\Omega})(2\lambda) \abs{\Pfaff(\lambda)}\,\dd \lambda ;
\]
in addition, Propositions~\ref{prop:58} and~\ref{prop:59} imply that 
\[
\Lc(\Delta_\Omega^{2 \vect s}\cdot \nu_{\Omega})(2\lambda)=4^{-\vect{s'}}\Gamma_\Omega(2 \vect s) \Delta_{\Omega'}^{-2\vect s}(\lambda)
\]
and 
\[
\abs{\Pfaff(\lambda)}=\abs{\Pfaff(e_{\Omega'})}\Delta_{\Omega'}^{-\vect b}(\lambda),
\] 
so that $\Pc$ is an isometry of $A^{2,2}_{\vect s}(D)$ into $\Lc^2_{\vect s}(\Omega')$. In addition, Proposition~\ref{prop:4} implies that $\Pc$ is onto.
For the second assertion, it suffices to observe that 
\[
2 \vect s-\vect b-\vect{d}\in \half \vect{m'}+  \vect m-\vect b+\vect 1_r+(\R_+^*)^r\subseteq \half \vect{m'}+(\R_+^*)^r
\] 
since $\vect m, -\vect b\in \R_+^r$, and to apply Corollary~\ref{cor:5}.
\end{proof}

In the following result we extend the reproducing properties of the weighted Bergman kernels to more general weighted Bergman spaces.

\begin{prop}\label{prop:9}
Take $p,q\in ]0,\infty]$, $\vect s,\vect {s'}\in \R^r$ and $f\in A^{p,q}_{\vect s}(D)$. Assume that the following hold:
\begin{itemize}
\item $\vect s\in \frac 1 p(\vect b+\vect d)+\frac{1}{2 q'}\vect{m'}+(\R_+^*)^r$; 

\item $ \vect s+\vect{s'}\in \frac{1}{\min(1,p)}(\vect b+\vect d)-\frac{1}{2 q'}\vect m-(\R_+^*)^r$ (resp.\ $ \vect s+\vect{s'}\in \frac{1}{\min(1,p)}(\vect b+\vect d)-\R_+^r$ if $ q'=\infty$);

\item $\vect{s'}\in \frac{1}{p'}(\vect b+\vect d)-\frac{1}{2p'} \vect{m'}-(\R_+^*)^r$. 
\end{itemize}
Then, for every $(\zeta,z)\in D$,
\[
f(\zeta,z)=  \int_D f(\zeta',z') K_{(\vect b+\vect d -\vect{s'})/2} ((\zeta,z),
(\zeta',z')) \Delta_\Omega^{-\vect{s'}} (\Im z'-\Phi(\zeta'))\,\dd \nu_D (\zeta',z'). 
\]
\end{prop}

\begin{oss}
Observe that, by Proposition~\ref{prop:6}, the
previous equality can be written as
\[
  f(\zeta,z)= c_{\vect{s'}} \int_D f(\zeta',z') B^{\vect{s'}} _{(\zeta',z')}
(\zeta,z)) \Delta_\Omega^{-\vect {s'}}
(\Im z'-\Phi(\zeta'))\,\dd \nu_D (\zeta',z'),
\] 
where
\[
c_{\vect{s'}}=\frac{\abs{\Pfaff(e_{\Omega'})}\Gamma_{\Omega'}(-\vect{s'})}{4^m
  \pi^{n+m} \Gamma_\Omega (\vect b+\vect d-\vect{s'}) } .
\]

Moreover, since the relations between the various parameters are simpler using the kernels $B^{\vect{s'}}$ rather than the Bergman kernels $K_{\vect s}$, we will privilege the former ones in the rest of the work. 

\end{oss}

The proof is based on~\cite[Lemma 4.4]{RicciTaibleson}, which deals with $\C_+$.
Notice that the conditions on $\vect s $ and $\vect{s'}$ simply ensure that $B^{\vect{s'}}_{(\zeta,z)}\in A^{p',q'}_{(\vect b  +\vect d)/\min(1,p)-\vect s-\vect{s'}}$.

\begin{proof}
Notice that we may assume that $\vect s\in \frac{1}{2 q}\vect m+(\R_+^*)^r$ if $q<\infty$ and that $\vect s\in \R_+^r$ if $q=\infty$, thanks to Proposition~\ref{prop:11}, so that $\vect{s'}\in \vect b+\vect d-\frac 1 2 \vect m-(\R_+^*)^r$.
Then, Proposition~\ref{prop:6} shows that the assertion holds if 
\[
f\in A^{p,q}_{\vect s}(D)\cap A^{2,2}_{(\vect {b}+\vect d- \vect{s'})/2}(D).
\]
In addition,  Proposition~\ref{prop:60} shows that $B^{\vect{s'}}_{(\zeta,z)}\in A^{p',q'}_{(\vect b  +\vect d)/\min(1,p)-\vect s-\vect{s'}}$, so that the assertion follows from Proposition~\ref{prop:8} by continuity. 
\end{proof}

Cf.~Corollary~\ref{cor:29} and Proposition~\ref{prop:30} for a much more general version of the following result. We provide a simple direct proof for the sake of completeness.

\begin{prop}\label{prop:7}
Take $\vect s\in \frac{1}{4}\vect{m}+ (\R_+^*)^r$ and $\vect{s'}\in \N_{\Omega'}$.
Then, convolution by $\sqrt{\frac{4^{ \vect{s'}}\Gamma_\Omega(2 \vect{s})}{ \Gamma_\Omega(2 \vect{s}+2\vect{s'})  } }I_\Omega^{-\vect{s'}}$ induces an isometry of $A^{2,2}_{\vect s}(D)$ onto $A^{2,2}_{\vect s+\vect{s'}}(D)$.
\end{prop}

\begin{proof}
Observe first that Lemma~\ref{lem:13} shows that $\Fc_F(I_\Omega^{-\vect{s'}})$ induces the homogeneous polynomial $i^{\vect{s'}}\Delta_{\Omega'}^{\vect{s'}}$ on $\Omega'$, so that 
\[
\dd \pi_\lambda(I_\Omega^{-\vect{s'}})=  i^{\vect{s'}} \Delta_{\Omega'}^{\vect{s'}}(\lambda) I_{H_\lambda}
\]
for every $\lambda\in F'\setminus W$.
Now, take $f\in A^{2,2}_{\vect s}(D)$, and observe that $f* I_\Omega^{-\vect{s'}}$ is a well-defined holomorphic function on $D$. In addition, if $\eta\in C^\infty_c(F)$, then  
\[
(f* I_\Omega^{-\vect{s'}}*\eta)_h= f_h*(I_\Omega^{-\vect{s'}}*\eta)\in L^2(\Nc)
\]
and 
\[
\pi_\lambda( (f* I_\Omega^{-\vect{s'}}*\eta)_h )= (\Fc_F \eta)(\lambda) i^{\vect{s'}}\Delta_{\Omega'}^{\vect{s'}}(\lambda) \pi_\lambda(f_h)
\]
for every almost $\lambda\in \Omega'$ and for every $h\in \Omega$. By the arbitrariness of $\eta$, this implies that $f*I_{\Omega}^{-\vect {s'}}\in A^{2,2}_{\vect s+\vect{s'}}(D)$ and that
\[
\norm{f*I_{\Omega}^{-\vect {s'}}}_{A^{2,2}_{\vect s+\vect{s'}}(D)}^2=\frac{ \Gamma_\Omega(2 \vect{s}+2\vect{s'})  }{4^{ \vect{s'}}\Gamma_\Omega(2 \vect{s})}\norm{f}_{A^{2,2}_{\vect s}(D)}^2
\]
by Proposition~\ref{prop:6}. Conversely, if $f\in A^{2,2}_{\vect s+\vect{s'}}(D)$, then there is a unique $\widetilde f\in A^{2,2}_{\vect s}(D)$ such that $\pi_\lambda(\widetilde f_h)=i^{-\vect{s'}}\Delta_{\Omega'}^{-\vect s'} \pi_\lambda(f_h)$ for almost every $\lambda\in \Omega'$ and for every $h\in \Omega$, so that  $f=\widetilde f * I_\Omega^{-\vect{s'}}$. 
\end{proof}

\begin{deff}\label{def:8}
Take $\vect s\in \R^r$ and  $\vect s'\in \N_{\Omega'}$ such that $\vect{s}+\vect{s'}\in \frac1 4 \vect{m}+(\R_+^*)^r$. Then, we define $\widehat A^{2,2}_{\vect s,\vect{s'}}(D)$ as  the Hausdorff locally convex space associated with the space of $f\in \Hol(D)$ such that $f*I^{-\vect{s'}}_{\Omega}\in A^{2,2}_{\vect s+\vect{s'}}(D)$, endowed with the corresponding topology. 
\end{deff}

In other words, the elements of $\widehat A^{2,2}_{\vect s,\vect{s'}}(D)$ are equivalence classes of $\Hol(D)$ modulo the elements of $\Hol(D)$ which are annihilated by convolution with $I^{-\vect{s'}}_{\Omega}$.

In the following result we provide an extension of Proposition~\ref{prop:6} to the spaces $\widehat A^{2,2}_{\vect s,\vect{s'}}(D)$. 
We shall see later (cf.~the remarks following Definition~\ref{def:6}) that $\widehat A^{2,2}_{\vect s,\vect{s'}}(D)$ may be canonically identified with a subspace (with a finer topology) of $A^{\infty,\infty}_{\vect s-(\vect b+\vect d)/2}(D)$ when $\vect s \in \frac 1 2(\vect b+\vect d)+\frac 1 4 \vect{m'}+(\R_+^*)^r$.

\begin{prop}\label{prop:54}
Take $\vect s\in\R^r$ and  $\vect s'\in \N_{\Omega'}$ such that $\vect{s}+\vect{s'}\in \frac1 4 \vect{m}+(\R_+^*)^r$. Then, there is a unique isomorphism 
\[
\Psi\colon \widehat A^{2,2}_{\vect s,\vect{s'}}(D)\to\Lc^2_{\vect s}(\Omega')
\]
such that
\[
\pi_\lambda(( f*I^{-\vect{s'}}_\Omega)_h   )= \ee^{-\langle \lambda,h\rangle}  i^{\vect{s'}}\Delta_{\Omega'}^{\vect{s'}}(\lambda)\Psi(f)(\lambda)
\]
for every $f\in \widehat A^{2,2}_{\vect s,\vect{s'}}(D)$, for almost every $\lambda\in \Omega'$, and for every $h\in \Omega$.\footnote{We identify $ f*I^{-\vect{s'}}_\Omega$ with its unique representative in $A^{2,2}_{\vect s+\vect{s'}}(D)$.}

In addition, $\sqrt{\frac{\abs{\Pfaff(e_{\Omega'}) \Gamma_\Omega(2 \vect s+ 2 \vect{s'})}}{4^{\vect s+\vect{s'}}2^{m-n} \pi^{n+m} }}\Psi$ is an isometry and  $\widehat A^{2,2}_{\vect{s},\vect{s'}}(D)$ is a hilbertian space.
\end{prop}

In particular, $\widehat A^{2,2}_{\vect 0,\vect{s'}}(D)$, for $\vect{s'}\in \N_{\Omega'}\cap \left(\frac 1 4 \vect m + (\R_+^*)^r\right)$, is canonically isomorphic to the Hardy space $A^{2,\infty}_{\vect 0}(D)$, thanks to Corollary~\ref{cor:1} and Proposition~\ref{prop:4}.

\begin{proof}
The property of the statement clearly  defines a unique continuous linear mapping $\Psi\colon  \widehat A^{2,2}_{\vect s,\vect{s'}}(D)\to \Lc^2_{\vect s}(\Omega')$ which is an isomorphism onto its image. Thus, it only remains to prove that $\Psi$ is onto. Then, take $\tau\in \Lc^2_{\vect s}(\Omega')$, and fix $\vect{s'}\in \N_{\Omega'}$ such that $\vect{s}+\vect{s'}\in \frac1 4 \vect{m}+(\R_+^*)^r$. 
Then, Proposition~\ref{prop:6} implies that there is a unique $f\in A^{2,2}_{\vect s+\vect{s'}}(D)$ such that 
\[
\pi_\lambda(f_h)= e^{\langle \lambda,h\rangle} i^{\vect{s'}} \Delta_{\Omega'}^{\vect{s'}}(\lambda) \tau(\lambda)
\]
for almost every $\lambda\in \Omega'$ and for every $h\in \Omega$. Since $D$ is a convex open subset of $E\times F_\C$,~\cite[Theorem 9.4]{Treves2} implies that there is  $g\in \Hol(D)$ such that $g* I^{-\vect{s'}}_\Omega=f$. Denoting by $\dot g$ the equivalence class of $g$ in $\widehat A^{2,2}_{\vect s,\vect{s'}}(D)$, it follows that $\Psi(\dot g)=\tau$, whence the conclusion.
\end{proof}

\begin{prop}\label{prop:34}
Take $\vect{s'}\in \N_{\Omega'}$ such that $\frac 1 2 (\vect b+\vect d)+\vect{s'}\in \frac 1 4 \vect m+(\R_+^*)^r$.
Then, the space $\widehat A^{2,2}_{(\vect b+\vect d)/2,\vect{s'}}(D)$ (and its hilbertian norm) is invariant under composition with the \emph{affine} automorphisms of $D$.
\end{prop}

The precise way in which the `composition' of the elements of $\widehat A^{2,2}_{(\vect b+\vect d)/2,\vect{s'}}(D)$ with the affine automorphisms of $D$ has to be intended will be defined in the proof.

Cf.~\cite[Theorem 5.5]{ArcozziMonguzziPelosoSalvatori} for the case $\Omega=\R_+^*$, in which case an explicit knowledge of the group of holomorphic automorphisms leads to a proof of the invariance of $\widehat A^{2,2}_{(\vect b+\vect d)/2,\vect{s'}}(D)$ under the group of all holomorphic automorphisms of $D$.

\begin{proof}
In order to simplify the proof, we shall denote by $\widehat A^{2,2}_{(\vect b+\vect d)/2}(D)$ the space $\widehat A^{2,2}_{(\vect b+\vect d)/2,\vect{s'}}(D)$ endowed with the norm induced by the mapping $\Psi$ of Proposition~\ref{prop:54}, which is a multiple of the norm of $\widehat A^{2,2}_{(\vect b+\vect d)/2,\vect{s'}}(D)$.

Observe first that, by~\cite[Propositions 2.1 and 2.2]{Murakami}, the group of affine automorphisms of $D$ is the semidirect product of the group $G_1$ of left translations by elements of $b D$ and the group
\[
G_2=\Set{g_1\times (g_2)_\C\colon g_1\in GL(E), g_2\in GL(F), g_2(\Omega)=\Omega, g_2 \Phi=\Phi\circ (g_1\times g_1)  }.
\]
Now, take $\tau\in \Lc^2_{\vect s}(\Omega')$ such that $\tau$ has compact support in $\Omega'$, and let $f$ be the unique element of $A^{2,\infty}_{\vect 0}(D)$ such that 
\[
\pi_\lambda(f_h)=\ee^{-\langle \lambda,h\rangle} \tau(\lambda)
\]
for almost every $\lambda\in \Omega'$ and for every $h\in \Omega$ (cf.~Proposition~\ref{prop:4}). Observe that $f$ induces an element of $\widehat A^{2,2}_{(\vect b+\vect d)/2}(D)$. Since they form a dense subspace of $\widehat A^{2,2}_{(\vect b+\vect d)/2 }(D)$, it will suffice to prove our assertion for such $f$. Then mapping $f\mapsto f\circ A$, for $A\in G_1G_2$, may then be extend by continuity to the whole of $\widehat A^{2,2}_{(\vect b+\vect d)/2}(D)$.

Observe that, if $A\in G_1$, then $\pi_\lambda(f_h\circ A)$ equals the composition of $\pi_\lambda(f_h)$ with a  unitary automorphism of $H_\lambda$, so that 
\[
\norm{\pi_\lambda(f_h\circ A)}_{L^2(H_\lambda)}=\norm{\pi_\lambda(f_h)}_{L^2(H_\lambda)}
\]
for almost every $\lambda\in \Omega'$ and for every $h\in \Omega$, whence 
\[
\norm{f_h\circ A}_{\widehat A^{2,2}_{(\vect b+\vect d)/2}(D)}=\norm{f_h}_{\widehat A^{2,2}_{(\vect b+\vect d)/2}(D)}.
\] 

Next, take $A\in G_2$, so that $A=g_1\times(g_2)_\C $ for some $g_1,g_2$ as above. Denote by $G(\Omega)$ the set of elements of $GL(F)$ which preserve $\Omega$, so that $G(\Omega)$ is a closed subgroup of $GL(F)$ which acts transitively on $\Omega$. Let $K$ be the stabilizer of $e_\Omega$ in $G(\Omega)$. Then, $K$ is a compact normal subgroup of $G(\Omega)$ (cf.~\cite[p.~31]{Murakami}). In addition, if $G_T$ denotes the subgroup of $G(\Omega)$ induced by the action of $T_+$ on $\Omega$, then $G(\Omega)$ is the semidirect product of $G_T$ and $K$. 
Since $K$ is compact, each of its elements has determinant $1$ in absolute value, so that there is $t\in T_+$ such that $\abs{\det(g_2)}=\Delta^{-\vect d}(t)$. 
In the same way, we see that $\abs{\det(g_1)}=\Delta^{-\vect b}(t)$ for a suitable choice of $t$. Then, $g_1\times g_2$ is an automorphism of $\Nc$ and $\abs{\det(g_1\times g_2)}=\Delta^{-\vect b- \vect d}(t) $, so that
\[
\begin{split}
\pi_\lambda((f\circ A)_h)&=\pi_\lambda(f_{g_2 h}\circ (g_1\times g_2))\\
&=\Delta^{\vect b+ \vect d}(t)\pi_{\lambda\circ g_2^{-1}}(f_{g_2 h})\\
&=\Delta^{\vect b+ \vect d}(t) \ee^{-\langle \lambda, h\rangle}\tau(\lambda\circ g_2^{-1}).
\end{split}
\] 
Therefore,
\[
\begin{split}
\norm{(f\circ A)}_{\widehat A^{2,2}_{(\vect b+\vect d)/2}(D)}^2&=  \int_{\Omega'} \norm{\tau(\lambda\circ g_2^{-1})}_2^2 \Delta^{2 (\vect b+\vect d)}(t)\Delta_{\Omega'}^{-2\vect b-\vect d}(\lambda)\,\dd \lambda\\
&= \int_{\Omega'} \norm{\tau(\lambda)}_2^2 \Delta^{2 \vect b+\vect d}(t)\Delta_{\Omega'}^{-2\vect b-\vect d}(\lambda\cdot t) \,\dd \lambda\\
&= \norm{f}_{\widehat A^{2,2}_{(\vect b+\vect d)/2}(D)}^2,
\end{split}
\]
whence the result.
\end{proof}

\begin{oss}
If  $D$ is an irreducible symmetric Siegel domain of type II and $\frac{m}{r}\in \N$, then a `generalized Dirichlet space' $\Dc$ has been defined in~\cite[pp.~223--224]{ArazyUpmeier}. Such space is  invariant (with its norm) under the (suitably defined) composition with all the elements of $G$, where $G$ is the component of the identity in the group of holomorphic automorphisms of $D$. In~\cite[Theorem 5.2]{Arazy} it is also proved that, if $D$ is irreducible and symmetric and $H$ is a complete prehilbertian space of holomorphic functions in which the constant functions are adherent to $\Set{0}$, which is invariant (with its norm) under the composition with all the elements of $G$, and for which the mapping $\phi\mapsto f\circ \phi$ is continuous on a maximal compact subgroup of $G$ for every $f\in H$, then $\frac m r\in \N$ and the Hausdorff space associated with $H$ is canonically isomorphic to $\Dc$.

In particular, if $\frac m r\not \in \N$, the space $\widehat A^{2,2}_{(\vect b+\vect d)/2,\vect{s'}}(D)$ cannot be $G$-invariant in such a natural way. It would be interesting to determine if $\widehat A^{2,2}_{(\vect b+\vect d)/2,\vect{s'}}(D)$ is isomorphic to $\Dc$ when $\frac m r\in \N$.
\end{oss}

\section{Sampling}\label{sec:6:3}

In this section we provide some sampling theorems. In particular, we show that, under suitable assumptions, the spaces $A^{p,q}_{\vect s}(D)$ can be effectively studied by means of a suitable discretization.

\begin{deff}\label{36}
Take $p,q\in ]0,\infty]$, and take two sets $J$ and $K$. Define
\[
\ell^{p,q}(J,K)\coloneqq \Set{\lambda\in \C^{J\times K}\colon ((\lambda_{j,k})_{j\in J})_{k\in K}\in \ell^q(K;\ell^p(J))},
\]
endowed with the corresponding quasi-norm, and define $\ell^{p,q}_0(J,K)$ as the closure of $\C^{(J\times K)}$ in $\ell^{p,q}(J,K)$.
\end{deff}

Then, $\ell^{p,q}(J,K)$ and $\ell^{p,q}_0(J,K)$ are locally bounded $F$-spaces. In addition, $\ell^{p,q}_0(J,K)$ is canonically isomorphic to $\ell^q_0(K;\ell^p_0(J))$ (cf.~Proposition~\ref{prop:50}).

\begin{teo}\label{teo:3}
Take $R_0>1$,  $p,q\in ]0,\infty]$ and $\vect s\in \frac{1}{2 q}\vect m+(\R_+^*)^r$ (resp.\ $\vect s\in \R_+^r$ if $q=\infty$). 
Then, there are $\delta_0>0$ and $C>0$ such that, for every $(\delta,R)$-lattice $(\zeta_{j,k},z_{j,k})_{j\in J,k\in K}$ on $D$ with $\delta\in ]0,\delta_0]$ and  $R\in ]1,R_0]$, defining $h_k\coloneqq \Im z_{j,k}-\Phi(\zeta_{j,k})$ for every $k\in K$ and for some (hence every) $j\in J$, the mapping
\[
S\colon\Hol(D)\ni f \mapsto \left(\Delta_\Omega^{\vect s-(\vect b+\vect d)/p}(h_k)f(\zeta_{j,k},z_{j,k})\right)_{j,k}\in \C^{J\times K}
\]
induces an isomorphism of $A^{p,q}_{\vect s,0}(D)$ onto a (closed) subspace of $\ell^{p,q}_0(J,K)$ (resp.\  of $A^{p,q}_{\vect s}(D)$ onto a (closed) subspace of $\ell^{p,q}(J,K)$) such that 
\[
\frac{1}{C}\norm{f}_{A^{p,q}_{\vect s}(D)}\meg \delta^{(2n+m)/p+m/q} \norm{S f}_{\ell^{p,q}(J,K)}\meg C \norm{f}_{A^{p,q}_{\vect s}(D)}
\]
\CR for every $f\in A^{p,q}_{\vect s,0}(D)$ (resp.\ for every $f\in A^{p,q}_{\vect s}(D)$).  In addition, 
\[
A^{\infty,\infty}_{\vect s-(\vect b+\vect d)/p}(D)\cap S^{-1}(\ell^{p,q}_0(J,K))\subseteq A^{p,q}_{\vect s,0}(D)\qquad \text{(resp.\ } A^{\infty,\infty}_{\vect s-(\vect b+\vect d)/p}(D)\cap S^{-1}(\ell^{p,q}(J,K))\subseteq A^{p,q}_{\vect s}(D)\text{).}
\]\CB
\end{teo}

Theorem~\ref{teo:3} is an immediate consequence of the following more general result.

\begin{teo}\label{teo:8}
\CR Take $\delta_+>0$, $R_0>1$,  $p,q\in ]0,\infty]$ and $\vect s\in \frac{1}{2 q}\vect m+(\R_+^*)^r$ (resp.\ $\vect s\in \R_+^r$ if $q=\infty$). 
Then, there are two constants $\delta_->0$ and $C>0$ such that, for every $\eps\in \Set{-,+}$, and for every $(\delta,R)$-lattice $(\zeta_{j,k},z_{j,k})_{j\in J,k\in K}$ on $D$ with $\delta\in ]0,\delta_\eps]$ and  $R\in ]1,R_0]$, defining
\[
S_+\colon\Hol(D)\ni f \mapsto \left(\Delta_\Omega^{\vect s-(\vect b+\vect d)/p}(h_k) \max_{\overline B((\zeta_{j,k},z_{j,k}),R\delta)} \abs{f}\right)_{j,k}\in \C^{J\times K}
\]
and
\[
S_-\colon\Hol(D)\ni f \mapsto \left(\Delta_\Omega^{\vect s-(\vect b+\vect d)/p}(h_k) \min_{\overline B((\zeta_{j,k},z_{j,k}),R\delta)} \abs{f}\right)_{j,k}\in \C^{J\times K},
\]
one has 
\[
\frac 1 C\norm{f}_{A^{p,q}_{\vect s}(D)}\meg \delta^{(2n+m)/p+m/q}\norm{S_\eps f}_{\ell^{p,q}(J,K)}\meg C \norm{f}_{A^{p,q}_{\vect s}(D)}
\]
for every $f\in A^{p,q}_{\vect s,0}(D)$ (resp.\ for every $f\in A^{p,q}_{\vect s}(D)$).
In addition,  
\[
 S^{-1}_+(\ell^{p,q}_0(J,K))\subseteq A^{p,q}_{\vect s,0}(D) \qquad\text{(resp.\ }  S^{-1}_+(\ell^{p,q}(J,K))\subseteq A^{p,q}_{\vect s}(D)\text{)}
 \]
 when $\eps=+$, while
 \[
 A^{\infty,\infty}_{\vect s-(\vect b+\vect d)/p}(D)\cap S^{-1}_-(\ell^{p,q}_0(J,K))\subseteq A^{p,q}_{\vect s,0}(D) \qquad \text{(resp.\ } A^{\infty,\infty}_{\vect s-(\vect b+\vect d)/p}(D)\cap S^{-1}_-(\ell^{p,q}(J,K))\subseteq A^{p,q}_{\vect s}(D)\text{)}
 \]
 when $\eps=-$.\CB
\end{teo}

The proof is based on~\cite[Lemma 6.3]{RicciTaibleson}, which deals with $S_+$ in the case of $\C_+$, and ~\cite[Theorem 5.6]{Bekolleetal}, which is Theorem~\ref{teo:3} in the case in which $p=q\in [1,\infty]$ and $D$ is an irreducible symmetric tube domain.

\CR With minor modifications, one may further prove that, for every $\vect{s'}\in \R^r$, there is $\delta_->0$ such that, if $f\in \Hol(D)$, $S_- f\in \ell^{p,q}_0(J,K)$ (resp.\ $S_- f\in \ell^{p,q}(J,K)$), and the mapping $(\zeta,z)\mapsto \Delta_\Omega^{\vect {s'}}(\Im z-\Phi(\zeta)) f(\zeta,z) \ee^{- \abs{\zeta}^{2\alpha}-\abs{\Re z}^\alpha-\abs{\Im z-\Phi(\zeta)}^\alpha}$ is bounded on $D$ for some $\alpha\in ]0,1/2[$, then $f\in A^{p,q}_{\vect s,0}(D)$ (resp.\ $f \in A^{p,q}_{\vect s}(D)$).\CB

We first need some simple lemmas.

\begin{lem}\label{lem:30}
There are $\rho_0>0$ and $C>0$ such that for every $p\in ]0,\infty[$, for every $\rho\in ]0,\rho_0]$, for every $f\in \Hol(D)$ and for every $(\zeta,z)\in D$,
\[
\abs{f(\zeta,z)}^p\meg C\dashint_{B((\zeta,z),\rho)} \abs{f}^p\,\dd \nu_D.
\]
\end{lem}

\begin{proof}
Observe that there are $\rho_0>0$ and $\gamma>0$ such that, for every $\rho\in ]0,\rho_0]$,
\[
B_{E\times F_\C}((0, i e_\Omega),\rho)\subseteq B((0, i e_\Omega), \rho\gamma).
\]
Now, Lemma~\ref{lem:76} shows that 
\[
\abs{f(0,i e_\Omega)}^{p}\meg  \dashint_{B_{E\times F_\C}((0,i e_\Omega),\rho/\gamma)} \abs{f(\zeta',z')}^{p}\,\dd (\zeta',z') 
\]
for every $\rho\in ]0, \rho_0]$ and for every $f\in \Hol(D)$, so that our assertion follows for $(\zeta,z)=(0,i e_\Omega)$, with 
\[
C\coloneqq \sup\limits_{0<\rho\meg \rho_0}\frac{\nu_D(B((0, i e_\Omega),\rho))}{\Hc^{2 n+2 m}(B_{E\times F_\C}((0, i e_\Omega),\rho/\gamma))} \sup\limits_{B_\Omega(e_\Omega,\rho_0)}\Delta_\Omega^{-\vect b-2\vect d}.
\]
The assertion for general $(\zeta,z)$ follows by homogeneity. 
\end{proof}

\begin{lem}\label{lem:28}
There are $\rho_0>0$ and $C>0$  such that, for every $\rho,\rho'\in ]0,\rho_0]$, for every $f\in \Hol(D)$, for every $(\zeta,z),(\zeta',z')\in D$ such that $\rho\coloneqq d((\zeta,z),(\zeta',z'))\meg \rho_0$,
\[
\abs{f(\zeta',z')-f(\zeta,z)}\meg \frac{ C \rho}{\rho' } \sup\limits_{B((\zeta,z),\rho+\rho')} \abs{f}.  
\]
\end{lem}

\begin{proof}
Observe first that there is $\rho_1>0$ such that $\exp_{(0, i e_\Omega)}$ induces diffeomorphisms
\[
B_{T_{(0, i e_\Omega)}(D)}(0,\rho)\to B((0, i e_\Omega),\rho) \quad \text{and} \quad \partial B_{T_{(0, i e_\Omega)}(D)}(0,\rho)\to\partial B((0, i e_\Omega),\rho)
\]
for every $\rho\in ]0, \rho_1]$ (cf.~\cite[Theorem 6.4]{Lang}). 
In addition, there are $\rho_2\in ]0, \rho_1]$ and $\gamma\Meg 1$ such that
\[
\frac{1}{\gamma}d((\zeta',z'),(\zeta'',z''))\meg  \abs{(\zeta',z')-(\zeta'',z'')}\meg \gamma d((\zeta',z'),(\zeta'',z''))
\]
for every $(\zeta',z'),(\zeta'',z'')\in B((0, i e_\Omega), \rho_2)\cup B_{E\times F_\C}((0, i e_\Omega),  \rho_2)$.
Finally, observe that there is $\gamma'$ such that 
\[
\abs{\partial_v \exp_{(0, i e_\Omega)}(t v) }\meg \gamma'
\]
for every $v \in \overline B_{T_{(0, i e_\Omega)}(D)}(0, 1)$ and for every $t\in [0, \rho_2]$.
Then, choose $\rho_0\coloneqq \rho_2/(2 \gamma)$.
In addition, take $\rho\meg \rho_0$, $f\in \Hol(D)$, $(\zeta',z')\in \partial B((0,i e_\Omega),\rho)$, and let $v\in T_{(0, i e_\Omega)}(D)$ be the unit vector such that 
\[
\exp_{(0, i e_\Omega)}(\rho  v)=(\zeta',z').
\]
Then,
\[
\begin{split}
\abs{f(\zeta',z')-f(0,i e_\Omega)}&\meg \int_0^{\rho} \abs*{\frac{\dd}{\dd t} f(\exp_{(0, i e_\Omega)}(t v))   }\,\dd t\\
&\meg \rho\sup\limits_{t\in ]0, \rho[}\abs*{f'\left(\exp_{(0, i e_\Omega)}(t v)\right) \cdot \left( \partial_v \exp_{(0, i e_\Omega)}\right) (t v)}\\
&\meg \gamma' \rho \sup\limits_{(\zeta'',z'')\in B((0, i e_\Omega),\rho) }\abs{f'(\zeta'',z'')}.
\end{split}
\]
In addition, Cauchy's integral formula implies that
\[
\abs{f'(\zeta'',z'')}\meg \frac{\gamma}{\rho'}\sup\limits_{\abs{v'}\meg 1}\dashint_{\partial B_\C(0,\rho'/\gamma)} \abs{f((\zeta'',z'')+i w v')}\,\dd w 
\]
for every $(\zeta'',z'')\in B((0, i e_\Omega),\rho)$ and for every $\rho'\in ]0,\rho_0]$.  Thus,
\[
\abs{f(\zeta',z')-f(0,i e_\Omega)}\meg \frac{\gamma\gamma'}{\rho'}\rho \sup\limits_{B((0,i e_\Omega), \rho+\rho')} \abs{f}.
\]
The general assertion follows by homogeneity.
\end{proof}

\begin{lem}\label{lem:36}
There are $\rho_0>0$ and a constant $C>0$ such that, for every $p,q\in ]0,\infty]$, for every $\rho\in ]0,\rho_0]$, for every $f\in \Hol(D)$ and for every $h\in \Omega$,
\[
\norm{f_h}_{L^p(\Nc)}\meg C^{1/\min(1,p,q)} \left( \dashint_{B_\Omega(h, \rho)}  \norm{f_{h'}}_{L^p(\Nc)}^q\,\dd \nu_\Omega(h)\right) ^{1/q}.
\]
(modification if $q=\infty$).
\end{lem}

\begin{proof}
Set $\ell\coloneqq \min(1,p,q)$ to simplify the notation.
By Lemma~\ref{lem:30}, there are $\rho_0>0$ and $C'>0$ such that
\[
\abs{f(\zeta,z)}^{\ell}\meg C' \dashint_{B((\zeta,z),\rho)} \abs{f}^{\ell}\,\dd \nu_D
\]
for every $f\in \Hol(D)$, for every $(\zeta,z)\in D$, and for every $\rho\in ]0,\rho_0]$.  Then, applying Minkowski's integral inequality (with exponent $\frac{p}{\ell}$) and Young's inequality,
\[
\begin{split}
\norm{f_h}_{L^p(\Nc)}^\ell&\meg C' C'_\rho  \dashint_{B_\Omega(h,\rho)} \norm*{\abs{f_{h'}}^{\ell}*[(\chi_{B((0, i h),\rho)})_{h'}]\check{\;}  }_{L^{p/\ell}(\Nc)}  \Delta_\Omega^{\vect b+\vect d}(h') \,\dd \nu_\Omega(h')\\
&\meg C''  \dashint_{B_\Omega(h,\rho)} \norm*{f}_{L^{p}(\Nc)}^\ell \frac{\Delta_\Omega^{\vect b+\vect d}(h') }{\Delta_\Omega^{\vect b+\vect d}(h) } \,\dd \nu_\Omega(h')
\end{split}
\]
for every $f\in \Hol(D)$ and for every $h\in \Omega$,
where 
\[
C'_\rho\coloneqq \frac{\nu_\Omega(B(e_\Omega,\rho))}{\nu_D(B((0, i e_\Omega),\rho))}
\]
and 
\[
C''\coloneqq C' \sup\limits_{0<\rho\meg \rho_0} \sup\limits_{h'\in \Omega} C'_\rho \norm*{\left(\chi_{B((0, i e_\Omega),\rho)}\right)_{h'}}_{L^1(\Nc)}.
\]
By Corollary~\ref{cor:34}, there is a constant $C>0$ such that, for every $f\in \Hol(D)$ and for every $h\in \Omega$,
\[
\begin{split}
\norm{f_h}_{L^p(\Nc)}^\ell&\meg C\dashint_{B_\Omega(h,\rho)} \norm*{f}_{L^{p}(\Nc)}^\ell \,\dd \nu_\Omega(h'),
\end{split}
\]
provided that $\rho_0$ is sufficiently small. Then, Jensen's inequality (with exponent $ \frac{q}{\ell}$) leads to the conclusion
\end{proof}

\begin{proof}[Proof of Theorem~\ref{teo:8}.]\CR
We leave to the reader the (purely formal) modifications needed to deal with the case $\max(p,q)=\infty$.

\textsc{Step I.} Define, for every $\rho>0$ and for every $h\in \Omega$, 
\[
M_\rho(h)\coloneqq \norm*{\left(\chi_{B((0,i e_\Omega),  \rho)}\right)_h}_{L^1(\Nc)}.
\]
Let us first show that $M_\rho$ is bounded. Observe that, for every $(\zeta,x)\in \Nc$ and for every $h\in \Omega$, 
\[
d((\zeta,x+ih +i\Phi(\zeta)),(0, i e_\Omega))=d((g^{-1}\zeta,t^{-1}\cdot x+i e_\Omega+i \Phi(g^{-1}\zeta)),(0, i t^{-1}\cdot e_\Omega)),
\]
where $t\in T_+$, $g\in GL(E)$, $h=t\cdot e_\Omega$, and $t\cdot \Phi=\Phi\circ (g\times g)$. 
Since 
\[
d((0, i t^{-1}\cdot e_\Omega),(0, i e_\Omega))=d_\Omega(t^{-1}\cdot e_\Omega, e_\Omega)=d_\Omega(e_\Omega,h)
\]
and since $(\chi_{B((0,i e_\Omega),R\delta)})_h=0$ if $h\not \in B_\Omega(e_\Omega, R\delta)$, thanks to Lemma~\ref{lem:45}, we then conclude that 
\[
\begin{split}
&\bigcup_{h\in B_\Omega(e_\Omega,R\delta)}\Set{(\zeta,x)\in \Nc\colon (\zeta,x+i\Phi(\zeta)+i h)\in B((0,i e_\Omega),\rho)}\\
	&\qquad\subseteq \bigcup_{(g,t)\in K'} (g\times t)\left( \Set{(\zeta,x)\in \Nc\colon (\zeta,x+i\Phi(\zeta)+i e_\Omega)\in B((0,i e_\Omega),2 \rho)}\right) ,
\end{split}
\]
where $K'$ is a suitable \emph{compact} subset of $GL(E)\times T_+$. Therefore, Proposition~\ref{prop:26} implies that there is a compact subset $K''$ of $\Nc$ such that $(\chi_{B((0,i e_\Omega),R\delta)})_h\meg \chi_{K''}$ for every $h\in \Omega$, so that $\norm{M_\rho}_{L^\infty(\Omega)}\meg \Hc^{2 n+m}(K'')$.

Then, for every $\ell\in ]0,\infty]$, for every $t'\in T_+ $, and for every $(\zeta',x')\in \Nc$,
\[
\norm*{\left(\chi_{B((\zeta',x'+i\Phi(\zeta)+i t'\cdot e_\Omega), \rho)}\right)_h}_{L^\ell(\Nc)}=\Delta_\Omega^{-(\vect b+ \vect d)/\ell}(h') M_\rho(t'^{-1}\cdot h)^{1/\ell},
\]
with the convention $0^0=0$. 
In particular, (cf.~Lemma~\ref{lem:45}),
\[
\norm*{\left(\chi_{B((\zeta_{j,k},z_{j,k}), \rho)}\right)_h}_{L^\ell(\Nc)}= \chi_{B_\Omega(h_k,\rho)}(h) 
\Delta_\Omega^{-(\vect b+ \vect d)/\ell}(h_k) M_\rho(t_k^{-1}\cdot h)^{1/\ell}
\]
for every $h\in \Omega$ and for every $k\in K$, where $t_k\in T_+$ and $h_k=t_k\cdot e_\Omega$. 
In addition,  
\[
\norm{M_\rho}_{L^\infty(\nu_\Omega)}\asymp \rho^{2 n+m} \text{ for $\rho\to 0^+$}.
\]
For every $h\in \Omega$, define 
\[
K_h\coloneqq \Set{k\in K\colon h\in B_\Omega(h_k,R\delta)},
\]
and observe that there is $N\in\N$ such that $ \card(K_h)\meg N$ for every $h\in\Omega$, provided that $R\meg R_0$ and $\delta\meg \delta_+$, thanks to Proposition~\ref{prop:56}. We may also assume that every $h\in \Omega$ is contained in at most $N$ balls $B_\Omega(h_k,2 R\delta)  $, $k\in K$, and that every $(\zeta,z)\in D$ is contained in at most $N$ balls $B((\zeta_{j,k},z_{j,k}),2 R\delta)$, $(j,k)\in J\times K$, provided that $R\meg R_0$ and $\delta\meg \delta_+$. 
Finally, set $\ell\coloneqq \min (1,p,q)$.

\textsc{Step II.} Let us prove that $S_+$ maps $ A^{p,q}_{\vect s}(D)$  into $\ell^{p,q}(J,K)$. Take $f\in A^{p,q}_{\vect s}(D)$ and define 
\[
C_{D,\rho}\coloneqq \nu_D(B((0,i e_\Omega),\rho)) \qquad \text{and} \qquad C_{\Omega,\rho}\coloneqq \nu_\Omega(B_\Omega(e_\Omega,\rho))
\] 
for every $\rho>0$ to simplify the notation.
Then, Lemma~\ref{lem:30} implies that there are $\rho_0\in ]0,1/2]$ and $C_1>0$ such that
\[
\max_{\overline B((\zeta_{j,k},z_{j,k}),R\delta)} \abs{f}^p\meg  \frac{C_1}{ C_{D,\rho_0 \delta}} \int_{B((\zeta_{j,k},z_{j,k}), (R+\rho_0)\delta)} \abs{f}^p\,\dd \nu_D
\]
for every $(j,k)\in J\times K$.
Therefore, Corollary~\ref{cor:34} implies that there is a constant $C_{2}>0$ such that
\[
(S_+ f)_{j,k}^p \meg \frac{C_{2}}{C_{D,\rho_0\delta}}\Delta_\Omega^{p\vect s}(h_k) \int_\Omega \int_\Nc \abs{(\chi_{B((\zeta_{j,k},z_{j,k}),(R+\rho_0)\delta)}f)_h(\zeta,x)}^p\,\dd (\zeta,x)\,\dd \nu_\Omega(h)
\]
for every $(j,k)\in J\times K$.
Hence, 
\[
\sum_{j\in J}(S_+ f)_{j,k}^p \meg \frac{C_{2}}{C_{D,\rho_0\delta}} N\Delta_\Omega^{p\vect s}(h_k) \int_{B_\Omega(h_k, (R+\rho_0)\delta)} \norm{f_h}_{L^p(\Nc)}^p\,\dd \nu_\Omega(h)
\]
for every $k\in K$.
Now, Lemma~\ref{lem:36} shows that there is a constant $C_{3}>0$ such that 
\[
\begin{split}
&\int_{B_\Omega(h_k, (R+\rho_0)\delta)} \norm{f_h}_{L^p(\Nc)}^p\,\dd \nu_\Omega(h)\\
&\qquad\meg C_3 \int_{B_\Omega(h_k, (R+\rho_0)\delta)}\left(  \dashint_{B_\Omega(h', \rho_0\delta)}  \norm{f_{h}}_{L^p(\Nc)}^q\,\dd \nu_\Omega(h)\right) ^{p/q}\,\dd \nu_\Omega(h')\\
&\qquad \meg C_3 \frac{C_{\Omega, (R+\rho_0)\delta}}{C_{\Omega,\rho_0\delta}^{p/q}  } \left(  \int_{B_\Omega(h_k, (R+2\rho_0)\delta)}  \norm{f_{h}}_{L^p(\Nc)}^q\,\dd \nu_\Omega(h)\right) ^{p/q}
\end{split}
\]
for every $k\in K$,  provided that $\rho_0\delta_+$ is sufficiently small.
Therefore, another application of Corollary~\ref{cor:34} shows that there is a constant $C_2'>0$ such that
\[
\norm{S_+ f}_{\ell^{p,q}(J,K)}\meg \frac{C_2' C_{\Omega,(R+\rho_0)\delta}^{1/p} }{C_{D,\rho_0\delta}^{1/p}C_{\Omega, \rho_0\delta}^{1/q} } N^{1/p+1/q} \norm{f}_{A^{p,q}_{\vect s}(D)}.
\]

Next, assume that $\vect{s}\in \frac{1}{2 q}\vect{m}+(\R_+^*)^r$, and let us prove that $S_+(A^{p,q}_{\vect s,0}(D))\subseteq \ell^{p,q}_0(J,K)$. Indeed, take $\tilde q\in ]0,q[$ so that $ \vect{s}\in \frac{1}{2 \tilde q}\vect{m}+(\R_+^*)^r$, and observe that the preceding computations show that $S_+(A^{p,q}_{\vect s,0}(D)\cap A^{p,\tilde q}_{\vect s}(D))\subseteq \ell^{p,q}(J,K)\cap \ell^{p,\tilde q}(J,K)$, so that the assertion follows by means of  Proposition~\ref{prop:8} if $p<\infty$, since in that case $\ell^{p,q}(J,K)\cap \ell^{p,\tilde q}(J,K)\subseteq \ell^{p,q}_0(J,K)$. 
If, otherwise, $p=\infty$, then it is clear that $\lim\limits_{j\to \infty}(S_+ f)_{j,k}=0$ for every $f\in A^{p,q}_{\vect s,0}(D)\cap A^{p,\tilde q}_{\vect s}(D)$ and for every $k\in K$, so that the assertion follows as before.

\textsc{Step III.} Conversely, assume that $S_+ f\in \ell^{p,q}_0(J,K)$ (resp.\ $S_+ f\in \ell^{p,q}(J,K)$), and let us prove that $f\in A^{p,q}_{\vect s}(D)$ (resp.\ $f\in A^{p,q}_{\vect s,0}(D)$).
Indeed, it is readily verified that
\[
\abs{f_h}\meg \sum_{(j,k)\in J\times K_h} \Delta_\Omega^{(\vect b+\vect d)/p-\vect s}(h_k) \left( \chi_{B((\zeta_{j,k},z_{j,k}),R\delta)} \right)_h (S_+ f)_{j,k}
\]
on $D$, for every $h\in \Omega$, so that $f_h\in L^p_0(\Nc)$ (resp.\ $f_h\in L^p(\Nc)$) for every $h\in \Omega$. In addition,
\[
\begin{split}
	\norm{f_h}_{L^p(\Nc)}&\meg N^{1/p'} \norm*{(\Delta_\Omega^{-\vect s}(h_k) M_{R\delta}(t_k^{-1}\cdot h )^{1/p}(S_+ f)_{j,k})_{j,k}}_{\ell^p(J\times K_h)}    \\
	&\meg   N^{1/p'+(1/p-1/q)_+}\norm{M_{R\delta}}_\infty^{1/p} \norm*{ \left( \Delta_\Omega^{-\vect s}(h_k)\norm*{  ((S_+ f)_{j,k})_j}_{\ell^p(J)}\right) _k}_{\ell^q(K_h)}  
\end{split}
\]
so that the assertion follows.

\textsc{Step IV.} Finally, take $f\in A^{\infty,\infty}_{\vect s-(\vect b+\vect d)/p}(D)$ such that $S_- f\in \ell^{p,q}_0(J,K)$ (resp.\ $S_- f\in \ell^{p,q}(J,K)$), and let us prove that $f\in A^{p,q}_{\vect s,0}(D)$ (resp.\ $f\in A^{p,q}_{\vect s}(D)$) and that
\[
\norm{f}_{A^{p,q}_{\vect s}(D)}\meg C \delta^{(2 n+m)/p+m/q} \norm{S_- f}_{\ell^{p,q}(J,K)}
\] 
for a suitable constant $C>0$, depending only on $\delta_-$ and $R_0$, provided that $\delta_-$ is sufficiently small. 
Observe that, for every $(j,k)\in J\times K$, we may find $(\zeta'_{j,k},z'_{j,k})\in \overline B((\zeta_{j,k},z_{j,k}),R\delta)$ such that
\[
\abs{f(\zeta'_{j,k},z'_{j,k})}=\min_{\overline B((\zeta_{j,k},z_{j,k}),R\delta)}\abs{f}.
\]
Now, Lemmas~\ref{lem:30} and~\ref{lem:28} imply that there are $\rho_1>0$ and  $C_3>0$ such that, for every $j\in J$, for every $k\in K_h$, and for every $(\zeta,x)\in \Nc$ such that $d((\zeta_{j,k},z_{j,k}),(\zeta,x+i\Phi(\zeta)+i h) )<R\delta$,
\[
\abs{f_h(\zeta,x)}\meg\abs{f(\zeta'_{j,k},z'_{j,k})}+C_{3} R \delta \norm{\chi_{B((\zeta,x+i\Phi(\zeta)+i h), 2 R \delta+\rho_1)} f }_{L^p(\nu_D)},
\]
provided that $R\delta \meg \rho_1$.
Then,
\[
\begin{split}
\norm{f_h}_{L^p(\Nc)}^p&\meg 2^{(p-1)_+}\norm{M_{R\delta}}_{L^\infty(\Omega)} \sum_{(j,k)\in J\times K_h}\Delta_\Omega^{-(\vect b+\vect d)}(h_k) \abs{f(\zeta'_{j,k},z'_{j,k}) }^p+2^{(p-1)_+} (C_3 R \delta)^p \Theta_1(h)  ,
\end{split}
\]
where
\[
\begin{split}
\Theta_1(h)&\coloneqq \sum_{(j,k)\in J\times K_h} \int_\Nc (\chi_{B((\zeta_{j,k},z_{j,k}),R\delta)})_h(\zeta,x)\int_D \chi_{B((\zeta,x+i\Phi(\zeta)+i h),2 R\delta+\rho_1)} \abs{f }^p\,\dd \nu_D \,\dd (\zeta,x) 
\end{split}
\]
Now, set $K^{(k')}\coloneqq \bigcup_{h\in  B_\Omega(h_{k'},R\delta)} K_h$ for every $k'\in K$, so that $\card\left(\Set{k'\in K\colon k\in K^{(k')}}\right)\meg N$ by~{step I}, provided that $\delta\meg \delta_+$. 
Then, Corollary~\ref{cor:34} implies that there are constants $C_4,C_4'>0$ such that, if $R\delta \meg \rho_1$,
\[
\begin{split}
&\int_\Omega \Delta_\Omega^{q\vect s}(h) \left(\sum_{(j,k)\in J\times K_h} \Delta_\Omega^{-(\vect b+\vect d)}(h_k) \abs{f(\zeta'_{j,k},z'_{j,k}) }^p\right)^{q/p}\,\dd \nu_\Omega(h) \\
&\qquad\meg C_4 C_{\Omega,R\delta} \sum_{k'\in K} \Delta_\Omega^{q \vect s}(h_{k'}) \left( \sum_{k\in K^{(k')}} \sum_{j\in J}   \Delta_\Omega^{-(\vect b+\vect d)}(h_k) \abs{f(\zeta'_{j,k},z'_{j,k}) }^p \right)^{q/p}\\
&\qquad \meg C_4' C_{\Omega, R\delta} N^{\max(1,q/p)} \norm{S_- f}_{\ell^{p,q}(J,K)}^q
\end{split}
\]
where the last inequality follows from the convexity or subadditivity of the mapping $x\mapsto x^{q/p}$ on $\R_+$.

Now, observe that
\[
\begin{split}
&\Theta_1(h)=\int_D  \abs{f(\zeta',z')}^p \Theta_2(\zeta',z',h)  \,\dd \nu_D(\zeta',z'),
\end{split}
\]
where
\[
\Theta_2(\zeta',z',h)\coloneqq\int_\Nc \sum_{(j,k)\in J\times K_h} (\chi_{B((\zeta_{j,k},z_{j,k}),R\delta)\cap B((\zeta',z'),2 R\delta+\rho_1)})_h(\zeta,x)   \,\dd (\zeta,x) .
\]
In addition, for every $(\zeta',z')\in D$ and for every $h\in \Omega$, setting $h'\coloneqq \Im z'-\Phi(\zeta')$, one has
\[
\begin{split}
\Theta_2(\zeta',z',h)&\meg N \norm{(\chi_{B((\zeta',z'),2 R\delta+\rho_1)})_h  }_{L^1(\Nc)}=N M_{2 R\delta+\rho_1}(h'^{-1}\cdot h) \Delta_\Omega^{-(\vect b+\vect d)}(h'),
\end{split}
\]
provided that $R\meg R_0$ and $\delta\meg \min(\rho_0,\delta_+)$.
Therefore, by~{step I} we see that
\[
\Theta_1(h) \meg N \norm{M_{2R\delta+\rho_1}}_\infty\int_{B_\Omega(h,2 R\delta+\rho_1)} \norm{f_{h'}}_{L^p(\Nc)}^p\,\dd \nu_\Omega(h'),
\]
provided that $R\meg R_0$ and $\delta\meg \min(\rho_0,\delta_+)$.
%Now, arguing as in~{step II} we see that Lemma~\ref{lem:36} implies that there are  $C_{5}, \rho_2>0$ such that
%\[
%\begin{split}
%&\int_{B_\Omega(h,2 R\delta+\rho_1)} \norm{f_{h'}}_{L^p(\Nc)}^p\,\dd \nu_\Omega(h')\meg C_5  \left(  \int_{B_\Omega(h, 2 R\delta+ \rho_1+\rho_2)}\norm{f_{h'}}_{L^p(\Nc)}^q\,\dd \nu_\Omega(h')   \right)^{p/q} .
%\end{split}
%\]
%Thus, there is a constant $C_6>0$ such that, for every $h\in \Omega$,
%\[
%\begin{split}
%\Theta_1(h) \meg C_6 \left(  \int_{B_\Omega(h, 2R\delta+ \rho_1+\rho_2)}\norm{f_{h'}}_{L^p(\Nc)}^q\,\dd \nu_\Omega(h')   \right)^{p/q}
%\end{split}
%\]
%provided that $R\meg R_0$ and $\delta\meg \rho_0$.
%Now, Corollary~\ref{cor:34} shows that there is a constant $C_7>0$ such that
%\[
%\begin{split}
%&\int_\Omega \Delta_\Omega^{q \vect s}(h)  \int_{B_\Omega(h, 2R\delta+ \rho_1+\rho_2)}\norm{f_{h'}}_{L^p(\Nc)}^q\,\dd \nu_\Omega(h') \,\dd \nu_\Omega(h) \\
%&\qquad \qquad\qquad \qquad\qquad \qquad\qquad \qquad\meg C_7 C_{\Omega,2 R\delta+\rho_1+\rho_2} \norm{f}_{A^{p,q}_{\vect s}(D)}^q.
%\end{split}
%\]

Now, take $(g^{(\eps)})_{\eps>0}$ as in Lemma~\ref{lem:34} for some $\alpha\in ]0,1/2[$. Then, it is readily seen that $f^{(\eps)}\coloneqq f g^{(\eps)}$ belongs to $A^{p,\infty}_{\vect s-(\vect b+\vect d)/p}(D)$, and that clearly $S_- f_\eps\meg S_- f$ for every $\eps>0$. In particular, the mapping $h\mapsto \norm{f^{(\eps)}_h}_{L^p(\Nc)}$ is (finite and) decreasing on $\Omega$, thanks to Corollary~\ref{cor:26}, for every $\eps>0$.
In addition, observe that we may take $\delta_1\in ]0,\min(\rho_0,\delta_+)]$ and $\rho_1$ sufficiently small such that $B_\Omega(e_\Omega, 2 R_0\delta_1+\rho_1)\subseteq e_\Omega/2+\Omega$, so that, by homogeneity,
\[
B_\Omega(h, 2 R_0\delta_1+\rho_1)\subseteq h/2+\Omega
\]
for every $h\in \Omega$. Then, the preceding estimates (applied to $f^{(\eps)}$) show that there is a constant $C_5>0$ such that
\[
\begin{split}
\norm*{f_h^{(\eps)}}_{L^p(\Nc)}&\meg C_5 \Delta_\Omega^{-\vect s}(h)\delta^{(2n+m)/p} \norm{S_- f}_{\ell^{p}(J\times K_h)}+\delta C_5 \norm*{h'\mapsto \chi_{B_\Omega(h,2 R\delta+\rho_1)}(h')\norm*{f^{(\eps)}_{h'}}_{L^p(\Nc)}  }_{L^p(\nu_\Omega)}\\
	&\meg C_5   \Delta_\Omega^{-\vect s}(h)\delta^{(2n+m)/p}\norm{S_- f}_{\ell^{p}(J\times K_h)}+\delta C_5  C_{\Omega,2 R_0\delta_1+\rho_1}^{1/p}	\norm*{f_{h/2}^{(\eps)}}_{L^p(\Nc)},
\end{split}
\]
for every $\eps>0$, provided that $\delta\meg \delta_1$ and $R\meg R_0$,
since the mapping $h\mapsto \norm*{f^{(\eps)}_h}_{L^p(\Nc)}$ is decreasing.
If we define $\chi_\ell\colon D\ni (\zeta,z)\mapsto \chi_{e_\Omega/2^\ell+\Omega }(\Im z-\Phi(\zeta))\in \R_+$ for every $\ell\in \N$, then there is a constant $C_6>0$ such that
\[
\norm*{\chi_\ell f^{(\eps)}}_{L^{p,q}_{\vect s}(D)}\meg C_6 \delta^{(2n+m)/p+m/q}\norm{S_- f}_{\ell^{p,q}(J,K)}+\delta C_6 \norm*{\chi_{\ell+1} f^{(\eps)}}_{L^{p,q}_{\vect s}(D)}
\]
for every $\eps>0$ and for every $\ell\in \N$, provided that $\delta\meg \delta_1$ and $R\meg R_0$. 
Now, observe that
\[
\begin{split}
	\norm*{f^{(\eps)} \chi_{\ell+1}}_{L^{p,q}_{\vect s}(D)}&\meg \norm{f}_{A^{\infty,\infty}_{\vect s-(\vect b+\vect d)/p}(D)} \norm*{g^{(\eps)}\chi_{\ell+1}}_{L^{p,q}_{(\vect b+\vect d)/p}(D)}\\
	&\meg C_{7,\eps} \norm{f}_{A^{\infty,\infty}_{\vect s-(\vect b+\vect d)/p}(D)} \norm*{\chi_{e_\Omega/2^{\ell+1}+\Omega} \Delta_\Omega^{(\vect b+\vect d)/p}\ee^{- C_{7,\eps}\abs{\,\cdot\,}^\alpha}}_{L^q(\nu_\Omega)}\\
	& \meg C_{8,\eps} 2^{((n+m)/p+m/q)(\ell+1)}\norm{f}_{A^{\infty,\infty}_{\vect s-(\vect b+\vect d)/p}(D)}
\end{split}
\]
for suitable constants $C_{7,\eps},C_{8,\eps}>0$. Then, choose $N'> (n+m)/p+m/q$ and choose $\delta_-\in ]0,\delta_1]$ so that $ C_6\delta_- \meg 2^{-N'}$, and observe that the preceding computations show that, if $\delta\in ]0,\delta_-]$,
\[
\begin{split}
\norm*{f^{(\eps)}\chi_{\ell}}_{L^{p,q}_{\vect s}(D)}&=\sum_{\ell\in\N} 2^{-\ell N'}\left(\norm*{f^{(\eps)} \chi_{\ell}}_{L^{p,q}_{\vect s}(D)}-\frac{1}{2^{N'}} \norm*{f^{(\eps)} \chi_{\ell+1}}_{L^{p,q}_{\vect s}(D)}\right)\\
	&\meg \frac{C_6}{1-2^{-N'}}\delta^{(2n +m)/p+m/q}\norm{S_- f}_{\ell^{p,q}(J,K)}
\end{split}
\]
for every $\eps>0$ and for every $\ell\in \N$. Passing to the limit for $\ell\to \infty$ we then infer that $f^{(\eps)}\in A^{p,q}_{\vect s}(D)$ and that
\[
\norm*{f^{(\eps)}}_{A^{p,q}_{\vect s}(D)}\meg \frac{C_6}{1-2^{-N'}}\delta^{(2n +m)/p+m/q}\norm{S_- f}_{\ell^{p,q}(J,K)}
\]
for every $\eps>0$.
Then, passing to the limit for $\eps\to 0^+$, we infer that $f\in A^{p,q}_{\vect s}(D)$ and that
\[
\norm{f}_{A^{p,q}_{\vect s}(D)}\meg \frac{C_6}{1-2^{-N'}}\delta^{(2n +m)/p+m/q}\norm{S_- f}_{\ell^{p,q}(J,K)}.
\]
It only remains to prove that, if $\vect s\in \frac{1}{2 q}\vect m+(\R_+^*)^r$ and $f\in A^{\infty,\infty}_{\vect s-(\vect b+\vect d)/p}(D)\cap S_-^{-1}(\ell^{p,q}_0(J,K))$, then $f\in A^{p,q}_{\vect s,0}(D)$.

Observe first that the preceding computations show that
\[
\norm{\chi_\ell (f-f^{(\eps)})}_{L^{p,q}_{\vect s}(D)}\meg  \frac{C_6}{1-2^{-N'}} \delta^{(2 n+m)/p+m/q} \norm{S_-(f-f^{(\eps)})}_{\ell^{p,q}(J,K)}.
\]
Since $S_-(f-f^{(\eps)})\meg S_- f \widetilde S_+(1-g^{(\eps)})$, where $\big[\widetilde S_+(1-g^{(\eps)})\big]_{j,k}\coloneqq\max_{B((\zeta_{j,k},z_{j,k}),R\delta)} \abs*{1-g^{(\eps)}}$ for every $(j,k)\in J\times K$, and since $1-g^{(\eps)}\to 0$ locally uniformly, it is readily seen that $\chi_\ell (f-f^{(\eps)})\to 0$ in $L^{p,q}_{\vect s}(D)$ for $\eps\to 0^+$. In particular, $f_h\in L^p_0(\Nc)$ for every $h\in \Omega$, and $h\mapsto \chi_{e_\Omega/2^\ell+\Omega}(h)\Delta^{\vect s}_\Omega(h)\norm{f_h}_{L^p(\Nc)} $ belongs to $L^q_0(\nu_\Omega)$ for every $\ell\in \N$. To conclude, it will essentially suffice to show that, if $q=\infty$, then $\Delta^{\vect s}_\Omega(h)\norm{f_h}_{L^p(\Nc)}\to 0$ as $h$ approaches the boundary of $\Omega$.
Observe that, by the preceding computations, there is a constant $C_9>0$ such that
\[
\norm{f_h}_{L^p(\Nc)}\meg C_9 \Delta_\Omega^{-\vect s}(h)\norm{S_- f}_{\ell^p(J\times K_h)}+\delta C_9 \norm{f_{h/2}}_{L^p(\Nc)}
\]
for every $h\in \Omega$. Observe that 
\[
\Delta_\Omega^{\vect s}(h) \norm{f_{h/2}}_{L^p(\Nc)}= 2^{\vect s} \Delta_\Omega^{\vect s}(h/2) \norm{f_{h/2}}_{L^p(\Nc)}\meg 2^{\vect s} \norm{f}_{A^{p,\infty}_{\vect s}(D)}
\]
for every $h\in \Omega$. Therefore, assuming that $\delta_-$ is so small that $\delta_- C_9 2^{\vect s}\meg 1/2 $,
\[
\begin{split}
\Delta_\Omega^{\vect s}(h)\norm{f_h}_{L^p(\Nc)}&= \sum_{\ell\in \N} 2^{-\ell}\left(\Delta_\Omega^{\vect s}(h/2^{\ell}) \norm{f_{h/2^{\ell}}}_{L^p(\Nc)}-\frac{1}{2 \cdot 2^{\vect s}}\Delta_\Omega^{\vect s}(h/2^{\ell}) \norm{f_{h/2^{\ell+1}}}_{L^p(\Nc)}\right)\\
	&\meg C_9  \sum_{\ell\in \N} 2^{-\ell}\norm{S_- f}_{\ell^p(J\times K_{h/2^\ell})}
\end{split}
\]
for every $h\in \Omega$. Now, observe that $\eta\coloneqq \min_{\ell\in \N} d_\Omega(e_\Omega,e_\Omega/2^{\ell+1})>0$, and that $\eta=\min_{\ell\in \N} d_\Omega(h,h/2^{\ell+1})$ for every $h\in\Omega$, by homogeneity. Therefore, if $\delta_-$ is so small that $2 R_0\delta_-<\eta$, then the sets $K_{h/2^\ell}$, as $\ell$ runs through $\N$, are pairwise disjoint for every $h\in \Omega$. Hence,
\[
\Delta_\Omega^{\vect s}(h)\norm{f_h}_{L^p(\Nc)}\meg 2 C_9 N^{1/p} \norm{S_- f}_{\ell^{p,\infty}(J, K'_h)}
\]
for every $h\in\Omega$, where $K'_h\coloneqq \bigcup_{\ell\in \N} K_{h/2^\ell}$. Since $K'_h$ is contained in the complement of every \emph{fixed} finite subset of $K$ if $h\in\Omega\setminus (e_\Omega/2^\ell+\Omega)$ and $\ell$ is sufficiently large, this and the preceding arguments prove that $\Delta_\Omega^{\vect s}(h) \norm{f_h}_{L^p(\Nc)}\to 0$ as $h\to \infty$ in $\Omega$, provided that $\delta_-$ is sufficiently small (independently of $f$). The proof is complete.
\CB
\end{proof}

\begin{cor}\label{cor:40}
	Take $p,q\in ]0,\infty ]$, $\vect s\in \frac{1}{2 q}\vect m+(\R_+^*)^r$ (resp.\ $\vect s\in \R_+^r$ if $q=\infty$) and $\vect{s'}\in \N_{\Omega'}$. Then, convolution by $I_\Omega^{-\vect{s'}}$ induces a continuous linear mapping $A^{p,q}_{\vect s,0}(D)\to A^{p,q}_{\vect s+\vect{s'},0}(D)$ (resp.\ $A^{p,q}_{\vect s}(D)\to A^{p,q}_{\vect s+\vect{s'}}(D)$).
\end{cor}

\begin{proof}
	Using Cauchy's estimates, we see that there is a constant $C>0$  such that
	\[
	\max_{B((0,i e_\Omega),1)}\abs{f*I_\Omega^{-\vect{s'}}}\meg C\max_{B((0,i e_\Omega),2)}\abs{f}
	\]
	for every $f\in \Hol(D)$. Therefore, by homogeneity we see that
	\[
	\Delta_\Omega^{\vect{s'}}( \Im z-\Phi(\zeta))\max_{B((\zeta,z),1)}\abs{f*I_\Omega^{-\vect{s'}}}\meg C\max_{B((\zeta,z),2)}\abs{f}
	\]
	for every $(\zeta,z)\in D$ and for every $f\in \Hol(D)$. Therefore, the assertion follows from Lemma~\ref{lem:32} and Theorem~\ref{teo:8}.
\end{proof}

\section{Atomic Decomposition}\label{sec:6:4}

In this section we deal with atomic decomposition\index{Atomic decomposition} for the spaces $A^{p,q}_{\vect s}(D)$. As mentioned earlier, along with property $\atomics^{p,q}_{\vect s}$, which is the properly called atomic decomposition, we consider also a stronger property $\atomic^{p,q}_{\vect s,+}$, which is somewhat easier to deal with.
As we shall see in Section~\ref{sec:6:7}, these properties are closely related with analogous statements concerning the Bergman projectors.

\begin{deff}\label{38}
Take $p,q\in ]0,\infty]$ and $\vect s,\vect{s'}\in \R^r$ such that the following conditions are satisfied:
\begin{itemize}
\item $\vect s\in \frac{1}{2 q}\vect m+(\R_+^*)^r$ (resp.\ $\vect s\in \R_+^r$ if $q=\infty$);

\item $\vect{s'}\in \frac 1 p(\vect b+\vect d)-\frac{1}{2 p}\vect{m'}- (\R_+^*)^r$ (resp.\ $\vect{s'}\in  -\R_+^r$ if $p=\infty$);

\item $\vect s+\vect{s'}\in \frac 1 p(\vect b+\vect d)-\frac{1}{2 q}\vect{m'}- (\R_+^*)^r$ (resp.\ $\vect{s'}\in  \frac 1 p(\vect b+\vect d)- \R_+^r$ if $q=\infty$).
\end{itemize}
Then, we say that property $\atomic^{p,q}_{\vect s,\vect{s'},0}$ (resp.\ $\atomic^{p,q}_{\vect s,\vect{s'}}$) holds if for every $\delta_0>0$ there is a $(\delta,4)$-lattice $(\zeta_{j,k},z_{j,k})_{j\in J, k\in K}$, with $\delta\in ]0,\delta_0]$,  such that, defining $h_k\coloneqq \Im z_{j,k}-\Phi(\zeta_{j,k})$ for every $k\in K$ and for some (hence every) $j\in J$, the mapping
\[
\Psi\colon\lambda \mapsto \sum_{j,k} \lambda_{j,k} B_{(\zeta_{j,k},z_{j,k})}^{\vect{s'}} \Delta_\Omega^{(\vect b+\vect d)/p-\vect s-\vect{s'}}(h_k)
\]
is well defined (with locally uniform convergence of the sum) and maps $\ell^{p,q}_0(J,K) $ into $A^{p,q}_{\vect s,0}(D)$ continuously (resp.\ maps $\ell^{p,q}(J,K) $ into $A^{p,q}_{\vect s}(D)$ continuously).

If we may take $(\zeta_{j,k},z_{j,k})_{j\in J, k\in K}$, for every $\delta_0>0$ as above, in such a way that the corresponding mapping $\Psi$ is onto, then we say that property $\atomics^{p,q}_{\vect s,\vect{s'},0}$ (resp.\ $\atomics^{p,q}_{\vect s,\vect{s'}}$) holds.
\end{deff}

Because of Theorems~\ref{teo:3} and~\ref{teo:8}, it is useful to know that atomic decomposition holds for sufficiently `fine' lattices. This is the reason why we formulated properties $\atomic^{p,q}_{\vect s,\vect{s'},0}$ and $\atomic^{p,q}_{\vect s,\vect{s'}}$ (and then $\atomics^{p,q}_{\vect s,\vect{s'},0}$ and $\atomics^{p,q}_{\vect s,\vect{s'}}$) is a somewhat strong way.

Notice that the conditions on $\vect s$ and $\vect{s'}$ are natural, since they ensure that $B_{(\zeta,z)}^{\vect{s'}}\in A^{p,q}_{\vect s,0}(D)$ (resp.\ $B_{(\zeta,z)}^{\vect{s'}}\in A^{p,q}_{\vect s}(D)$) for every $(\zeta,z)\in D$.
In the following result we provide further necessary conditions for the validity of property $\atomic^{p,q}_{\vect s,\vect{s'},0}$ (resp.\ $\atomic^{p,q}_{\vect s,\vect{s'}}$).

\begin{lem}\label{lem:46}
Take $p,q\in ]0,\infty]$ and $\vect s,\vect{s'}\in \R^r$ such that  property $\atomic^{p,q}_{\vect s,\vect{s'},0}$ (resp.\ $\atomic^{p,q}_{\vect s,\vect{s'}}$) holds.
Then, the following hold:
\begin{itemize}
\item $\vect s\in \sup\left( \frac{1}{2 q}\vect m , \frac 1 p (\vect b+\vect d) +\frac{1}{2 q'} \vect{m'}\right)+(\R_+^*)^r$ (resp.\  $\vect s\in \R_+^r$ if $q=\infty$);

\item $\vect{s'}\in \frac{1}{\min(p,p')}(\vect b+\vect d)-\frac{1}{2 \min(p,p')} \vect{m'}-(\R_+^*)^r$;

\item $\vect s+\vect{s'}\in \inf\left(\frac{1}{\min(1,p)}(\vect b+\vect d)-\frac{1}{2 q'}\vect m, \frac 1 p (\vect b+\vect d)-\frac{1}{2 q} \vect{m'} \right)-(\R_+^*)^r$ if $q'<\infty$ and $\vect s+\vect{s'}\in \Big(\frac{1}{\min(1,p)}(\vect b+\vect d)-\R_+^r\Big) \cap\Big( \frac 1 p (\vect b+\vect d)-\frac{1}{2 q} \vect{m'}-(\R_+^*)^r  \Big)$ if $q'=\infty$.
\end{itemize}
\end{lem}

\begin{proof}
By Proposition~\ref{prop:60}, it will suffice to show that 
\[
B^{\vect{s'}}_{(0, i e_\Omega)}\in A^{p',q'}_{(\vect b+\vect d)/\min(1,p)-\vect s-\vect{s'}}(D).
\]
Take  $\delta_0>0$. Then, there is a $(\delta,4)$-lattice $(\zeta_{j,k},z_{j,k})_{j\in J, k\in K}$, with $\delta\meg \delta_0$, such that the mapping
\[
\ell^{p,q}_0(J,K)\ni \lambda \mapsto \sum_{j,k} \lambda_{j,k} B^{\vect{s'}}_{(\zeta_{j,k},z_{j,k})} \Delta^{(\vect b+\vect d)-\vect s-\vect{s'}}_\Omega(h_k)\in A^{p,q}_{\vect s}(D)
\]
is well defined and continuous, where $h_k\coloneqq \Im z_{j,k}-\Phi(\zeta_{j,k})$ for every $k\in K$ and for some (hence every) $j\in J$.
Observe that the continuity of the mapping $f\mapsto f(0,i e_\Omega)$ on $A^{p,q}_{\vect s}(D)$ implies that there is a constant $C_1>0$ such that
\[
\abs*{\sum_{j,k} \lambda_{j,k} B^{\vect{s'}}_{(\zeta_{j,k},z_{j,k})}(0,i e_\Omega) \Delta^{(\vect b+\vect d)-\vect s-\vect{s'}}_\Omega(h_k)}\meg C_1 \norm{\lambda}_{\ell^{p,q}_0(J,K)}
\]
for every $\lambda\in \ell^{p,q}_0(J,K)$. 
Therefore, Proposition~\ref{prop:52} implies that
\[
\Big(B^{\vect{s'}}_{(\zeta_{j,k},z_{j,k})}(0,i e_\Omega) \Delta^{(\vect b+\vect d)-\vect s-\vect{s'}}_\Omega(h_k)\Big)\in \ell^{p',q'}(J,K),
\]
so that the conclusion follows from \CR Theorems~\ref{teo:2} and~\ref{teo:8}.\CB
\end{proof}

\begin{deff}\label{37}
Take $p,q\in ]0,\infty]$ and $\vect s,\vect{s'}\in \R^r$ such that the following conditions are satisfied:
\begin{itemize}
\item $\vect s\in \frac{1}{2 q}\vect m+(\R_+^*)^r$ (resp.\ $\vect s\in \R_+^r$ if $q=\infty$);

\item $\vect{s'}\in \frac 1 p(\vect b+\vect d)-\frac{1}{2 p}\vect{m'}- (\R_+^*)^r$ (resp.\ $\vect{s'}\in  -\R_+^r$ if $p=\infty$);

\item $\vect s+\vect{s'}\in \frac 1 p(\vect b+\vect d)-\frac{1}{2 q}\vect{m'}- (\R_+^*)^r$ (resp.\ $\vect{s'}\in  \frac 1 p(\vect b+\vect d)- \R_+^r$ if $q=\infty$).
\end{itemize}
Then, we say that property $\atomic^{p,q}_{\vect s,\vect{s'},0,+}$ (resp.\ $\atomic^{p,q}_{\vect s,\vect{s'},+}$) holds if there are a $(\delta,R)$-lattice $(\zeta_{j,k},z_{j,k})_{j\in J, k\in K}$, for some $\delta>0$ and $R>1$, and a constant $C>0$ such that, defining $h_k\coloneqq \Im z_{j,k}-\Phi(\zeta_{j,k})$ for every $k\in K$ and for some (hence every) $j\in J$ and
\[
\Psi_+(\lambda)(\zeta,z)\coloneqq \sum_{j,k} \abs*{\lambda_{j,k} B_{(\zeta_{j,k},z_{j,k})}^{\vect{s'}}(\zeta,z)} \Delta_\Omega^{(\vect b+\vect d)/p-\vect s-\vect{s'}}(h_k)
\]
for every $\lambda\in \C^{J\times K}$ and for every $(\zeta,z)\in D$, one has
\[
\norm{\Psi_+(\lambda)}_{L^{p,q}_{\vect s}(D)}\meg C  \norm{\lambda}_{\ell^{p,q}(J,K)}
\]
for every $\lambda\in \ell^{p,q}_0(J,K)$ (resp.\ $\lambda \in \ell^{p,q}(J,K)$).
\end{deff}

Notice that, unlike in Definition~\ref{38},  we do not need to define properties $\atomic^{p,q}_{\vect s,\vect{s'},0,+}$ and $\atomic^{p,q}_{\vect s,\vect{s'},+}$ requiring some condition to hold for several lattices, thanks to the following result.

\begin{prop}\label{prop:16}
Take $p,q\in ]0,\infty]$ and $\vect s,\vect{s'}\in \R^r$, and assume that property $\atomic^{p,q}_{\vect s,\vect{s'},0,+}$ (resp.\ $\atomic^{p,q}_{\vect s,\vect{s'},+}$) holds. Then, for every $\delta_0>0$ and for every $R_0>0$ there is a constant $C>0$ such that, for every $\vect{s_0}\in \vect{s}+\R_+^r$, for every $\vect{s_0'}\in \vect{s'}+\vect s-\vect{s_0}-\R_+^r$ and for every $(\delta,R)$-lattice $(\zeta_{j,k},z_{j,k})_{j\in J,k\in K}$ such that $\delta\in ]0,\delta_0]$ and $R\in ]1,R_0]$,  defining $h_k\coloneqq \Im z_{j,k}-\Phi(\zeta_{j,k})$ for every $k\in K$ and
\[
\Psi_+(\lambda)(\zeta,z)\coloneqq \sum_{j,k} \abs*{\lambda_{j,k} B_{(\zeta_{j,k},z_{j,k})}^{\vect{s'_0}}(\zeta,z)} \Delta_\Omega^{(\vect b+\vect d)/p-\vect{s_0}-\vect{s'_0}}(h_k)
\]
for every $\lambda\in \C^{J\times K}$ and for every $(\zeta,z)\in D$, one has
\[
\norm*{\Psi_+(\lambda)}_{L^{p,q}_{\vect{s_0}}(D)}\meg \frac{C }{\delta^{(2 n+2 m)/p'+(1/p-1/q)_+ m}}\norm{\lambda}_{\ell^{p,q}(J,K)}
\]
for every $\lambda\in \ell^{p,q}_0(J,K)$ (resp.\ $\lambda \in \ell^{p,q}(J,K)$).
\end{prop}

In particular, property $\atomic^{p,q}_{\vect s,\vect{s'},0,+}$ (resp.\ $\atomic^{p,q}_{\vect s,\vect{s'},+}$)  implies property $\atomic^{p,q}_{\vect{s_0},\vect{s'_0},0,+}$ (resp.\ $\atomic^{p,q}_{\vect{s_0},\vect{s'_0},+}$). 

\begin{proof}
By assumption, the are a $(\delta',R')$-lattice $(\zeta'_{j',k'},z_{j',k'})_{j'\in J',k'\in K'}$ on $D$, with $\delta'>0$ and $R'>1$, and a constant $C'>0$ such that, defining $h'_{k'}\coloneqq \Im z'_{j',k'}-\Phi(\zeta'_{j',k'})$ for every $k'\in K'$
and
\[
\Psi'_+(\lambda')(\zeta,z)\coloneqq\sum_{j',k'} \abs*{\lambda'_{j',k'} B_{(\zeta'_{j',k'},z'_{j',k'})}^{\vect{s'}}(\zeta,z)} \Delta_\Omega^{(\vect b+\vect d)/p-\vect s-\vect{s'}}(h'_{k'})
\]
for every $\lambda'\in \C^{J\times K}$ and for every $(\zeta,z)\in D$, one has
\[
\norm*{\Psi'_+(\lambda')}_{L^{p,q}_{\vect s}(D)}\meg C'\norm{\lambda'}_{\ell^{p,q}(J',K')}
\]
for every $\lambda'\in \ell^{p,q}_0(J',K')$ (resp.\ $\lambda' \in \ell^{p,q}(J',K')$). 
Then, define mapping $\rho\colon J\times K\to J'\times K'$ in such a way that 
\[
d((\zeta'_{\rho(j,k)},z'_{\rho(j,k)}),(\zeta_{j,k},z_{j,k}))< R'\delta'
\]
for every $j\in J$ and for every $k\in K$, and observe that, for every $(j',k')\in J'\times K'$,
\[
\card(\rho^{-1}(j',k'))\meg \frac{\nu_D(B((0,i e_\Omega), R'\delta'+\delta))}{\nu_D(B((0,i e_\Omega)),\delta)}
\]
since the balls $B((\zeta_{j,k},z_{j,k}),\delta)$, as  $(j,k)$ runs through $ \rho^{-1}(j',k')$, are pairwise disjoint and contained in $B((\zeta'_{j',k'},z'_{j',k'}),R'\delta'+\delta)$. In addition, define
\[
K_{k'}\coloneqq \Set{k\in K\colon \exists (j,j')\in J\times J' \quad \rho(j,k)=(j',k')},
\] 
and observe that the balls $B(h_k,\delta)$, for $k\in K_{k'}$, are pairwise disjoint and contained in $B(h'_{k'},R'\delta'+\delta)$ thanks to Lemma~\ref{lem:45}. Hence,
\[
\card(K_{k'})\meg \frac{\nu_\Omega(B(e_\Omega,R'\delta'+\delta) )}{\nu_\Omega(B(e_\Omega,\delta))}.
\]
Thus, there is $C''>0$ such that 
\[
\card(\rho^{-1}(j',k'))\meg C'' \delta^{-(2n+2m)} \qquad \text{and} \qquad\card(K_{k'})\meg C''\delta^{-m}
\]
for every $j'\in J'$ and for every $k'\in K'$, provided that $\delta\meg \delta_0$ and $R\meg R_0$.
Now, take $\lambda\in \C^{J\times K}$ and define $\lambda'_{j',k'}\coloneqq \sum_{\rho(j,k)=(j',k')}\abs{\lambda_{j,k}}$, so that 
\[
\norm{\lambda'}_{\ell^{p,q}(J',K')}\meg C''^{1/p'+(1/p-1/q)_+}\delta^{-(2 n+2 m)/p'-(1/p-1/q)_+ m}\norm{\lambda}_{\ell^{p,q}(J,K)}
\]
and $\lambda'\in \ell^{p,q}_0(J',K')$ if $\lambda\in \ell^{p,q}_0(J,K)$. 
Now, Theorem~\ref{teo:2} implies that there is a constant $C'''>0$ such that
\[
\sum_{j,k} \abs*{\lambda_{j,k} B_{(\zeta_{j,k},z_{j,k})}^{\vect{s'}}} \Delta_\Omega^{(\vect b+\vect d)/p-\vect{s}-\vect{s'}}(h_k)\meg C'''\Psi'_+(\lambda') .
\] 
Therefore, it will suffice to choose $C\coloneqq C' C''^{1/p'+(1/p-1/q)_+} C'''$, since 
\[
\Delta_\Omega^{\vect{s_0}-\vect{s}  }(\Im z-\Phi(\zeta))\Delta^{\vect s+\vect{s'}-\vect{s_0}-\vect{s'_0}}_\Omega(h_k) \abs{B^{\vect{s'_0}-\vect{s'}}_{(\zeta_{j,k},z_{j,k})}(\zeta,z) }\meg1
\]
for every $(j,k)\in J\times K$ and for every $(\zeta,z)\in D$, thanks to Proposition~\ref{prop:60}.
\end{proof}

In the next result we show that property $\atomic^{p,q}_{\vect s,\vect{s'},0,+}$ (resp.\ $\atomic^{p,q}_{\vect s,\vect{s'},+}$) implies property $\atomic^{p,q}_{\vect s,\vect{s'},0}$ (resp.\ $\atomic^{p,q}_{\vect s,\vect{s'}}$). As we shall see in Theorem~\ref{teo:4}, it actually implies property $\atomics^{p,q}_{\vect s,\vect{s'},0}$ (resp.\ $\atomics^{p,q}_{\vect s,\vect{s'}}$).

\begin{prop}\label{prop:17}
Take $p,q\in ]0,\infty]$ and $\vect s,\vect{s'}\in \R^r$ such that  property $\atomic^{p,q}_{\vect s,\vect{s'},0,+}$ (resp.\ $\atomic^{p,q}_{\vect s,\vect{s'},+}$) holds. Let $(\zeta_{j,k},z_{j,k})_{j\in J,k\in K}$ be a $(\delta,R)$-lattice for some $\delta>0$ and some $R>1$. Define $h_k\coloneqq \Im z_{j,k}-\Phi(\zeta_{j,k})$ for every $k\in K$ and for some (hence every) $j\in J$. Then, the mapping
\[
\Psi\colon \lambda \mapsto \sum_{j,k}  \lambda_{j,k} B_{(\zeta_{j,k},z_{j,k})}^{\vect{s'}} \Delta_\Omega^{(\vect b+\vect d)/p-\vect s-\vect{s'}}(h_k)
\]
induces a continuous linear mapping from $\ell^{p,q}_0(J,K)$ (resp.\ from $\ell^{p,q}(J,K)$) into $A^{p,q}_{\vect s,0}(D)$ (resp.\ $A^{p,q}_{\vect s}(D)$).
\end{prop}

\begin{proof}
The assertion follows immediately from Propositions~\ref{prop:60} and~\ref{prop:16} when property $\atomic^{p,q}_{\vect s,\vect{s'},0,+}$ holds. 
Then, assume that property $\atomic^{p,q}_{\vect s,\vect{s'},+} $ holds and let us prove that the sum defining $\Psi(\lambda)$ converges locally uniformly for every $\lambda\in \ell^{p,q}(J,K)$.
Observe that Theorem~\ref{teo:2} implies that there is a constant $C>0$ such that
\[
\abs{B_{(\zeta_{j,k},z_{j,k})}^{\vect{s'}}(\zeta',z')}\meg C\abs{B_{(\zeta_{j,k},z_{j,k})}^{\vect{s'}}(\zeta,z)}
\]
for every $(\zeta,z),(\zeta',z')\in D$ such that $d((\zeta,z),(\zeta',z'))\meg 1$, and for every $(j,k)\in J\times K$. 
In addition,  there is a negligible subset $N$ of $D$ such that, for every $(\zeta,z)\in D\setminus N$,
\[
\sum_{j,k} \abs{ \lambda_{j,k} B_{(\zeta_{j,k},z_{j,k})}^{\vect{s'}}(\zeta',z')} \Delta_\Omega^{(\vect b+\vect d)/p-\vect s-\vect{s'}}(h_k)<\infty.
\]
Then, the sum defining $\Psi(\lambda)$ converges uniformly on $B((\zeta,z),1)$ for every $(\zeta,z)\in D\setminus N$, so that $\Psi(\lambda)$ is pointwise well defined and holomorphic. It is then easily verified that $\Psi(\lambda)\in A^{p,q}_{\vect s}(D)$, and that $\Psi\colon \ell^{p,q}(J,K)\to A^{p,q}_{\vect s}(D)$ is continuous.
\end{proof}

\begin{teo}\label{teo:4}
Take $\vect s,\vect{s'}\in \R^r$, $p,q\in ]0,\infty]$ such that property $\atomic^{p,q}_{\vect s,\vect{s'},0,+}$ (resp.\ $\atomic^{p,q}_{\vect s,\vect{s'},+}$) holds.
Fix $R_0>1$. Then, there is $\delta_0>0$ such that, if $\delta\in ]0,\delta_0]$ and $R\in ]1,R_0]$, then the mapping $\Psi$ defined in Proposition~\ref{prop:17} has a continuous linear section $\Psi' \colon A^{p,q}_{\vect s,0}(D)\to \ell^{p,q}_0(J,K)$ (resp.\  $\Psi'\colon A^{p,q}_{\vect s}(D)\to \ell^{p,q}(J,K)$). 
\end{teo}

Thus,  property $\atomic^{p,q}_{\vect s,\vect{s'},0,+}$ implies property $\atomics^{p,q}_{\vect s,\vect{s'},0}$, while property $\atomic^{p,q}_{\vect s,\vect{s'},+}$ implies property $\atomics^{p,q}_{\vect s,\vect{s'}}$.
Since clearly property $\atomic^{p,q}_{\vect s,\vect{s'},0,+}$ implies property $\atomic^{p,q}_{\vect s,\vect{s'},+}$, it also implies property $\atomics^{p,q}_{\vect s,\vect{s'}}$.

The proof is based on~\cite[Theorem 1.5]{RicciTaibleson}, which deals with $\C_+$.

\begin{proof}
Put a well-ordering on $J\times K$ and define 
\[
 U_{j,k} \coloneqq B((\zeta_{j,k},z_{j,k}),R\delta)\setminus \left(\bigcup_{(j',k')<(j,k)} B((\zeta_{j',k'},z_{j',k'}),R\delta) \right)
\]
for every $(j,k)\in J\times K$, so that $(U_{j,k})_{(j,k)\in J\times K}$ is a Borel measurable partition of $D$ (since $J$ and $K$ are countable). 
In addition, define $c_{j,k}\coloneqq \nu_D(U_{j,k})$ for every $(j,k)\in J\times K$, so that 
\[
\nu_D(B((0,i e_\Omega),\delta))\meg c_{j,k}\meg \nu_D(B((0,i e_\Omega),R\delta))
\]
for every $(j,k)\in J\times K$.
Then, define 
\[
S\colon A^{p,q}_{\vect s}(D)\ni f \mapsto \left(c_{j,k}\Delta_\Omega^{\vect s-(\vect b+\vect d)/p}(h_k) f(\zeta_{j,k},z_{j,k}) \right)\in \ell^{p,q}(J,K),
\]
so that Theorem~\ref{teo:8} shows that $S$ is well defined and continuous, and maps $A^{p,q}_{\vect s,0}(D)$ into $\ell^{p,q}_0(J,K)$ under the finer assumptions. In addition, define $S'\coloneqq \Psi S$ and $h_k\coloneqq \Im z_{j,k}-\Phi(\zeta_{j,k})$ for some (hence every) $j\in J$ and for every $k\in K$. 
Then, Proposition~\ref{prop:9} and Lemma~\ref{lem:46} imply that, for every $ f\in A^{p,q}_{\vect s}(D) $,
\[
\begin{split}
(f- S' f)(\zeta,z)&= \sum_{j,k}\int_{ U_{j,k} } \Big(f(\zeta',z') B^{\vect{s'}}_{(\zeta',z')}(\zeta,z) \Delta_\Omega^{-\vect{s'}}(\Im z'-\Phi(\zeta')) \\
&\qquad\qquad\qquad -f(\zeta_{j,k},z_{j,k}) B^{\vect{s'}}_{(\zeta_{j,k},z_{j,k})}(\zeta,z) \Delta_\Omega^{-\vect{s'}}(h_k) \Big) \,\dd \nu_D(\zeta',z').
\end{split}
\]
Hence, Theorem~\ref{teo:2}, Corollary~\ref{cor:34}, and Lemma~\ref{lem:28} imply that there are $\rho_0>0$ and $C_1>0$ such that
\[
\begin{split}
&\abs{(f-S' f)(\zeta,z)}\\
&\quad\meg C_1 R\delta\sum_{j,k} c_{j,k}\sup\limits_{(\zeta',z')\in B((\zeta_{j,k},z_{j,k}), R\delta+\rho_0)} \abs{f(\zeta',z')} \abs{B^{\vect{s'}}_{(\zeta_{j,k},z_{j,k})}(\zeta,z)} \Delta_\Omega^{-\vect{s'}}(h_k) 
\end{split}
\]
for every $(\zeta,z)\in D$. 
Now, let $(\zeta'_{j',k'},z'_{j',k'})_{j'\in J', k'\in K'}$ be a $(1,4)$-lattice on $D$ (cf.~Lemma~\ref{lem:32}), and observe that the proof of  Proposition~\ref{prop:56}, together with Theorem~\ref{teo:2} and Corollary~\ref{cor:34} again, implies that there is a constant $C_2>0$ such that
\[
\begin{split}
&\abs{(f-S' f)(\zeta,z)}\\
&\meg C_2 R \delta\sum_{j',k'}\sup\limits_{(\zeta',z')\in B\big((\zeta'_{j',k'},z'_{j',k'}), R\delta+\rho_0+4\big)} \abs{f(\zeta',z')} \abs{B^{\vect{s'}}_{(\zeta'_{j',k'},z'_{j',k'})}(\zeta,z)} \Delta_\Omega^{-\vect{s'}}(h'_{k'}) 
\end{split}
\]
for every $(\zeta,z)\in D$, where $h'_{k'}= \Im z'_{j',k'}-\Phi(\zeta'_{j',k'})$ for every $k'\in K'$ and for some (hence every) $j'\in J'$.
Hence, Theorem~\ref{teo:8} and Proposition~\ref{prop:16} show that there is a constant $C_3>0$ such that, if $R\meg R_0$ and $\delta\meg 1$, then
\[
\norm{f-S' f}_{A^{p,q}_{\vect s}(D)}\meg C_3 \delta \norm{f}_{A^{p,q}_{\vect s}(D)}.
\]
Take $\delta_0>0$ so that $C_3 \delta_0\meg \frac 1 2$. Then, 
\[
\norm*{\sum_{j\Meg k}   (I-S')^j f }_{A^{p,q}_{\vect s}(D)}^{\min(1,p,q)}\meg \sum_{j\Meg k} 2^{-\min(1,p,q) j} \norm{f}_{A^{p,q}_{\vect s}(D)}^{\min(1,p,q)}
\]
for every $k\in \N$, so that $\sum_{j\in \N} (I- S')^j$ is a well defined endomorphism of $A^{p,q}_{\vect s,0}(D) $ (resp.\ $A^{p,q}_{\vect s}(D)$), and is the inverse of $S'$. 
Hence, 
\[
\Psi'\coloneqq S\sum_{j\in \N} (I-S')^j
\]
is a well defined and continuous linear mapping from $A^{p,q}_{\vect s,0}(D)$ into $\ell^{p,q}_0(J,K)$ (resp.\ from $A^{p,q}_{\vect s}(D)$ into $\ell^{p,q}(J,K)$), and $\Psi \Psi'= S'\sum_{j\in \N}(I-S')^j=I$.
\end{proof}

In the next result we shall provide some sufficient conditions for properties  $\atomic^{p,q}_{\vect s,\vect{s'},0,+}$ and $\atomic^{p,q}_{\vect s,\vect{s'},+}$.

\begin{teo}\label{teo:5}
Take $p,q\in ]0,\infty]$ and $\vect s, \vect{s'}\in \R^r$ such that the following hold:
\begin{itemize}
\item $\vect s\in \frac{1}{2 q}\vect m+\left(\frac{1}{2 \min(1,p)}-\frac{1}{2 q} \right)_+\vect{m'}+(\R_+^*)$;

\item $\vect{s'}\in \frac{1}{\min(1,p)} (\vect b+\vect d)-\frac{1}{2 \min(1,p)} \vect{m'}-(\R_+^*)^r$;

\item $\vect s +\vect{s'}\in \frac{1}{\min(1,p)} (\vect b+\vect d)-\frac{1}{2 q} \vect{m'}-\left(\frac{1}{2 \min(1,p)}-\frac{1}{2 q} \right)_+\vect{m}-(\R_+^*)^r$;
\end{itemize}
Then, properties $\atomic^{p,q}_{\vect s,\vect{s'},0,+}$ and $\atomic^{p,q}_{\vect s,\vect{s'},+}$ hold.
\end{teo}

This result is optimal when $q\meg p\meg 1$ and gives atomic decomposition for all spaces $A^{p,q}_{\vect s}(D)$ such that $q\meg \min(1,p)$.

The proof of the case $p,q\in [1,\infty]$ is based on~\cite[Theorem 4.10]{Bekolleetal}, which deals with irreducible symmetric tube domains. 
The strategy of~\cite[Lemma 5.1]{RicciTaibleson}, which gives optimal results for $\C_+$, leads to worse results in the general case.

We shall prepare the proof of Theorem~\ref{teo:5} by means of some simple consequences of Schur's lemma.

\begin{lem}\label{lem:79}
Take $q\in [1,\infty]$ and $\vect{s_1},\vect{s_2}\in \R^r$ such that the following conditions hold:
\begin{enumerate}
\item[\em(1)] $\vect{s_2}\in \frac{1}{2 q'}\vect{m}+\frac{1}{2 q}\vect{m'}+(\R_+^*)^r$;

\item[\em(2)] $\vect{s_1}+\vect{s_2}\in-\frac{1}{2 q}\vect{m}-\frac{1}{2 q'}\vect{m'}-(\R_+^*)^r $.
\end{enumerate}
Then, the mapping
\[
T\colon C_c(\Omega)\ni f \mapsto \Delta_\Omega^{-\vect{s_1}-\vect{s_2}} \int_\Omega f(h) \Delta_\Omega^{\vect{s_1}}(\,\cdot\,+h)\Delta_\Omega^{\vect{s_2}}(h)\,\dd \nu_\Omega(h)\in C(\Omega)
\]
induces endomorphisms of $L^q_0(\nu_\Omega)$ and of $L^q(\nu_\Omega)$.
\end{lem}

The proof is based on~\cite[Theorem 4.5]{Bekolleetal}, which deals with the case in which $\Omega$ is irreducible and symmetric, and $\vect {s_1},\vect{s_2}\in \R \vect{1}_r$.

\begin{proof}
Assume first that $q=1$. Then, for every $f\in C_c(\Omega)$,
\[
\norm{T f}_{L^1(\nu_\Omega)}\meg \int_{\Omega\times \Omega} \abs{f(h)} \Delta_\Omega^{-\vect{s_1}-\vect{s_2}}(h') \Delta_\Omega^{\vect{s_1}}(h+h')\Delta_\Omega^{\vect{s_2}}(h)\,\dd (\nu_\Omega \times \nu_\Omega)(h,h').
\]
Therefore, the assertion follows from Corollary~\ref{cor:10}.

Next, assume that $q=\infty$. Then, for every  $f\in L^\infty(\nu_\Omega)$, by an abuse of notation,
\[
\norm*{T f}_{L^\infty(\nu_\Omega)}\meg  \norm{f}_{L^\infty(\nu_\Omega )} \sup\limits_{h'\in \Omega} \Delta_\Omega^{-\vect{s_1}-\vect{s_2}}(h') \int_\Omega  \Delta_\Omega^{\vect{s_1}}(h+h')\Delta_\Omega^{\vect{s_2}}(h)\,\dd \nu_\Omega(h),
\]
so that by means of Proposition~\ref{prop:58} we see that $T$ induces an endomorphism of $L^\infty(\nu_\Omega)$ if conditions~{(1)} and~{(2)} are satisfied. Further, by means of Lemma~\ref{lem:23} we see that $ \Delta_\Omega^{-\vect{s_1}-\vect{s_2}} \Delta_\Omega^{\vect{s_1}}(\,\cdot\,+h)\in C_0(\Omega)$ for every $h\in \Omega$, so that $T$ induces also an endomorphism of $C_0(\Omega)$ under the same assumptions. 

Finally, assume that $q\in ]1,\infty[$. 
Define 
\[
T'\colon C_c(\Omega)\ni g \mapsto \Delta_\Omega^{\vect{s_2}} \int_\Omega g(h') \Delta_\Omega^{-\vect{s_1}-\vect{s_2}}(h') \Delta_\Omega^{\vect{s_1}}(\,\cdot\,+h')\,\dd \nu_\Omega(h')\in C(\Omega),
\]
so that $\trasp T (g \cdot \nu_\Omega)= (T' g)\cdot \nu_\Omega$ for every $g\in C_c(\Omega)$.
Take $\vect{s'}\in \R^r$, and observe that Corollary~\ref{cor:10} implies that there are constants $C_1,C_2>0$ such that
\[
T(\Delta_\Omega^{q' \vect{s'}} )=C_1 \Delta_\Omega^{q'\vect{s'}} \qquad \text{and}\qquad T'(\Delta_\Omega^{q \vect{s'}} )=C_2 \Delta_\Omega^{q\vect{s'}} ,
\]
provided that the following conditions hold:
\begin{enumerate}
\item[(i)] $q'\vect{s'}+\vect{s_2}\in \frac 1 2 \vect m+(\R_+^*)^r$;

\item[(ii)] $q'\vect{s'}+\vect{s_1}+\vect{s_2}\in- \frac 1 2 \vect{m'}-(\R_+^*)^r$;

\item[(iii)] $q\vect{s'}-\vect{s_1}-\vect{s_2}\in \frac 1 2 \vect m+(\R_+^*)^r$;

\item[(iv)] $q\vect{s'}-\vect{s_2}\in- \frac 1 2 \vect{m'}-(\R_+^*)^r$;
\end{enumerate}
It is then clear that we may find  $\vect{s'}$ satisfying conditions~{(i)} to~{(iv)} if and only if conditions~{(1)} and~{(2)} are satisfied, in which case~\cite[Lemma I.2]{Grafakos} implies that $T$ induces an endomorphism of $L^q(\nu_\Omega)$.
\end{proof}

\begin{cor}\label{cor:33}
Take $q\in [1,\infty]$ and $\vect{s_1},\vect{s_2}\in \R^r$ such that the following conditions hold:
\begin{enumerate}
\item[\em(1)] $\vect{s_2}\in \frac{1}{2 q'}\vect{m}+\frac{1}{2 q}\vect{m'}+(\R_+^*)^r$;

\item[\em(2)] $\vect{s_1}+\vect{s_2}\in-\frac{1}{2 q}\vect{m}-\frac{1}{2 q'}\vect{m'}-(\R_+^*)^r $.
\end{enumerate}
Let $(h_k)_{k\in K}$ be a $(\delta,R)$-lattice on $\Omega$, with $\delta>0$ and $R>1$.
Then, the mapping
\[
T\colon \C^{(K)}\ni \lambda \mapsto  \Delta_\Omega^{-\vect{s_1}-\vect{s_2}} \sum_{k\in K} \lambda_k \Delta_\Omega^{\vect{s_1}}(\,\cdot\,+h_k)\Delta_\Omega^{\vect{s_2}}(h_k) \in C(\Omega)
\]
induces continuous linear mappings $\ell^q_0(K)\to L^q_0(\nu_\Omega)$ and $\ell^q(K)\to L^q(\nu_\Omega)$.
\end{cor}

\begin{proof}
By Lemma~\ref{lem:79}, let $T'$ be the endomorphism of $L^q(\nu_\Omega)$ such that, for every $f\in L^q(\nu_\Omega)$,
\[
(T' f)(h)=\Delta_\Omega^{-\vect{s_1}-\vect{s_2}} (h)\int_\Omega f(h) \Delta_\Omega^{\vect{s_1}}(h+h')\Delta_\Omega^{\vect{s_2}}(h')\,\dd \nu_\Omega(h').
\]
In addition, endow $K$ with a well-ordering and define 
\[
U_k \coloneqq B(h_k,R\delta)\setminus \left( \bigcup_{k'<k} B(h_{k'},R\delta)\right) 
\] for every $k\in K$, so that $(U_k )$ is a measurable partition of $\Omega$ and $B(h_k,\delta)\subseteq U_k  \subseteq B(h_k,R\delta)$ for every $k\in K$. Therefore, by means of Corollary~\ref{cor:34} we see that there is a constant $C>0$ such that the linear mapping
\[
\Psi\colon \ell^q(K)\ni \lambda \mapsto \sum_{k\in K} \lambda_k \chi_{ U_k  }\in L^q(\nu_\Omega)
\]
is well-defined and satisfies 
\[
\frac{1}{C} \norm{\lambda}_{\ell^q(K)}\meg \norm{\Psi(\lambda)}_{L^q(\nu_\Omega)}\meg C\norm{\lambda}_{\ell^q(K)}
\]
for every $\lambda\in \ell^q(K)$.
In addition, Corollary~\ref{cor:34} implies that there is a constant $C'>0$ such that, for every positive $\lambda\in \ell^q(K)$,
\[
\frac{1}{C'}(T'\Psi) (\lambda)\meg T(\lambda)\meg C' (T'\Psi)(\lambda),
\]
so that $T$ induces a continuous linear mappings $\ell^q(K)\to L^q(\nu_\Omega)$. To see that $T(c_0(K))\subseteq C_0(\nu_\Omega)$, it suffices to observe that $\Delta_\Omega^{-\vect{s_1}-\vect{s_2}}  \Delta_\Omega^{\vect{s_1}}(\,\cdot\,+h_k)\in C_0(\Omega)$ for every $k\in K$, thanks to Lemma~\ref{lem:23}.
\end{proof}

\begin{proof}[Proof of Theorem~\ref{teo:5}.]
Take a $(\delta,R)$-lattice $(\zeta_{j,k},z_{j,k})_{j\in J,k\in K}$ on $D$ for some $\delta>0$ and some $R>1$ (cf.~\ref{lem:32}).
Define $h_k\coloneqq \Im z_{j,k}-\Phi(\zeta_{j,k})$ for every $k\in K$ and for some (hence every) $j\in J$.
In addition, define 
\[
B_{j,k}^{\vect{s''}}\coloneqq B^{\vect{s''}}_{(\zeta_{j,k},z_{j,k})} 
\]
for every $\vect{s''}\in \R^r$ and for every $(j,k)\in J\times K$, to simplify the notation. 
Further, for every $\lambda\in \ell^{p,q}(J,K)$ define
\[
\Psi(\lambda)\coloneqq \sum_{j,k} \abs{\lambda_{j,k} B^{\vect {s'}}_{j,k} } \Delta_\Omega^{(\vect b+\vect d)/p-\vect s-\vect{s'}}(h_k)\in [0,\infty].
\]

\textsc{Step I.} Assume first that $q\meg p\meg 1$. Then, for every $h\in \Omega$,
\[
\norm{\Psi(\lambda)_h}_{L^p(\Nc)}^p\meg\sum_{j,k}\abs{\lambda_{j,k}}^p   \Delta_\Omega^{\vect b+\vect d-p(\vect s+\vect{s'})}(h_k)  \norm{(B^{\vect {s'}}_{j,k}) _h }_{L^p(\Nc)}^p.
\]
In addition, Lemma~\ref{lem:21} shows that there is a constant $C_1>0$ such that 
\[
\norm{(B^{\vect {s'}}_{j,k}) _h }_{L^p(\Nc)}^p =C_1 \Delta_\Omega^{p\vect{s'}-(\vect b+\vect d)}(h+h_k) 
\]
for every $(j,k)\in J\times K$ and for every $h\in \Omega$.
Therefore, using the subadditivity of the mapping $x\mapsto x^{q/p}$ on $\R_+$,
\[
\begin{split}
\norm{\Psi(\lambda)}_{L^{p,q}_{\vect s}(D)}^q&\meg C_1^{q/p} \sum_{k}  \left( \sum_j \abs{\lambda_{j,k}}^p\right) ^{q/p} \Delta_\Omega^{ q((\vect b+\vect d)/p-(\vect s+\vect{s'}))}(h_k) \times \\
&\qquad \times\int_\Omega  \left( \Delta_\Omega^{(\vect{s'}-(\vect b+\vect d)/p)}(h+h_k) \Delta^{\vect s}_\Omega(h)\right) ^q\,\dd \nu_\Omega(h).
\end{split}
\]
Now, Corollary~\ref{cor:10} implies that  there is a constant $C_2>0$ such that
\[
\int_\Omega  \left( \Delta_\Omega^{(\vect{s'}-(\vect b+\vect d)/p)}(h+h_k) \Delta^{\vect s}_\Omega(h)\right) ^q\,\dd \nu_\Omega(h)= C_2 \Delta_\Omega^{ q(\vect s+\vect{s'}-(\vect b+\vect d)/p)}(h_k)
\]
for every $k\in K$. Hence,
\[
\norm{\Psi(\lambda)_h}_{L^{p,q}_{\vect s}(D)}\meg C_1^{1/p} C_2^{1/q} \norm{\lambda }_{\ell^{p,k}(J,K)},
\]
whence the result in this case.

\textsc{Step II.} Assume, now, that $q\Meg p\meg 1$, and define
\[
T\colon  \lambda' \mapsto \Delta_\Omega^{p\vect s} \sum_{k\in K} \lambda'_k \Delta_\Omega^{(\vect b+\vect d)-p(\vect s+\vect{s'})}(h_k) \Delta_\Omega^{p\vect{s'}-(\vect b+\vect d)}(\,\cdot\,+h_k) ,
\]
so that $T$ maps continuously $\ell^{q/p}_0(K)$ and  $\ell^{q/p}(K)$ into $L^{q/p}_0(\nu_\Omega)$ and $ L^{q/p}(\nu_\Omega)$, respectively, by Corollary~\ref{cor:33}.
Arguing as in~{step I}, we then see that 
\[
\norm{\Psi(\lambda)}_{L^{p,q}_{\vect s}(D)}\meg C_1 \norm*{ T\left(\left( \sum_j \abs{\lambda_{j,k}}^p\right) _k\right) }_{L^{q/p}(\nu_\Omega)}^{1/p} 
\]
for every $\lambda\in \ell^{p,q}(J,K)$, whence the assertion in this case.

\textsc{Step III.} Assume, now, that $q\meg 1 \meg p$. For every $(j,k)\in J\times K$, choose $\tau_{j,k}\in C_c(\Omega)$ so that $\chi_{B((\zeta_{j,k},z_{j,k}),\delta/2)}\meg \tau_{j,k}\meg \chi_{B((\zeta_{j,k},z_{j,k}),\delta)}$, and define $C_3\coloneqq \sup\limits_{h\in \Omega} \int_\Nc \left( \chi_{B((0,i e_\Omega),\delta)}\right)_h\,\dd \Hc^{2 n+m}$. Define $\varrho  \colon D\ni (\zeta,z)\mapsto \Im z-\Phi(\zeta)\in \Omega$ and
\[
\Psi'\colon \ell^{p,q}(J,K)\ni \lambda \mapsto \sum_{j,k} \lambda_{j,k} \tau_{j,k} \Delta_\Omega^{(\vect b+\vect d)/p-\vect s}\circ \varrho  \in C(D),
\]
and let us prove that $\Psi'$ maps continuously $\ell^{p,q}_0(J,K)$ and $\ell^{p,q}(J,K)$ into $L^{p,q}_{\vect s,0}(D)$ and $L^{p,q}_{\vect s}(D)$, respectively. Indeed, take $\lambda\in \ell^{p,q}(J,K)$, and observe that
\[
\begin{split}
\norm{(\Psi'( \lambda))_h}_{L^p(\Nc)}&\meg \Delta_\Omega^{(\vect b+\vect d)/p-\vect s}(h)\sum_{k\in K} \norm*{\sum_{j\in J} \abs{\lambda_{j,k}}\left( \chi_{B((\zeta_{j,k},z_{j,k}),\delta)} \right)_h  }_{L^p(\Nc)}\\
&\meg C_3^{1/p}\sum_{k\in K}\frac{\Delta_\Omega^{(\vect b+\vect d)/p-\vect s}(h)}{\Delta_\Omega^{-(\vect b+\vect d)/p}(h_k)}\chi_{B(h_k,\delta)}(h) \norm*{(\lambda_{j,k})_j}_{\ell^p(J)}
\end{split}
\]
for every $h\in \Omega$, so that Corollary~\ref{cor:34} there is a constant $C_3'>0$ such that
\[
\begin{split}
\norm{\Psi'(\lambda)}_{L^{p,q}_{\vect s}(D)}&\meg  C_3' \norm*{\sum_{k\in K}\chi_{B(h_k,\delta)} \norm*{(\lambda_{j,k})_j}_{\ell^p(J)}  }_{L^q(\nu_\Omega)}\\
&= C_3' \nu_\Omega(B_\Omega(e_\Omega,\delta))^{1/q} \norm{\lambda}_{\ell^{p,q}(J,K)}.
\end{split}
\]
Thus, $\Psi$ induces a continuous linear mapping $\ell^{p,q}(J,K)\to L^{p,q}_{\vect s}(D)$.
In addition, since   $\Psi'(\C^{(J\times K)})\subseteq C_c(D)$, we also see that $\Psi$ induces a continuous linear mapping $\ell^{p,q}_0(J,K)\to L^{p,q}_{\vect s,0}(D)$.

Now, observe that Theorem~\ref{teo:2} and Corollary~\ref{cor:34} imply that there is a constant $C_4>0$ such that
\[
\begin{split}
&\abs{\Psi(\lambda)_h(\zeta,x)}\\
&\qquad\meg  C_4 \int_\Omega \int_\Nc \Psi'(\abs{\lambda})_{h'}(\zeta',x') \abs*{\left(B_{(\zeta',z')}^{\vect{s'}} \right)_{h}(\zeta,x) } \,\dd (\zeta',x')\Delta_\Omega^{\vect b+\vect d-\vect{s'}}(h')\,\dd \nu_\Omega(h')  \\
&\qquad= C_4 \int_\Omega \Psi'(\abs{\lambda})_{h'}* \abs*{\left(B_{(0, i h')}^{\vect{s'}} \right)_{h}(\zeta,x) } \Delta_\Omega^{\vect b+\vect d-\vect{s'}}(h')\,\dd \nu_\Omega(h')  \\
\end{split}
\]
for every $\lambda\in \ell^{p,q}(J,K)$, for every $h\in \Omega$, and for every $(\zeta,x)\in \Nc$.
Therefore, Young's inequality and Lemma~\ref{lem:21} show that there is a constant $C_5>0$ such that
\[
\norm{\Psi(\lambda)_h}_{L^p(\Nc)}\meg C_5 \int_\Omega \norm{\Psi'(\abs{\lambda})_{h'}}_{L^p(\Nc)}\Delta_\Omega^{\vect{s'}-(\vect b+\vect d)}(h+h')\Delta_\Omega^{\vect b+\vect d-\vect{s'}}(h') \,\dd \nu_\Omega(h')
\]
for every $\lambda\in \ell^{p,q}(J,K)$, and for every \CR$h\in \Omega$, \CB since $\vect{s'}\in \vect b+\vect d-\frac 1 2 \vect{m'}-(\R_+^*)^r$.
Therefore, the preceding arguments and Corollary~\ref{cor:34} imply that there is a constant $C_6>0$ such that
\[
\norm{\Psi(\lambda)_h}_{L^p(\Nc)}\meg C_6\sum_{k\in K}\norm{(\lambda_{j,k})_j}_{\ell^p(J)} \Delta_\Omega^{\vect{s'}-(\vect b+\vect d)}(h+h_k)\Delta_\Omega^{\vect b+\vect d-\vect s-\vect{s'}}(h_k),
\]
so that, by the subadditivity of the mapping $x\mapsto x^q$ on $\R_+$,
\[
\begin{split}
\norm{\Psi(\lambda)}_{L^{p,q}_{\vect s}(D)}^q&\meg C_6^q \sum_k \norm{(\lambda_{j,k})_j}_{\ell^p(J)}^q\Delta_\Omega^{q(\vect b+\vect d-\vect s-\vect{s'})}(h_k)\times\\
&\qquad\times \int_\Omega \left( \Delta_\Omega^{\vect s}(h) \Delta_\Omega^{\vect{s'}-(\vect b+\vect d)}(h+h_k)\right) ^q \,\dd\nu_\Omega(h).
\end{split}
\]
Hence, Corollary~\ref{cor:10}  implies that there is a constant $C_7>0$ such that 
\[
\norm{\Psi(\lambda)}_{L^{p,q}_{\vect s}(D)}\meg C_7 \norm*{\lambda }_{\ell^{p,q}(J,K)}.
\]
In order to conclude, it suffice to observe that $\Psi(\C^{(J\times K)})\subseteq L^{p,q}_{\vect s,0}(D)$.

\textsc{Step IV.} Finally, assume that $p,q\Meg 1$. Define $\Psi'$ as in~{step III}, and define 
\[
T'\colon  f \mapsto \Delta_\Omega^{\vect s} \int_\Omega f(h')\Delta_\Omega^{\vect{s'}-(\vect b+\vect d)}(\,\cdot\,+h')\Delta_\Omega^{\vect b+\vect d-\vect s-\vect{s'}}(h') \,\dd \nu_\Omega(h').
\]
so that $T'$ induces endomorphisms of $L^q_0(\nu_\Omega)$ and $L^q(\nu_\Omega)$ by Lemma~\ref{lem:79}. Therefore, taking~{step III} into account,
\[
\begin{split}
\norm{\Psi(\lambda)}_{L^{p,q}_{\vect s}(D)}&\meg C_5 \norm{T'}_{\Lin(L^q(\nu_\Omega))} \norm{\Psi'(\abs{\lambda})}_{L^{p,q}_{\vect s}(D)}\\
&\meg C_5 \norm{T'}_{\Lin(L^q(\nu_\Omega))} C_3' \nu_\Omega(B_\Omega(e_\Omega,\delta))^{1/q} \norm{\lambda}_{\ell^{p,q}(J,K)}.
\end{split}
\] 
Since clearly $\Psi(\C^{(J\times K)})\subseteq L^{p,q}_{\vect s,0}(D)$, the assertion follows.
\end{proof}

\section{Duality}\label{sec:6:5}

In this section we shall explore the duality induced by the sesquilinear forms of Corollary~\ref{cor:13}. As we shall see, this topic is closely related to that of atomic decomposition. On the one hand, atomic decomposition on $A^{p,q}_{\vect s,0}(D)$ allows us to identify its dual with a space of holomorphic functions. On the other hand, the dual of $A^{p,q}_{\vect s,0}(D)$ always enjoys an analogue of the atomic decomposition studied in Section~\ref{sec:6:4}.

We first show how atomic decomposition implies a simple characterization of the dual of weighted Bergman spaces.
Notice that, as the proof shows, the assumptions that property $\atomics^{p,q}_{\vect s, (\vect b+\vect d)/\min(1,p) -\vect s- \vect{s'},0}$ holds can be slightly weakened. Indeed, it suffices to assume that property $\atomic^{p,q}_{\vect s, (\vect b+\vect d)/\min(1,p) -\vect s- \vect{s'},0}$ and that the $B^{(\vect b+\vect d)/\min(1,p)-\vect s-\vect{s'}}_{(\zeta,z)}$, as $(\zeta,z)$ runs through $D$, form a total subset of $A^{p,q}_{\vect s,0}(D)$ (i.e., generate a dense vector subspace of $A^{p,q}_{\vect s,0}(D)$).

\begin{prop}\label{prop:29}
Take $\vect{s},\vect{s'}\in \R^r$ and $p,q\in ]0,\infty]$ such that property 
\[
\atomics^{p,q}_{\vect s, (\vect b+\vect d)/\min(1,p) -\vect s- \vect{s'},0}
\]
holds. 
Then, the sesquilinear form on $A^{p,q}_{\vect s,0}(D)\times A^{p',q'}_{\vect{s'}}(D)$
\[
(f,g)\mapsto \int_D f(\zeta,z) \overline{g(\zeta,z)}\Delta_\Omega^{\vect s+\vect{s'}-(\vect b+\vect d)/\min(1,p)}(\Im z-\Phi(\zeta))\,\dd\nu_D (\zeta,z)
\]
induces a antilinear isomorphism of $A^{p',q'}_{\vect{s'}}(D)$ onto $A^{p,q}_{\vect s,0}(D)'$.
\end{prop}

The proof is based on~\cite[Theorem 8.2]{RicciTaibleson}, which deals with the case $D=\C_+$. We shall prepare the proof by means of the following lemma, which is of independent interest.

\begin{lem}\label{lem:35}
Take $p,q\in ]0,\infty[$ and $\vect{s}, \vect{s'}\in \R^r$ such that property $\atomic^{p,q}_{\vect s,\vect{s'},0}$ holds. 
For every  $L\in(A^{p,q}_{\vect s,0}(D))'$, define 
\[
F(L)\colon (\zeta,z)\mapsto \overline{\big\langle L, B^{\vect{s'}}_{(\zeta,z)}\big\rangle}.
\]
Then, the antilinear mapping
\[
F\colon (A^{p,q}_{\vect s,0}(D))'\to A^{p',q'}_{(\vect b+\vect d)/\min(1,p)-\vect s-\vect{s'}}(D)
\]
is continuous.
\end{lem}

The proof is based on~\cite[Lemma 8.1]{RicciTaibleson}, which deals with the case $D=\C_+$.

\begin{proof}
Let us first prove that $F(L)$ is holomorphic. Observe that Theorem~\ref{teo:2} implies that  there are two constants $R,C>0$ such that
\[
\frac{\abs{B^{\vect{s'}}_{(\zeta',z')}(\zeta,z) - B^{\vect{s'}}_{(\zeta'',z'')}(\zeta,z) }}{d((\zeta',z'),(\zeta'',z''))}\meg C\abs{B^{\vect{s'}}_{(\zeta'',z'')}(\zeta,z) }
\]
for every $(\zeta,z),(\zeta',z'),(\zeta'',z'')\in D$ such that $d((\zeta',z'),(\zeta'',z''))\meg R$. 
Taking into account the fact that $d((\zeta',z'),(\zeta'',z'')) \asymp \norm{(\zeta',z')-(\zeta'',z'')} $ as $(\zeta',z')\to (\zeta'',z'')$,  Proposition~\ref{prop:60} and the dominated convergence theorem imply that the mapping  $F(L)$ is holomorphic on $D$.
Now, by assumption, for every $\delta_0>0$ there is a $(\delta,4)$-lattice $(\zeta_{j,k},z_{j,k})_{j\in J, k\in K}$  on $D$ for some $\delta\in ]0,\delta_0]$ such that the mapping
\[
\Psi\colon\ell^{p,q}_0(J,K)\ni\lambda\mapsto  \sum_{j,k}\lambda_{j,k}\Delta_\Omega^{(\vect b+\vect d)/p-\vect s-\vect{s'}}(h_k) B^{\vect{s'}}_{(\zeta_{j,k}, z_{j,k})}\in A^{p,q}_{\vect{s},0}(D)
\]
is well defined and continuous. \CR Observe that $F(L)\in A^{\infty,\infty}_{\vect{s''}-(\vect b+\vect d)/p'}(D)$ thanks to Proposition~\ref{prop:60} and Lemma~\ref{lem:46}.

Then, Theorem~\ref{teo:3} shows that, if we choose $\delta_0$ sufficiently small, then there is a constant $C>0$ such that, setting $\vect{s''}\coloneqq (\vect b+ \vect d)/\min(1,p)-\vect s-\vect{s'}$,
\[
\norm{F(L)}_{A^{p',q'}_{\vect{s''}}(D)}\meg C \norm*{\left( \Delta_\Omega^{(\vect b+ \vect d)/p-\vect s-\vect{s'}}(h_k)F(L)(\zeta_{j,k},z_{j,k})\right) _{j,k}}_{\ell^{p',q'}(J,K)},
\] 
where $h_k\coloneqq \Im z_{j,k}-\Phi(\zeta_{j,k})$ for some (hence every) $j\in J$ and for every $k\in K$. \CB
Then,
\[
\norm{F(L)}_{A^{p',q'}_{\vect{s''}}(D)}\meg C\sup\limits_{\norm{\lambda}_{\ell_0^{p,q}(J,K)}\meg 1} \abs*{\sum_{j,k}\lambda_{j,k}\Delta_\Omega^{(\vect b+ \vect d)/p-\vect s-\vect{s'}}(h_k) F(L)(\zeta_{j,k}, z_{j,k})},
\]
thanks to Proposition~\ref{prop:52}, so that
\[
\begin{split}
\norm{F(L)}_{A^{p',q'}_{\vect{s''}}(D)}&\meg C\sup\limits_{\norm{\lambda}_{\ell^{p,q}_0(J,K)}\meg 1} \abs*{\bigg\langle L,\sum_{j,k}\lambda_{j,k}\Delta_\Omega^{(\vect b+ \vect d)/p-\vect s-\vect{s'}}(h_k)B^{\vect{s'}}_{(\zeta_{j,k}, z_{j,k})}\bigg\rangle}\\
&\meg C \norm{\Psi}_{\Lin(\ell^{p,q}_0(J,K); A^{p,q}_{\vect s,0}(D))}\norm{L}_{A^{p,q}_{\vect s,0}(D)'}.
\end{split}
\]
Therefore, $F(L)\in A^{p',q'}_{\vect{s''}}(D)$ and the asserted continuity holds.
\end{proof}

\begin{proof}[Proof of Proposition~\ref{prop:29}.]
The given sesquilinear pairing is continuous by Corollary~\ref{cor:13}, so that it induces a continuous antilinear mapping of $A^{p',q'}_{\vect{s'}}(D)$ into $A^{p,q}_{\vect s,0}(D)'$. Then, take $L\in A^{p,q}_{\vect s,0}(D)'$, and observe that Lemma~\ref{lem:35} implies that the mapping $F(L)\colon(\zeta,z)\mapsto \overline{\Big\langle L, B^{  (\vect b+\vect d)/\min(1,p)-\vect s-\vect{s'}}_{(\zeta,z)}\Big\rangle}$ belongs to $A^{p',q'}_{\vect{s'} }(D)$. Then, Proposition~\ref{prop:9} and Lemma~\ref{lem:46} imply that, for every $(\zeta,z)\in D$,
\[
\begin{split}
\Big\langle L, B^{  (\vect b+\vect d)/\min(1,p)-\vect s-\vect{s'}}_{(\zeta,z)}\Big\rangle&=\overline{F(L)(\zeta,z)}\\
&=\Big\langle B^{  (\vect b+\vect d)/\min(1,p)-\vect s-\vect{s'}}_{(\zeta,z)}\Big\vert F(L)\Big\rangle.
\end{split}
\]
Now, by assumption the set of $B^{  (\vect b+\vect d)/\min(1,p)-\vect s-\vect{s'}}_{(\zeta,z)}$, for $(\zeta,z)\in D$, is total in $A^{p,q}_{\vect s,0}(D)$, so that $\langle L, f\rangle=\langle f \vert F(L)\rangle$ for every $f\in A^{p,q}_{\vect s,0}(D)$ by continuity. The assertion follows.
\end{proof}

We now show that the dual of a weighted Bergman space has suitable atomic decompositions. This is a consequence of \CR Theorem~\ref{teo:3}, by transposition.\CB

\begin{prop}\label{prop:13}
Take $p,q\in ]0,\infty]$, and  $\vect{s},\vect{s'}\in \R^r$ such that the following conditions hold:
\begin{itemize}
\item $\vect s\in \sup\left( \frac{1}{2 q}\vect m, \frac 1 p (\vect b+ \vect d) +\frac{1}{2 q'} \vect{m'} \right)+(\R_+^*)^r$;

\item $\vect{s'}\in \inf\left( \frac{1}{p'}(\vect b+\vect d)-\frac{1}{2 p'}\vect{m'}, \vect b+\vect d-\frac 1 2 \vect m\right) -(\R_+^*)^r$;

\item $\vect s+\vect{s'}\in \frac{1}{\min(1,p)}(\vect b+\vect d)-\frac{1}{2 q'} \vect m-(\R_+^*)^r$.
\end{itemize}
Denote by $\iota\colon A^{p',q'}_{(\vect b+\vect d)/\min(1,p)-\vect s-\vect{s'}}(D)\to A^{p,q}_{\vect s,0}(D)'$ the continuous antilinear mapping induced by the sesquilinear form on $A^{p,q}_{\vect s,0}(D)\times A^{p',q'}_{(\vect b+\vect d)/\min(1,p)-\vect s-\vect{s'}}(D)$
\[
(f,g)\mapsto \int_D f(\zeta,z) \overline{g(\zeta,z)}\Delta_\Omega^{-\vect{s'}}(\Im z-\Phi(\zeta))\,\dd\nu_D (\zeta,z),
\]
and take a $(\delta,R)$-lattice $(\zeta_{j,k}, z_{j,k})_{j\in J, k\in K}$ on $D$ for some $\delta>0$ and $R>1$.
Then, the following hold:
\begin{enumerate}
\item[\em(1)] for every $\lambda\in \ell_0^{p',q'}(J,K)$ (resp.\ $\ell^{p',q'}(J,K)$), the sum 
\[
\Psi(\lambda)\coloneqq \sum_{j,k} \lambda_{j,k}\, \iota\big(B^{\vect{s'}}_{(\zeta_{j,k},z_{j,k})} \big)\Delta_\Omega^{\vect{s'}-(\vect b+\vect d)/p}(h_k)
\]
converges strongly (resp.\ weakly) in $A^{p,q}_{\vect s,0}(D)'$;

\item[\em(2)] for every $R_0>1$ and $\delta_0>0$ there is a constant $C>0$ such that 
\[
\norm{\Psi(\lambda)}_{A^{p,q}_{\vect s,0}(D)'}\meg C\delta^{-(2 n+m)/p-m/q} \norm{\lambda}_{\ell^{p',q'}(J,K)}
\]
whenever $\delta \meg \delta_0$ and $R\meg R_0$; 

\item[\em(3)] for every $R_1>1$ there is $\delta_1>0$ such that $\Psi$ is a strict morphism of $\ell^{p',q'}(J,K)$ onto $A^{p,q}_{\vect s,0}(D)'$ whenever $\delta\meg \delta_1$ and $R\meg R_1$.
\end{enumerate}
\end{prop}

In particular, when $A^{p',q'}_{\vect b+\vect d-\vect s-\vect{s'}}(D)$ can be identified with $A^{p,q}_{\vect s,0}(D)'$, this result gives atomic decompositions for $A^{p',q'}_{\vect b+\vect d-\vect s-\vect{s'}}(D)$. In addition, this result shows that the image of $\iota$ is weakly dense in $A^{p,q}_{\vect s,0}(D)'$, and also dense when $p',q'<\infty$.

The proof is based on~\cite[Theorem 5.7]{Bekolleetal}, which deals with the case in which $p=q$, $\vect s\in \R\vect 1_r$, and $D$ is an irreducible symmetric tube domain.

\begin{proof}
Define
\[
S\colon A^{p,q}_{\vect s,0}(D)\ni f \mapsto (\Delta_\Omega^{\vect s-(\vect b+\vect d)/p}(h_k))\in \ell^{p,q}_0(J,K),
\]
so that Theorem~\ref{teo:3} shows that $S$  is continuous and also an injective strict morphism if $\delta$ is small enough. 
Hence, $\trasp S \colon \ell^{p',q'}(J,K)\to A^{p,q}_{\vect s,0}(D)'$ is continuous and is a surjective strict morphism if $\delta$ is small enough.
Let us prove that 
\[
\trasp S(\lambda)=\sum_{j,k}\lambda_{j,k} \,\iota\big(B^{\vect{s'}}_{(\zeta_{j,k},z_{j,k})}\big) \Delta_\Omega^{\vect{s'}-(\vect b+\vect d)/p}(h_k)
\]
for every $\lambda\in \ell^{p',q'}(J,K)$. 
Indeed, assume first that $\lambda$ has finite support and take $f\in A^{p,q}_{\vect s,0}(D)$. Then, Proposition~\ref{prop:9} shows that
\[
\begin{split}
\langle f, \trasp S(\lambda)\rangle&= \sum_{j,k} \lambda_{j,k} f(\zeta_{j,k},z_{j,k}) \Delta_\Omega^{\vect s-(\vect b+\vect d)/p}(h_k)\\
&=\sum_{j,k} \lambda_{j,k} \Delta_\Omega^{\vect s-(\vect b+\vect d)/p}(h_k)\int_D f(\zeta',z') B_{(\zeta',z')}^{\vect{s'}}(\zeta_{j,k}, z_{j,k})\times \\
&\qquad \times\Delta^{-\vect{s'}}_\Omega(\Im z'-\Phi(\zeta'))\,\dd \nu_D(\zeta',z')\\
&=\bigg\langle f,\sum_{j,k}   \lambda_{j,k} \,\iota\big(B_{(\zeta_{j,k}, z_{j,k})}^{\vect{s'}}\big)\Delta_\Omega^{\vect s-(\vect b+\vect d)/p}(h_k)\bigg\rangle,
\end{split}
\]
whence the assertion in this case. 
Thus, 
\[
\trasp S(\lambda)=\sum_{j,k} \lambda_{j,k} \iota\big(B^{\vect{s'}}_{(\zeta_{j,k},z_{j,k})}\big) \Delta_\Omega^{\vect{s'}-(\vect b+\vect d)/p}(h_k) 
\]
with convergence in the strong topology of $A^{p,q}_{\vect s,0}(D)'$ when $\lambda\in \ell^{p',q'}_0(J,K)$, and with convergence in the weak topology $\sigma(A^{p,q}_{\vect s,0}(D)',A^{p,q}_{\vect s,0}(D))$  when $\lambda \in \ell^{p',q'}(J,K)$.
The assertion follows.
\end{proof}

\begin{cor}
Take $p,q\in ]0,\infty]$, and  $\vect{s},\vect{s'}\in \R^r$ such that property $\atomics^{p,q}_{\vect s,\vect{s'},0}$ holds.
Then, property $\atomics^{p',q'}_{(\vect b+\vect d)/\min(1,p)-\vect s-\vect{s'},\vect{s'}}$ holds.
\end{cor}

\begin{proof}
This follows from Lemma~\ref{lem:46} and Propositions~\ref{prop:29} and~\ref{prop:13}.
\end{proof}

\section{Notes and Further Results}

\paragraph{3.6.1} \label{sec:3:6:1}  As noted in Proposition~\ref{prop:6}, the function
\[
\widetilde K_{\vect s}\colon D\times D \ni ((\zeta,z),(\zeta',z'))\mapsto B^{\vect b+\vect d-2\vect s}_{(\zeta',z')}(\zeta,z)
\]
is (up to a constant) the reproducing kernel of $A^{2,2}_{\vect s}(D)$ for every $\vect s\in \frac 1 4 \vect m\in (\R_+^*)^r$. As noted in~\cite{VergneRossi} in the case in which $D$ is an irreducible symmetric Siegel domain of type II, there are other values of $\vect s \in \R\vect 1_r$ for which $\widetilde K_{\vect s}$ is the reproducing kernel of a hilbertian space $\Ac_{\vect s}$ of holomorphic functions. 
More precisely, $\Ac_{\vect s}$ is the Hausdorff completion of the subspace of $\Hol(D)$ generated by the $B^{\vect b+\vect d-2\vect s}_{(\zeta,z)}$, for $(\zeta,z)\in D$, endowed with the semi-norm
\[
\norm*{\sum_{j=1}^k c_j B^{\vect b+\vect d-2 \vect s}_{(\zeta_j,z_j)}}^{ 2}_{ \Ac_{\vect s}}= \sum_{j,j'=1}^k c_j\overline{c_{j'}}B^{\vect b+\vect d-2 \vect s}_{(\zeta_{j'},z_{j'})}(\zeta_j,z_j)
\]
for every $k\in \N$, for every $c_j\in \C$ and for every  $(\zeta_j,z_j)\in D$, $j=1,\dots,k$. 

Thanks to~\cite[Proposition 3.1.5]{VergneRossi} this happens if and only if $\vect s\in \frac{\vect b+\vect d+\Gc(\Omega')}{2}$, where $\Gc(\Omega')$ denotes the Gindikin--Wallach set associated with $\Omega'$ (cf.~2.6.2).
When $\vect s\in \frac{\vect b+\vect d}{2}+\frac 1 4 \vect{m'}+(\R_+^*)^r$, the space $\Ac_{\vect s}$ embeds canonically into $\widehat A^{2,2}_{\vect s,\vect{s'}}$ for every $\vect{s'}\in \N_{\Omega'}$ such that $\vect s+\vect{s'}\in \frac 1 4 \vect m+(\R_+^*)^r$, and is actually equal to the space $\widetilde A^{2,2}_{\vect s}(D)$ to be defined in Section~\ref{sec:6:6} (cf.~Proposition~\ref{prop:39}).

When $\vect s\in \frac{\vect b+\vect d}{2}+\frac 1 2 \big[\Gc(\Omega')\setminus \big(\frac 1 2 \vect{m'}+(\R_+^*)^r\big)\big]$, however, the space $\Ac_{\vect s}$ is somewhat singular and has no relationship with the spaces  $\widehat A^{2,2}_{\vect s,\vect{s'}}$. For example, $\Ac_{(\vect b+\vect d)/2}$ is the space of constant functions on $D$, which has nothing in common with the generalized Dirichlet space $\widehat A^{2,2}_{\vect s,\vect{s'}}$.
 More generally, by means of Proposition~\ref{prop:39} it is readily seen that $\Ac_{\vect s}*I^{-\vect{s'}}_\Omega=\Set{0}$ for every $\vect{s'}\in \N_{\Omega'}$ such that $\vect s+\vect{s'}\in \frac 1 4 \vect m+(\R_+^*)^r$.
These singularities have been addressed when $D$ is an irreducible symmetric  domain and $\vect s\in \R \vect 1_r$ in~\cite{ArazyFisher,Arazy,ArazyUpmeier}, where some new `invariant' spaces are defined, \emph{in some situations}, as suitable (possibly higher order) `residues'. It would be interesting to compare such spaces with the corresponding spaces $\widehat A^{2,2}_{\vect s,\vect{s'}}(D)$.
We refer the reader to~\cite{Arazy,ArazyUpmeier} for further details.

We wish to return to this kind of questions in a future work.

\paragraph{3.6.2}  When $\vect s\in \frac 1 2 \Gc(\Omega)$, the space \CR$\Ac_{\vect s}=\widetilde A^{2,2}_{\vect s}(D)$ \CB (cf.~3.6.1) can be  characterized as the space of $f\in \Hol(D)$ such that
\[
\CR\sup_{h\in \Omega} \int_{\overline \Omega} \norm{f_{h+h'}}_{L^2(\Nc)}^2\,\dd I_{\Omega}^{2\vect s}(h')<\infty.\CB
\]
In addition, there is a constant $c>0$ such that
\[
\norm{f}_{\Ac_{\vect s}}=\CR c \sup_{h\in \Omega} \left( \int_{\overline \Omega} \norm{f_{h+h'}}_{L^2(\Nc)}^2\,\dd I_{\Omega}^{2\vect s}(h')\right) ^{1/2}\CB
\]
for every $f\in \Ac_{\vect s}$. See~\cite{VergneRossi} for the case in which $D$ is an irreducible symmetric Siegel domain of type II. The details are left to the reader.

Thus, the spaces $\Ac_{\vect s}$ for $\vect s\in \frac 1 2 \Gc(\Omega)$ can be considered as `generalized Hardy spaces.' One may therefore consider the space $H^{p,q}_{\vect s}(D)$ of $f\in \Hol(D)$ such that
\[
\sup_{h\in \Omega} \norm*{ h'\mapsto \norm{f_{h+h'}}_{L^p(\Nc)}}_{L^q(I_{\Omega}^{\vect s})}<\infty
\]
for $\vect s\in \Gc(\Omega)$ and $p,q\in ]0,\infty]$, endowed with the corresponding norm.\footnote{Actually, with a finer analysis of the orbit of $T_+$ in $\overline \Omega$ on which  $I^{\vect s}_\Omega$ is concentrated, one may define $H^{p,q}_{q\vect s}(D)$ in complete analogy with $A^{p,q}_{\vect s}(D)$, thus getting more interesting spaces when $q=\infty$.} When $\vect s\in \frac 1 2 \vect m+(\R_+^*)^r$, one simply recovers the space $A^{p,q}_{\vect s/q}(D) $.
Such spaces have been extensively studied in~\cite{Garrigos} when $D$ is an irreducible symmetric tube domain and $p=q$. 
We do not know if the resulting spaces fit into one of the scales of spaces $A^{p_1,q_1}_{\vect {s_1}}(D)$, $\widetilde A^{p_1,q_1}_{\vect {s_1}}(D)$, or $\widehat A^{p_1,q_1}_{\vect {s_1},\vect{s_2}}(D)$ (cf.~Chapter~\ref{sec:7}) when $(p,q)\neq (2,2)$ or $\vect s\not \in \Set{\vect 0} \cup\big(\frac{1}{2 p}\vect m+(\R_+^*)^r\big)$.

\paragraph{3.6.3} With the notation of~3.6.1, we remark that, when $r=1$ and $n>0$, so that $D$ is biholomorphic to the unit ball of $\C^{n+1}$, the space $\Ac_{\vect b/2}$ may be canonically identified with the Drury--Arveson space (cf.~\cite{ArcozziMonguzziPelosoSalvatori}).

\chapter{Besov Spaces of Analytic Type}\label{sec:5}

In this chapter, we introduce the Besov spaces $B^{\vect s}_{p,q}(\Nc,\Omega)$. In comparison with our main reference~\cite{BekolleBonamiGarrigosRicci}, two new difficulties arise: on the one hand, the group $\Nc$ is not necessarily commutative, so that the associated Fourier transform is far less manageable. On the other hand, we shall consider the full range of exponents $p,q\in ]0,\infty]$ instead of dealing only with $p,q\in [1,\infty]$.

In order to deal with the general case $p,q\in ]0,\infty]$, the classical techniques presented, e.g., in~\cite{Triebel} can be effectively extended to our situation. On the other hand, dealing with the case in which $\Nc$ is not necessarily commutative provides new kinds of issues, which are basically related to the fact that the Fourier transform of the Schwartz space $\Sc(\Nc)$ is not easy to manage. 
In order to deal with this inconvenience, we shall introduce three spaces: $\Sc_\Omega(\Nc)$, $\Sc_{\Omega,L}(\Nc)$, and $\widetilde \Sc_\Omega(\Nc)$.
Let us briefly explain the role played by each one of them.
First of all, because of the results of Section~\ref{sec:1:4}, it is convenient to consider only functions $f\in \Sc(\Nc)$ such that 
\[
\pi_\lambda(f)= \chi_{\Omega'}(\lambda) \pi_\lambda(f) P_{\lambda,0}
\]
for every $\lambda\in F'\setminus W$. Nonetheless, the space of such $f$ is too big for our purposes. In addition to that, the description of the image of $\Sc(\Nc)$ under the Fourier transform provided in~\cite{Geller} at least when $\Nc$ is a Heisenberg group, is not easy to work with.
For this reason, it is convenient to consider the space 
\[
\Sc_\Omega(\Nc)\coloneqq \big\{f\in \Sc(\Nc)\colon \exists \varphi\in C^\infty_c(\Omega')\; \forall\lambda\in F'\setminus W\quad \pi_\lambda(f)= \varphi(\lambda) P_{\lambda,0} \big\}.
\]
Since the Fourier transform of the elements of $\Sc_\Omega(\Nc)$ is essentially scalar-valued, this choice implies, in particular, that $\Sc_\Omega(\Nc)$ is commutative under convolution, a fact that will prove very useful in the computations. In addition, by means of the classical Paley--Wiener theorems, it is not hard to prove that the Fourier transform maps $\Sc_\Omega(\Nc)$ isomorphically onto $C^\infty_c(\Omega')$ (times the field $\lambda \mapsto P_{\lambda,0}$). 

Even though the space $\Sc_\Omega(\Nc)$ has several important
properties and is relatively easy to work with, it is far too small to
be dense in any of the  Besov  spaces  $B^{\vect s}_{p,q}(\Nc,\Omega)$ to be defined. Moreover, it is not stable under left translations. 
Nonetheless, the `left-invariant completion' $\Sc_{\Omega,L}(\Nc)$ of
$\Sc_\Omega(\Nc)$ has the properties that $\Sc_\Omega(\Nc)$ is lacking
in order to complete the study of the spaces $B^{\vect s}_{p,q}(\Nc,\Omega)$. On the one hand, the spaces $B^{\vect s}_{p,q}(\Nc,\Omega)$ embed canonically into $\Sc_{\Omega,L}'(\Nc)$. On the other hand, $\Sc_{\Omega,L}(\Nc)$ embeds canonically as a dense subspace of $B^{\vect s}_{p,q}(\Nc,\Omega)$ (at least when $p,q<\infty$).

One last issue occurs: it is unclear if any of the spaces $\Sc_\Omega(\Nc)$ and $\Sc_{\Omega,L}(\Nc)$ is stable under pointwise multiplication. Even though this fact is of minor importance when dealing with the case $p,q\in [1,\infty]$, in order to extend to the general case the techniques presented in~\cite{Triebel}, we need to be able to multiply functions with a control on their Fourier transforms. For this reason, we shall introduce an auxiliary space $\widetilde \Sc_\Omega(\Nc)$, which is defined with the aid of the Euclidean Fourier transform on the space $F$. 

With these three spaces in hand, we can develop the basic theory of the Besov spaces $B^{\vect s}_{p,q}(\Nc,\Omega)$ following the classical case with only minor modifications.

Concerning the results which follow, we only remark that Theorem~\ref{teo:10}, which characterizes the dual of $B^{\vect s}_{p,q}(\Nc,\Omega)$, only covers the case $p\Meg 1$. The case $p<1$ will be established later on (Corollary~\ref{cor:30}) with the aid of the theory of weighted Bergman spaces.

\section{Spaces of Test Functions}\label{sec:5:1}

This section deals with the study of several subspaces of $\Sc(\Nc)$ which are necessary to deal with the Besov spaces to be defined in Section~\ref{sec:5:2}.
Since these spaces are essentially instrumental, we shall content ourselves with a few basic results.

\begin{deff}\label{28}
Define
\[
\Sc_\Omega(\Nc)\coloneqq\big\{f\in \Sc(\Nc)\colon \exists \varphi\in C^\infty_c(\Omega')\; \forall\lambda\in F'\setminus W\quad \pi_\lambda(f)= \varphi(\lambda) P_{\lambda,0} \big\}.
\]
In addition, for every compact subset $K$ of $\Omega'$, define $\Sc_\Omega(\Nc,K)$ as the set of $f\in \Sc_\Omega(\Nc)$ such that $\pi_\lambda(f)=0$ for every $\lambda\in \Omega'\setminus K$, endowed with the topology induced by $\Sc(\Nc)$.
We endow $\Sc_\Omega(\Nc)$ with the finest locally convex topology for which the inclusions $\Sc_\Omega(\Nc,K)\to \Sc_\Omega(\Nc)$ are continuous.
\end{deff}

 Thus, $\Sc_\Omega(\Nc)$ is a complete locally convex space and embeds
continuously into $\Sc(\Nc)$. In addition, a subset $\Kc$   of
$\Sc_\Omega(\Nc)$ is bounded if and only if $\Kc$ is contained and bounded in $ \Sc_\Omega(\Nc,K) $ for some compact subset
$K$ of $\Omega'$.

The following result collects the main properties of the space $\Sc_\Omega(\Nc)$. 

\begin{prop}\label{prop:40}
The following hold:
\begin{enumerate}
\item[\em(1)] the mapping $\Fc_\Nc\colon \varphi\mapsto [\lambda \mapsto \tr(\pi_\lambda(\varphi) ) ] $ induces an isomorphism of $\Sc_\Omega(\Nc)$ onto $C^\infty_c(\Omega')$;

\item[\em(2)] for every $\psi\in C^\infty_c(\Omega')$ and for every $(\zeta,x)\in \Nc$,
\[
(\Fc_\Nc^{-1} \psi)(\zeta,x)=\frac{2^{n-m} \abs{\Pfaff(e_{\Omega'})}}{\pi^{n+m}}\int_{\Omega'} \psi(\lambda)\Delta^{-\vect b}_{\Omega'}(\lambda) \ee^{i\langle\lambda,x\rangle-\langle \lambda, \Phi(\zeta)\rangle}\,\dd \lambda;
\]

\item[\em(3)] for every $\varphi_1,\varphi_2\in \Sc_\Omega(\Nc)$, $\varphi_1*\varphi_2\in \Sc_\Omega(\Nc)$ and 
\[
\Fc_\Nc(\varphi_1*\varphi_2)=(\Fc_\Nc \varphi_1)(\Fc_\Nc \varphi_2);
\]

\item[\em(4)] if $t\in T_+$, $g\in GL(E)$, $\varphi \in \Sc_\Omega(\Nc)$, and $t\cdot
  \Phi=\Phi\circ (g\times g)$,  then  $(g\times t)_*\varphi\in \Sc_\Omega(\Nc)$ and
\[
\Fc_\Nc((g\times t)_*\varphi )= (\Fc_\Nc \varphi)(\,\cdot\,t).
\]
\end{enumerate} 
\end{prop}

In particular,  convolution is commutative on $\Sc_\Omega(\Nc)$.
Observe that the mapping $\Fc_\Nc$ is essentially the (non-commutative) Fourier transform on $\Nc$, since
\[
\pi_\lambda(\varphi)= \Fc_\Nc(\varphi)(\lambda) P_{\lambda,0}
\]
for every $\lambda\in F'\setminus W$ and for every $\varphi \in \Sc_\Omega(\Nc)$.
The characterization of the image of $\Sc_\Omega(\Nc)$ under $\Fc_\Nc$ will be of particular importance in the study of Besov spaces.

\smallskip

In the statement, we wrote $(g\times t)_*\varphi$ instead of $\Delta^{\vect b+\vect d}_\Omega(t)\, \varphi\circ (g\times t)^{-1}$.\label{45} 

\begin{proof}
{ (1)--(2)} Take $\varphi \in \Sc_\Omega(\Nc)$ and observe that, by Proposition~\ref{prop:36},
\[
(\Fc_\Nc\varphi)(\lambda)=\tr(\pi_\lambda(\varphi))= \langle \pi_\lambda(\varphi) e_{\lambda,0}\vert e_{\lambda,0}\rangle=\int_\Nc \varphi(\zeta,x) \ee^{-i \langle \lambda,x\rangle-\langle \lambda, \Phi(\zeta)\rangle }\,\dd (\zeta,x)
\]
for every $\lambda\in\Omega'$. It then follows that $\Fc_\Nc\varphi\in C^\infty_c(\Omega')$.

Conversely, take $\psi\in C^\infty_c(\Omega')$, and define
\[
(\Gc\psi)(\zeta,x)\coloneqq \frac{2^{n-m} \abs{\Pfaff(e_{\Omega'})}}{\pi^{n+m}}\int_{\Omega'} \psi(\lambda)\Delta^{-\vect b}_{\Omega'}(\lambda) \ee^{i\langle\lambda,x\rangle-\langle \lambda, \Phi(\zeta)\rangle}\,\dd \lambda,
\]
for every $(\zeta,x)\in \Nc$. Then, Proposition~\ref{prop:36} and Corollary~\ref{cor:7} imply that $\Gc\psi\in L^2(\Nc)$ and that $\pi_\lambda(\Gc\psi)=\psi(\lambda) P_{\lambda,0}$ for almost every $\lambda\in F'\setminus W$. 
Now, take $\alpha_1,\alpha_2\in \N$ and $\alpha_3\in \N^m$, and observe that Faà di Bruno's formula and some integrations by parts show that
\[
\begin{split}
&\partial_E^{\alpha_1}\partial_F^{\alpha_2}(\Gc \psi)(\zeta,x)=\frac{2^{n-m} \abs{\Pfaff(e_{\Omega'})}}{\pi^{n+m}}\int_{\Omega'} \psi(\lambda) \ee^{i\langle\lambda,x\rangle-\langle \lambda, \Phi(\zeta)\rangle} \Theta_{\alpha_2}(\lambda) \Delta^{-\vect b}_{\Omega'}(\lambda) \,\dd \lambda\\
&\qquad \qquad=\frac{(-1)^{\alpha_3}2^{n-m} \abs{\Pfaff(e_{\Omega'})}}{\pi^{n+m} (i x-\Phi(\zeta)  )^{\alpha_3}}\int_{\Omega'}  \ee^{i\langle\lambda,x\rangle-\langle \lambda, \Phi(\zeta)\rangle}\partial^{\alpha_3} (\psi\, \Theta_{\alpha_2}\Delta_{\Omega'}^{-\vect b}  )(\lambda)  \,\dd \lambda,
\end{split}
\] 
where
\[
\Theta_{\alpha_2}(\lambda)\coloneqq\sum_{\gamma_1+2\gamma_2=\alpha_1} \frac{\alpha_1!}{\gamma_1!\gamma_2!} (i\lambda')^{\alpha_2}\cdot\left(-2\langle \lambda' ,\Re\Phi(\zeta, \,\cdot\,)\rangle  \right)^{\gamma_1} \cdot\left( -\langle \lambda',\Re \Phi(\,\cdot\,,\,\cdot\,)\rangle \right) ^{\gamma_2}
\]
for every $\lambda\in \Omega'$.\footnote{Here, $\partial_E^{\alpha_1}$ denotes the partial derivative of order $\alpha_1$ with respect to the subspace $E$. In other words, if $f\in C^{\alpha_1}(\Nc)$, then $\partial_E^{\alpha_1} f(\zeta,x)=(f(\,\cdot\,,x))^{(\alpha_1)}(\zeta)$ for every $(\zeta,x)\in \Nc$. Similar remarks apply to $\partial_F^{\alpha_2}$. Finally, the $\cdot$ appearing in the definition of $\Theta_{\alpha_2}(\lambda)$ stand for symmetrized tensor products. In other words, if $T$ and $S$ are symmetric $k$- and $h$-multilinear mappings on $E\times F$, then $(T\cdot S)(v_1,\dots, v_{h+k})= \frac{1}{(k+h)!}\sum_{\sigma\in \Sf_{k+h}} T(v_{\sigma(1)},\dots, v_{\sigma(k)}) S(v_{\sigma(k+1)},\dots, v_{\sigma(k+h)}) $ for every $v_1,\dots, v_{h+k}\in E\times F$, where $\Sf_{k+h}$ denotes the set of permutations on $\Set{1,\dots, k+h}$. }
Now, observe that the mapping $\zeta\mapsto \Phi(\zeta)$ is absolutely homogeneous of degree $2$ and vanishes only at $0$. 
Therefore, for every $N\in \N$ there is a constant $C'_N>0$ such that
\[
\abs{\partial^{\alpha}(\Gc \psi)(\zeta,x)}\meg \frac{C_N}{(1+\abs{\zeta}^2+\abs{x})^N}
\]
for every $(\zeta,x)\in \Nc$ and for every $\alpha\in \N^{2 n+m} $ with length at most $N$ (identifying $E$ with $\R^{2 n}$). In addition, since the mapping $\lambda\mapsto J'_\lambda=\abs{J_\lambda}^{-1} J_\lambda$ is (semialgebraic and) clearly analytic on $F'\setminus W$, it is easily seen that $\pi_\lambda(\Gc\psi)=\psi(\lambda)P_{\lambda,0}$ for \emph{every} $\lambda\in F'\setminus W$.
Thus, $\Gc$ maps $C^\infty_c(F')$ into $\Sc_\Omega(\Nc)$, and clearly $\Fc_\Nc$ and $\Gc$ are inverses of one another. 
 Finally, it is readily seen that $\Fc_\Nc$ and $\Gc$ induce homeomorphisms between $\Sc_\Omega(\Nc,K)$ and $\Set{\psi\in C^\infty_c(\Omega')\colon \Supp{\psi}\subseteq K}$ for every compact subset $K$ of $\Omega'$, so that $\Fc_\Nc$ and $\Gc$ are homeomorphisms between $\Sc_\Omega(\Nc)$ and $C^\infty_c(\Omega')$. 

{(3)} Take $\varphi_1,\varphi_2\in \Sc_\Nc(\Omega)$, and observe that $\pi_\lambda(\varphi_1*\varphi_2)=\pi_\lambda(\varphi_1)\pi_\lambda(\varphi_2)$ for every $\lambda\not \in W$. Since $\pi_\lambda(\varphi_j)= \Fc_\Nc(\varphi_j) P_{\lambda,0}$ for every $\lambda\not \in W$ and for every $j=1,2$, this implies that $\varphi_1*\varphi_2\in \Sc_\Omega(\Nc)$ and that
\[
\Fc_\Nc(\varphi_1*\varphi_2)=(\Fc_\Nc \varphi_1)(\Fc_\Nc \varphi_2),
\]
as we wished to show. 

{(4)} Finally,  define $\psi\coloneqq \Fc_\Nc \varphi$, and observe that~{(2)} and~{(3)} imply that $\psi(\,\cdot\,t)\in C^\infty_c(\Omega')$ and that 
\[
\Fc_\Nc^{-1} (\psi(\,\cdot\,t))=\Delta_\Omega^{-\vect b-\vect d}(t^{-1}) (\Fc_\Nc^{-1} \psi)(g^{-1}\times t^{-1}\,\cdot\,)=(g\times t)_* \varphi,
\]
whence the result.
\end{proof}

\begin{deff}
We define $\Fc_\Nc\colon \Sc_\Omega(\Nc)\to C^\infty_c(F')$ as in Proposition~\ref{prop:40}. \label{30}
\end{deff}

As we mentioned above, the space $\Sc_\Omega(\Nc)$ is too small for some of our purposes. In order to get invariance under left translations we introduce $\Sc_{\Omega,L}(\Nc)$. 
The auxiliary space $\widetilde \Sc_\Omega(\Nc)$, in turn, is introduced to get control over pointwise multiplication.

\begin{deff}\label{31}
Define $\widetilde \Sc_\Omega(\Nc)$\label{32} as the space of $\varphi\in \Sc(\Nc)$ such that $\Fc_F(\varphi(\zeta,\,\cdot\,))$ is supported in $ \overline{\Omega'}$ for every $\zeta\in E$ and such that $\pi_\lambda(\varphi)= \pi_\lambda(\varphi)P_{\lambda,0} $ for every $\lambda\in \Omega'$. Endow $\widetilde \Sc_\Omega(\Nc)$ with the topology induced by $\Sc(\Nc)$.

In addition, for every compact subset $K$ of $\Omega'$, define $\Sc_{\Omega,L}(\Nc,K)$ as the space of $\varphi \in \widetilde \Sc_\Omega(\Nc)$ such that $\pi_\lambda(\varphi)=0$ for every $\lambda\in F'\setminus (W\cup K)$, endowed with the topology induced by $\Sc(\Nc)$. 
Define $\Sc_{\Omega,L}(\Nc)$ as the inductive limit of the spaces $\Sc_{\Omega,L}(\Nc,K)$, endowed with the corresponding locally convex topology.
\end{deff}

The following result collects some of the most important results concerning the spaces $\Sc_{\Omega,L}(\Nc)$ and $\widetilde \Sc_\Omega(\Nc)$. In particular, it shows that $\widetilde \Sc_\Omega(\Nc)$ is an algebra under pointwise multiplication. 

Recall that a Montel space is a Hausdorff barrelled space whose bounded subsets are relatively compact.

\begin{prop}\label{prop:41}
The following hold:
\begin{enumerate}
\item[\em(1)] $\widetilde \Sc_\Omega(\Nc)$ is a  Fréchet Montel  space;

\item[\em(2)] $\Sc_{\Omega,L}(\Nc)$ is a Hausdorff, complete,
  bornological Montel space;
  
  \item[\em(3)] $\Sc_{\Omega,L}(\Nc)$ embeds canonically into $\widetilde
  \Sc_\Omega(\Nc)$ and  induces   on $\Sc_{\Omega,L}(\Nc,K)$ its topology for every compact subset $K$ of $\Omega'$;
  
 \item[\em(4)]  a subset of $\Sc_{\Omega,L}(\Nc)$ is bounded if and only if it is contained and bounded in $\Sc_{\Omega,L}(\Nc,K)$ for some compact subset $K$ of $\Omega'$; 

\item[\em(5)] $\pi_\lambda(\varphi)=0$ for every $\varphi\in \widetilde \Sc_\Omega(\Nc)$ and for every $\lambda \in F'\setminus(\Omega'\cup W)$;

\item[\em(6)] $\widetilde \Sc_\Omega(\Nc)$ and $\Sc_{\Omega,L}(\Nc)$ are left-invariant;

\item[\em(7)] $\Sc_\Omega(\Nc)$ embeds canonically into
  $\Sc_{\Omega,L}(\Nc)$ and   the left-invariant subspace of $\Sc(\Nc)$ generated by $\Sc_\Omega(\Nc)$ is dense in $\Sc_{\Omega,L}(\Nc)$;

\item[\em(8)] $\widetilde \Sc_\Omega(\Nc)\widetilde \Sc_\Omega(\Nc)\subseteq \widetilde \Sc_\Omega(\Nc)$.
\end{enumerate}
\end{prop}

Thus, $\Sc_{\Omega,L}(\Nc)$ can be interpreted as the completion of the left-invariant subspace generated by $\Sc_\Omega(\Nc)$, with respect to a suitable topology. This remark should justify the `$L$' appearing the symbol $\Sc_{\Omega,L}(\Nc)$.

\begin{proof}
{(1)} It is clear that $\widetilde \Sc_\Omega(\Nc)$ is a closed subspace of $\Sc(\Nc)$, so that the assertion follows.

{(2)--(4)} Since clearly $\Sc_{\Omega,L}(\Nc,K)$ is closed in $\Sc(\Nc)$ for every compact subset $K$ of $\Omega'$, and since $\Omega'$ is locally compact and $\sigma$-compact,~\cite[Proposition 9 of Chapter II, \S 4, No.\ 6]{BourbakiTVS} implies that $\Sc_{\Omega,L}(\Nc)$ is Hausdorff and complete, and induces on $\Sc_{\Omega,L}(\Nc,K)$ its topology for every compact subset $K$ of $\Omega'$. In addition,~\cite[Example 3 of Chpater III, \S 2 and Corollary 3 of Chapter III, \S 4, No.\ 1]{BourbakiTVS} imply that $\Sc_{\Omega,L}(\Nc)$ is bornological and barrelled. Further,~\cite[Proposition 6 of Chapter III, \S 1, No.\ 4]{BourbakiTVS} implies that a subset of $\Sc_{\Omega,L}(\Nc)$ is bounded if and only if it is contained and bounded in $\Sc_{\Omega,L}(\Nc,K)$ for some compact subset $K$ of $\Omega'$. Since $\Sc(\Nc)$ is a Montel space, it then follows that $\Sc_{\Omega,L}(\Nc)$ is a Montel space. Finally, it is clear that $\Sc_{\Omega,L}(\Nc,K)$ embeds canonically into $\widetilde \Sc_{\Omega}(\Nc)$ for every compact subset  $K$ of $\Omega'$, so that also $\Sc_{\Omega,L}(\Nc)$ embeds canonically into $\widetilde \Sc_\Omega(\Nc)$.

{(5)} Indeed, take $\varphi\in \widetilde \Sc_\Omega(\Nc)$ and  $\lambda \in F'\setminus(\Omega'\cup W)$. Then,
\[
\pi_\lambda(\varphi)=\int_E \Fc_F(\varphi(\zeta,\cdot\,))(\lambda)\pi_\lambda(\zeta,0)\,\dd \zeta=0.
\]

{(6)} It suffices to prove that $\widetilde \Sc_\Omega(\Nc)$ is left-invariant. Indeed, take $\varphi \in \widetilde \Sc_\Omega(\Nc)$ and $(\zeta,x)\in \Nc$. Then 
\[
\pi_\lambda(L_{(\zeta,x)}\varphi)=\pi_\lambda(\zeta,x) \pi_\lambda(\varphi)=\pi_\lambda(\zeta,x) \pi_\lambda(\varphi) P_{\lambda,0}= \pi_\lambda(L_{(\zeta,x)}\varphi) P_{\lambda,0}
\]
for every $\lambda \in \Omega'$. In addition,
\[
\begin{split}
\Fc_F((L_{(\zeta,x)}\varphi )(\zeta',\,\cdot\,))&=\Fc_F( \varphi(\zeta'-\zeta, \,\cdot\,-x-2\Im \Phi(\zeta,\zeta')))\\
&= \ee^{ -i\langle\,\cdot\, ,x+2 \Im \Phi(\zeta,\zeta')\rangle }\Fc_F(\varphi(\zeta'-\zeta, \,\cdot ))\in C^\infty_c(\Omega')
\end{split}
\]
for every $\zeta'\in E$. Therefore, $L_{(\zeta,x)}\varphi \in \widetilde \Sc_\Omega(\Nc)$.

{(7)} Let us first prove that $\Sc_\Omega(\Nc)$ embeds canonically into $\Sc_{\Omega,L}(\Nc,K)$. Indeed, take $\varphi \in \Sc_\Omega(\Nc)$ and observe that clearly $\pi_\lambda(\varphi)=\pi_\lambda(\varphi) P_{\lambda,0}$ for every $\lambda\in \Omega'$. In addition, Proposition~\ref{prop:40} implies that
\[
\varphi(\zeta,x)=\frac{2^{n-m}\abs{\Pfaff(e_{\Omega'})}}{\pi^{n+m}} \Lc(\psi\, \Delta_{\Omega'}^{-\vect b} )(\Phi(\zeta)- i x)
\]
for every $(\zeta,x)\in \Nc$,where $\psi\coloneqq \Fc_\Nc \varphi\in C^\infty_c(\Omega')$. Therefore,
\[
\Fc_F(\varphi(\zeta,\,\cdot\,))=\frac{2^{n}\abs{\Pfaff(e_{\Omega'})}}{\pi^{n}}  \psi\, \Delta_{\Omega'}^{-\vect b}  \ee^{-\langle\,\cdot\, ,\Phi(\zeta)\rangle}\in C^\infty_c(\Omega')
\]
for every $\zeta\in E$. Thus, $\varphi\in \Sc_{\Omega,L}(\Nc)$. Continuity is then easily established.

Conversely, take $\varphi\in \Sc_{\Omega,L}(\Nc)$, and let $K$ be a compact subset of $\Omega'$ such that $\pi_\lambda(\varphi)=0 $ for every $\lambda\in F'\setminus (K\cup W)$. 
Then, Proposition~\ref{prop:40} implies that there is $\tau\in \Sc_\Omega(\Nc)$ such that $\Fc_\Nc \tau$ equals $1$ on $K$, so that $\varphi=\varphi*\tau$. 
In addition, the space $M$ of measures with finite support is dense in the space $\Oc'_{C,L}(\Nc)$ of left convolutors of the space $\Sc(\Nc)$ into itself, so that there is a filter $\Ff$ on $M$ which converges to $\varphi$ in $\Oc'_{C,L}(\Nc)$.
Now,  $M*\tau$ is contained in the left-invariant subspace of $\Sc(\Nc)$ generated by $\Sc_\Omega(\Nc)$, and clearly $\Ff*\tau$ converges to $\varphi=\varphi*\tau$ in $\Sc_{\Omega,L}(\Nc,\Supp{\Fc_\Nc\tau})$, hence in $\Sc_{\Omega,L}(\Nc)$. The assertion follows.

{(8)} Take $\varphi_1,\varphi_2\in \widetilde \Sc_\Omega(\Nc)$. Then, 
\[
\Fc_F( (\varphi_1 \varphi_2)(\zeta,\,\cdot\,))=(2 \pi)^{-m} \Fc_F(\varphi_1(\zeta,\,\cdot\,))*\Fc_F(\varphi_2(\zeta,\,\cdot\,))
\]
is supported in $\overline{\Omega'}$ for every $\zeta\in E$, since $\overline {\Omega'}$ is a closed convex cone. 
Thus, it remains to prove that $\pi_\lambda(\varphi_1\varphi_2)=\pi_\lambda(\varphi_1\varphi_2) P_{\lambda,0}$ for every $\lambda\in \Omega'$.
To this aim, we define
\[
L^2_\Omega(\Nc)\coloneqq \Set{f\in L^2(\Nc)\colon \pi_\lambda(f)= \chi_{\Omega'}(\lambda)\pi_\lambda(f) P_{\lambda,0} \text{ for almost every $\lambda\in F'\setminus W$} },
\]
and an operator $\Ec\colon L^2_\Omega(\Nc)\to A^{2,\infty}_{\nu_\Omega}(D)$ (with the notation of Section~\ref{sec:1:3}) such that
\[
\pi_\lambda((\Ec f)_h)= \ee^{-\langle \lambda,h\rangle}\pi_\lambda(f) 
\]
for every $f\in L^2_\Omega(\Nc)$, for every $h\in \Omega$, and for almost every $\lambda\in F'\setminus W$. 
Then, by means of Corollary~\ref{cor:1} and Proposition~\ref{prop:4} we see that $\Ec$ is well defined and maps $L^2_\Omega(\Nc)$ isometrically onto $A^{2,\infty}_{\nu_\Omega}(D)$. 
In addition, $\Ec(f)_h\to  f $ in $L^2(\Nc)$ as $h\to 0$ for every $f\in L^2_\Omega(\Nc)$.
Now, it is clear that $\widetilde \Sc_\Omega(\Nc)$ embeds continuously into $L^2_\Omega(\Nc)$, so that $\varphi_1 \varphi_2$ is the limit of $(\Ec \varphi_1)_h (\Ec \varphi_2)_h$ in $L^1(\Nc)$. In particular,
\[
\begin{split}
\pi_\lambda(\varphi_1 \varphi_2)&=\lim_{h\to 0} \pi_\lambda( (\Ec \varphi_1)_h (\Ec \varphi_2)_h )\\
&=\lim_{h\to 0} \pi_\lambda( (\Ec \varphi_1)_h (\Ec \varphi_2)_h ) P_{\lambda,0}\\
&=\pi_\lambda(\varphi_1 \varphi_2) P_{\lambda,0}
\end{split}
\]
for every $\lambda\in \Omega'$,where the second equality follows from Proposition~\ref{prop:2}. The assertion follows.
\end{proof}

From the following result until Corollary~\ref{cor:24}, we shall extend to this setting some of the tools needed to deal with the Besov spaces $B^{\vect s}_{p,q}$ for the full range of exponents $p,q\in ]0,\infty]$. Corollaries~\ref{cor:23} and~\ref{cor:24}, in particular, are trivial consequences of Young's inequality when $p\Meg 1$.

\begin{cor}\label{cor:21}
The following hold:
\begin{enumerate}
\item[\em(1)] the trilinear mapping
\[
\Sc'(\Nc)\times \widetilde\Sc_\Omega(\Nc)\times \widetilde\Sc_\Omega(\Nc)\ni (T,\varphi, \eta)\mapsto (T*\varphi)\eta\in \widetilde\Sc_{\Omega}(\Nc)
\]
is well-defined and  separately continuous;

\item[\em(2)] $\pi_\lambda((T*\varphi)\eta)=0$ for every $\lambda \in \Omega'\setminus ( \Supp{\Fc_\Nc \varphi}+\Supp{\Fc_\Nc \eta}) $, for every $T\in \Sc'(\Nc)$ and for every $\varphi,\eta\in \Sc_{\Omega}(\Nc)$. 
\end{enumerate}
\end{cor}

\begin{proof}
{(1)} The considered trilinear mapping is clearly separately continuous with values in $\Sc(\Nc)$. 
Now, observe that the set of $\delta_{(\zeta,x)}$, for $(\zeta,x)\in \Nc$, is total in $\Sc'(\Nc)$. To conclude, it will then suffice to show that $(\delta_{(\zeta,x)}*\varphi)\eta\in \widetilde \Sc_\Omega(\Nc)$ for every $\varphi,\eta\in \widetilde \Sc_\Omega(\Nc)$. However, since $\delta_{(\zeta,x)}*\varphi=L_{(\zeta,x)}\varphi$, the assertion follows from Proposition~\ref{prop:41}.

{(2)} As in~{(1)}, we may reduce to the case in which $T=\delta_{(\zeta,x)} $. In this case, then, 
\[
\begin{split}
\Fc_F((L_{(\zeta,x)}\varphi)(\zeta',\,\cdot\,))&=\ee^{ -i\langle\,\cdot\, ,x+2 \Im \Phi(\zeta,\zeta')\rangle }\Fc_F( \varphi(\zeta'-\zeta,\,\cdot\,) )\\
&=\frac{2^n\abs{\Pfaff(e_{\Omega'})}}{\pi^n} \ee^{ -i\langle\,\cdot\, ,x+2 \Im \Phi(\zeta,\zeta')\rangle } (\Fc_\Nc \varphi)\Delta_{\Omega'}^{-\vect b} \ee^{-\langle \,\cdot\,,\Phi(\zeta'-\zeta)\rangle}
\end{split}
\]
and
\[
\Fc_F(\eta(\zeta',\,\cdot\,))= \frac{2^n\abs{\Pfaff(e_{\Omega'})}}{\pi^n} (\Fc_\Nc \eta)\Delta_{\Omega'}^{-\vect b} \ee^{-\langle \,\cdot\,,\Phi(\zeta')\rangle}
\]
for every $\zeta'\in E$ (cf.~the computations of the proof of Proposition~\ref{prop:41}). Since
\[
\Fc_F(((L_{(\zeta,x)}\varphi) \eta)(\zeta',\,\cdot\,)  )=(2 \pi)^{-m} \Fc_F((L_{(\zeta,x)}\varphi)(\zeta',\,\cdot\,))*\Fc_F(\eta(\zeta',\,\cdot\,))
\]
and since
\[
\pi_\lambda((L_{(\zeta,x)}\varphi) \eta)=\int_E \Fc_F(((L_{(\zeta,x)}\varphi) \eta)(\zeta',\,\cdot\,))(\lambda) \pi_\lambda(\zeta',0)\,\dd \zeta'
\]
for every $\lambda\in F'\setminus W$, the  assertion follows. 
\end{proof}

\begin{cor}\label{cor:22}
Take $\varphi \in \Sc_\Omega(\Nc)$. Then, for every $p,q\in ]0,\infty]$ such that $p\meg q$  there is a constant $C>0$ such that
\[
\norm{T}_{L^q(\Nc)}\meg C \norm{T}_{L^p(\Nc)}
\]
for every $T\in \Sc'(\Nc)$ such that $T=T*\varphi$.
\end{cor}

The proof is based on~\cite[Theorem 1.4.1]{Triebel}, which deals with the abelian case.

\begin{proof}
\textsc{Step I.} Assume first that $T\in \widetilde \Sc_\Omega(\Nc)$. It will suffice to prove the assertion for $q=\infty$.
Then,
\[
\begin{split}
\abs{T(\zeta,x)}&\meg \int_{\Nc} \abs{T(\zeta',x') \varphi((\zeta',x')^{-1}(\zeta,x)  )  }\,\dd (\zeta',x')\\
&\meg \norm{\varphi}_{L^{p'}(\Nc)} \norm{T}_{L^\infty(\Nc)}^{1-\min(1,p)}\norm{T}_{L^p(\Nc)}^{\min(1,p)}
\end{split}
\]
for every $(\zeta,x)\in \Nc$,
so that
\[
\norm{T}_{L^\infty(\Nc)}\meg  \norm{\varphi}_{L^{p'}(\Nc)}^{1/\min(1,p)} \norm{T}_{L^p(\Nc)}
\]
since $\norm{T}_{L^\infty(\Nc)}$ is finite.

\textsc{Step II.} Take $\eta\in \Sc_\Omega(\Nc)$ such that  $\eta(e)=1$, that is,  such that 
\[
\int_{\Omega'}(\Fc_\Nc \varphi)(\lambda) \Delta^{-\vect b}_{\Omega'}(\lambda)\,\dd \lambda= \frac{\pi^{n+m}}{2^{n-m} \abs{\Pfaff(e_{\Omega'})}}
\]
(cf.~Proposition~\ref{prop:40}). Define 
\[
\eta_\rho\colon (\zeta,x)\mapsto \eta(\rho \zeta,\rho^2 x)
\]
for every $\rho>0$.
Let $K$ be a compact convex neighbourhood of $0$ in $\overline{\Omega'}$ which contains $\Supp{\Fc_\Nc\eta}$, so that $K$ contains the support of 
\[
\Fc_\Nc(\eta_\rho)=\rho^{-2 n-2m} (\Fc_\Nc \eta)(\rho^{-2}\,\cdot\,)
\]
for every $\rho\in ]0,1]$.
Observe that there is $\varphi'\in \Sc_\Omega(\Nc)$ such that $\Fc_\Nc\varphi'$ equals $1$ on the compact subset $K+\Supp{\Fc_\Nc\varphi}$ of $\Omega'$. Therefore, Corollary~\ref{cor:21} shows that $T\eta_\rho\in \widetilde \Sc_\Omega(\Nc)$ and that $T \eta_\rho=(T \eta_\rho)*\varphi'$ for every $\rho\in ]0,1]$. Thus, applying~{step I} with $\varphi'$ in place of $\varphi$, we see that there is a constant $C>0$ such that
\[
\norm{T\eta_\rho}_{L^q(\Nc)}\meg C \norm{T\eta_\rho}_{L^p(\Nc)}
\]
for every $T\in \Sc'(\Nc)$ such that $T=T*\varphi$ and for every $\rho\in ]0,1]$.
Now, passing to the limit for $\rho\to 0^+$ we see that
\[
\norm{T}_{L^q(\Nc)}\meg C \liminf_{\rho\to 0^+}\norm{T\eta_\rho}_{L^p(\Nc)} \meg C\norm{T}_{L^p(\Nc)},
\]
whence the result
\end{proof}

\begin{cor}\label{cor:25}
Take $p,q\in ]0,\infty]$ such that $p\meg q$ and let $B$ be a bounded subset of $\Sc_\Omega(\Nc)$. For every $\varphi \in B$ and for every $t\in T_+$, define $\varphi_t\in \Sc_\Omega(\Nc)$ so that $\Fc_\Nc \varphi_t= (\Fc_\Nc \varphi)(\,\cdot\,t^{-1})$.
Then, there is a constant $C>0$ such that, for every $T\in \Sc'(\Nc)$, for every $\varphi\in B$, and for every $t\in T_+$,
\[
\norm{T*\varphi_t}_{L^q(\Nc)}\meg C \Delta^{(1/p-1/q)(\vect b+\vect d) }(t) \norm{T*\varphi_t}_{L^p(\Nc)}.
\]
In addition, if $T*\varphi_t\in L^p_0(\Nc)$, then $T*\varphi_t\in L^q_0(\Nc)$.
\end{cor}

\begin{proof}
For every $t\in T_+$, choose $g_t\in GL(E)$ such that $t\cdot \Phi=\Phi\circ (g_t\times g_t)$, and observe that Proposition~\ref{prop:40} and the fact that $(g_t \times t)$ is an automorphism of $\Nc$ imply that
\[
(g_t\times t)_*(T*\varphi_t)= [(g_t\times t)_* T] *\varphi
\]
for every $T\in \Sc'(\Nc)$, for every $\varphi\in B$, and for every $t\in T_+$.  
Now, observe that there is a compact subset $K$ of $\Omega'$ such that $\Fc_\Nc \varphi$ is supported in $K$ for every $\varphi\in B$. Then, fix $\psi\in \Sc_\Omega(\Nc)$  such  that $\Fc_\Nc \psi$  equals $1$ on $K$, so that $(T*\varphi)*\psi= T*\varphi$ for every $T\in \Sc'(\Nc)$ and for every $\varphi\in B$. Then, the first assertion follows from Lemmas~\ref{lem:70} and~\ref{lem:19},  and Corollary~\ref{cor:22}.

For what concerns the second assertion, it will suffice to show that, if $p<\infty$ and $T*\varphi_t\in L^p(\Nc)$, then $T*\varphi_t\in C_0(\Nc)$. Observe that, applying the preceding remarks with $q=\max(1,p)$, we see that $T*\varphi_t\in L^{\max(1,p)}(\Nc)$. Therefore, $T*\varphi_t=T*\varphi_t*\psi_t\in C_0(\Nc)$, whence the result.
\end{proof}

\begin{cor}\label{cor:23}
Take $\varphi\in \Sc_\Omega(\Nc)$, and fix a compact subset $K$ of $\Omega'$. Then, for every $p\in ]0,\infty]$ there is a constant $C>0$ such that
\[
\norm{T*\varphi'}_{L^p(\Nc)}\meg C \norm{T}_{L^p(\Nc)} \norm{\varphi'}_{L^{\min(1,p)}(\Nc)}
\]
for every $T\in \Sc'(\Nc)$ such that $T=T*\varphi$ and for every $\varphi '\in \Sc_\Omega(\Nc)$ such that $\Fc_\Nc \varphi'$ is supported in $K$.
\end{cor}

The proof is based on~\cite[Proposition 1.5.1]{Triebel}, which deals with the classical case.

\begin{proof}
The assertion follows by Young's inequality (with $C= 1$) when $p\Meg 1$, so that we may assume that $p<1$. Observe that $\varphi'^*\in \Sc_\Omega(\Nc)$, so that $L_{(\zeta,x)}^{-1}\varphi'^*\in \widetilde \Sc_\Omega(\Nc)$ by Proposition~\ref{prop:41} for every $(\zeta,x)\in \N$ and for every $\varphi '\in \Sc_\Omega(\Nc)$. 
Therefore, Proposition~\ref{prop:41} again shows that $T L_{(\zeta,x)}^{-1}\varphi'^*\in \widetilde \Sc_\Omega(\Nc)$. In addition, by means of Corollary~\ref{cor:21} we see that, if $\Fc_\Nc \varphi'$ is supported in $K$, then
\[
T L_{(\zeta,x)}^{-1}\varphi'^* =(T L_{(\zeta,x)}^{-1}\varphi'^*)*\varphi'',
\]
where $\varphi''\in \Sc_\Omega(\Nc)$ and $ \Fc_\Nc \varphi''$ equals $1$ on $\Supp{\Fc_\Nc \varphi}+K $. Therefore, Corollary~\ref{cor:22} shows that there is a constant $C'>0$ such that 
\[
\begin{split}
\norm{T*\varphi'}_{L^p(\Nc)}^p&\meg\int_\Nc \norm{T L_{(\zeta,x)}^{-1}\check\varphi'  }^p_{L^1(\Nc)} \,\dd (\zeta,x)\\
&=\int_\Nc \norm{T L_{(\zeta,x)}^{-1}\varphi'^*  }_{L^1(\Nc)}^p\,\dd (\zeta,x)\\
&\meg C'\int_\Nc \int_\Nc \abs{T(\zeta',x') (L_{(\zeta,x)}^{-1}\varphi'^*)(\zeta',x') }^p\,\dd (\zeta',x')\,\dd (\zeta,x)\\
&= C'\norm{T}_{L^p(\Nc)}^p \norm{\varphi'}_{L^{p}(\Nc)}^p
\end{split}
\]
for every $T\in \Sc'(\Nc)$ such that $T=T*\varphi$ and for every $\varphi '\in \Sc_\Omega(\Nc)$ such that $\Fc_\Nc \varphi'$ is supported in $K$.
\end{proof}

\begin{cor}\label{cor:24}
Take $p\in ]0,\infty]$ and let $B$ be a bounded subset of $\Sc_\Omega(\Nc)$. For every $\varphi \in B$ and for every $t\in T_+$, define $\varphi_t\in \Sc_\Omega(\Nc)$ so that $\Fc_\Nc \varphi_t= (\Fc_\Nc \varphi)(\,\cdot\,t^{-1})$.
Then, there is a constant $C>0$ such that, for every $T\in \Sc'(\Nc)$, for every $\varphi,\varphi'\in B$, and for every $t,t'\in T_+$,
\[
\norm{T*\varphi_t*\varphi'_{t'}}_{L^p(\Nc)}\meg C \norm{T*\varphi_t}_{L^p(\Nc)}.
\]
\end{cor}

\begin{proof}
By assumption, there is a compact subset $K$ of $T_+$ such that $\Fc_\Nc \varphi$ is supported in $e_{\Omega'}\cdot K$ for every $\varphi \in B$. Take $\varphi,\varphi'\in B$ and $t,t'\in T_+$. If $\varphi_t*\varphi'_{t'}\neq 0$, then $(K t)\cap (K t')\neq \emptyset$, so that $t' t^{-1}\in K K^{-1} $. 
Next, for every $t\in T_+$, choose $g_t\in GL(E)$ such that $t\cdot \Phi=\Phi\circ (g_t\times g_t)$, and observe that Proposition~\ref{prop:40} and the fact that $(g_t \times t)$ is an automorphism of $\Nc$ imply that
\[
(g_t\times t)_*(T*\varphi_t*\varphi'_{t'})= [(g_t\times t)_* T] *\varphi*\varphi'_{t' t^{-1}}
\]
and
\[
(g_t\times t)_*(T*\varphi_t)= [(g_t\times t)_* T] *\varphi
\]
for every $T\in \Sc'(\Nc)$, for every $\varphi,\varphi'\in B$, and for every $t,t'\in T_+$.  Thus, we may reduce to proving that
\[
\norm{T*\varphi*\varphi'_{t''}}_{L^p(\Nc)}\meg C_1 \norm{T*\varphi}_{L^p(\Nc)}
\]
for every $T\in \Sc'(\Nc)$, for every $\varphi,\varphi'\in B$, and for every $t''\in K K^{-1}$.  Now, take $\psi\in \Sc_\Omega(\Nc)$ so that $\Fc_\Nc \psi$  equals $1$ on $K$, and observe that $(T*\varphi)*\psi= T*\varphi$ for every $T\in \Sc'(\Nc)$ and for every $\varphi\in B$. Thus, Corollary~\ref{cor:23} implies that there is a constant $C_1>0$ such that
\[
\norm{T*\varphi*\varphi'_{t''}}_{L^p(\Nc)}\meg C_1 \norm{T*\varphi}_{L^p(\Nc)} \norm{\varphi'_{t''}}_{L^{\min(1,p)}(\Nc)}
\]
for every $T\in \Sc'(\Nc)$, for every $\varphi,\varphi'\in B$, and for every $t''\in K K^{-1}$. 
To conclude, observe that Lemma~\ref{lem:50} implies that the $\Fc_\Nc \varphi'_{t''}$, as $\varphi'$ runs through $B$ and $t''$ runs through $K K^{-1}$, stay in a bounded subset of $C^\infty_c(\Omega')$, so that the $\varphi'_{t''}$ stay in a bounded subset of $\Sc_\Omega(\Nc)$ by Proposition~\ref{prop:40}. The assertion follows.
\end{proof}

We now describe the interactions of Riemann--Liouville operators with the spaces $\Sc_\Omega(\Nc)$, $\Sc_{\Omega,L}(\Nc)$, and $\widetilde \Sc_\Omega(\Nc)$.

\begin{prop}\label{prop:55}
Take $p,q\in ]0,\infty]$ and $\vect s\in \C^r$. Then, the following hold:
\begin{enumerate}
\item[\em(1)] for every $\varphi\in \widetilde\Sc_{\Omega}(\Nc)$ and for every $\lambda\in F'\setminus W$, 
\[
\dd \pi_\lambda(\varphi*I^{\vect s}_{\Omega}  )= i^{-\vect s}\Delta_{\Omega'}^{-\vect s}(\lambda) \pi_\lambda(\varphi);
\]

\item[\em(2)] the mapping $\varphi\mapsto \varphi*I^{\vect s}_{\Omega} $ induces automorphisms of $\Sc_\Omega(\Nc)$, $\Sc_{\Omega,L}(\Nc)$, and $\widetilde \Sc_\Omega(\Nc)$;

\item[\em(3)] $\varphi*(i^{\vect s}I_{\Omega}^{\vect s})^*=\varphi* (i^{\overline{\vect s}}I_{\Omega}^{\overline {\vect s}})$ for every $\varphi\in \widetilde \Sc_\Omega(\Nc)$.
\end{enumerate}
\end{prop}

Notice that $(i^{\vect s} I_{\Omega}^{\vect s})^*\neq i^{\overline{\vect s}}I_{\Omega}^{\overline {\vect s}}$, unless $\vect s\in - \N_{\Omega'}$.

\begin{proof}
{(1)} Take $\varphi \in \widetilde \Sc_\Omega(\Nc)$, and observe that $(\varphi* I^{\vect s}_{\Omega})(\zeta,\,\cdot\,)=\varphi(\zeta,\,\cdot\,)*I^{\vect s}_{\Omega}$ for every $\zeta\in E$. Since $\Fc_F(I^{\vect s}_{\Omega'})= i^{-\vect s}\Delta^{-\vect s}_{\Omega'}$ on $\Omega'$ thanks to Lemma~\ref{lem:13}, this implies that $\varphi(\zeta,\,\cdot\,)*I^{\vect s}_{\Omega}\in \Sc(F)$. 
In addition,
\[
\pi_\lambda(\varphi*I^{\vect s}_\Omega) = \int_{E}  \Fc_F(\varphi(\zeta,\,\cdot\,)* I^{\vect s}_\Omega)(\lambda) \pi_\lambda(\zeta,0)\,\dd \zeta,
\]
for every $\lambda\in \Omega'$, 
so that 
\[
\pi_\lambda(\varphi*I^{\vect s}_\Omega) = i^{-\vect s}\Delta^{-\vect s}_{\Omega'}(\lambda) \pi_\lambda(\varphi)
\]
for every $\lambda\in F'\setminus W$, whence the assertion.

{(2)}  By~{(1)}, it is clear that convolution by $I^{\vect s}_{\Omega}$ maps $\widetilde \Sc_\Omega(\Nc)$ into itself, so that it induces an endomorphism $\Ic^{\vect s}$ of $\widetilde \Sc_\Omega(\Nc)$ by the closed graph theorem. Since $\Ic^{\vect s}\Ic^{-\vect s}=\Ic^{-\vect s}\Ic^{\vect s}=I_{\widetilde \Sc_\Omega(\Nc)}$, $\Ic^{\vect s}$ is an automorphism of $\widetilde \Sc_\Omega(\Nc)$. 
From~{(1)} it also follows that $\Ic^{\vect s}$ maps $\Sc_\Omega(\Nc)$ into itself, so that it induces an automorphism of $\Sc_\Omega(\Nc)$. Finally,  it is clear that $\Ic^{\vect s}$ induces an automorphism of $\Sc_{\Omega,L}(\Nc)$.

{(3)} Take $\varphi \in \widetilde \Sc_\Omega(\Nc)$. Then,
\[
\begin{split}
\Fc_F((\varphi* (i^{\vect s}I^{\vect s}_{\Omega})^*)(\zeta,\,\cdot\,))&= \Fc_F(\varphi(\zeta,\,\cdot\,)*(I^{\vect s}_{\Omega})^*)\\
&=\Fc_F(\varphi(\zeta,\,\cdot\,)) \overline{\Fc_F(i^{\vect s} I^{\vect s}_{\Omega})}\\
&= \Fc_F (\varphi(\zeta,\,\cdot\,)) \Fc_F(i^{\overline{\vect s}} I^{\overline {\vect s}}_\Omega)\\
&=\Fc_F((\varphi* (i^{\overline{\vect s}} I^{\overline{\vect s}}_{\Omega}))(\zeta,\,\cdot\,))
\end{split}
\]
for every $\zeta\in E$, since $\overline{\Fc_F ( i^{{\vect s}}I^{\vect s}_{\Omega})}=\overline{\Delta^{-\vect s}_{\Omega'} }= \Delta^{-\overline {\vect s}}_{\Omega'}=\Fc_F (i^{\overline{\vect s}}I^{\overline{\vect s}}_\Omega)$ on $\Omega'$ by Lemma~\ref{lem:13}. The assertion follows.
\end{proof}

We conclude this section with some remarks on the dual of $\Sc_{\Omega,L}(\Nc)$.

\begin{deff}\label{48}
We denote by $\Sc_{\Omega,L}'(\Nc,K)$ and $\Sc_{\Omega,L}'(\Nc)$ the strong duals of $\Sc_{\Omega,L}(\Nc,K)$ and $\Sc_{\Omega,L}(\Nc)$, respectively, for every compact subset $K$ of $\Omega'$.
We also define $\langle T\vert \varphi\rangle\coloneqq \langle T, \overline \varphi\rangle$ for every $(T,\varphi)\in \Sc_{\Omega,L}'(\Nc,K)\times \Sc_{\Omega,L}(\Nc,K)$ and for every $(T,\varphi)\in \Sc_{\Omega,L}'(\Nc)\times \Sc_{\Omega,L}(\Nc)$.
\end{deff}

\begin{prop}\label{prop:63}
The following hold:
\begin{enumerate}
\item[\em(1)] $\Sc_{\Omega,L}'(\Nc,K)$ is canonically isomorphic to the quotient of $\Sc'(\Nc)$ by the polar of $\Sc_{\Omega,L}(\Nc,K)$, for every compact subset $K$ of $\Omega'$;

\item[\em(2)] $\Sc_{\Omega,L}'(\Nc)$ is canonically isomorphic to the projective limit of the $\Sc_{\Omega,L}'(\Nc,K)$, where $K$ runs through the set of compact subsets of $\Omega'$;

\item[\em(3)] $\Sc_{\Omega,L}'(\Nc,K)$ is a complete Montel space;

\item[\em(4)] for every $T\in \Sc_{\Omega,L}'(\Nc)$ and for every $\varphi\in \Sc_{\Omega}(\Nc)$, we may define
\[
(T*\varphi)(\zeta,x)\coloneqq \langle T\vert L_{(\zeta,x)}\varphi^*\rangle
\]
for every $(\zeta,x)\in \Nc$. In addition,  $T*\varphi=T'*\varphi$ for every $T'\in \Sc'(\Nc)$ such that the canonical images of $T$ and $T'$ in $\Sc_{\Omega,L}(\Nc,\Supp{\Fc_\Nc \varphi})$ are equal.
\end{enumerate}
\end{prop}

In particular, Corollaries~\ref{cor:22} to~\ref{cor:24} apply with $\Sc'(\Nc)$ replaced by $\Sc_{\Omega,L}'(\Nc)$.

\begin{proof}
{(1)} Obvious.

{(2)} This follows from the characterization of the bounded subsets of $\Sc_{\Omega,L}(\Nc)$ provided in Proposition~\ref{prop:41}.

{(3)} This follows from Proposition~\ref{prop:41} and~\cite[Corollary 1 to Proposition 12 of Chapter III, \S 3, No.\ 8, and Proposition 9 of Chapter IV, \S 2, No.\ 5]{BourbakiTVS}.

{(4)} Since clearly $\varphi^*\in \Sc_\Omega(\Nc)$, $T*\varphi$ is well defined thanks to Proposition~\ref{prop:41}. The second assertion is clear, since $L_{(\zeta,x)}\varphi^*\in \Sc_{\Omega,L}(\Nc, \Supp{\Fc_\Nc\varphi})$ for every $(\zeta,x)\in \Nc$.
\end{proof}

\section{Besov Spaces}\label{sec:5:2}

In this section we develop the basics of the theory of the Besov spaces $B^{\vect s}_{p,q}(\Nc,\Omega)$.
As we shall see in Theorem~\ref{teo:11} (cf.~also Theorem~\ref{teo:9} and the proof of Proposition~\ref{prop:54}), for every $T\in B^{\vect s}_{p,q}(\Nc,\Omega)$ there is a holomorphic function $f$ on $D$ such that
\[
T=\lim_{h\to 0} f_h
\]
in $\Sc'_{\Omega,L}(\Nc)$, even though $f$ is not uniquely determined by this property. 
For this reason, one may say that the $B^{\vect s}_{p,q}(\Nc,\Omega)$ are Besov spaces `of analytic type.'

\smallskip

\emph{In the following results, by an abuse of notation we shall often write $\norm{a_k}_{\ell^q(K)}$ instead of $\norm{(a_k)}_{\ell^q(K)}$, for $(a_k)\in \ell^q(K)$.}

\medskip

Before we define the spaces $B^{\vect s}_{p,q}(\Nc,\Omega)$, we prove the equivalence of several (generalized) quasi-norms on $\Sc_{\Omega,L}'(\Nc)$.

\begin{lem}\label{lem:51}\label{lem:53}
Take $\delta,\delta'>0$, $R,R'>1$, $p,q\in ]0,\infty]$, and $\vect s\in \R^r$.
Let $(\lambda_k)_{k\in K}$ and $(\lambda'_{k'})_{k'\in K'}$ be a $(\delta,R)$- and a $(\delta',R')$-lattice on $\Omega$, respectively, and  define $t_k, t'_{k'}\in T_+$ so that 
\[
\lambda_k=e_{\Omega'}\cdot t_k \qquad \text{and} \qquad \lambda'_{k'}=e_{\Omega'}\cdot t'_{k'}
\]
for every $k\in K$ and for every $k'\in K'$. 
Let $(\varphi_k)_{k\in K}$ and $(\varphi'_{k'})_{k'\in K'}$ be two bounded families of positive elements of $C^\infty_c(\Omega')$ such that 
\[
\sum_{k\in K} \varphi_k(\,\cdot\,t_k^{-1})\Meg 1 \qquad \text{and} \qquad \sum_{k'\in K'} \varphi_{k'}(\,\cdot\,t_{k'}'^{-1})\Meg 1 
\]
on $\Omega'$, and define $\psi_k, \psi'_{k'}\in \Sc_\Omega(\Nc)$ so that $\Fc_\Nc \psi_k=\varphi_k(\,\cdot\, t_k^{-1})$ and $\Fc_\Nc\psi'_{k'}=\varphi'_{k'}(\,\cdot\,t'^{-1}_{k'})$ for every $k\in K$, and for every $k'\in K'$.
Then, there is a constant $C>0$ such that, for every $T\in \Sc_{\Omega,L}'(\Nc)$,
\[
\begin{split}
\frac 1 C \norm*{ \Delta_{\Omega'}^{\vect s}(h_{k'}) \norm{T*\psi'_{k'} }_{L^p(\Nc)}  }_{L^q(K')}&\meg \norm*{ \Delta_{\Omega'}^{\vect s}(h_k) \norm{T*\psi_k }_{L^p(\Nc)}  }_{L^q(K)}\\
&\meg C \norm*{ \Delta_{\Omega'}^{\vect s}(h_{k'}) \norm{T*\psi'_{k'} }_{L^p(\Nc)}  }_{L^q(K')}.
\end{split}
\] 
\end{lem}

The argument is classical. Cf., for example,~\cite[Lemma 3.8]{BekolleBonamiGarrigosRicci}.

\begin{proof}
It will suffice to prove the first inequality.
Define
\[
\widetilde \varphi\coloneqq \sum_{k\in K} \varphi_k(\,\cdot\,t_k^{-1}),
\]
and observe that, since $(\varphi_k)$ is bounded in $C^\infty_c(\Omega')$, Proposition~\ref{prop:56} implies that the sum defining $\widetilde \varphi$ is locally finite on $\Omega'$. Therefore, $\widetilde \varphi$ is of class $C^\infty$ on $\Omega'$. 
For every $k'\in K'$, let $K_{k'}$ be the set of $k\in K$ such that $\varphi'_{k'}(\,\cdot\,t_{k'}^{-1})\varphi_k(\,\cdot\,t_k^{-1})\neq 0 $. By Proposition~\ref{prop:56}, there is $N\in \N$ such that $\card(K_{k'})\meg N$ for  every $k'\in K'$ and such that each $k\in K$ is contained in at most $N$ of the sets $K_{k'}$, for $k'\in K'$.
Then, for every $k'\in K'$,
\[
\varphi'_{k'}(\,\cdot\,t'^{-1}_{k'})=\sum_{k\in K_{k'}} \frac{\varphi_{k'}'(\,\cdot\,t'^{-1}_{k'}) \varphi_k (\,\cdot\,t^{-1}_{k})}{\widetilde \varphi}= \sum_{k\in K_{k'}} \widetilde \varphi'_{k'}(\,\cdot\,t'^{-1}_{k'}) \varphi_k (\,\cdot\,t^{-1}_{k}),
\]
where
\[
\widetilde \varphi'_{k'}\coloneqq \frac{\varphi_{k'}'}{\widetilde \varphi(\,\cdot\, t_{k'}')}= \frac{\varphi'_{k'}}{\sum_{k\in K_{k'}} \varphi_k(\,\cdot\, (t'_{k'} t_k^{-1} )) }.
\]
Observe that $\sum_{k\in K_{k'}} \varphi_k(\,\cdot\, (t'_{k'} t_k^{-1} ))\Meg 1$ on the support of $\varphi'_{k'}$ and that the $t'_{k'} t_k^{-1} $ stay in a compact subset of $T_+$ (independent of $k'$) as $k$ runs through $K_{k'}$, so thatLemma~\ref{lem:50} and the preceding remarks imply that the family $(\widetilde \varphi'_{k'})_{k'\in K'}$ is bounded in $C^\infty_c(\Omega')$.
Then, choose $\widetilde \psi'_{k'}\in \Sc_\Omega(\Nc)$ so that $\Fc_\Nc \widetilde \psi'_{k'}=\widetilde \varphi'_{k'}(\,\cdot\,t_{k'}'^{-1})$ for every $k'\in K'$.
Then, Corollary~\ref{cor:24}  and Proposition~\ref{prop:63} imply that there is a constant $C_1>0$ such that, for every $T\in \Sc_{\Omega,L}'(\Nc)$, for every $k'\in K'$, and for every $k\in K$,
\[
\norm{(T*\psi_{k})*\widetilde \psi'_{k'}}_{L^p(\Nc)}\meg C_1\norm{T*\psi_k}_{L^p(\Nc)}.
\]
Now, for every $k'\in K'$,
\[
T*\psi'_{k'}=\sum_{k\in K_{k'}}(T*\psi_k)*\widetilde \psi'_{k'}
\]
by the associativity of convolution on $\Sc_{\Omega,L}'(\Nc)\times \Sc(\Nc)\times \Sc(\Nc)$, so that
\[
\norm{T*\psi'_{k'}}_{L^p(\Nc)}\meg C_1 N^{(1/p-1)_+}\sum_{k\in K_{k'}}\norm{T*\psi_k}_{L^p(\Nc)} 
\]
for every $k'\in K'$. In addition, Corollary~\ref{cor:34} shows that there is a constant $C_2>0$ such that
\[
\Delta^{\vect s}_{\Omega'}(\lambda_{k'}')\norm{T*\psi'_{k'}}_{L^p(\Nc)}\meg C_2\sum_{k\in K_{k'}} \Delta^{\vect s}_{\Omega'}(\lambda_{k})\norm{T*\psi_k}_{L^p(\Nc)}
\]
for every $k'\in K'$.

Therefore, 
\[
\norm*{ \Delta^{\vect s}_{\Omega'}(\lambda_{k'}')\norm{T*\psi'_{k'}}_{L^p(\Nc)} }_{L^q(K')}\meg C_2 N^{ 1/q+1/q' }\norm*{ \Delta^{\vect s}_{\Omega'}(\lambda_{k})\norm{T*\psi_{k}}_{L^p(\Nc)} }_{L^q(K)},
\]
whence the result.
\end{proof}

\begin{deff}\label{def:4}\index{Besov spaces}
Take $\vect s\in \R^r$ and $p,q\in ]0,\infty]$. In addition, take $(\lambda_k)_{k\in K}$, $(\varphi_k)$, and $(\psi_k)$ as in Lemma~\ref{lem:51}.
Then, we define $B^{\vect s}_{p,q}(\Nc,\Omega)$ as the space of $T\in \Sc_{\Omega,L}'(\Nc)$ such that 
\[
(\Delta^{\vect s}_{\Omega'}(\lambda_k) (T*\psi_k))_k\in L^q(K;L^p(\Nc)),
\]
endowed with the corresponding topology.\label{33}

We denote by  $\mathring{B}^{\vect s}_{p,q}(\Nc,\Omega)$ the closed subspace of $B^{\vect s}_{p,q}(\Nc,\Omega)$ consisting of the $T\in \Sc_{\Omega,L}'(\Nc)$ such that 
\[
(\Delta^{\vect s}_{\Omega'}(\lambda_k) (T*\psi_k))_k\in L^q_0(K;L^p_0(\Nc)).
\]
\end{deff}

As we shall see later,  $\Sc_{\Omega,L}(\Nc)$ is (canonically embedded and) dense in $\mathring B^{\vect s}_{p,q}(\Nc,\Omega)$ (cf.~Theorem~\ref{teo:10}), whence the notation.

\begin{prop}\label{prop:20}
Take $p,q\in ]0,\infty]$ and $\vect s\in \R^r$.  Then, $B^{\vect s}_{p,q}(\Nc,\Omega)$ and $\mathring{B}^{\vect s}_{p,q}(\Nc,\Omega)$ are locally bounded $F$-spaces. In addition, the inclusion $B^{\vect s}_{p,q}(\Nc,\Omega)\subseteq \Sc_{\Omega,L}'(\Nc)$ is continuous.
\end{prop}

In particular, $B^{\vect s}_{p,q}(\Nc,\Omega)$ and $\mathring{B}^{\vect s}_{p,q}(\Nc,\Omega)$ are complete. 

\begin{proof}
Observe first that, with the notation of Definition~\ref{def:4}, the mapping
\[
N\colon \Sc_{\Omega,L}'(\Nc)\ni T\mapsto \norm*{\Delta_{\Omega'}^{\vect s}(\lambda_k )\norm{T*\psi_k}_{L^p(\Nc)}}_{\ell^q(K)}\in \R_+
\]
is lower semi-continuous. In addition, $N$ is finite exactly on $B^{\vect s}_{p,q}(\Nc,\Omega) $, and the mapping $(T_1,T_2)\mapsto N(T_1-T_2)^{\min(1,p,q)}$ is a distance on $B^{\vect s}_{p,q}(\Nc,\Omega)$ which is compatible with its topology.  Then, let $(T_j)$ be a Cauchy sequence on $B^{\vect s}_{p,q}(\Nc,\Omega)$, so that $(T_j)$ is a Cauchy sequence on $\Sc_{\Omega,L}'(\Nc)$, hence converges to some $T$ in $\Sc_{\Omega,L}'(\Nc)$ by Proposition~\ref{prop:63}. Since $N$ is lower semi-continuous on $\Sc_{\Omega,L}'(\Nc)$, it then follows that $T\in B^{\vect s}_{p,q}(\Nc,\Omega)$ and that $(T_j)$ converges to $T$ in $B^{\vect s}_{p,q}(\Nc,\Omega)$. The assertion follows.
\end{proof}

In the following result, we show how the various affine automorphisms of $D$ interact with the spaces $B^{\vect s}_{p,q}(\Nc,\Omega)$. 

\begin{prop}\label{prop:42}
Take $\vect s\in \R^r$ and $p,q\in ]0,\infty]$. Denote by $G$ the set of automorphisms of $\Nc$ of the form $g\times t$, where $t\in T_+$, $g\in GL(E)$, and $t\cdot \Phi=\Phi\circ (g\times g)$.
Then, the following hold:
\begin{enumerate}
\item[\em(1)] the mappings $T\mapsto L_{(\zeta,x)} T$, for $(\zeta,x)\in \Nc$, induce equicontinuous automorphisms of $B^{\vect s}_{p,q}(\Nc,\Omega)$ and $\mathring B^{\vect s}_{p,q}(\Nc,\Omega)$;

\item[\em(2)] for every $T\in \mathring B^{\vect s}_{p,q}(\Nc,\Omega)$, the mapping $\Nc\ni (\zeta,x)\mapsto L_{(\zeta,x)} T\in \mathring B^{\vect s}_{p,q}(\Nc,\Omega)$ is continuous;

\item[\em(3)]  the mappings \CR $T\mapsto \Delta^{\vect s-(1-1/p)(\vect b+\vect d)}(t)(g\times t)_* T$, \CB for $g\times t\in G$, induce equicontinuous automorphisms of $B^{\vect s}_{p,q}(\Nc,\Omega)$ and $\mathring B^{\vect s}_{p,q}(\Nc,\Omega)$;

\item[\em(4)] for every $T\in \mathring B^{\vect s}_{p,q}(\Nc,\Omega)$, the mapping $G\ni g\times t \mapsto  (g\times t)_* T\in \mathring B^{\vect s}_{p,q}(\Nc,\Omega)$ is continuous.
\end{enumerate}
\end{prop}

\begin{proof}
Fix a $(\delta,R)$-lattice $(\lambda_k)_{k\in K}$ on $\Omega'$ for some $\delta>0$ and some $R>1$ (cf.~Lemma~\ref{lem:32}), and define $t_k\in T_+$ so that $ \lambda_k=e_{\Omega'}\cdot t_k$ for every $k\in K$. 
Fix a positive $\varphi\in C^\infty_c(\Omega')$ such that $\sum_{k\in K}\varphi(\,\cdot\, t_k^{-1})\Meg 1$ on $\Omega'$, and define $\psi_k\in \Sc_\Omega(\Nc)$ so that $\Fc_\Nc \psi_k=  \varphi(\,\cdot\,t_k^{-1})$ for every $k\in K$. By Lemma~\ref{lem:53}, we may choose the norm on $B_{p,q}^{\vect s}(\Nc,\Omega)$ induced by the mapping
\[
T\mapsto \norm*{\Delta^{\vect s}_{\Omega'}(\lambda_k) \norm{T*\psi_k}_{L^p(\Nc)}}_{\ell^q(K)}.
\]
Then, it is immediately seen that the mapping $T\mapsto L_{(\zeta,x)} T$ induces an isometry of $B_{p,q}^{\vect s}(\Nc,\Omega)$ and $\mathring B^{\vect s}_{p,q}(\Nc,\Omega)$ for every $(\zeta,x)\in \Nc$, whence~{(1)}. 
For what concerns~{(2)}, fix $T\in \mathring B_{p,q}^{\vect s}(\Nc,\Omega)$ and observe that 
\[
\lim_{(\zeta,x)\to (0,0)}(L_{(\zeta,x)} T)*\psi_k=\lim_{(\zeta,x)\to (0,0)} L_{(\zeta,x)} (T*\psi_k)= T*\psi_k
\]
in $L^p_0(\Nc)$, for   every $k\in K$.
Thus, $L_{(\zeta,x)} T$ converges to $T$ in $B_{p,q}^{\vect s}(\Nc,\Omega)$ as $(\zeta,x)\to (0,0)$ (by dominated convergence, if $q<\infty$. By equicontinuity at the point at infinity of $K$, if $q=\infty$). This proves continuity at $(0,0)$. Continuity on $\Nc$ follows since $\Nc$ is a group.

Now, take $(g\times t)\in G$. Then, 
\[
\begin{split}
&\norm*{\Delta^{\vect s}_{\Omega'}(\lambda_k) \norm*{[(g\times t)_*T]*\psi_k}_{L^p(\Nc)}}_{\ell^q(K)}\\&\qquad \qquad=\norm*{\Delta^{\vect s}_{\Omega'}(\lambda_k) \norm*{(g\times t)_*[T*(g\times t)^*\psi_k]}_{L^p(\Nc)}}_{\ell^q(K)}\\
&\qquad \qquad= \Delta^{(\vect b+\vect d)(1-1/p)}(t)\norm*{ \Delta^{\vect s}_{\Omega'}(\lambda_k) \norm*{T*(g\times t)^*\psi_k}_{L^p(\Nc)}}_{\ell^q(K)}
\end{split}
\]
for every $T\in \Sc_{\Omega,L}'(\Nc)$. In addition, Proposition~\ref{prop:40} shows that $(g\times t)^*\psi_k\in\Sc_\Omega(\Nc)$ and that
\[
\Fc_\Nc((g\times t)^*\psi_k)=\varphi( \,\cdot\,(t_k t)^{-1}),
\] 
for every $k\in K$. 
Since $(\lambda_k\cdot t)$ is still a $(\delta,R)$-lattice, Lemma~\ref{lem:51} shows that the mapping 
\[
T\mapsto \norm*{ \Delta^{\vect s}_{\Omega'}(\lambda_k) \norm*{T*(g\times t)^*\psi_k}_{L^p(\Nc)}}_{\ell^q(K)}
\]
induces a norm on $B^{\vect s}_{p,q}(\Nc,\Omega)$ which defines its topology, whence~{(3)}. Equicontinuity follows from a refinement of the proof of Lemma~\ref{lem:51}. Finally,~{(4)} is proved as~{(2)}.
\end{proof}

\begin{lem}\label{lem:55}
The canonical embedding $\Sc_{\Omega,L}(\Nc)\to \Sc_{\Omega,L}'(\Nc)$ induces a  canonical embedding $\Sc_{\Omega,L}(\Nc)\to \mathring{B}^{\vect s}_{p,q}(\Nc,\Omega)$.
\end{lem}

\begin{proof}
This follows immediately from the definitions of the spaces $\Sc_{\Omega,L}(\Nc)$ and $\mathring{B}^{\vect s}_{p,q}(\Nc,\Omega)$, and from the fact that $\Sc(\Nc)\subseteq L^p_0(\Nc)$.
\end{proof}

In the following result we prove the analogues of the classical Sobolev embeddings between homogeneous Besov spaces on $\R$.

\begin{prop}\label{prop:27}
Take $p_1,p_2,q_1,q_2\in ]0,\infty]$ and $\vect{s_1},\vect{s_2}\in \R^r$ such that 
\[
p_1\meg p_2, \qquad q_1\meg q_2, \qquad \text{and} \qquad
\vect{s_2}=\vect{s_1}+\left(\frac{1}{p_1}-\frac{1}{p_2}\right)(\vect b+\vect d).
\]
Then,  we have continuous inclusions
\[
B^{\vect{s_1}}_{p_1,q_1}(\Nc,\Omega)\subseteq B^{\vect{s_2}}_{p_2,q_2}(\Nc,\Omega) \qquad \text{and} \qquad\mathring{B}^{\vect{s_1}}_{p_1,q_1}(\Nc,\Omega)\subseteq \mathring{B}^{\vect{s_2}}_{p_2,q_2}(\Nc,\Omega).
\]
\end{prop}

\begin{proof}
The assertions follow from Corollary~\ref{cor:25} and the continuous inclusions $\ell^{q_1}(K)\subseteq \ell^{q_2}(K)$ and $\ell^{q_1}_0(K)\subseteq \ell^{q_2}_0(K)$.
\end{proof}

In the following, we shall denote by $\langle \,\cdot\,\vert \,\cdot\,\rangle$ the sesquilinear pairing between $L^p(\Nc)$ and $L^{p'}(\Nc)$ defined by 
\[
\langle f_1\vert f_2\rangle= \int_\Nc f_1(\zeta,x)\overline{f_2(\zeta,x)}\,\dd (\zeta,x),
\]
for every $p\in [1,\infty]$.

Our aim is to extend this sesquilinear pairing to a sesquilinear form on 
\[
B^{\vect s}_{p,q}(\Nc,\Omega)\times B^{-\vect s-(1/p-1)_+(\vect b+\vect d)}_{p',q'}(\Nc,\Omega).
\]
If $p,q\in ]0,\infty[$ (resp.\ if $p,q\in ]1,\infty]$), then $\Sc_{\Omega,L}(\Nc)$ is dense in $ B^{\vect s}_{p,q}(\Nc,\Omega) $ (resp.\ $B^{-\vect s-(1/p-1)_+(\vect b+\vect d)}_{p',q'}(\Nc,\Omega)$) by Theorem~\ref{teo:10}, so that this form is uniquely determined by its restriction to $\Sc_{\Omega,L}(\Nc)\times \Sc_{\Omega,L}'(\Nc)$ (resp.\ $\Sc_{\Omega,L}'(\Nc)\times \Sc_{\Omega,L}(\Nc)$).
Nonetheless, in the general  case  we can only give a direct definition and show its independence from the most reasonable constructions.

\begin{prop}\label{prop:18}
Take $p,q\in ]0,\infty]$ and $\vect s\in \R^r$.
Let $(\lambda_k)_{k\in K}$ and $(\lambda'_{k'})_{k'\in K'}$ be a $(\delta,R)$- and a $(\delta',R')$-lattice on $\Omega'$ for some $\delta,\delta'>0$ and some $R,R'>1$. Define $t_k,t'_{k'}\in T_+$ so that 
\[
\lambda_k=e_{\Omega'}\cdot t_k\qquad \text{and} \qquad \lambda'_{k'}=e_{\Omega'}\cdot t'_{k'}
\]
for every $k\in K$ and for every $k'\in K'$.
Let $(\varphi_k)_{k\in K}$ and $(\varphi'_{k'})_{k'\in K'}$ be two bounded families of positive elements of $C^\infty_c(\Omega')$ such that 
\[
\sum_{k\in K} \varphi_k( \,\cdot \,t_k^{-1})^2=\sum_{k'\in K'} \varphi'_{k'}( \,\cdot \,t_{k'}'^{-1})^2=1 
\]
on $ \Omega'$, and define $\psi_k,\psi'_{k'}\in \Sc_\Omega(\Nc)$ so that $\Fc_\Nc \psi_k=\varphi_k(\,\cdot\, t_k^{-1})$ and $\Fc_\Nc \psi'_{k'}=\varphi'_{k'}(\,\cdot\, t'^{-1}_{k'})$ for every $k\in K$ and for every $k'\in K'$.
Then, the following hold:
\begin{enumerate}
\item[\em(1)] the sesquilinear forms on $B^{\vect s}_{p,q}(\Nc,\Omega)\times B^{-\vect s-(1/p-1)_+(\vect b+\vect d)}_{p',q'}(\Nc,\Omega)$
\[
(T_1,T_2)\mapsto \sum_{k\in K} \langle T_1*\psi_k\vert T_2*\psi_k\rangle
\]
and
\[
(T_1,T_2)\mapsto \sum_{k'\in K'} \langle T_1*\psi'_{k'}\vert T_2*\psi'_{k'}\rangle
\]
are well defined and continuous;

\item[\em(2)] for every $T_1\in B^{\vect s}_{p,q}(\Nc,\Omega)$ and for every $ T_2\in B^{-\vect s-(1/p-1)_+(\vect b+\vect d)}_{p',q'}(\Nc,\Omega)$,
\[
\sum_{k\in K} \langle T_1*\psi_k\vert T_2*\psi_k\rangle=\sum_{k'\in K'} \langle T_1*\psi'_{k'}\vert T_2*\psi'_{k'}\rangle.
\]
\end{enumerate}
\end{prop}

\begin{proof}
{(1)} This follows easily from  a twofold application of H\"older's inequality, since Proposition~\ref{prop:27} shows that the space $B^{\vect s}_{p,q}(\Nc,\Omega)$ embeds continuously into the space $B^{\vect s+(1/p-1)_+(\vect b+\vect d)}_{\max(1,p),\max(1,q)}(\Nc,\Omega)$. 

{(2)} Fix $T_1\in B^{\vect s}_{p,q}(\Nc,\Omega)$ and $T_2\in B^{-\vect s-(1/p-1)_+(\vect b+\vect d)}_{p',q'}(\Nc,\Omega)$, and observe that clearly
\[
\psi_k=\sum_{k'\in K'} \psi_k*\psi'_{k'}*\psi'_{k'}
\]
for every $k\in K$ (the sum being finite).  Therefore, 
\[
T_1*\psi_k= \sum_{k'\in K'} T_1*\psi_k*\psi'_{k'}*\psi'_{k'}
\]
for every $k\in K$ (the sum being finite).
Therefore,
\[
\langle T_1*\psi_k\vert T_2*\psi_k\rangle= \sum_{k'\in K'}\langle T_1*\psi_k*\psi'_{k'}*\psi'_{k'} \vert  T_2*\psi_k\rangle
\]
for every $k\in K$. 
Now, observe that $\psi_k=\psi_k^*$ and that $\psi'_{k'}= \psi_{k'}'^*$ since $\varphi_k$ and $\varphi'_{k'}$ are (positive, hence) real, for every $k\in K$ and for every $k'\in K'$. In addition, convolution is commutative on $\Sc_\Omega(\Nc)$ by Proposition~\ref{prop:40}, so that
\[
\sum_{k'\in K'}\langle T_1*\psi_k*\psi'_{k'}*\psi'_{k'} \vert  T_2*\psi_k\rangle=\sum_{k'\in K'}\langle T_1*\psi'_{k'}*\psi_k*\psi_k\vert  T_2*\psi'_{k'}\rangle
\] 
for every $k\in K$.
Now, clearly 
\[
\begin{split}
&\sum_{k\in K}\sum_{k'\in K'}\abs{\langle T_1*\psi'_{k'}*\psi_k*\psi_k\vert  T_2*\psi'_{k'}\rangle}\\
&\qquad\meg \sum_{k\in K}\sum_{k'\in K'} \norm{T_1*\psi'_{k'}*\psi_k*\psi_k}_{L^{\max(1,p)}(\Nc)} \norm{T_2*\psi'_{k'}}_{L^{p'}(\Nc)}\\
&\qquad \meg  N \sup\limits_{k\in K} \norm{\psi_k}_{L^1(\Nc)}^2 \sum_{k'\in K'} \norm{T_1*\psi'_{k'}}_{L^{\max(1,p)}(\Nc)} \norm{T_2*\psi'_{k'}}_{L^{p'}(\Nc)}<\infty,
\end{split}
\]
where $N\coloneqq \sup\limits_{k'\in K'} \card(\Set{k\in K \colon \psi_k*\psi'_{k'}\neq 0 })<\infty$ (cf.~Proposition~\ref{prop:56}).
Therefore, arguing as above we see that
\[
\begin{split}
\sum_{k\in K}\langle T_1*\psi_k\vert T_2*\psi_k\rangle&=\sum_{k'\in K'}\sum_{k\in K} \langle T_1*\psi'_{k'}*\psi_k*\psi_k\vert  T_2*\psi'_{k'}\rangle\\
&=\sum_{k'\in K'} \langle T_1*\psi'_{k'}\vert  T_2*\psi'_{k'}\rangle,
\end{split}
\]
whence the result.
\end{proof}

\begin{deff}
Take $p,q\in ]0,\infty]$ and $\vect s\in \R^r$. Then, we define a sesquilinear form $\langle \,\cdot\,\vert \,\cdot\,\rangle\colon B^{\vect s}_{p,q}(\Nc,\Omega)\times B^{-\vect s-(1/p-1)_+(\vect b+\vect d)}_{p',q'}(\Nc,\Omega)\to \C$ in any  one  of the equivalent ways of Proposition~\ref{prop:18}.

We denote by $\sigma^{\vect s}_{p,q}$ the corresponding weak topology (cf.~\cite[Definition 2 of Chapter II, \S 6, No.\ 2]{BourbakiTVS})\label{44} \index{Weak topology}
\[
\sigma\Big(B^{\vect s}_{p,q}(\Nc,\Omega), \mathring B^{-\vect s-(1/p-1)_+(\vect b+\vect d)}_{p',q'}(\Nc,\Omega)\Big).
\]
\end{deff}

As we shall see later, (cf.~Theorem~\ref{teo:10} and Corollary~\ref{cor:30}), the sesquilinear form $\langle \,\cdot\,\vert \,\cdot\,\rangle$ induces a antilinear isomorphism of $B^{\vect s}_{p,q}(\Nc,\Omega)$ onto the dual of  $\mathring B^{-\vect s}_{p',q'}(\Nc,\Omega)$, provided that $p,q\in [1,\infty]$. Therefore, in this case the weak topology $\sigma^{\vect s}_{p,q}$ is simply the weak dual topology. 
Even though in the general situation there is no such interpretation, the weak topology $\sigma^{\vect s}_{p,q}$ still has some of the most important properties of the weak dual topology. In particular, it is Hausdorff (cf.~ Proposition~\ref{prop:20}) and the bounded subsets of $B^{\vect s}_{p,q}(\Nc,\Omega)$ are relatively compact for $\sigma^{\vect s}_{p,q}$ (cf.~Corollary~\ref{cor:28} below).

The next result shows that the  sum $\sum_k \,\cdot\,*\psi_k$ (with the $\psi_k$ chosen as below) enjoys nice `reproducing properties' in the spaces $\mathring B^{\vect s}_{p,q}(\Nc,\Omega)$ and $  B^{\vect s}_{p,q}(\Nc,\Omega) $.

\begin{lem}\label{lem:59}
Take $p,q\in ]0,\infty]$ and $\vect s\in \R^r$. In addition, take $(\lambda_k)_{k\in K}$, $(\varphi_k)$, $(t_k)$, and $(\psi_k)$ as in Lemma~\ref{lem:51}, and assume further that
\[
\sum_{k\in K} \varphi_k(\lambda\cdot t_k^{-1})=1
\]
for every $\lambda\in \Omega'$.
Then, for every $T\in \mathring B^{\vect s}_{p,q}(\Nc,\Omega)$ (resp.\ $T\in  B^{\vect s}_{p,q}(\Nc,\Omega)$), 
\[
T=\sum_{k\in K} T*\psi_k,
\]
with convergence in $B^{\vect s}_{p,q}(\Nc,\Omega)$ (resp.\ in the weak topology $\sigma_{p,q}^{\vect s}$).
\end{lem}

\begin{proof}
By transposition, it suffices to show the assertion when $T\in  \mathring B^{\vect s}_{p,q}(\Nc,\Omega)$.
Let $K'$ be a finite subset of $K$ and define $T_{K'}\coloneqq \sum_{k\in K'} T*\psi_k$. 
Let $V$ be the (symmetric) set of $(k,k')\in K\times K$ such that $\psi_k*\psi_{k'}\neq 0 $. By means of Proposition~\ref{prop:56} we then see that there is $N\in\N$ such that $\card(V(k))\meg N$ for every $k\in K$, where $V(k)\coloneqq \Set{k'\in K\colon (k,k')\in V}$.
In addition, Corollary~\ref{cor:24} shows that  there is a constant $C_1>0$ such that 
\[
\norm{T'*\psi_k*\psi_{k'}}_{L^p(\Nc)}\meg C_1 \norm{T'*\psi_k}_{L^p(\Nc)}
\]
for every $T'\in \Sc_{\Omega,L}'(\Nc)$ and for every $k,k'\in K$.
Then,
\[
\norm{T_{K'}*\psi_k}_{L^p(\Nc)}\meg C_1 N^{\max(1,1/p)} \chi_{V(K')}(k) \norm{T*\psi_{k}}_{L^p(\Nc)},
\] 
so that
\[
\begin{split}
&\norm*{ \Delta^{\vect s}_{\Omega'}(\lambda_k) \norm*{T_{K'}*\psi_k}_{L^p(\Nc)}}_{\ell^q(K)}\\
&\qquad \qquad\meg C_1 N^{\max(1,1/p)} \norm*{\chi_{V(K')}(k) \Delta^{\vect s}_{\Omega'}(\lambda_k) \norm*{T*\psi_k}_{L^p(\Nc)}}_{\ell^q(K)}.
\end{split}
\]
By Cauchy's criterion and Proposition~\ref{prop:20}, the sum $\sum_{k\in K} T*\psi_k$ converges in $B^{\vect s}_{p,q}(\Nc,\Omega)$. To prove that $T=\sum_{k\in K} T*\psi_k$, it suffices to show that $T*\psi_{k'}=\left(\sum_{k\in K} T*\psi_{k}\right)*\psi_{k'}$ for every $k'\in K$, but this follows from Proposition~\ref{prop:40}.
\end{proof}

\CR 
%Even though it is possible to give a direct (but somewhat more technical) proof of the following result also for $p\in ]0,\infty]$, we prefer to wait until  Corollary~\ref{cor:30}, where a different interpretation of $\langle\,\cdot\,\vert \,\cdot\,\rangle$ is provided.\CB \textbf{COMMENTO:} sono tentato di mettere qui la dimostrazione diretta e di eliminare il Corollario~\ref{cor:30}. Tuttavia, è pur vero che il Corollario~\ref{cor:30} fornisce una interpretazione naturale della dualità tra gli spazi $\widetilde A$, che verrebbe altrimenti a mancare. Forse un buon compromesso sarebbe mettere qui una dimostrazione diretta e trasformare il Corollario~\ref{cor:30} in un risultato sulla dualità tra gli spazi $\widetilde A$.

%We shall extend the second assertion of the following result to the case of general $p\in ]0,\infty]$ in Corollary~\ref{cor:30}. \CR Nonetheless, it is possible to give a direct proof of this result.\CB

\begin{teo}\label{teo:10}
Take $p,q\in ]0,\infty]$ and $\vect s\in \R^r$. Then, the following hold:
\begin{enumerate}
\item[\em(1)] $\Sc_{\Omega,L}(\Nc)$ is dense in $\mathring{B}^{\vect s}_{p,q}(\Nc,\Omega) $;

\item[\em(2)]  the sesquilinear form 
\[
\langle \,\cdot\,\vert \,\cdot\,\rangle\colon \mathring{B}^{\vect s}_{p,q}(\Nc,\Omega)\times B^{-\vect s-(1/p-1)_+(\vect b+\vect d)}_{p',q'}(\Nc,\Omega)\to \C
\]
induces a antilinear isomorphism of 
\[
B^{-\vect s-(1/p-1)_+(\vect b+\vect d)}_{p',q'}(\Nc,\Omega) \text{ onto } \mathring{B}^{\vect s}_{p,q}(\Nc,\Omega)'.
\]
\end{enumerate}
\end{teo}

Notice that the space $\Sc_{\Omega,L}(\Nc)$ is always dense in $B^{\vect s}_{p,q}(\Nc,\Omega)$ for the (Hausdorff) weak topology $\sigma^{\vect s}_{p,q}$, thanks to Proposition~\ref{prop:20}.

The proof is based on~\cite[Proposition 3.27]{BekolleBonamiGarrigosRicci}, which deals with the case in which $p,q\Meg 1$ and $D$ is a symmetric tube domain.
\CR Before we pass to the proof, we need a technical lemma.

\begin{lem}\label{lem:90}
	Take $\psi\in \Sc_\Omega(\Nc)$ and $p\in ]0,1]$. Define $\psi_t\coloneqq (g\times t)_* \psi$ for every $t\in T_+$, where $g\in GL(E)$ is such that $t\cdot \Phi=\Phi\circ (g\times g)$.\footnote{Notice that the definition of $\psi_t$ does not depend on $g$ by Proposition~\ref{prop:40}.} 
	Then, there is a constant $C>0$ such that, for every $T\in \Sc'(\Nc)$ and for every $t\in T_+$,
	\[
	\norm{T*\psi_t}_{L^\infty(\Nc)}\meg C\Delta^{(\vect b+\vect d)(1-1/p)}(t)\sup_{\substack{f\in \Sc_{\Omega,L}(\Nc)\\ \norm{f*\psi_t}_{L^p(\Nc)}=1}  } \abs*{\langle f\vert T*\psi_t\rangle}.
	\]
\end{lem}

\begin{proof}
	Let us first show that we may reduce to proving the statement for $t=e_{T_+}$, where $e_{T_+}$ is the identity element of $T_+$. Indeed, assume that the statement holds for $t=e_{T_+}$, in which case $\psi_{e_{T_+}}=\psi$. Observe that
	\[
	\langle f\vert T*\psi_t\rangle=\langle f\circ (g\times t)\vert [(g\times t)^* T]*\psi\rangle
	\]
	and that
	\[
	\norm{f*\psi_t}_{L^p(\Nc)}=\norm*{[[f\circ (g\times t)]*\psi]\circ (g^{-1}\times t^{-1})}_{L^p(\Nc)}=\Delta^{-(\vect b+\vect d)/p}(t)\norm{[f\circ (g\times t)]*\psi}_{L^p(\Nc)}
	\]
	for every $f\in \Sc_{\Omega,L}(\Nc)$.
	Then, 
	\[
	\sup_{\substack{f\in \Sc_{\Omega,L}(\Nc)\\ \norm{f*\psi_t}_{L^p(\Nc)}=1}  } \abs*{\langle f\vert T*\psi_t\rangle}=\Delta^{(\vect b+\vect d)/p}(t)\sup_{\substack{f\in \Sc_{\Omega,L}(\Nc)\\ \norm{f*\psi}_{L^p(\Nc)}=1}  } \abs*{\langle f\vert [(g\times t)^* T]*\psi\rangle},
	\]
	so that
	\[
	\norm{T*\psi_t}_{L^\infty(\Nc)}=\Delta^{\vect b+\vect d}(t)\norm{[(g\times t)^* T]*\psi}_{L^\infty(\Nc)}\meg C\Delta^{(\vect b+\vect d)(1-1/p)}(t)\sup_{\substack{f\in \Sc_{\Omega,L}(\Nc)\\ \norm{f*\psi_t}_{L^p(\Nc)}=1}  } \abs*{\langle f\vert T*\psi_t\rangle}.
	\]
	Now, take $\psi'\in \Sc_\Omega(\Nc)$ so that $\psi=\psi*\psi'^*=\psi*\psi'$, and fix a non-zero positive $f\in C^\infty_c(\Nc)$ with integral $1$; define $f_j\coloneqq 2^{-2(n+m)} f(2^{-j}\,\cdot\,)$ for every $j\in \N$, so that $f_j\to \delta_e$ vaguely. Then, $\delta_{(\zeta,x)}*f_j*\psi'\in \Sc_{\Omega,L}(\Nc)$,
	\[
	\lim_{j\to \infty}\abs*{\langle \delta_{(\zeta,x)}*f_j*\psi'\vert T*\psi\rangle }=\lim_{j\to \infty}\abs*{\langle \delta_{(\zeta,x)}*f_j\vert T*\psi\rangle }=\abs{(T*\psi)(\zeta,x)}
	\]
	and 
	\[
	C\coloneqq \sup_{j\in \N}\norm{\delta_{(\zeta,x)}*f_j*\psi'*\psi}_{L^p(\Nc)}=\sup_{j\in \N}\norm{f_j*\psi}_{L^p(\Nc)}<\infty
	\]
	for every $T\in \Sc'(\Nc)$, for every $(\zeta,x)\in \Nc$ and for every $j\in \N$. Therefore,
	\[
	\norm{T*\psi}_{L^\infty(\Nc)}\meg C\sup_{\substack{f\in \Sc_{\Omega,L}(\Nc)\\ \norm{f*\psi}_{L^p(\Nc)}=1}  } \abs*{\langle f\vert T*\psi\rangle}
	\]
	for every $T\in \Sc'(\Nc)$.
\end{proof}
\CB

\begin{proof}[\CR Proof of Theorem~\ref{teo:10}.]
{(1)} Take $T\in \mathring{B}^{\vect s}_{p,q}(\Nc,\Omega)$, and let us prove that $T$ may be approximated by elements of $\Sc_{\Omega,L}(\Nc)$. 
Indeed,  take $(\lambda_k)_{k\in K}$, $(\varphi_k)$, and $(\psi_k)$ as in Lemma~\ref{lem:59}, so that
\[
T=\sum_{k\in K} T*\psi_k
\]
in $\mathring{B}^{\vect s}_{p,q}(\Nc,\Omega)$. Hence, we may reduce to proving that $T*\psi_k$ may be approximated by elements of $ \Sc_{\Omega,L}(\Nc)$ in $B^{\vect s}_{p,q}(\Nc,\Omega)$ for every $k\in K$. 

Then, let $H$ be a convex compact neighbourhood of $0$ in $\overline{\Omega'}$, and take $\eta\in \Sc_\Omega(\Nc)$ such that $\eta(e)=1$ and such that $\Fc_\Nc\eta$ is supported in $H$. Then,
\[
\eta_\rho\colon (\zeta,x)\mapsto \eta(\rho\zeta,\rho ^2 x) 
\]
is an element of $\Sc_\Omega(\Nc)$  and 
\[
\Supp{\Fc_\Nc \eta_\rho}\subseteq H
\]
for every $\rho\in ]0,1]$. 
In addition, 
\[
(T*\psi_k)\eta_\rho\in\Sc_{\Omega,L}(\Nc)
\]
for every $\rho\in ]0,1]$ thanks to Corollary~\ref{cor:21}. 
Hence, it will suffice to show that $(T*\psi_k)\eta_\rho$ converges to $T*\psi_k$ in $B^{\vect s}_{p,q}(\Nc,\Omega)$  as $\rho\to 0^+$. Observe that there is $\psi'_k\in \Sc_\Omega(\Nc)$ such that $\Fc_\Nc\psi_k$ equals $1$ on $\Supp{\Fc_\Nc \psi_k}+ K $, so that Corollary~\ref{cor:21} shows that 
\[
(T*\psi_k)\eta_\rho=[(T*\psi_k)\eta_\rho]*\psi'_k 
\]
for every $\rho\in [0,1]$. Here, $\eta_0(\zeta,x)=\eta(e)=1$  for every $(\zeta,x)\in \Nc$. 
Since $(T*\psi_k)\eta_\rho$ converges to $T*\psi_k$ in $L^p(\Nc)$ as $\rho\to 0^+$, by dominated convergence, we have
\[
\lim_{\rho\to 0^+} [(T*\psi_k)\eta_\rho]*\psi'_k =(T*\psi_k)*\psi'_k =T*\psi_k
\]
in $L^p(\Nc)$. It is then easily seen that $(T*\psi_k)\eta_\rho$ converges to $T*\psi_k$ in $B^{\vect s}_{p,q}(\Nc,\Omega)$, whence the result.

\CR
{(2)} Take $L\in \mathring{B}^{\vect s}_{p,q}(\Nc,\Omega)'$, and observe that the preceding remarks show that there is a unique $T\in \Sc_{\Omega,L}'(\Nc)$ such that
\[
\langle L, \eta\rangle= \langle\eta\vert T\rangle
\]
for every $\eta\in \Sc_{\Omega,L}(\Nc)$. Keep the notation of (1) and observe that, if $K'$ is a finite subset of $K$ and $(f_k)\in \Sc_{\Omega,L}(\Nc)^{K'}$, then
\[
\sum_{k\in K'}\langle f_k\vert T*\psi_k\rangle= \Big\langle \sum_{k\in K' } f_k*\psi_k\Big\vert T\Big\rangle= \Big\langle L, \sum_{k\in K' } f_k*\psi_k\Big\rangle,
\]
so that
\[
\abs*{\sum_{k\in K'}\langle f_k\vert T*\psi_k\rangle}\meg \norm{L}_{\mathring{B}^{\vect s}_{p,q}(\Nc,\Omega)'} \norm*{\sum_{k\in K' } f_k*\psi_k}_{\mathring{B}^{\vect s}_{p,q}(\Nc,\Omega)}.
\]
For every $k\in K$, define $K_k$ as the set of $k'\in K$ such that $\psi_k*\psi_{k'}\neq 0$, and observe that Proposition~\ref{prop:56} implies that there is  $N\in \N$ such that $\card(K_k)\meg N$ for every $k\in K$.
Now, assume that $p\Meg 1$,   set $C_1\coloneqq \sup\limits_{k\in K}\norm{\psi_k}_{L^1(\Nc)}$, and observe that  Young's inequality implies that
\[
\norm*{\sum_{k\in K' } f_k*\psi_k*\psi_{k'}}_{L^p(\Nc)}\meg C_1^2 \sum_{k\in K'\cap K_{k'}}\norm{f_k}_{L^p(\Nc)}
\]
for every $k'\in K$. Taking Corollary~\ref{cor:34} into account, we then see that there is a constant $C_2>0$ such that
\[
\abs*{\sum_{k\in K'}\langle f_k\vert T*\psi_k\rangle}\meg C_1^2 C_2 N^{1/q+1/q'} \norm{L}_{B^{\vect s}_{p,q}(\Nc,\Omega)'} \norm*{\Delta^{\vect s}_{\Omega'}(\lambda_k) \norm{f_k}_{L^p(\Nc)}}_{\ell^q(K')}.
\]
By the arbitrariness of $(f_k)$ and $K'$,~\cite[Theorem 1 of \S 2]{BenedekPanzone} shows that $(\Delta^{-\vect s}_{\Omega'}(\lambda_k) T*\psi_k)\in L^{q'}(K; L^{p'}(\Nc))$ when $q\Meg 1$. The case $q<1$ is easier: it suffice to consider the case $\card(K')=1$ and to apply the duality between $L^p(\Nc)$ and $L^{p'}(\Nc)$. 

Conversely, assume that $p<1$, and observe that Corollary~\ref{cor:24} shows that there is a constant $C_3>0$ such that
\[
\norm*{\sum_{k\in K' } f_k*\psi_k*\psi_{k'}}_{L^p(\Nc)}\meg C_3 \sum_{k\in K'\cap K_{k'}}\norm{f_k*\psi_k}_{L^p(\Nc)}
\]
for every $k'\in K$. Taking Corollary~\ref{cor:34} into account, we then see that there is a constant $C_4>0$ such that
\[
\abs*{\sum_{k\in K'}\langle f_k\vert T*\psi_k\rangle}\meg C_3 C_4 N^{1/q+1/q'} \norm{L}_{B^{\vect s}_{p,q}(\Nc,\Omega)'} \norm*{\Delta^{\vect s}_{\Omega'}(\lambda_k) \norm{f_k*\psi_k}_{L^p(\Nc)}}_{\ell^q(K')}.
\]
If $q\meg 1$, take $K$ so that $\card(K')=1$ and apply Lemma~\ref{lem:90} to obtain that  $\big(\Delta^{-\vect s-(1/p-1)_+(\vect b+\vect d)}_{\Omega'}(\lambda_k) T*\psi_k\big)\in L^{\infty}(K; L^{\infty}(\Nc))$. If, otherwise, $q>1$, then use Lemma~\ref{lem:90} and argue as in the proof of~\cite[Theorem 1 of \S 2]{BenedekPanzone} to obtain $\big(\Delta^{-\vect s-(1/p-1)_+(\vect b+\vect d)}_{\Omega'}(\lambda_k) T*\psi_k\big)\in L^{q'}(K; L^{\infty}(\Nc))$.

Thus, we have proved that $\big(\Delta^{-\vect s-(1/p-1)_+(\vect b+\vect d)}_{\Omega'}(\lambda_k) T*\psi_k\big)\in L^{q'}(K; L^{p'}(\Nc))$, that is, 
\[
T\in B^{-\vect s-(1/p-1)_+(\vect b+\vect d)}_{p',q'}(\Nc,\Omega).
\]
By Proposition~\ref{prop:18} and the density of $\Sc_{\Omega,L}(\Nc)$ in $\mathring{B}^{\vect s}_{p,q}(\Nc,\Omega)$, we then see that
\[
\langle L, T'\rangle= \langle T'\vert T\rangle
\]
for every $T'\in \mathring{B}^{\vect s}_{p,q}(\Nc,\Omega)$.\CB
\end{proof}

\begin{cor}\label{cor:28}
Take $p,q\in ]0,\infty]$ and $\vect s\in \R^r$. Then, the bounded subsets of $B^{\vect s}_{p,q}(\Nc,\Omega)$ are relatively compact for the weak topology $\sigma^{\vect s}_{p,q}$.
\end{cor}

\begin{proof}
Take $(\lambda_k)_{k\in K}$, $(\varphi_k)$, and $(\psi_k)$ as in Lemma~\ref{lem:51}.
It will suffice to prove that the set $U$ of  $T\in B^{\vect s}_{p,q}(\Nc,\Omega)$ such that
\[
\norm*{\Delta^{\vect s}_{\Omega'}(\lambda_k) \norm{T*\psi_k}_{L^p(\Nc)}}_{\ell^q(K)}\meg 1
\]
is compact for the weak topology $\sigma^{\vect s}_{p,q}$. Then, let $\Uf$ be an ultrafilter on $U$. By Proposition~\ref{prop:20} (and its proof) $U$ is closed and bounded in $\Sc_{\Omega,L}'(\Nc)$, hence  compact thanks to Proposition~\ref{prop:63}. Therefore, $\Uf$ converges to some $T_0\in U$ in $ \Sc_{\Omega,L}'(\Nc)$. Since $\Sc_{\Omega,L}(\Nc)$ is dense in $\mathring B_{p',q'}^{-\vect s-(1/p-1)_+(\vect b+\vect d)}(\Nc,\Omega)$, it is then easily verified that $\Uf$ converges to $T_0$ in the weak topology $\sigma^{\vect s}_{p,q}$.
\end{proof}

The following result shows that the Riemann--Liouville operators play for the spaces $B^{\vect s}_{p,q}(\Nc,\Omega)$ the same role that the fractional powers of the Laplacian play for the classical homogeneous Besov spaces.

\begin{teo}\label{teo:9}
Take $p,q\in ]0,\infty]$ and $\vect s\in \R^r$ and $\vect{s'}\in \C^r$. Then, convolution by $I^{\vect{s'}}_\Omega$ on $\Sc_{\Omega,L}(\Nc)$  induces a unique isomorphism of 
\[
\mathring{B}^{\vect s}_{p,q}(\Nc,\Omega)\text{ onto }\mathring{B}^{\vect s+\Re\vect{s'}}_{p,q}(\Nc,\Omega)
\]
and a unique isomorphism of 
\[
B^{\vect s}_{p,q}(\Nc,\Omega) \text{ onto } B^{\vect s+\Re\vect{s'}}_{p,q}(\Nc,\Omega)
\]
which is also continuous for the weak topologies $\sigma^{\vect s}_{p,q}$ and $\sigma^{\vect s+\Re \vect{s'}}_{p,q}$. 
\end{teo}

\begin{proof}
Take $(\lambda_k)_{k\in K}$,  $(\varphi_k)$, and  $(\psi_k)$ as in Lemma~\ref{lem:51}.
Observe that Proposition~\ref{prop:55} implies that $ \Delta_{\Omega'}^{\vect{s'}}(\lambda_k)\psi_k* I^{\vect{s'}}_\Omega \in \Sc_\Omega(\Nc)$ and that 
\[
\Fc_\Nc(\Delta_{\Omega'}^{\vect{s'}}(\lambda_k)\psi_k* I^{\vect{s'}}_\Omega)= i^{-\vect{s'}}(\Delta_{\Omega'}^{-\vect{s'}}\varphi_k) (\,\cdot\,t_k^{-1}),
\]
for every $k\in K$. 

Then, choose $\varphi''\in C^\infty_c(\Omega')$ so that $\varphi_k=\varphi_k \varphi''$ for every $k\in K$, and define $\psi'_k$ and $\psi''_k$ as the elements of $\Sc_\Omega(\Nc)$  such that  
\[
\Fc_\Nc \psi'_k=(\Delta_{\Omega'}^{-\Re\vect{s'}}\varphi_k) (\,\cdot\,t_k^{-1}) \qquad \text{and} \qquad \Fc_\Nc \psi''_k= i^{-\vect{s'}}(\Delta_{\Omega'}^{-\Im\vect{s'}}\varphi'') (\,\cdot\,t_k^{-1}),
\]
so that
\[
\psi_k*I^{\vect{s'}}_\Omega=\Delta_{\Omega'}^{-\vect{s'}}(\lambda_k) \psi'_k* \psi''_k
\]
for every $k\in K$. 
Now, Corollary~\ref{cor:24} implies that there is a constant $C>0$ such that
\[
\norm{\eta*I^{\vect{s'}}_\Omega*\psi_k }_{L^p(\Nc)}\meg C \Delta_{\Omega'}^{-\Re\vect{s'}}(\lambda_k) \norm{ \eta*\psi'_k }_{L^p(\Nc)}
\]
for every $\eta\in \Sc_{\Omega,L}(\Nc)$ and for every $k\in K$, whence
\[
\norm*{\Delta_{\Omega'}^{\vect s+\Re\vect{s'}}(\lambda_k)\norm{  \eta*I^{\vect{s'}}_\Omega*\psi_k}_{L^p(\Nc)}  }_{\ell^q(K)}\meg C \norm*{ \Delta_{\Omega'}^{\vect s}(\lambda_k)\norm{  \eta*\psi'_k}_{L^p(\Nc)}  }_{\ell^q(K)}
\]
for every $\eta\in \Sc_{\Omega,L}(\Nc)$. Thus, Lemma~\ref{lem:51} and Theorem~\ref{teo:10} imply that the automorphism $\eta\mapsto \eta* (i^{\vect s}I_{\Omega}^{\vect s})$ of $\Sc_{\Omega,L}(\Nc)$ extends to an continuous linear mapping 
\[
\Ic^{\vect s,\vect{s'}}_{p,q}\colon\mathring{B}^{\vect s}_{p,q}(\Nc,\Omega)\to\mathring{B}^{\vect s+\Re\vect{s'}}_{p,q}(\Nc,\Omega).
\] 
It is then easily that $\Ic^{\vect s,\vect{s'}}_{p,q}$ is actually an isomorphism with inverse $\Ic^{\vect s+\vect{s'},-\vect {s'}}_{p,q}$.

Now, the transpose of $\Ic^{-\vect s-(1/p-1)_+(\vect b+\vect d),\overline{\vect{s'}}}_{p',q'}$ with respect to the sesquilinear forms 
\[
\langle\,\cdot\,\vert \,\cdot\,\rangle\colon {B}^{\vect s+\Re\vect{s'}}_{p,q}(\Nc,\Omega)\times B^{-\vect s-\Re\vect{s'}-(1/p-1)_+(\vect b+\vect d)}_{p',q'}(\Nc,\Omega)
\]
and
\[
\langle\,\cdot\,\vert \,\cdot\,\rangle\colon B^{\vect s}_{p,q}(\Nc,\Omega)\times \mathring B^{-\vect s-(1/p-1)_+(\vect b+\vect d)}_{p',q'}(\Nc,\Omega)
\]
equals  $ \Ic^{\vect s,\vect{s'}}_{p,q}$ on $\Sc_{\Omega,L}(\Nc)$ by Proposition~\ref{prop:55}, whence the result since $ \Sc_{\Omega,L}(\Nc)$ is  clearly  dense in $B^{\vect s}_{p,q}(\Nc,\Omega)$ for the weak topology $\sigma^{\vect s}_{p,q}$. 
\end{proof}

\section{Notes and Further Results}

\paragraph{4.3.1}  When $p,q\in [1,\infty]$, it is possible to give a `continuous' characterization of the spaces $B^{\vect s}_{p,q}(\Nc,\Omega)$. Indeed, if $\Omega'\ni\lambda \mapsto \varphi_\lambda\in C^\infty_c(\Omega')$ is a \emph{bounded} measurable mapping such that each $\varphi_\lambda$ is positive and, choosing $t_\lambda\in T_+$ so that $\lambda= e_{\Omega'}\cdot t_\lambda$,
\[
\int_{\Omega'} \varphi_{\lambda}(\lambda'\cdot t_\lambda^{-1})\,\dd \nu_{\Omega'}(\lambda)\Meg 1
\] 
for every $\lambda'\in \Omega'$, and if $\psi_{\lambda}\in \Sc_\Omega(\Nc)$ is defined so that $\Fc_\Nc(\psi_{\lambda})= \varphi_\lambda(\,\cdot\, t_\lambda^{-1})$, then $T\in B^{\vect s}_{p,q}(\Nc,\Omega)$ if and only if
\[
\left( \int_{\Omega'} \left(\Delta_{\Omega'}^{\vect s}(\lambda ) \norm{T*\psi_{ \lambda}}_{L^p(\Nc)}  \right)^q\,\dd \nu_{\Omega'}(\lambda)\right) ^{1/q}<\infty
\]
(modification if $q=\infty$). The preceding expression then gives an equivalent norm on $B^{\vect s}_{p,q}(\Nc,\Omega)$. 
Several results of this chapter can be modified in this terms.

\paragraph{4.3.2}  Several aspects of the classical theory of Besov spaces have not been treated in this chapter.
Some of them can be (at least partially) addressed by means of the results of Chapter~\ref{sec:6}.
For example, thanks to Theorem~\ref{teo:9} and Corollary~\ref{cor:29}, we see that, if $p_1,p_2,p_3\in [1,\infty]$, $q_1,q_2,q_3\in ]0,\infty]$, $\vect{s_1}, \vect{s_2}\in \R^r$,
\[
\frac{1}{p_1'}+\frac{1}{p_2'}=\frac{1}{p_3'}, \qquad \text{and} \qquad  \frac{1}{q_1}+\frac{1}{q_2}=\frac{1}{q_3},
\]
then (by means of Young's inequality) it is possible to give a reasonable definition of convolution which induces a continuous bilinear mapping
\[
B^{\vect{s_1}}_{p_1,q_1}(\Nc,\Omega)\times B^{\vect{s_2}}_{p_2,q_2}(\Nc,\Omega)\to B^{\vect{s_1}+\vect{s_2}}_{p_3,q_3}(\Nc,\Omega).
\]
One may also argue directly and extend the preceding discussion to the case of general $p_1,p_2,p_3$.

\paragraph{4.3.3}  Using Corollary~\ref{cor:29},  one may  prove that, if $p_1,p_2,p_3,q_1,q_2,q_3\in ]0,\infty]$, $\vect{s_1}, \vect{s_2}\in \R^r$,
\[
\frac{1}{p_1}+\frac{1}{p_2}=\frac{1}{p_3}, \qquad \text{and} \qquad  \frac{1}{q_1}+\frac{1}{q_2}=\frac{1}{q_3},
\]
and $\vect{s_1},\vect{s_2}$ belong to a suitable translate of $-\R_+^r$, then it is possible to give a reasonable definition of multiplication which induces a continuous bilinear mapping
\[
B^{\vect{s_1}}_{p_1,q_1}(\Nc,\Omega)\times B^{\vect{s_2}}_{p_2,q_2}(\Nc,\Omega)\to B^{\vect{s_1}+\vect{s_2}}_{p_3,q_3}(\Nc,\Omega).
\]
In this case, though, one cannot make use of Corollary~\ref{cor:29} to extend the range of $\vect{s_1}$ and $\vect{s_2}$ for which the preceding assertion holds. Nonetheless, arguing directly more precise results can be obtained.

\paragraph{4.3.4}  As in the classical case, one may investigate the validity of multiplier theorems for Besov spaces. 
Cf.~\cite[Proposition 3.34 and Remark 3.37]{BekolleBonamiGarrigosRicci} for the case of symmetric tube domains (and $p,q\in [1,\infty[$), where such multiplier theorems can be extended in a similar form. 
The case of (homogeneous) Siegel domains of type II is more subtle. On the one hand, one may investigate the spectral multipliers associated with finite (necessarily commutative) families of Riemann--Liouville operators, and find results which are analogous to those presented in~\cite[Proposition 3.34 and Remark 3.37]{BekolleBonamiGarrigosRicci}. 
On the other hand, general Fourier multipliers can be investigated. Nonetheless, since in this case the Fourier transform is vector-valued, the situation is  more delicate.

\paragraph{4.3.5}  It is well known that classical Besov spaces enjoy very useful (real and complex) interpolation properties. On the one hand, complex interpolation (which is meaningful only when $p,q\in [1,\infty]$) can be treated as in the classical case, since it is not hard to prove that the space $B^{\vect s}_{p,q}(\Nc,\Omega)$ is a retract of the space $(\Delta_{\Omega'}^{\vect s}(\lambda_k)) L^q(K;L^p(\Nc))$ (the retractions being independent of $\vect s$, $p$, and $q$). Cf.~\cite[Theorem 6.4.5]{BerghLofstrom} for the classical case (which applies to the case $p,q\in ]1,\infty[$).
One may also consider the variant of the complex method considered, e.g., in~\cite{Triebel} to treat the general case $p,q\in]0,\infty]$, and investigate whether it can be generalized to the spaces $B^{\vect s}_{p,q}(\Nc,\Omega)$.

On the other hand, real interpolation seems to be more problematic. Indeed, the classical case heavily relies on the interplay between integration on $\R_+^*$ with respect to the invariant measure $\nu_{\R_+^*}$ (related to the real method) and a dyadic summation over $\N$ or $\Z$ (related to the definition of non-homogeneous and homogeneous Besov spaces, respectively). When the rank $r$ of the cone $\Omega$ is $>1$, though, the situation becomes more delicate and the na\"ive generalization of the proof of, e.g.,~\cite[Theorem 2.4.2]{Triebel} fails. 
On the contrary, it seems that a suitable modification of the real method, where the role of $\R_+^*$ is suitably replaced by $\Omega'$, should prove more fruitful to extend the results of the classical case to this context.

\paragraph{4.3.6}  Classical Besov spaces enjoy several equivalent definitions, e.g., by means of finite differences or atomic decomposition. In connection to these remarks, we mention that in~\cite{RicciTaibleson} the Besov spaces $B^s_{p,q}(\R,\R_+)$ were defined by means of their atomic decomposition. Cf.~\cite{Christensen1} for a study of atomic decomposition on the spaces $B^{s\vect{1}_r}_{p,q}(\Nc,\Omega)$ in the case in which $p,q\in [1,\infty[$, $E=\Set{0}$, and $\Omega$ is symmetric.

\paragraph{4.3.7}  The Besov spaces considered in this chapter were not proved to embed canonically into a suitable quotient of the space $\Sc'(\Nc)$ of tempered distributions, but rather in the projective limit $\Sc_{\Omega,L}'(\Nc)$.
In order to address this problem, denote by $\widetilde \Sc_{\Omega,L}'(\Nc)$ the strong dual of the closure of $\Sc_{\Omega,L}(\Nc)$ in $\Sc(\Nc)$ (endowed with the topology induced by $\Sc(\Nc)$).
Observe that $B^{\vect 0}_{2,2}(\Nc,\Omega)$ can be canonically identified with the space of $f\in L^2(\Nc)$ such that $\pi_\lambda(f)=\chi_{\Omega'}(\lambda) \pi_\lambda(f) P_{\lambda,0}$ for almost every $\lambda\in F'\setminus W$, thanks to Corollary~\ref{cor:7}. Therefore, $B^{\vect 0}_{2,2}(\Nc,\Omega)$ embeds canonically into  $\widetilde \Sc_{\Omega,L}'(\Nc)$. Since convolution by the Riemann--Liouville operators induce automorphisms of $\widetilde \Sc_{\Omega,L}'(\Nc)$ by Proposition~\ref{prop:55}, by means of  Proposition~\ref{prop:27} and Theorem~\ref{teo:9} we see that $B^{\vect s}_{p,q}(\Nc,\Omega)$ embeds canonically into $\widetilde \Sc_{\Omega,L}(\Nc)$ for every $p,q\in ]0,2]$ and for every $\vect s\in \R^r$. In~\cite{BekolleBonamiGarrigosRicci}, the same assertion is proved for $p,q\in [1,\infty[$ and for $\vect s\in \R 1_r $, provided that $E=\Set{0}$ and $\Omega$ is an irreducible symmetric cone.
We shall not pursue this investigation any further.

\chapter[Weighted Bergman Spaces II]{Weighted Bergman Spaces: Boundary Values and Bergman Projectors}\label{sec:7}

In this chapter we develop further the theory of weighted Bergman spaces on homogeneous Siegel domains of type II. We recall that the provisional notation of Sections~\ref{sec:1:3} and~\ref{sec:1:4} will no longer be used.

In Section~\ref{sec:6:6}, we define an extension operator 
\[
\Ec \colon B^{-\vect s}_{p,q}(\Nc,\Omega)\to A^{\infty,\infty}_{\vect s-(\vect b+\vect d)/p}(D)
\]
(Theorem~\ref{teo:11}) and prove that $A^{p,q}_{\vect s}(D)\subseteq \Ec(B^{-\vect s}_{p,q}(\Nc,\Omega))$ (Proposition~\ref{prop:23}). We then provide sufficient conditions for the equality $A^{p,q}_{\vect s}(D)= \Ec(B^{-\vect s}_{p,q}(\Nc,\Omega))$ (Theorem~\ref{teo:6} and Corollary~\ref{cor:29}).

By means of these results we are able to \CR provide a new interpretation of the sesquilinear dual pairing between $B^{\vect s}_{p,q}(\Nc,\Omega)$ and $B^{-\vect s-(1/p-1)_+(\vect b+\vect d)}_{p',q'}(\Nc,\Omega)$  (Proposition~\ref{cor:30}) \CB and to prove the equivalence of several properties concerning the characterization of the boundary values of $A^{p,q}_{\vect s}(D)$, the validity of properties $\atomic^{p,q}_{\vect s}$ and $\atomics^{p,q}_{\vect s}$, and the fact that the Riemann--Liouville operators $I^{-\vect{s'}}_\Omega$ induce isomorphisms between various weighted Bergman spaces (Corollary~\ref{cor:23}).

Finally, in Section~\ref{sec:6:7} we shall study the boundedness of the Bergman projectors $P_{\vect s}$. As for atomic decomposition, it is sometimes convenient to study also the boundedness of the operator $P_{\vect s,+}$, whose kernel is the absolute value of the kernel of $P_{\vect s}$. Since the boundedness of $P_{\vect{s'},+}$ on $A^{p,q}_{\vect s}(D)$, for $p,q\in [1,\infty]$, turns out to be equivalent to property $\atomic^{p,q}_{\vect s,+}$ (Proposition~\ref{prop:31}), no new sufficient conditions are provided. Instead, we prove that, if $\Ec$ induces isomorphisms of $\mathring B^{-\vect s}_{p,q}(\Nc,\Omega)$ onto  $A^{p,q}_{\vect s,0}(D)$ and of $B^{\vect s+\vect{s'}-(\vect b+\vect d)}_{p',q'}(\Nc,\Omega)$ onto $A^{p',q'}_{\vect b+\vect d-\vect s-\vect {s'}}(D)$, then $P_{\vect{s'}}$ induces an endomorphism of $L^{p,q}_{\vect s,0}(D)$ (cf.~Definition~\ref{26}) with image $A^{p,q}_{\vect s,0}(D)$.

\medskip

We keep the hypotheses and the notation of Chapters~\ref{sec:6} and~\ref{sec:5}.

\section{Boundary Values}\label{sec:6:6}

This section deals with the limits $\lim\limits_{h\to 0} f_h$, for $f\in A^{p,q}_{\vect s}(D)$. As we shall see (cf.~Proposition~\ref{prop:23}), under very general assumptions these limits exist and belong to the Besov space $B^{-\vect s}_{p,q}(\Nc,\Omega)$. Conversely, under the same assumptions every element of $B^{-\vect s}_{p,q}(\Nc,\Omega)$ is the `boundary value' of a unique element of $A^{p,q}_{\vect s-(\vect b+\vect d)/p}(D)$. 
We shall then define $\widetilde A^{p,q}_{\vect s}(D)$ as the space of such holomorphic extensions, and provide sufficient conditions for the equality $A^{p,q}_{\vect s}(D)=\widetilde A^{p,q}_{\vect s}(D)$ (cf.~Theorem~\ref{teo:6} and Corollary~\ref{cor:29}). 

We shall then draw some consequences of the interplay between weighted Bergman spaces and Besov spaces.
\CR %On the one hand, we shall use the results of Section~\ref{sec:6:5} to characterize the dual of $B^{\vect s}_{p,q}(\Nc,\Omega)$ in full generality (cf.~Corollary~\ref{cor:30}).
In particular, we shall use the results of Chapter~\ref{sec:5}, translated in terms of the spaces $\widetilde A^{p,q}_{\vect s}(D)$ (cf.~Proposition~\ref{prop:30}), to show the equivalence (in a somewhat weak sense) of several relevant properties of weighted Bergman spaces (cf.~Corollary~\ref{cor:32}).\CB

\medskip

Notice that, unlike in the case of tube domains (cf.~\cite{BekolleBonamiGarrigosRicci}), we may not define the holomorphic extension of the elements of $B^{-\vect s}_{p,q}(\Nc,\Omega)$ composing the Laplace transform with the inverse Fourier transform, since the inverse Fourier transform is not  readily  available for general tempered distributions. Nonetheless, it is possible to give a direct pointwise characterization of such extension by means of the `boundary value' of the Cauchy--Szeg\H o kernel. In order to do this, we need the following lemma.

\begin{lem}\label{lem:61}
Take $p,q\in [1,\infty]$, and $\vect s\in \R^r$ such that $\vect s\in\frac{1}{p'}(\vect b+\vect d)+\frac{1}{2 q}\vect{m'}+(\R_+^*)^r$. Define, for every $(\zeta,z)\in D$ and for every $(\zeta',x')\in \Nc$,
\[
\begin{split}
\Szego_{(\zeta,z)}(\zeta',x')&\coloneqq \frac{ \abs{\Pfaff(e_{\Omega'})}\Gamma_{\Omega'}(-\vect b-\vect d)  }{4^m\pi^{n+m}}\left( B^{\vect b+\vect d}_{(\zeta,z)}\right)_0(\zeta',x').
\end{split}
\]
Then, the following hold:
\begin{enumerate}
\item[\em(1)] $\pi_\lambda(\Szego_{(\zeta,z)})=\chi_{\Omega'}(\lambda)\ee^{-\langle \lambda, \Im z-\Phi(\zeta)\rangle} \pi_\lambda(\zeta,\Re z)  P_{\lambda,0}$ for every $\lambda\in F'\setminus W$;

\item[\em(2)] the family $\big(\Delta_\Omega^{\vect s-(\vect b+\vect d)/p'}(\Im z-\Phi(\zeta))\Szego_{(\zeta,z)}\big)_{(\zeta,z)\in D}$ is bounded in $\mathring{B}^{\vect s}_{p,q}(\Nc,\Omega)$;

\item[\em(3)] the mapping $D\ni (\zeta,z)\mapsto \langle T\vert \Szego_{(\zeta,z)}\rangle \in \C$ is holomorphic for every $T\in  {B}^{-\vect s}_{p',q'}(\Nc,\Omega)$.
\end{enumerate}
\end{lem}

The proof is based on~\cite[Proposition 3.43]{BekolleBonamiGarrigosRicci}, which deals with the case in which $p,q\Meg 1$, $\vect s\in \R\vect 1_r$, and $D$ is a symmetric tube domain.

Observe that $ \Szego_{(\zeta,z)}$ is the `boundary value' of the Cauchy--Szeg\H o kernel (cf.~Corollary~\ref{cor:4} and Propositions~\ref{prop:58} and~\ref{prop:59}). Then, Corollary~\ref{cor:1} shows that, for every $f\in A^{2,\infty}_{\vect 0}(D)$, 
\[
f(\zeta,z)=\langle f_0\vert \Szego_{(\zeta,z)}\rangle
\]
for every $(\zeta,z)\in D$, where $f_0=\lim\limits_{h\to 0} f_h$ in $L^2(\Nc)$. 
As we shall see in Proposition~\ref{prop:23}, a similar statement holds for more general weighted Bergman spaces.

\begin{proof}
{(1)} Observe that $\Szego_{(\zeta,z)}\in L^2(\Nc)$ by Lemma~\ref{lem:21} and that clearly the measurable field 
\[
\lambda \mapsto \chi_{\Omega'}(\lambda) \ee^{-\langle \lambda, \Im z-\Phi(\zeta)\rangle} \pi_\lambda(\zeta,\Re z) P_{\lambda,0}
\]
belongs to $\frac{2^{n-m}}{ \pi^{n+m}} \int_{F'\setminus W}^\oplus \Lin^2(H_\lambda) \abs{\Pfaff(\lambda)}\dd \lambda$, with the notation of Corollary~\ref{cor:7}. 
Since clearly 
\[
\int_{\Omega'} \norm{\pi_\lambda(\zeta,\Re z) P_{\lambda,0}}_{\Lin^1(H_\lambda)} \ee^{-\langle \lambda, \Im z- \Phi(\zeta)\rangle} \Delta^{-\vect b}_{\Omega'}(\lambda)\,\dd \lambda=\int_{\Omega'} \ee^{-\langle \lambda, \Im z- \Phi(\zeta)\rangle} \Delta^{-\vect b}_{\Omega'}(\lambda)\,\dd \lambda
\]
is finite,  it will suffice to show that
\[
\begin{split}
&\Szego_{(\zeta,z)}(\zeta',x')\\
&\qquad=  \frac{2^{n-m} \abs{\Pfaff(e_{\Omega'})}}{\pi^{n+m}}\int_{\Omega'} \tr(\pi_\lambda(\zeta,\Re z) P_{\lambda,0} \pi_\lambda(-\zeta',-x')  )  \ee^{-\langle \lambda, \Im z-\Phi(\zeta)\rangle}\Delta_{\Omega'}^{-\vect b}(\lambda)\,\dd \lambda
\end{split}
\]
for every $(\zeta,z)\in D$. However, this follows from Propositions~\ref{prop:36}  and~\ref{prop:58}. 

{(2)} Choose $g\in GL(E)$ and $t\in T_+$ such that $t\cdot e_\Omega=\Im z-\Phi(\zeta)$ and $t\cdot \Phi=\Phi\circ (g\times g)$, and observe that
\[
\Szego_{(\zeta,z)}=L_{(\zeta,\Re z)}(g\times t)_* \Szego_{(0,i e_\Omega)},
\]
so that it will suffice to show that $ \Szego_{(0,i e_\Omega)}\in B^{\vect s}_{p,q}(\Nc,\Omega)$, thanks to Proposition~\ref{prop:42}.
Take $(\lambda_k)_{k\in K}$, $(t_k)$, $(\varphi_k)$, $(\psi_k)$ as in Lemma~\ref{lem:51}.
Then,~{(1)} implies that $\Szego_{(0,i e_\Omega)}*\psi_k\in\Sc_\Omega(\Nc)$  and 
\[
\Fc_\Nc(\Szego_{(0,i e_\Omega)}*\psi_k)=\varphi_k(\,\cdot\,t_k^{-1}) \ee^{-\langle \,\cdot\,,e_{\Omega}\rangle} 
\] for every $k\in K$. 
Now, Lemma~\ref{lem:48} shows that there is a constant $C\Meg 1$ such that
\[
\ee^{-C\langle \lambda_k, e_\Omega\rangle}\meg \ee^{-\langle  \lambda\cdot t_k,e_\Omega\rangle }\meg \ee^{-\langle \lambda_k, e_\Omega\rangle/C}
\]
for every $\lambda\in \Supp{\varphi_k}$ and for every $k\in K$, so that the family 
\[
\left( \ee^{\langle \lambda_k, e_\Omega\rangle/C}\varphi_k \ee^{-\langle \,\cdot\, t_k,e_{\Omega}\rangle}  \right)_{k\in K} 
\]
is bounded in $C^\infty_c(\Omega')$.
Hence, the family 
\[
\left(\ee^{\langle \lambda_k, e_\Omega\rangle/C} \Delta_{\Omega'}^{(\vect b+\vect d)/p'}(\lambda_k)\Szego_{(0,i e_\Omega)}*\psi_k  \right)_{k\in K}  
\]
is bounded in $L^p_0(\Nc)$, whence $\Szego_{(0, i e_\Omega)}\in \mathring{B}^{\vect s}_{p,q}(\Nc,\Omega)$ thanks to Proposition~\ref{prop:58}, Lemmas~\ref{lem:18} and~\ref{lem:48}, and Corollary~\ref{cor:34}.

{(3)} Take $T\in B^{-\vect s}_{p',q'}(\Nc,\Omega)$, keep the notation of~{(2)}, and assume further that 
\[
\sum_{k\in K}\varphi_k(\,\cdot\, t_k^{-1})^2=1
\]
on $\Omega'$. 
Let us first prove that the mapping 
\[
D\ni (\zeta,z) \mapsto \langle T*\psi_k\vert\Szego_{(\zeta,z)}*\psi_k \rangle\in \C
\]
is holomorphic for every $k\in K$.
Indeed, arguing as in~{(1)} we see that
\[
\overline{(\Szego_{(\zeta,z)}*\psi_k)(\zeta',x')}= \frac{2^{n-m} \abs{\Pfaff(e_{\Omega'})}  }{\pi^{n+m}}
\Fc_{\Omega'}(\varphi_k(\,\cdot\,t_k^{-1})\Delta_{\Omega'}^{-\vect b})(x'-z+2i \Phi(\zeta,\zeta') )
\]
for every $(\zeta,z)\in D$ and for every $(\zeta',x')\in \Nc$. 
The mapping 
\[
(\zeta,z)\mapsto \overline{(\Szego_{(\zeta,z)}*\psi_k)(\zeta',x')}
\]
is therefore holomorphic for every $(\zeta',x')\in \Nc$. In addition, observe as in the beginning of~{(2)} that
\[
\Szego_{(\zeta,z)}*\psi_k=L_{(\zeta,\Re z)} [[\Szego_{(0,i e_\Omega)} *(\psi_k\circ (g\times t))]\circ (g^{-1}\times t^{-1})],
\]
so that $ \Szego_{(\zeta,z)}*\psi_k$ stays in a bounded subset of $\Sc(\Nc)$  as long as $(\zeta,z)$ stays in a compact subset of $D$. In particular, as long as $(\zeta,z)$ stays in a compact subset of $D$, the function $\abs{\Szego_{(\zeta,z)}*\psi_k}$  is uniformly bounded by an element of $L^p(\Nc)$. Hence, the dominated convergence theorem shows that the mapping 
\[
D\ni (\zeta,z)\mapsto\overline{(\Szego_{(\zeta,z)}*\psi_k)} \in L^p(\Nc)
\]
is continuous. Then, by means of Morera's and Fubini's theorem, we see that the mapping 
\[
D\ni (\zeta,z)\mapsto\langle T*\psi_k\vert \Szego_{(\zeta,z)}*\psi_k\rangle \in L^p(\Nc)
\]
is actually holomorphic.

To conclude, it will suffice to show that the sum 
\[
\sum_{k\in K} \langle T*\psi_k\vert \Szego_{(\zeta,z)}*\psi_k\rangle
\]
converges (to $\langle T\vert \Szego_{(\zeta,z)}\rangle$) locally uniformly  in $(\zeta,z)\in D$.
However, this follows from the fact that the  families
\[
\big( \Delta_{\Omega'}^{\vect{s}}(\lambda_k) \norm{\Szego_{(\zeta,z)}*\psi_k}_{L^p(\Nc)} \big)_{k\in K} 
\]
are uniformly bounded by an element of  $\ell^q_0(K)$  as long as $(\zeta,z)$ stays in a compact subset of $D$,  thanks to~(2).
\end{proof}

We are now able to define the announced extension operator.

\begin{teo}\label{teo:11}
Take $p,q\in]0,\infty]$, and $\vect s\in-\frac{1}{p}(\vect b+\vect d)-\frac{1}{2 q'}\vect{m'}-(\R_+^*)^r$.
Define $\Szego_{(\zeta,z)}$, for every $(\zeta,z)\in D$, as in Lemma~\ref{lem:61}, and define
\[
\Ec\colon B^{\vect s}_{p,q}(\Nc,\Omega)\ni T \mapsto \big[(\zeta,z)\mapsto\langle T\vert \Szego_{(\zeta,z)}\rangle\big]\in A^{\infty,\infty}_{-\vect s-(\vect b+\vect d)/p}(D).
\]
Set $(\Ec T)_0\coloneqq T$ for every $T\in B^{\vect s}_{p,q}(\Nc,\Omega)$.
Then, the following hold:
\begin{enumerate}
\item[\em(1)] the linear mappings $T\mapsto (\Ec T)_h$, as $h$ runs through $\Omega$, induce equicontinuous endomorphisms of $B^{\vect s}_{p,q}(\Nc,\Omega)$;

\item[\em(2)] if $T\in \mathring{B}^{\vect s}_{p,q}(\Nc,\Omega)$, then the mapping 
\[
\Omega\cup \Set{0}\ni h\mapsto (\Ec T)_h\in\mathring{B}^{\vect s}_{p,q}(\Nc,\Omega)
\]
is continuous;

\item[\em(3)] if $T\in B^{\vect s}_{p,q}(\Nc,\Omega)$, then the mapping 
\[
\Omega\cup \Set{0}\ni h\mapsto (\Ec T)_h\in B^{\vect s}_{p,q}(\Nc,\Omega)
\]
is continuous for the weak topology $\sigma^{\vect s}_{p,q}$.
\end{enumerate}
\end{teo}

The proof is based on~\cite[Proposition 3.43]{BekolleBonamiGarrigosRicci}, which deals with the case in which $p,q\Meg 1$, $\vect s\in \R \vect 1_r$, and $D$ is an irreducible symmetric tube domain.

\begin{proof}
{(1)} Take $(\lambda_k)_{k\in K}$, $(t_k)$, $(\varphi_k)$, and $(\psi_k)$ as in Lemma~\ref{lem:59}.
Then, Lemmas~\ref{lem:59} and~\ref{lem:61} imply that
\[
\Szego_{(\zeta,z)}=\sum_{k\in K} \Szego_{(\zeta,z)}*\psi_k
\]
for every $(\zeta,z)\in D$, with convergence in $\mathring{B}^{-\vect s-(1/p-1)_+(\vect b+\vect d)}_{p',q'}(\Nc,\Omega)$. Therefore, Proposition~\ref{prop:18} shows that
\[
\langle T\vert\Szego_{(\zeta,z)}\rangle=\sum_{k\in K}\langle T\vert  \Szego_{(\zeta,z)}*\psi_k\rangle
\]
for every $(\zeta,z)\in D$. In addition, the sum converges locally uniformly  if the $\Szego_{(\zeta,z)}$ stay in a compact subset of $\mathring{B}^{-\vect s-(1/p-1)_+(\vect b+\vect d)}_{p',q'}(\Nc,\Omega)$, which is the case  if  $(\zeta,z)$ stays in a compact subset of $D$.
Furthermore, setting $h\coloneqq \Im z-\Phi(\zeta)$ to simplify the notation,  
\[
\Szego_{(\zeta,z)}*\psi_k=L_{(\zeta,\Re z)} (\Szego_{(0,i h)}*\psi_k),
\]
and 
\[
\Szego_{(0,i h)}*\psi_k\in \Sc_\Omega(\Nc) 
\]
for every $k\in K$, thanks to Lemma~\ref{lem:61}. Now, $\overline{\Szego_{(0,i h)}*\psi_k  }=(\Szego_{(0,i h)}*\psi_k){\check\;}$, so that
\[
\langle T\vert  \Szego_{(\zeta,z)}*\psi_k\rangle= (T*(\Szego_{(0,i h)}*\psi_k ))(\zeta,\Re z)
\]
for every $k\in K$. Next, for every $k'\in K$, define 
\[
K_{k'}\coloneqq \Set{k\in K\colon \psi_k*\psi_{k'}\neq 0},
\]
so that Proposition~\ref{prop:56} implies that there is $N\in\N$ such that $\card(K_{k'})\meg N$ for every $k'\in K$.
By the previous remarks, for every $k'\in K$ and for every $h\in \Omega$,
\[
\begin{split}
(\Ec T)_h*\psi_{k'}&=\sum_{k\in K} (T*(\Szego_{(0,i h)}*\psi_k ))*\psi_{k'}\\
&=\sum_{k\in K} T*((\Szego_{(0,i h)}*\psi_k )*\psi_{k'} )\\
&=\sum_{k\in K_{k'}} T*(\psi_{k'}*(\Szego_{(0,i h)}*\psi_k ) )\\
&=\sum_{k\in K_{k'}} (T*\psi_{k'})*(\Szego_{(0,i h)}*\psi_k ) .
\end{split}
\]
Now, $\Szego_{(0,i h)}*\psi_k\in \Sc_\Omega(\Nc)$ and  
\[
\Fc_\Nc(\Szego_{(0,i h)}*\psi_k)=\varphi_k(\,\cdot\,t_k^{-1}) \ee^{-\langle \,\cdot\,,  h\rangle} .
\]
In addition, it is easily seen that the family 
\[
\left(\ee^{-\langle \,\cdot\,,t_k\cdot h\rangle}\varphi_k \right)_{h\in \Omega}
\]
is bounded in $C^\infty_c(\Omega')$, so that Corollary~\ref{cor:24} implies that there is a constant $C_1>0$ such that
\[
\norm{T'*\psi_{k'}*\Szego_{(0,i h)}*\psi_k}_{L^p(\Nc)}\meg C_1 \norm{T'*\psi_{k'}}_{L^p(\Nc)}
\]
for every $h\in \Omega$, for every $k,k'\in K$, and for every $T'\in \Sc'(\Nc)$. Thus,
\[
\norm{(\Ec T)_h*\psi_{k'}}_{L^p(\Nc)}\meg N^{1/\min(1,p)} C_1 \norm{T*\psi_{k'}}_{L^p(\Nc)}
\]
for every $k'\in K$ and for every $h\in \Omega$, whence~{(1)}.

{(2)} By~{(1)}, Proposition~\ref{prop:42}, and Theorem~\ref{teo:10}, it will suffice to prove the assertion when $T\in \Sc_\Omega(\Nc)$. 
In this case, $(\Ec T)_h= T*\Szego_{(0,i h)}\in \Sc_\Omega(\Nc)$ and $\Fc_\Nc(T*\Szego_{(0,i h)})=(\Fc_\Nc T) \ee^{-\langle \,\cdot\,,h\rangle} $. 
Since clearly the mapping 
\[
\Omega'\cup \Set{0}\ni h\mapsto (\Fc_\Nc T)\ee^{-\langle \,\cdot\,,h\rangle}  \in C^\infty_c(\Omega')
\]
is continuous, the assertion follows from Proposition~\ref{prop:40} and Lemma~\ref{lem:55}.

{(3)} Observe first that, by~{(1)}, Propositions~\ref{prop:41} and~\ref{prop:42}, and Theorem~\ref{teo:10}, it will suffice to prove that, for every $T\in B^{\vect s}_{p,q}(\Nc,\Omega)$ and for every $\eta\in \Sc_\Omega(\Nc)$, the mapping
\[
\Omega \cup\Set{0}\ni h\mapsto \langle (\Ec T)_h\vert \eta \rangle\in \C
\]
is continuous. Now, the arguments of~{(1)} imply that
\[
\langle (\Ec T)_h\vert \eta \rangle= \sum_{k,k'\in K} \langle T\vert \eta* \Szego_{(0,i h)}*\psi_k *\psi_{k'}\rangle= \langle T \vert \eta * \Szego_{(0, i h)}\rangle
\]
for every $h\in \Set{0}\cup \Omega$ (defining  $\Szego_{(0,0)}\coloneqq \delta_e$). Thus, the assertion follows from the arguments of~{(2)}.
\end{proof}

\begin{deff}\label{39}
Take $p,q\in]0,\infty]$, and $\vect s\in \R^r$ such that $\vect s\in \frac{1}{p}(\vect b+\vect d)+\frac{1}{2 q'}\vect{m'}+(\R_+^*)^r$. Then, we define 
\[
\Ec\colon B^{-\vect s}_{p,q}(\Nc,\Omega)\to A^{\infty,\infty}_{\vect s-(\vect b+\vect d)/p}(D)
\]
as in Theorem~\ref{teo:11}.
In addition, we define 
\[
\widetilde A^{p,q}_{\vect s}(D)\coloneqq \Ec(B^{-\vect s}_{p,q}(\Nc,\Omega))\qquad \text{and} \qquad\widetilde A^{p,q}_{\vect s,0}(D)\coloneqq \Ec(\mathring B^{-\vect s}_{p,q}(\Nc,\Omega)),
\]
endowed with the corresponding (direct image) topology.
\end{deff}

\emph{As in Section~\ref{sec:5:2}, from now until the end of this section, in order to simplify the notation, we shall sometimes write $\norm{a_k}_{\ell^q(K)}$ instead of $\norm{(a_k)}_{\ell^q(K)}$, for $(a_k)\in \ell^q(K)$.}

\begin{prop}\label{prop:23}
Take $p,q\in ]0,\infty]$, and $\vect s\in \R^r$ such that the following hold:
\begin{itemize}
\item $\vect s\in \frac{1}{2 q}\vect m+(\R_+^*)^r$ (resp.\ $\vect s\in \R_+^r$ if $q=\infty$);

\item $\vect s\in\frac{1}{p}(\vect b+\vect d)+\frac{1}{2 q'}\vect{m'}+(\R_+^*)^r$.
\end{itemize}
Then, there are continuous inclusions
\[
\Ec(\Sc_{\Omega,L}(\Nc))\subseteq A^{p,q}_{\vect s,0}(D) \subseteq \widetilde A^{p,q}_{\vect s,0}(D)
\]
(resp.\
\[
\Ec(\Sc_{\Omega,L}(\Nc))\subseteq A^{p,q}_{\vect s}(D) \subseteq \widetilde A^{p,q}_{\vect s}(D)).
\]
\end{prop}

In particular,  $A^{p,q}_{\vect s,0}(D)$ is dense in $\widetilde A^{p,q}_{\vect s,0}(D)$  (resp.\  $\Ec^{-1}(A^{p,q}_{\vect s}(D))$ is dense in $B_{p,q}^{-\vect s}(\Nc,\Omega)$ for the weak topology $\sigma_{p,q}^{-\vect s}$).

The proof is based on~\cite[Theorem 1.7]{BekolleBonamiGarrigosRicci}, which deals with the case in which $p,q\Meg 1$, $\vect s\in \R\vect 1_r$, and $D$ is a symmetric tube domain.

\begin{proof}
We shall prove that there is a continuous linear mapping 
\[
\Bc\colon A^{p,q}_{\vect s,0}(D)\to \mathring{B}^{-\vect s}_{p,q} (\Nc,\Omega) \qquad \text{(resp.\ $\Bc\colon A^{p,q}_{\vect s}(D)\to B^{-\vect s}_{p,q}(\Nc,\Omega)$)}
\]
such that $\Ec \Bc=I$, and that  
\[
\Ec(\Sc_{\Omega,L}(\Nc))\subseteq A^{p,q}_{\vect s,0}(D)\qquad \text{(resp.\ $\Ec(\Sc_{\Omega,L}(\Nc))\subseteq A^{p,q}_{\vect s}(D)$)}.
\]
This will prove all assertions.

\textsc{Step I.}  We define first an isomorphism $\Bc\colon A^{2,\infty}_{\vect 0}(D)\to B^{\vect 0}_{2,2} (\Nc,\Omega) $ such that $\Ec \Bc=I$.  By Corollary~\ref{cor:1}, the linear mapping 
\[
\Bc\colon A^{2,\infty}_{\vect 0}(D)\ni f \mapsto \lim_{h\to 0} f_h\in L^2(\Nc)
\]
is well defined and continuous. 
 By Proposition~\ref{prop:4}, $\Bc$ is an isomorphism of $A^{2,\infty}_{\vect 0}(D)$ onto 
\[
L^2_\Omega(\Nc)\coloneqq \Set{f\in L^2(\Nc)\colon \pi_\lambda(f)=\chi_{\Omega'}(\lambda) \pi_\lambda(f) P_{\lambda,0} \text{ for almost every $\lambda\in F'\setminus W$}}.
\] 
Since it is easily seen that the canonical mapping $L^2_\Omega(\Nc)\to \Sc'_{\Omega,L}(\Nc)$ induces an isomorphism of   $L^2_\Omega(\Nc)$ onto $ B^{\vect 0}_{2,2}  (\Nc,\Omega)  $, it follows that $\Bc$ is an isomorphism of $A^{2,\infty}_{\vect 0}(D)$ onto $B^{\vect 0}_{2,2}  (\Nc,\Omega) $. Finally, it is clear that $\Ec \Bc = I$.

\textsc{Step II.} We now prove the existence of $\Bc$ on $A^{p,q}_{\vect s,0}(D)$. By Proposition~\ref{prop:8} and~{step I}, it will suffice to prove that $\Bc$ maps $A^{p,q}_{\vect s,0}(D)\cap A^{2,\infty}_{\vect 0}(D)$ into $\mathring{B}^{-\vect s}_{p,q}(\Nc,\Omega)$ continuously.
Then, take $(\lambda_k)_{k\in K}$, $(t_k)$, $(\varphi_k)$, and $(\psi_k)$ as in Lemma~\ref{lem:59}.
In addition, define $\widetilde \psi_k\in \Sc_\Omega(\Nc)$ so that 
\[
\Fc_\Nc\widetilde \psi_k= e^{-\langle \,\cdot\,, t_k^{-1}\cdot e_\Omega\rangle}\varphi_k(\,\cdot\, t_k^{-1})
\]
for every $k\in K$, thanks to Proposition~\ref{prop:40}.
Then, take $f\in A^{p,q}_{\vect s,0}(D)\cap A^{2,\infty}_{\vect 0}(D)$, and observe that
\[
\Bc f= \sum_{k\in K} (\Bc f)*\psi_k
\]
in $L^2(\Nc)$, thanks to~{step I} and Corollary~\ref{cor:7}. In addition, observe that~{step I} and Proposition~\ref{prop:2} imply that 
\[
\pi_\lambda\Big(f_{t_k^{-1}\cdot e_\Omega}\Big)= e^{-\langle \lambda, t_k^{-1}\cdot e_\Omega\rangle}\pi_\lambda(\Bc f),
\]
so that 
\[
(\Bc f)*\psi_k=f_{t_k^{-1}\cdot e_\Omega}*\widetilde \psi_k
\]
for every $k\in K$.
Now, for every $k\in K$ define 
\[
K_k\coloneqq \Set{k'\in K\colon d_{\Omega'}(\lambda_k, \lambda_{k'})\meg 2R\delta},
\]
so that Proposition~\ref{prop:56} implies that there is $N\in \N$ such that $\card(K_k)\meg N$ for every $k\in K$. In addition, Corollary~\ref{cor:24} implies that there is a constant $C_1>0$ such that
\[
\norm*{T*\psi_k*\psi_{k'}}_{L^p(\Nc)}\meg C_1\norm*{T*\psi_k}_{L^p(\Nc)}
\]
for every $T\in \Sc'(\Nc)$ and for every $k,k'\in K$. Thus,
\[
\norm{(\Bc f)*\psi_{k'}}_{L^p(\Nc)}^{\min(1,p)}\meg C_1 \sum_{k\in K_{k'}} \norm{f_{t_k^{-1}\cdot e_\Omega} *\widetilde \psi_k}_{L^p(\Nc)}^{\min(1,p)}
\]
for every $k'\in K$.
Now,  by Lemma~\ref{lem:32} (and its proof) we may assume that there is a $(\delta,R)$-lattice $(\zeta_{j,k},z_{j,k})_{j\in J, k\in K}$ on $D$ such that 
\[
h_k\coloneqq\Im z_{j,k}-\Phi(\zeta_{j,k})= t_k^{-1}\cdot e_\Omega
\]
for every $j\in J$ and for every $k\in K$, and such that
\[
b D+(0, i h_k)\subseteq \bigcup_{j\in J} B((\zeta_{j,k},z_{j,k}),R\delta).
\]
Then, define $B_{j,k}$ as the (relatively compact open) subset of $\Nc$ such that
\[
\chi_{B_{j,k}}=\left(\chi_{B((\zeta_{j,k},z_{j,k}),R\delta)}  \right)_{h_k}
\] 
for every $(j,k)\in J\times K$.

Now, assume that $p\Meg 1$. Then, by means of Young's inequality, we see that there is a constant $C_2>0$ such that 
\[
\norm{(\Bc f)*\psi_{k'}}_{L^p(\Nc)}\meg C_2 \sum_{k\in K_{k'}} \norm{f_{h_k}}_{L^p(\Nc)},
\]
for every $k'\in K$,
so that by means of Corollary~\ref{cor:34} we see that there is a constant $C_3>0$ such that
\[
\norm*{\Delta_\Omega^{-\vect s}(\lambda_{k'})\norm{(\Bc f)*\psi_{k'}}_{L^p(\Nc)}  }_{\ell^q(K)}\meg C_3 N^{1/q+1/q'} \norm*{\Delta_\Omega^{\vect s}(h_k)\norm{f_{h_k}}_{L^p(\Nc)}  }_{\ell^q(K)}.
\]
Since clearly there is a constant $C_4>0$ such that
\[
\norm{f_{h_k}}_{L^p(\Nc)} \meg C_4 \Delta_\Omega^{-(\vect b+\vect d)/p}(h_k)\norm*{\bigg(\sup\limits_{B_{j,k}} \abs{f}  \bigg)_{j}}_{\ell^p(J)}
\]
\CR for every $k\in K$, the assertion \CB follows from Theorem~\ref{teo:8}, provided that $\delta$ is sufficiently small.

Then, assume that $p<1$, and observe that
\[
\begin{split}
&\int_\Nc \abs{(f_{h_k}*\widetilde \psi_k)(\zeta,x)}^p\,\dd (\zeta,x)\meg C_5
\int_\Nc \sum_{j\in J} \sup\limits_{B_{j,k}} \abs{f}^p\times \\
&\qquad \qquad \qquad \qquad \times \bigg( \int_{B_{j,k}} \abs{\widetilde \psi_k((\zeta',x')^{-1}(\zeta,x))  }\,\dd (\zeta',x')\bigg)^p \,\dd (\zeta,x)\\   
\end{split}
\]
for every $k\in K$. In addition, by homogeneity we see that there is a constant $C_6>0$ such that
\[
\int_\Nc \bigg( \int_{B_{j,k}} \abs{\widetilde \psi_k((\zeta',x')^{-1}(\zeta,x))  }\,\dd (\zeta',x')\bigg)^p \,\dd (\zeta,x)\meg C_6\Delta_\Omega^{-(\vect b+\vect d)}(h_k)
\]
for every $j\in J$ and for every $k\in K$. Thus, by means of Corollary~\ref{cor:34} we see that there is a constant $C_7>0$ such that
\[
\norm{(\Bc f)*\psi_{k'}}_{L^p(\Nc)}^{p}\meg C_7 \Delta_\Omega^{-(\vect b+\vect d)}(h_k) \sum_{k\in K_{k'}}\sum_{j\in J} \sup\limits_{B_{j,k}} \abs{f} ^p.
\]
Therefore, the first assertion follows as above by means of Theorem~\ref{teo:8}, provided that $\delta$ is sufficiently small.

 Finally, take $\varphi \in \Sc_{\Omega,L}(\Nc)$ and choose $\varphi'\in \Sc_\Omega(\Nc)$ such that $\Fc_\Nc \varphi'(\lambda)=1$ for every $\lambda\in \Omega'$ such that $\pi_\lambda(\varphi)\neq 0$.  Then, $\varphi=\varphi*\varphi'$ and
\[
(\Ec \varphi)_h=\varphi*\big(\varphi'*\Szego_{(0,i h)}\big)\in \Sc_{\Omega,L}(\Nc)
\]
for every $h\in \Omega$. It then follows easily that $\Ec \varphi\in A^{p,q}_{\vect s,0}(D)$.

\textsc{Step III.} Take $f\in A^{p,q}_{\vect s}(D)$, and choose $(g^{(\eps)})_{\eps>0}$ as in Lemma~\ref{lem:34}. Observe that 
\[
f^{(\eps,h)}\coloneqq f(\,\cdot\,+i h) g^{(\eps)}\in A^{p,q}_{\vect s}(D)\cap A^{2,\infty}_{\vect 0}(D)
\]
for every $h\in \Omega$ and for every $\eps>0$. In addition,
\[
\norm{f^{(\eps,h)}}_{A^{p,q}_{\vect s}(D)}\meg  \norm{f}_{A^{p,q}_{\vect s}(D)}
\]
for every $h\in \Omega$ and for every $\eps>0$ (cf.~Corollary~\ref{cor:26}). 
Furthermore, arguing as in~{step II}, we see that there is a constant $C_8>0$ such that 
\[
\norm{\Bc f^{(\eps,h)}}_{B^{-\vect s}_{p,q}(\Nc,\Omega)}\meg C_8 \norm{f^{(\eps,h)}}_{A^{p,q}_{\vect s}(D)}
\]
for every $h\in \Omega$ and for every $\eps>0$.
Therefore, Corollary~\ref{cor:28} implies that the $\Bc(f^{(\eps,h)})$ stay in a relatively compact subset of $B^{-\vect s}_{p,q}(\Nc,\Omega)$ for the weak topology $\sigma^{-\vect s}_{p,q}$. 
Let $T$ be the limit of $\Bc f^{(\eps,h)}$ (for the weak topology  $\sigma^{-\vect s}_{p,q}$) along an ultrafilter $\Uf$ which is finer than the filter `$\eps\to 0^+$ and $h\to 0$.' 
Then 
\[
f=\lim_{(\eps,h), \Uf} f^{(\eps,h)}=\lim_{(\eps,h), \Uf} \Ec(\Bc f^{(\eps,h)})=\Ec T
\]
pointwise.
Since $\Ec$ is one-to-one on $B^{-\vect s}_{p,q}(\Nc,\Omega)$ by Theorem~\ref{teo:11}, this implies that $T$ does not depend on $\Uf$, so that $\Bc(f^{(\eps,h)})$ converges to $T$ in  the weak topology  $\sigma^{-\vect s}_{p,q}$. 
If we define $\Bc f\coloneqq T$, it is then easily verified that $\Bc$ induces a continuous linear mapping of $A^{p,q}_{\vect s}(D)$ into $B^{-\vect s}_{p,q}(\Nc,\Omega)$ such that $\Ec\Bc=I$.

The inclusion $\Ec(\Sc_{\Omega,L}(\Nc))\subseteq A^{p,q}_{\vect s}(D)$ is proved as in~{step II}.
\end{proof}

We shall now introduce an auxiliary property which will allow us to give sufficient conditions for the equality $A^{p,q}_{\vect s}(D)=\widetilde A^{p,q}_{\vect s}(D)$. Cf.~\cite{BekolleBonamiGarrigosRicci} for a discussion of its significance in other areas of mathematics.

\begin{deff}\label{def:2}
Take $\vect s\in \R^r$ and $p,q\in ]0,\infty]$. We say that property $(D)^{\vect s,0}_{p,q}$ (resp.\ $(D)^{\vect s}_{p,q}$) holds if there are a $(\delta,R)$-lattice $(\lambda_k)_{k\in K}$, for some $\delta>0$ and some $R>0$, a bounded family $(\varphi_k)_{k\in K}$ of elements of $C^\infty_c(\Omega')$ such that 
\[
\sum_{k\in K} \varphi_k(\,\cdot\, t_k^{-1})\Meg 1
\]
on $\Omega'$ (where $t_k\in T_+$ and $\lambda_k=e_{\Omega'}\cdot t_k$ for every $k\in K$), and two constants $c,C>0$ such that
\[
\norm*{\sum_{k\in K} T_k*\psi_k}_{L^p(\Nc)}\meg C \norm*{\Delta_{\Omega'}^{\vect s}(\lambda_k) \ee^{c\langle \lambda_k,e_\Omega\rangle} \norm{T_k*\psi_{k}}_{L^p(\Nc)}  }_{\ell^q(K)}
\]
for every  $(T_k)\in \Sc(\Nc)^{(K)}$ (resp.\ for every $(T_k)\in \Sc'(\Nc)^{(K)}$), where $\psi_k\in\Sc_\Omega(\Nc)$ and $\Fc_\Nc\psi_k=\varphi_k(\,\cdot\,t_k^{-1})$.
\end{deff}

Observe that property $(D)^{\vect s,0}_{p,q}$ (resp.\  $(D)^{\vect s}_{p,q}$) implies property $(D)^{\tilde{\vect s},0}_{p,\tilde q}$ (resp.\  properties $(D)^{\tilde{\vect s},0}_{p,\tilde q}$  and  $(D)^{\tilde{\vect s}}_{p,\tilde q}$) for every $\tilde q\in ]0,q]$ and for every $\tilde{\vect s}\in \vect s-\R_+^r $.

\begin{lem}\label{lem:66}
Take $p,q\in ]0,\infty]$ such that $q\meg \min(p,p')$. Then, property $(D)^{\vect 0}_{p,q}$ holds.
\end{lem}

The proof is based on~\cite[Lemma 4.8]{BekolleBonamiGarrigosRicci}, which deals with the case $p,q\Meg 1$.

\begin{proof}
Observe first that property $(D)^{\vect 0}_{\ell,\min(1,\ell)} $ and $(D)^{\vect 0}_{2,2}$ clearly hold for every $\ell\in ]0,\infty]$. The general assertion follows by interpolation, as in the proof of~\cite[Lemma 4.8]{BekolleBonamiGarrigosRicci}).
\end{proof}

\begin{lem}\label{lem:64}
Take $\vect s\in \R^r$ and $p,q\in ]0,\infty]$, and assume that property $(D)^{\vect s,0}_{p,q}$ (resp.\ $(D)^{\vect s}_{p,q}$) holds.
Take $(\lambda_k)_{k\in K}$, $(\varphi_k)$, and $(\psi_k)$ as in Lemma~\ref{lem:51}.
Then, there are two constants $c,C>0$ such that
\[
\norm*{\sum_{k\in K} T_k*\psi_k}_{L^p(\Nc)}\meg C \norm*{\Delta_{\Omega'}^{\vect s}(\lambda_k) \ee^{c\langle \lambda_k,e_\Omega\rangle} \norm{T_k*\psi_k}_{L^p(\Nc)}  }_{\ell^q(K)}
\]
for every  $(T_k)\in \Sc(\Nc)^{(K)}$ (resp.\   $(T_k)\in \Sc'(\Nc)^{(K)}$).
\end{lem}

\begin{proof}
By assumption, there are a $(\delta',R')$-lattice $(\lambda'_{k'})_{k'\in K'}$ on $\Omega'$ for some $\delta'>0$ and some $R'>1$, a bounded family $(\varphi'_{k'})_{k'\in K'}$ of elements of $C^\infty_c(\Omega')$ such that 
\[
\sum_{k'\in K'}\varphi'_{k'}(\,\cdot\, t_{k'}'^{-1})\Meg 1
\]
on $\Omega'$, where $t'_{k'}\in T_+$ and 
\[
\lambda'_{k'}=e_{\Omega'}\cdot t'_{k'}
\]
for every $k'\in K'$, and two constants $c',C'>0$ such that
\[
\norm*{\sum_{k'\in K'} T'_{k'}*\psi'_{k'}}_{L^p(\Nc)}\meg C' \norm*{\Delta_{\Omega'}^{\vect s}(\lambda'_{k'}) \ee^{c'\langle \lambda'_{k'},e_\Omega\rangle} \norm{T'_{k'}*\psi'_{k'}}_{L^p(\Nc)}  }_{\ell^q(K')}
\]
for every  $(T'_{k'})\in \Sc(\Nc)^{(K')}$ (resp.\  $(T'_{k'})\in \Sc'(\Nc)^{(K')}$), where $\psi'_{k'}\in\Sc_\Omega(\Nc)$  and $\Fc_\Nc \psi'_{k'}=\varphi'_k(\,\cdot\,t_{k'}'^{-1})$ for every $k'\in K'$. Choose $t_k\in T_+$ so that
\[
\lambda_k=e_{\Omega'}\cdot t_k
\]
for every $k\in K$.

For every $k'\in K'$, define
\[
K_{k'}\coloneqq \Set{k\in K\colon \psi_k*\psi'_{k'}\neq 0} \quad \text{ and} \quad K'_{k'}\coloneqq\Set{k''\in K'\colon \psi'_{k''}*\psi'_{k'}\neq 0},
\]
and observe that Proposition~\ref{prop:56} shows that there is $N\in \N$ such that 
\[
\card(K_{k'}), \card(K'_{k'})\meg N
\]
for every $k'\in K'$ and such that each $k\in K$ belongs to at most $N$ of the sets $K_{k'}$.
Define 
\[
\widetilde \varphi'\coloneqq \sum_{k'\in K'} \varphi'_{k'}(\,\cdot\, t_{k'}'^{-1}),
\]
so that $\widetilde \varphi'$ is well-defined, of class $C^\infty$, and $\Meg 1$ on $\Omega'$. Define, in addition,
\[
\widetilde \varphi_k\coloneqq \frac{\varphi_k}{\widetilde \varphi'(\,\cdot\,t_k)} \qquad \text{and} \qquad \widetilde \varphi'_{k'}\coloneqq \frac{\varphi'_{k'}}{\widetilde \varphi'(\,\cdot\,t'_{k'})} 
\]
for every $k\in K$ and for every $k'\in K'$, so that
\begin{align*}
\varphi_k(\,\cdot\,t_k^{-1})&=\sum_{k\in K_{k'}} \widetilde \varphi_k(\,\cdot\, t_k^{-1}) \varphi'_{k'}(\,\cdot\, t'^{-1}_{k'}),\\
\widetilde \varphi_k&=\frac{\varphi_k}{\sum_{k\in K_{k'}} \varphi'_{k'}(\,\cdot\, (t'_{k'} t_k^{-1})^{-1} )}
\end{align*}
and
\[
\widetilde \varphi'_{k'}=\frac{\varphi'_{k'}}{\sum_{k''\in K'_{k'}} \varphi'_{k''}(\,\cdot\, (t'_{k'} t_{k''}'^{-1})^{-1} )}
\]
for every $k\in K$ and for every $k'\in K'$. 
By means of Lemma~\ref{lem:50} and the preceding arguments, we see that the families $(\widetilde \varphi_k)$ and $(\widetilde \varphi'_{k'})$ are bounded in $C^\infty_c(\Omega')$. Then, define $\widetilde \psi_k, \widetilde \psi'_{k'}\in\Sc_\Omega(\Nc)$ so that 
\[
\Fc_\Nc \widetilde \psi_k=\widetilde \varphi_k(\,\cdot\, t_k^{-1})\qquad \text{and} \qquad\Fc_\Nc \widetilde \psi'_{k'}=\widetilde \varphi'_{k'}(\,\cdot\, t'^{-1}_{k'})
\]
for every $k\in K$ and for every $k'\in K'$.

Fix $(T_k)\in \Sc(\Nc)^{(K)}$ (resp.\ $(T_k)\in \Sc'(\Nc)^{(K)}$), and define 
\[
T'_{k'}\coloneqq \sum_{k\in K_{k'}} T_k*\widetilde \psi_k
\]
for every $k'\in K'$, so that $(T'_{k'})\in \Sc(\Nc)^{(K')}$ (resp.\ $(T'_{k'})\in \Sc'(\Nc)^{(K')}$) and
\[
\sum_{k\in K} T_k*\psi_k=\sum_{k'\in K'} T'_{k'}*\psi'_{k'}.
\]
Now, Corollary~\ref{cor:24} implies that there is a constant $C_1>0$ such that
\[
\norm{T'*\psi'_{k'}*\widetilde \psi_k}_{L^p(\Nc)}=\norm{T'*\widetilde\psi'_{k'}* \psi_k}_{L^p(\Nc)}\meg C_1 \norm{T'* \psi_k}_{L^p(\Nc)}
\]
for every $T'\in \Sc'(\Nc)$, for every $k'\in K'$, and for every $k\in K$, so that
\[
\norm{T'_{k'}*\psi'_{k'}}_{L^p(\Nc)}\meg C_1 N^{(1/p-1)_+} \sum_{k\in K_{k'}} \norm{T_k*\psi_k}_{L^p(\Nc)}.
\]
In addition, Lemma~\ref{lem:48} shows that there is a constant $c>0$ such that $\ee^{c'\langle \lambda'_{k'},e_\Omega\rangle}\meg \ee^{c\langle \lambda_k,e_\Omega\rangle}$ for every $k'\in K'$ and for every $k\in K_{k'}$. Analogously, Corollary~\ref{cor:34} shows that there is a constant $C_2>0$ such that $\Delta_{\Omega'}^{\vect s}(\lambda'_{k'})\meg C_2 \Delta_{\Omega'}^{\vect s}(\lambda_k)$ for every $k'\in K'$ and for every $k\in K_{k'}$.
Therefore,
\[
\norm*{\sum_{k\in K} T_k*\psi_k}_{L^p(\Nc)}\meg C_3  \norm*{\Delta_{\Omega'}^{\vect s}(\lambda_k) \ee^{c\langle \lambda_k,e_\Omega\rangle} \norm{T_k*\psi_k}_{L^p(\Nc)}  }_{\ell^q(K)},
\]
where $C_3\coloneqq C' C_1 N^{(1/p-1)_+ +\max(1,1/q)}C_2$. The assertion follows.
\end{proof}

The next lemma, which is the first step in the proof of Theorem~\ref{teo:6}, characterizes properties $(D)^{\vect s}_{p,q}$ and $(D)^{\vect s,0}_{p,q}$ in terms of the extension operator $\Ec$. See~\cite[Proposition 4.16]{BekolleBonamiGarrigosRicci} for another equivalent formulation of property $(D)^{\vect s}_{p,q}$ in the case in which $p,q\in [1,\infty[$, $\vect s\in \R\vect1_r$, and $D$ is  an irreducible symmetric tube domain.

\begin{prop}\label{prop:45}
Keep the hypotheses and the notation of Theorem~\ref{teo:11}. Then, the following conditions are equivalent:
\begin{enumerate}
\item[\em(1)] property  $(D)^{\vect s}_{p,q}$ (resp.\  $(D)^{\vect s,0}_{p,q}$) holds;

\item[\em(2)] $\Ec$ induces a continuous linear mapping from $B^{\vect s}_{p,q}(\Nc,\Omega)$ (resp.\ $\mathring{B}^{\vect s}_{p,q}(\Nc,\Omega)$) into $ A^{p,\infty}_{-\vect s}(D)$.
\end{enumerate}
In addition, the following conditions are equivalent:
\begin{enumerate}
\item[\em(1\ensuremath{'})] property $(D)^{\vect{s'},0}_{p,q}$ holds for every $\vect{s'}\in \vect s-(\R_+^*)^r$;

\item[\em(2\ensuremath{'})] $\Ec$ induces a continuous linear mapping from $\mathring{B}^{\vect{s'}}_{p,q}(\Nc,\Omega)$ into $ A^{p,\infty}_{-\vect {s'},0}(D)$ for every $\vect{s'}\in \vect s-(\R_+^*)^r$.
\end{enumerate}
\end{prop}

The proof is based on~\cite[Proposition 4.16]{BekolleBonamiGarrigosRicci}, which deals with the case in which $p,q\in [1,\infty[$, $\vect s\in \R\vect1_r$, and $D$ is  an irreducible symmetric tube domain.

\begin{proof}
{(1) $\implies$ (2).} Take $T\in B^{\vect s}_{p,q}(\Nc,\Omega)$ (resp.\ $T\in \Sc_{\Omega,L}(\Nc)$). 
By homogeneity, it will suffice to prove that 
\[
(\Ec T)_h \in L^p(\Nc) \quad \text{for some $h\in \Omega$}.
\] 
Take $(\lambda_k)_{k\in K}$, $(\varphi_k)_{k\in K}$, $(t_k)_{k\in K}$, $(\psi_k)_{k\in K}$, $C$ and $c$ as in Definition~\ref{def:2}. 
Observe that, by Lemma~\ref{lem:64}, we may assume that 
\[
\sum_{k\in K} \varphi_k(\,\cdot\,t_k^{-1})=1
\]
on $\Omega'$, so that 
\[
T=\sum_{k\in K}  T*\psi_k
\]
in  the  weak topology $\sigma^{\vect s}_{p,q}$  thanks to Lemma~\ref{lem:59}. Hence, 
\[
(\Ec T)_h=\sum_{k\in K} (\Ec T)_h*\psi_k
\]
pointwise for every $h\in \Omega$.
For every $k\in K$, define 
\[
K_k\coloneqq \Set{k'\in K\colon \psi_k*\psi_{k'}\neq 0},
\]
and observe that,  by Proposition~\ref{prop:56}, there is $N\in \N$ such that $\card(K_k)\meg N$ for every $k\in K$.
By the proof of Theorem~\ref{teo:11}, we know that
\[
(\Ec T)_h*\psi_{k'}=\sum_{k\in K_{k'}} (T*\psi_{k'})*(\Szego_{(0, i h)}*\psi_k)
\]
for every $h\in \Omega$ and for every $k'\in K$. 
Then, define
\[
T_{K',k}\coloneqq \sum_{k'\in K'\cap K_{k}} (T*\psi_{k'})*\Szego_{(0, i h)}
\]
for every finite subset $K'$ of $K$ and for every $k\in K$, so that
\[
\sum_{k'\in K'}(\Ec T)_h*\psi_{k'}  = \sum_{k\in K} T_{K',k} *\psi_k.
\] 
Therefore, 
\[
\norm*{\sum_{k'\in K'}(\Ec T)_h*\psi_{k'}   }_{L^p(\Nc)}\meg C\norm*{\Delta^{\vect s}_{\Omega'}(\lambda_k) e^{c \langle \lambda_k, e_\Omega\rangle}\norm*{T_{K',k}*\psi_k}_{L^p(\Nc)}  }_{\ell^q(K)}.
\]
Now, observe that Lemma~\ref{lem:61} shows that $\Szego_{(0,i h)}*\psi_{k}\in \Sc_\Omega(\Nc)$  and that 
\[
\Fc_\Nc(\Szego_{(0,i h)}*\psi_{k})=\varphi_{k}(\,\cdot\, t_{k}^{-1}) \ee^{-\langle \,\cdot\,,h\rangle} .
\]
By means of Lemma~\ref{lem:48}, we see that there is a constant $C_1>0$ such that the family 
\[
\left(\ee^{\langle \lambda_k,h\rangle/C_1}\ee^{-\langle \,\cdot\,t_k, h\rangle} \varphi_{k}\right)_{h\in \Omega, k\in K}
\]
is bounded in $C^\infty_c(\Omega')$, so that Corollary~\ref{cor:24} implies that there is a constant $C_2>0$ such that
\[
\norm{T'*\psi_{k'}*(\Szego_{(0,i h)}*\psi_{k})}_{L^p(\Nc)}\meg C_2 \ee^{-\langle \lambda_{k'},h\rangle/C_1}\norm{T'*\psi_{k}}_{L^p(\Nc)}
\]
for every $T'\in \Sc'(\Nc)$, for every $k,k'\in K$, and for every $h\in \Omega$.
Hence,
\[
\norm*{T_{K',k}*\psi_k }_{L^p(\Nc)}\meg C_2 N^{(1/p-1)_+}\norm{T*\psi_{k}}_{L^p(\Nc)} \sum_{k'\in K_k}\ee^{-\langle \lambda_{k'},h\rangle/C_1} 
\]
for every $k\in K$ and for every $h\in \Omega$. By means of Lemma~\ref{lem:48} again we see that there is a constant $C_3>0$ such that $e^{-\langle \lambda_{k'},h\rangle/C_1}\meg e^{-\langle \lambda_{k},h\rangle/C_3}$ for every $k\in K$, for every $k'\in K_{k}$, and for every $h\in \Omega$. Therefore,
\[
\norm*{T_{K',k}*\psi_k }_{L^p(\Nc)}\meg C_2 N^{\max(1,1/p)} \norm{T*\psi_{k}}_{L^p(\Nc)} \ee^{-\langle \lambda_{k},h\rangle/C_3} 
\]
for every $k'\in K$ and for every $h\in \Omega$. Now, take $h\in \Omega$ so that $c e_\Omega-h/C_3\in -\overline \Omega$ (e.g.,\ $h= c C_3 e_\Omega$), and observe that
\[
\norm*{\sum_{k'\in K'}(\Ec T)_h*\psi_{k'}   }_{L^p(\Nc)}\meg C C_2 N^{\max(1,1/p)} \norm*{\Delta_{\Omega'}^{\vect s}(\lambda_k) \norm{T*\psi_{k}}_{L^p(\Nc)}  }_{\ell^q(K)}
\]
so that, by the arbitrariness of $K'$, 
\[
\norm*{(\Ec T)_h }_{L^p(\Nc)}\meg C C_2 N^{\max(1,1/p)} \norm*{\Delta_{\Omega'}^{\vect s}(\lambda_k) \norm{T*\psi_{k}}_{L^p(\Nc)}  }_{\ell^q(K)},
\]
whence~{(2)}.

{(2) $\implies$ (1).} Take $(\lambda_k)_{k\in K}$, $(t_k)$, $(\varphi_k)$, and $(\psi_k)$ as in Lemma~\ref{lem:59}.
In addition, take $(T_k)\in \Sc'(\Nc)^{(K)}$ (resp.\ $(T_k)\in \Sc(\Nc)^{(K)}$), fix $h\in \Omega$, and define
\[
T'_h\coloneqq \sum_{k\in K} T_k*\psi_{k,h},
\]
where $\psi_{k,h}\in\Sc_\Omega(\Nc)$  and 
\[
\Fc_\Nc \psi_{k,h}=\varphi_k(\,\cdot\, t_k^{-1}) \ee^{\langle\,\cdot\,,h\rangle}
\]
for every $k\in K$. 
For every $k\in K$, define
\[
K_k\coloneqq \Set{k'\in K\colon \psi_k*\psi_{k'}\neq 0 },
\]
and observe that Proposition~\ref{prop:56} implies that there is $N\in \N$ such that $\card(K_k)\meg N$ for every $k\in K$.
Now, by means of Lemma~\ref{lem:48} we see that there is a constant $C_4>0$ such that the family 
\[
\left(\ee^{-C_4\langle \lambda_k, h\rangle}\varphi_k \ee^{\langle\,\cdot\,t_k,h\rangle}\right)_{k\in K,h\in \Omega}
\]
is bounded in $C^\infty_c(\Omega')$, so that Corollary~\ref{cor:24} shows that there is a constant $C_5>0$ such that
\[
\norm{T'*\psi_{k,h}*\psi_{k'}}_{L^p(\Nc)}=\norm{T'*\psi_k*\psi_{k',h}}_{L^p(\Nc)}\meg C_5 \ee^{C_4\langle \lambda_{k'},h\rangle}\norm{T'*\psi_{k}}_{L^p(\Nc)}
\]
for every $T'\in \Sc'(\Nc)$, for every $k,k'\in K$, and for every $h\in \Omega$. Therefore,
\[
\begin{split}
& \norm{T'_h*\psi_k}_{L^p(\Nc)}\meg C_5 N^{(1/p-1)_+} \sum_{k'\in K_k} \ee^{C_4\langle \lambda_{k'},h\rangle} \norm{T_{k'}*\psi_{k'}}_{L^p(\Nc)}
\end{split}
\]
for every $h\in\Omega$ and for every $k\in K$.  
In particular, since  $T'_h\in B^{\vect s}_{p,q}(\Nc,\Omega)$ (resp.\ $T'_h\in \Sc_{\Omega,L}(\Nc)$) for every $h\in \Omega$, there is a constant $C_6>0$ such that
\[
\norm{(\Ec T'_h)_{h'}}_{L^p(\Nc)}\meg C_6 \Delta_\Omega^{\vect s}(h')\norm*{ \Delta_{\Omega'}^{\vect s}(\lambda_k) \sum_{k'\in K_k} \ee^{C_4\langle \lambda_{k'},h\rangle} \norm{T_{k'}*\psi_{k'}}_{L^p(\Nc)} }_{\ell^q(K)}
\]
for every $h,h'\in \Omega$. 
Let us prove that $(\Ec T'_h)_h=\sum_{k\in K}  f_k*\psi_{k}$ for every $h\in \Omega$. Now, by the proof of  Theorem~\ref{teo:11} (and observing that the considered sums are finite),
\[
\begin{split}
(\Ec T'_h)_h&=\sum_{k,k'\in K} (T'_h*(\Szego_{(0,i h)}*\psi_k)  )*\psi_{k'}\\
&= \sum_{k,k',k''\in K} ( T_{k''}*\psi_{k'',h}*(\Szego_{(0,i h)}*\psi_k) )*\psi_{k'}\\
&=\sum_{k,k',k''\in K}  T_{k''}*\psi_{k''}*\psi_k *\psi_{k'}\\
&=\sum_{k\in K}  T_k*\psi_k .
\end{split}
\]
Therefore, by means of Corollary~\ref{cor:34} we see that there is a constant $C_7>0$ such that
\[
\norm*{\sum_{k\in K}  T_k*\psi_k}_{L^p(\Nc)}\meg C_7  \Delta_\Omega^{\vect s}(h) \norm*{  \Delta_{\Omega'}^{\vect s}(\lambda_k) \ee^{C_4 \langle \lambda_{k'},h\rangle} \norm{T_k*\psi_k}_{L^p(\Nc)}}_{\ell^q(K)}
\]
Thus,~{(1)} follows choosing $h=e_\Omega$.

{(1\ensuremath{'}) $\implies$ (2\ensuremath{'}).} Take $\vect{s'}\in \vect s-(\R_+^*)^r$, and take $\vect{s''}\in \vect s-(\R_+^*)^r$ such that $ \vect{s'}\in \vect{s''}-(\R_+^*)^r$. By the implication  {(1) $\implies$ (2)}, we see that $\Ec$  induces two continuous linear mappings from $\mathring{B}^{\vect{s''}}_{p,q}(\Nc,\Omega)$ and  $\mathring{B}^{\vect{s'}}_{p,q}(\Nc,\Omega)$ into $A^{p,\infty}_{-\vect{s''}}(D)$ and $A^{p,\infty}_{-\vect{s'}}(D)$, respectively. 
Therefore, 
\[
\Ec(\Sc_{\Omega,L}(\Nc)) \subseteq A^{p,\infty}_{-\vect{s'}}(D)\cap A^{p,\infty}_{-\vect{s''}}(D).
\]
Since clearly $(\Ec f)_h\in L^p_0(\Nc)$ for every $f\in\Sc_{\Omega,L}(\Nc)$, this implies that    
\[
\Ec(\Sc_{\Omega,L}(\Nc)) \subseteq A^{p,\infty}_{-\vect{s'},0}(D)
\]
by Corollary~\ref{cor:14}.
The assertion follows by means of Theorem~\ref{teo:10}.

{(2\ensuremath{'}) $\implies$ (1\ensuremath{'}).} This follows from the implication {(2) $\implies$ (1)}.
\end{proof}

\begin{cor}
Take $q\in ]0,\infty]$ and $\vect s\in -\frac{1}{2 q'}\vect{m'}-(\R_+^*)^r$. Then, property $(D)^{\vect s}_{\infty,q} $ holds.
\end{cor}

\begin{proof}
This is a consequence of Theorem~\ref{teo:11} and Proposition~\ref{prop:45}.
\end{proof}

In the next result we give sufficient conditions for the equality $A^{p,q}_{\vect s}(D)=\widetilde A^{p,q}_{\vect s}(D)$. Since it is not known, in general, when property $(D)^{\vect{s'},0}_{p,q}$ holds, we give here a general result, and provide more explicit conditions in Corollary~\ref{cor:29}.

\begin{teo}\label{teo:6}
Take $\vect s,\vect{s'}\in \R^r$ and $p,q,\tilde q\in ]0,\infty]$ such that the following hold:
\begin{itemize}
\item property $(D)^{\vect{s'},0}_{p,\tilde q}$  (resp.\ $(D)^{\vect{s'}}_{p,\tilde q}$) holds;

\item $\vect s\in \frac{1}{2 q}\vect m-\vect{s'}+\left( \frac{1}{2 \tilde q}-\frac{1}{2 q}\right)_+ \vect{m'}+(\R_+^*)^r $ (resp.\ $\vect s\in- \vect{s'}+\R_+^r $ if $q=\tilde q=\infty$);

\item $\vect s\in \frac 1 p(\vect b+\vect d)+\frac{1}{2 q'}\vect{m'}+(\R_+^*)^r$.
\end{itemize}
Then, $A^{p,q}_{\vect s,0}(D)=\widetilde A^{p,q}_{\vect s,0}(D)$ (resp.\ of  $A^{p,q}_{\vect s}(D)=\widetilde A^{p,q}_{\vect s}(D)$).
\end{teo}

The proof is based on~\cite[Theorem 4.11]{BekolleBonamiGarrigosRicci}, which deals with the case in which $p,q\in [1,\infty[$, $\vect s\in \R \vect 1_r$, and $D$ is  an irreducible  symmetric tube domain.

\begin{proof}
Take $T\in  \mathring{B}^{-\vect s}_{p,q}(\Nc,\Omega)$ (resp.\  $T\in B^{-\vect s}_{p,q}(\Nc,\Omega)$). 
In addition, take  $(\lambda_k)_{k\in K}$, $(t_k)$, $(\varphi_k)$, and $(\psi_k)$ as in Lemma~\ref{lem:59}
Let us first prove that 
\[
\Ec((\Ec T)_h)=(\Ec T)(\,\cdot\,+i h)
\]
for every $h\in \Omega$. Indeed, the proof of Theorem~\ref{teo:11} shows that
\[
\begin{split}
(\Ec((\Ec T)_h))_{h'}*\psi_k &=\sum_{k'\in K} ((\Ec T)_h*\psi_k  )*(\Szego_{(0, i h')}*\psi_{k'})\\
&= \sum_{k',k''\in K}  (T*\psi_k  )*(\Szego_{(0, i h)}*\psi_{k''})*(\Szego_{(0, i h')}*\psi_{k'})\\
&= \sum_{k'\in K}  (T*\psi_k  )*(\Szego_{(0, i (h+h'))}*\psi_{k'})\\
&= (\Ec T)_{h+h'}*\psi_k
\end{split}
\]
for every $k\in K$ and for every $h,h'\in \Omega$, since clearly 
\[
\Szego_{(0, i h')}*\Szego_{(0, i h)}=\Szego_{(0, i (h+h'))}.
\]

Now, Proposition~\ref{prop:45}  shows that there is a constant $C_1>0$ such that $(\Ec T)_{2 h}\in L^p_0(\Nc)$ (resp.\ $(\Ec T)_{2 h}\in L^p(\Nc)$) and
\[
\Delta_\Omega^{-\vect{s'}}(h)\norm{(\Ec T)_{2 h}}_{L^p(\Nc)}\meg C_1 \norm*{\Delta_{\Omega'}^{\vect{s'}}(\lambda_k)\norm{(\Ec T)_h* \psi_k}_{L^p(\Nc)}  }_{\ell^{\tilde q}(K)}
\]
for every $h\in \Omega$.
Further, observe that 
\[
 (\Ec T)_h* \psi_k=\sum_{k'\in K_k}( T*\psi_k)*(\Szego_{(0,i h)}*\psi_{k'}),
 \]
 so that, by means of Lemma~\ref{lem:48} and Corollary~\ref{cor:24}, we see that there is a constant $C_2>0$ such that
\[
\norm{(\Ec T)_h* \psi_k}_{L^p(\Nc)} \meg C_2 \ee^{-\langle \lambda_k,h\rangle/C_2} \norm{T* \psi_k}_{L^p(\Nc)} 
\]
for every $h\in \Omega$ and for every $k\in K$. Then, take $\vect{s''}\in \R^r $, and observe that 
\[
\begin{split}
&\norm*{\Delta_{\Omega'}^{\vect{s'}}(\lambda_k)\ee^{-\langle \lambda_k,h\rangle/C_2} \norm{T* \psi_k}_{L^p(\Nc)} }_{\ell^{\tilde q}(K)}\meg\norm*{\Delta_{\Omega'}^{-\vect{s''}}(\lambda_k)\ee^{-\langle \lambda_k,h\rangle/(2 C_2)} }_{\ell^{\hat q}(K)}\times \\
&\qquad \qquad\qquad\qquad\qquad\qquad\times  \norm*{ \Delta_{\Omega'}^{\vect{s''}+\vect{s'}}(\lambda_k)\ee^{-\langle \lambda_k,h\rangle/(2C_2)} \norm{T* \psi_k}_{L^p(\Nc)} }_{\ell^q(K)} ,
\end{split}
\]
where $\hat q\in [1,\infty]$ and $\frac{1}{\hat q}=\left(\frac{1}{\tilde q}-\frac 1 q\right)_+$. Observe that Proposition~\ref{prop:58}, Lemma~\ref{lem:18}, and Corollary~\ref{cor:34} show that there is a constant $C_3>0$ such that
\[
\norm*{\Delta_{\Omega'}^{-\vect{s''}}(\lambda_k)\ee^{-\langle \lambda_k, h\rangle/(2C_2)} }_{\ell^{\hat q}(K)}\meg C_3 \Delta^{\vect{s''}}_{\Omega}(h)
\]
for every $h\in \Omega$, provided that 
\[
\vect{s''}\in -\frac{1}{2 \hat q}\vect{m'}-(\R_+^*)^r \qquad (\text{resp.\ } \vect{s''}\in -\R_+^r \text{if $\hat q=\infty$}).
\]
In addition, Proposition~\ref{prop:58} and Lemma~\ref{lem:18} again show that there is a constant $C_4>0$ such that
\[
\norm*{\Delta_\Omega^{\vect s+\vect{s'} +\vect{s''}}\, \ee^{-\langle \lambda_k,\,\cdot\,\rangle/(4 C_2)} }_{L^q(\nu_\Omega)}=C_4\Delta_{\Omega'}^{ -\vect s-\vect{s'}-\vect{s''} }(\lambda_k)
\]
provided that
\[
\vect{s''}\in-\vect{s'}-\vect s +\frac{1}{2 q}\vect{m}+(\R_+^*)^r \qquad (\text{resp.\ } \vect{s''}\in -\vect{s'}-\vect s+\R_+^r \text{if $ q=\infty$}).
\]
 Observe that our assumptions imply that we may find $\vect{s''}$ satisfying the preceding requirements, so that, setting $C_5\coloneqq 2^{-\vect{s'}-\vect{s''}}C_1 C_2 C_3 C_4 $,
\[
\norm{\Ec T}_{A^{p,q}_{\vect s}(D)}\meg C_5\norm*{\Delta_{\Omega'}^{-\vect s}(\lambda_k)\norm{T*\psi_k}_{L^p(\Nc)}  }_{\ell^q(K)}.
\]
Thus, $\Ec$ maps $\mathring{B}^{-\vect s}_{p,q}(\Nc,\Omega)$ into $A^{p,q}_{\vect s}(D)$ (resp.\  $B^{-\vect s}_{p,q}(\Nc,\Omega)$ into $A^{p,q}_{\vect s}(D)$) continuously. 
To conclude, by Proposition~\ref{prop:23} and Corollary~\ref{cor:14}, it will suffice to show that, if $q=\infty$ and $T\in \Sc_{\Omega,L}(\Nc)$, then the mapping 
\[
h\mapsto\Delta^{\vect s}_{\Omega}(h) \norm{(\Ec T)_h}_{L^p(\Nc)}
\]
belongs to $C_0(\Omega)$. Since we may further reduce to the case $T\in \Sc_\Omega(\Nc)$ by left-invariance, the assertion follows easily from Proposition~\ref{prop:40}, since $\vect s\in (\R_+^*)^r$ under the first set of assumptions.
\end{proof}

\begin{cor}\label{cor:29}
Take $p,q\in ]0,\infty]$ and 
\[
\vect s\in\sup\left(\frac{1}{2 q}\vect m+\left( \frac{1}{2 \min(p,p')}-\frac{1}{2q} \right)_+\vect{m'} , \frac{1}{p}(\vect b+\vect d)+\frac{1}{2 q'}\vect{m'}\right) +(\R_+^*)^r.
\]
Then, $A^{p,q}_{\vect s,0}(D)=\widetilde A^{p,q}_{\vect s,0}(D)$ and $A^{p,q}_{\vect s}(D)=\widetilde A^{p,q}_{\vect s}(D)$.
\end{cor}

This result is optimal when $q\meg \min(p,p')$.

\begin{proof}
This follows from Lemma~\ref{lem:66} and Theorem~\ref{teo:6}.
\end{proof}

We are now able to extend the second assertion of Theorem~\ref{teo:10} to the general case $p,q\in ]0,\infty]$.

\CR
In the following result we both provide a (sesquilinear) dual pairing between the spaces $\widetilde A^{p,q}_{\vect s,0}(D)$ and $\widetilde A^{p',q'}_{\vect{s'}}(D)$, and a new interpretation of the sesquilinear form 
\[
\langle \,\cdot\,\vert \,\cdot\,\rangle\colon \mathring B^{\vect s}_{p,q}(\Nc,\Omega)\times B^{-\vect s-(1/p-1)_+(\vect b+\vect d)}_{p',q'}(\Nc,\Omega).
\]

\begin{prop}\label{cor:30}
Take $p,q\in ]0,\infty]$, $\vect s\in \frac 1 p (\vect b+\vect d)+\frac{1}{2 q'}\vect{m'}+(\R_+^*)^r$, and $\vect{s'}\in \frac{1}{p'}(\vect b+\vect d)+\frac{1}{2 \max(1,q)}\vect{m'}+(\R_+^*)^r$. Define  $\vect{s''}\coloneqq \vect s+\vect{s'}-(1/p-1)_+(\vect b+\vect d)$, and assume that $\vect{s''}\in \frac 1 2 \vect m+(\R_+^*)^r$.
Then, the sesquilinear form
\[
\langle\,\cdot\,\vert \,\cdot\,\rangle_{\vect s,\vect{s'}}\colon\Ec(\Sc_{\Omega,L}(\Nc))\times \Ec(\Sc_{\Omega,L}(\Nc))\ni(f,g)\mapsto \int_D f(\zeta,z)\overline{g(\zeta,z)}\Delta_\Omega^{\vect{s''}}(\Im z-\Phi(\zeta))\,\dd \nu_D(\zeta,z)
\]
is well-defined and extends to a unique continuous sequilinear form on $\widetilde A^{p,q}_{\vect s,0}(D)\times\widetilde A^{p',q'}_{\vect{s'}}(D)$ which is continuous on the second factor with respect to the topology induced by the weak topology $\sigma^{p',q'}_{\vect s-(1/p-1)_+(\vect b+\vect d)}$ through $\Ec$.\footnote{This latter requirement is needed to ensure uniqueness when $\max(p',q')=\infty$.} In addition, this extension of $\langle\,\cdot\,\vert \,\cdot\,\rangle_{\vect s,\vect{s'}}$ induces an antilinear isomorphism of $\widetilde A^{p',q'}_{\vect{s'}}(D)$ onto $\widetilde A^{p,q}_{\vect s,0}(D)'$, and
\[
\langle \Ec T\vert \Ec (T'* I^{-\vect{s''}}_\Omega)\rangle_{\vect s,\vect{s'}}= \frac{2^{n-m}i^{\vect{s''}}\Gamma_{\Omega}(\vect{s''})}{  2^{\vect {s''}} \pi^{n+m}} \big\langle T\big\vert T'\big\rangle
\]
for every $T,T'\in \Sc_{\Omega,L}(\Nc)$.
\end{prop}

\begin{proof}
%By Theorem~\ref{teo:9} and Corollary~\ref{cor:29}, we may assume that $\Ec$ induces an isomorphism of $\mathring B^{\vect s}_{p,q}(\Nc,\Omega)$ onto $A^{p,q}_{-\vect s,0}(D)$.
%In addition, by Theorems~\ref{teo:4} and~\ref{teo:5} and Corollary~\ref{cor:29}, we may take $\vect{s'}\in \R^r$ so that the following conditions hold:
%\begin{enumerate}
%\item[(1)] $\Ec$ induces an isomorphism of $ B^{\vect {s'}}_{p',q'}(\Nc,\Omega)$ onto $A^{p',q'}_{-\vect {s'}}(D)$;
%
%\item[(2)] $-\vect s-\vect{s'}-(1/p-1)_+(\vect b+\vect d)\in \frac 1 2 \vect{m'}+(\R_+^*)^r$;
%
%\item[(3)] property $\atomics^{p,q}_{-\vect s, (\vect b+\vect d)/\min(1,p)+\vect s+\vect{s'},0}$ holds.
%\end{enumerate}
%Then,~{(3)} and Proposition~\ref{prop:29} imply that the sesquilinear form on $A^{p,q}_{-\vect s,0}(D)\times A^{p',q'}_{-\vect{s'}}(D)$
%\[
%(f,g)\mapsto \int_D f(\zeta,z)\overline{g(\zeta,z)}\Delta_\Omega^{-\vect s-\vect{s'}-(\vect b+\vect d)/\min(1,p)}(\Im z-\Phi(\zeta))\,\dd \nu_D(\zeta,z)
%\]
%induces an antilinear isomorphism of $A^{p',q'}_{-\vect{s'}}(D)$ onto $A^{p,q}_{-\vect s,0}(D)'$.
Take $T,T'\in \Sc_{\Omega,L}(\Nc)$, and observe that
\[
\begin{split}
&\int_D (\Ec T)(\zeta,z)\overline{(\Ec (T'* I^{-\vect{s''}}_\Omega))(\zeta,z)}\Delta_\Omega^{\vect{s''}}(\Im z-\Phi(\zeta))\,\dd \nu_D(\zeta,z)=\\
&\qquad = \frac{2^{n-m}i^{\vect{s''}}}{\pi^{n+m}}\int_\Omega \int_{\Omega'} \tr(\pi_\lambda(T)\pi_\lambda(T')^*) \Delta^{\vect{s''}}_{\Omega'}(\lambda)e^{-2\langle \lambda,h\rangle}\abs{\Pfaff(\lambda)}\,\dd \lambda\, \Delta_\Omega^{\vect{s''}}(h)\,\dd \nu_\Omega(h)\\
&\qquad= \frac{2^{n-m}i^{\vect{s''}}\Gamma_{\Omega}(\vect{s''})}{ 2^{\vect {s''}} \pi^{n+m}} \int_{\Omega'} \tr(\pi_\lambda(T)\pi_\lambda(T')^*) \abs{\Pfaff(\lambda)}\,\dd \lambda\\
&\qquad= \frac{2^{n-m}i^{\vect{s''}}\Gamma_{\Omega}(\vect{s''})}{ 2^{\vect {s''}} \pi^{n+m}} \langle T\vert T'\rangle
\end{split}
\]
where the first equality follows from Corollary~\ref{cor:7} and Proposition~\ref{prop:55}, the second equality follows from Proposition~\ref{prop:58}, while the last equality follows from Corollary~\ref{cor:7}. 
The assertion follows easily.\CB
\end{proof}
\CB

The following result simply translates Theorem~\ref{teo:9} in terms of the spaces $\widetilde A^{p,q}_{\vect s}(D)$. It still has relevant consequences.

\begin{prop}\label{prop:30}
Take $p,q\in ]0,\infty]$ and $\vect s,\vect{s'}\in \R^r$ such that $\vect s,\vect s+\vect{s'}\in \frac 1 p (\vect b+\vect d)+\frac{1}{2 q'}\vect{m'}+(\R_+^*)$.
Let 
\[
 \Ic^{-\vect{s'}}\colon B^{-\vect s}_{p,q}(\Nc,\Omega)\to B^{-\vect s-\vect{s'}}_{p,q}(\Nc,\Omega)
\]
be the unique isomorphism which is continuous for the weak topologies $\sigma^{-\vect s}_{p,q}$ and $\sigma^{-\vect s-\vect{s'}}_{p,q}$ and which  induces the automorphism $f\mapsto f* I^{-\vect{s'}}_\Omega$ of $\Sc_{\Omega,L}(\Nc)$.
In addition, define a mapping 
\[
\widetilde\Ic^{-\vect{s'}}\colon \widetilde A^{p,q}_{\vect s}(D)\to \widetilde A^{p,q}_{\vect s+\vect{s'}}(D)
\]
so that 
\[
(\widetilde\Ic^{-\vect{s'}} f)_h=\Ic^{-\vect{s'}}(f_h)
\]
for every $f\in A^{p,q}_{\vect s}(D)$ and for every $h\in \Omega$. Then, $\widetilde\Ic^{-\vect{s'}}$ induces isomorphisms
\[
\widetilde A^{p,q}_{\vect s}(D) \to \widetilde A^{p,q}_{\vect s+\vect{s'}}(D) \qquad \text{and} \qquad\widetilde A^{p,q}_{\vect s,0}(D)\to \widetilde A^{p,q}_{\vect s+\vect{s'},0}(D).
\]
\end{prop}

\begin{proof}
The assertion follows from Theorem~\ref{teo:9}, since clearly $ \Ic^{-\vect{s'}}(\Ec T)_h=\Ec(\Ic^{-\vect{s'}} T)_h$ for every $T\in B^{-\vect s}_{p,q}(\Nc,\Omega)$ and for every $h\in \Omega$ (argue as in the proof of Theorem~\ref{teo:11}).
\end{proof}

By means of Proposition~\ref{prop:30}, several properties of \emph{some} weighted Bergman spaces $A^{p,q}_{\vect s}(D)$ `propagate' to all the spaces $\widetilde A^{p,q}_{\vect s}(D)$.
In particular, all spaces $\widetilde A^{p,q}_{\vect s}(D)$ enjoy an analogue of the atomic decomposition studied in Section~\ref{sec:6:4}.

\begin{cor}\label{cor:31}
Take $p,q\in ]0,\infty]$ and $\vect s\in \frac 1 p (\vect b+\vect d)+\frac{1}{2 q'}\vect{m'}+(\R_+^*)$. Take $R_0>1$ and  $\vect{s'}\in \R^r$ such that
\[
\vect s+\vect{s'}\in \frac{1}{\min(1,p)}(\vect b+\vect d)-\frac{1}{2 q} \vect{m'}-\left( \frac{1}{2 \min(1,p)}-\frac{1}{2 q} \right)_+\vect m-(\R_+^*)^r.
\]
Then, there is $\delta_0>0$ such that, for every $(\delta,R)$-lattice $(\zeta_{j,k},z_{j,k})_{j\in J,k\in K}$ on $D$, with $\delta\in ]0,\delta_0]$ and $R\in ]0,R_0]$, defining $h_k\coloneqq \Im z_{j,k}-\Phi(\zeta_{j,k})$ for every $k\in K$ and for some (hence every) $j\in J$, the mapping
\[
\Szego\colon \ell^{p,q}(J,K)\ni \lambda \mapsto \sum_{j,k} \lambda_{j,k} B^{\vect{s'}}_{(\zeta_{j,k},z_{j,k})} \Delta^{(\vect b+\vect d)/p-\vect s-\vect{s'}}_\Omega(h_k)\in \widetilde A^{p,q}_{\vect s}(D)
\]
is a surjective strict morphism and induces a strict morphism $\Szego_0$ of $\ell^{p,q}_0(J,K)$ onto $\widetilde A^{p,q}_{\vect s,0}(D)$. Further, both $\Szego$ and $\Szego_0$ have a continuous linear section.
\end{cor}

To prove this result, we need a lemma.

\begin{lem}\label{lem:80}
Take $p,q\in ]0,\infty]$, $\vect s\in \frac 1 p (\vect b+\vect d)+\frac{1}{2 q'}\vect{m'}+(\R_+^*)$, and $\vect{s'}\in\C^r$. Then, the following conditions are equivalent:
\begin{enumerate}
\item[\em(1)] $B^{\vect{s'}}_{(\zeta,z)}\in \widetilde A^{p,q}_{\vect s,0}(D)$ (resp.\ $B^{\vect{s'}}_{(\zeta,z)}\in \widetilde A^{p,q}_{\vect s}(D)$) for some $(\zeta,z)\in D$;

\item[\em(2)] $B^{\vect{s'}}_{(\zeta,z)}\in \widetilde A^{p,q}_{\vect s,0}(D)$ (resp.\ $B^{\vect{s'}}_{(\zeta,z)}\in \widetilde A^{p,q}_{\vect s}(D)$)  for every $(\zeta,z)\in D$;

\item[\em(3)] $\vect s+\Re\vect{s'}\in\frac 1 p (\vect b+\vect d)-\frac{1}{2 q}\vect{m'}-(\R_+^*)^r $ (resp.\ $\vect s+\Re\vect{s'}\in\frac 1 p (\vect b+\vect d)-\R_+^r $ if $q=\infty$).
\end{enumerate}
\end{lem}

\begin{proof}
Observe first that, by Corollary~\ref{cor:29}, there is $\vect{s''}\in \N_{\Omega'}$ such that 
\[
A^{p,q}_{\vect s+\vect{s''},0}(D)=\widetilde A^{p,q}_{\vect s+\vect{s''},0}(D)\qquad \text{and}\qquad A^{p,q}_{\vect s+\vect{s''}}(D)=\widetilde A^{p,q}_{\vect s+\vect{s''}}(D) .
\]
We define $\widetilde\Ic^{-\vect{s''}}$ as in Proposition~\ref{prop:30}.

{(1) $\implies$ (3).} 
Observe that  Proposition~\ref{prop:30} implies that
\[
\widetilde\Ic^{-\vect{s''}}(B^{\vect{s'}}_{(\zeta,z)})\in A^{p,q}_{\vect s+\vect{s''},0}(D)\qquad \text{(resp.\ $\widetilde\Ic^{-\vect{s''}}(B^{\vect{s'}}_{(\zeta,z)})\in A^{p,q}_{\vect s+\vect{s''}}(D)$)}.
\]
Now, Proposition~\ref{prop:39} shows that
\[
B^{\vect{s'}}_{(\zeta,z)}* I^{-\vect{s''}}_\Omega= (2 i)^{\vect{s''}} \left( \vect{s'}+\textstyle{\frac12} \vect{m'} \right)_{\vect{s''}} B^{\vect{s'}-\vect{s''}}_{(\zeta,z)}.
\]
Observe that $\left( \vect{s'}+\frac 1 2 \vect{m'} \right)_{\vect{s''}}\neq 0$ since $\widetilde\Ic^{-\vect{s''}}$ is injective on $\widetilde A^{p,q}_{\vect s}(D)$,  so that
\[
B^{\vect{s'}-\vect{s''}}_{(\zeta,z)}\in A^{p,q}_{\vect s+\vect{s''},0}(D) \qquad \text{(resp.\ $B^{\vect{s'}-\vect{s''}}_{(\zeta,z)}\in A^{p,q}_{\vect s+\vect{s''}}(D)$)}.
\]
Therefore, Proposition~\ref{prop:60} implies that~{(3)} holds.

{(3) $\implies$ (2).} Fix $(\zeta,z)\in D$. Observe first that Proposition~\ref{prop:39} shows that
\[
B^{\vect{s'}}_{(\zeta,z)}* I^{-\vect{s''}}_\Omega= (2 i)^{\vect{s''}} \left( \vect{s'}+\textstyle{\frac12} \vect{m'} \right)_{\vect{s''}} B^{\vect{s'}-\vect{s''}}_{(\zeta,z)}.
\]
In addition, since $\Re\vect{s'},\vect s+\Re\vect{s'}\in \frac 1 p (\vect b+\vect d) -\R_+^r$, 
\[
B^{\vect{s'}}_{(\zeta,z)}\in A^{\infty,\infty}_{\vect s-(\vect b+\vect d)/p} \qquad \text{and} \qquad B^{\vect{s'}-\vect{s''}}_{(\zeta,z)}\in A^{\infty,\infty}_{\vect s+\vect{s''}-(\vect b+\vect d)/p}
\]
by Proposition~\ref{prop:60}. Further, if $\vect{s''}$ is large enough, then Proposition~\ref{prop:60} also shows that $B^{\vect{s'}-\vect{s''}}_{(\zeta,z)}\in A^{p,q}_{\vect s+\vect{s''},0}(D)$ (resp.\ $B^{\vect{s'}-\vect{s''}}_{(\zeta,z)}\in A^{p,q}_{\vect s+\vect{s''}}(D)$).
 Now, $\widetilde\Ic^{-\vect{s''}}$ induces  a continuous linear mapping
\[
A^{\infty,\infty}_{\vect s-(\vect b+\vect d)/p}(D)\to A^{\infty,\infty}_{\vect s+\vect{s''}-(\vect b+\vect d)/p}(D),
\]
which extends the isomorphism
\[
\widetilde A^{p,q}_{\vect s,0}(D)\to \widetilde A^{p,q}_{\vect s+\vect{s''},0}(D) \qquad \text{(resp.\ $\widetilde A^{p,q}_{\vect s}(D)\to \widetilde A^{p,q}_{\vect s+\vect{s''}}(D)$)},
\]
thanks to Corollary~\ref{cor:40}, Theorem~\ref{teo:11}, and Proposition~\ref{prop:30}. 
Therefore, $B^{\vect{s'}}_{(\zeta,z)}\in \widetilde A^{p,q}_{\vect s,0}(D)$ (resp.\ $B^{\vect{s'}}_{(\zeta,z)}\in \widetilde A^{p,q}_{\vect s}(D)$), whence~{(2)} by the arbitrariness of $(\zeta,z)\in D$.
\end{proof}

\begin{proof}
By Corollary~\ref{cor:29}, there is $\vect{s''}\in \N_{\Omega'}$ such that
\[
A^{p,q}_{\vect s+\vect{s''}}(D)=\widetilde A^{p,q}_{\vect s+\vect{s''}}(D)\qquad \text{and}\qquad A^{p,q}_{\vect s+\vect{s''},0}(D)=\widetilde A^{p,q}_{\vect s+\vect{s''},0}(D).
\]
In addition, Theorems~\ref{teo:4} and~\ref{teo:5} imply that, if $\vect{s''}$ is large enough, then for every $\vect{s'}\in \R^r$ such that
\[
\vect s+\vect{s'}\in \frac{1}{\min(1,p)}(\vect b+\vect d)-\frac{1}{2 q} \vect{m'}-\left( \frac{1}{2 \min(1,p)}-\frac{1}{2 q} \right)_+\vect m-(\R_+^*)^r
\]
and for every $R_0>1$, there is $\delta_0>0$ such that, for every $\delta\in ]0,\delta_0]$, for every $R\in ]0,R_0]$, and for every  $(\delta,R)$-lattice $(\zeta_{j,k},z_{j,k})_{j\in J,k\in K}$ on $D$, defining $h_k\coloneqq \Im z_{j,k}-\Phi(\zeta_{j,k})$ for every $k\in K$ and for some (hence every) $j\in J$, the mapping
\[
\Szego_{\vect s+\vect{s''}}\colon\ell^{p,q}(J,K)\ni \lambda \mapsto \sum_{j,k} \lambda_{j,k} B^{\vect{s'}-\vect{s''}}_{(\zeta_{j,k},z_{j,k})} \Delta^{(\vect b+\vect d)/p-\vect s-\vect{s'}}_\Omega(h_k)\in  A^{p,q}_{\vect s+\vect{s''}}(D)
\]
is a surjective strict morphism and induces a strict morphism $\Szego_{\vect s+\vect{s''},0}$ of $\ell^{p,q}_0(J,K)$ onto $A^{p,q}_{\vect s+\vect{s''},0}(D)$. 
Further, both $\Szego_{\vect s+\vect{s''}}$ and $\Szego_{\vect s+\vect{s''},0}$ have a continuous linear section.
With the notation of Proposition~\ref{prop:30}, define  
\[
\Szego_{\vect s}\coloneqq \Ic^{\vect{s''}}\Szego_{\vect s+\vect{s''}}
\]
and
\[
\Szego_{\vect s,0}\coloneqq \Ic^{\vect{s''}}\Szego_{\vect s+\vect{s''},0}
\]
so that $\Szego_{\vect s}\colon \ell^{p,q}(J,K)\to\widetilde A^{p,q}_{\vect s}(D)$ and $\Szego_{\vect s,0}\colon \ell^{p,q}_0(J,K)\to\widetilde A^{p,q}_{\vect s,0}(D)$ are surjective strict morphisms and have a continuous linear section. 
In order to conclude, it will suffice to show that 
\[
\Ic^{\vect{s''}}(B^{\vect{s'}-\vect{s''}}_{(\zeta,z)})= c_{\vect{s'},\vect{s''}} B^{\vect{s'}}_{(\zeta,z)}  
\]
for a suitable $c_{\vect{s'},\vect{s''}}\neq 0$. This follows from Proposition~\ref{prop:39} and Lemma~\ref{lem:80}.
\end{proof}

In the following result, we characterize the equality $A^{p,q}_{\vect s,0}(D)=\widetilde A^{p,q}_{\vect s,0}(D)$ by means of atomic decomposition and the properties of Riemann--Liouville operators. 
Since we do not know, in general, if, e.g., property $\atomic^{p,q}_{\vect s, \vect{s'}}$ implies property $\atomic^{p,q}_{\vect s, \vect{s''}}$ for all $\vect{s''}\in \vect{s'}-\R_+^r$, this characterization is somewhat weaker than one may hope. 

\begin{cor}\label{cor:32}
Take $p,q\in ]0,\infty]$ and $\vect s\in \frac 1 p (\vect b+\vect d)+\frac{1}{2 q'}\vect{m'}+(\R_+^*)$ such that $\vect s\in \frac{1}{2 q}\vect m+(\R_+^*)$ (resp.\ $ \vect s\in \R_+^r$ if $q=\infty$). Then, the following conditions are equivalent:
\begin{enumerate}
\item[\em(1)] $A^{p,q}_{\vect s,0}(D)=\widetilde A^{p,q}_{\vect s,0}(D)$ (resp.\  $A^{p,q}_{\vect s}(D)=\widetilde A^{p,q}_{\vect s}(D)$);

\item[\em(2)] property $\atomics^{p,q}_{\vect s, \vect{s'},0}$ (resp.\ $\atomics^{p,q}_{\vect s, \vect{s'}}$) holds for every $\vect{s'}$ in some translate of $- \R_+^r$;

\item[\em(3)] property $\atomic^{p,q}_{\vect s, \vect{s'},0}$ (resp.\ $\atomic^{p,q}_{\vect s, \vect{s'}}$) holds for every $\vect{s'}$ in some translate of $- \R_+^r$;

\item[\em(4)] convolution by $I^{-\vect{s'}}_\Omega$ induces an isomorphism of $ A^{p,q}_{\vect s,0}(D)$ onto $A^{p,q}_{\vect s+\vect{s'},0}(D)$ (resp.\ of $ A^{p,q}_{\vect s}(D)$ onto $A^{p,q}_{\vect s+\vect{s'}}(D)$) for every $\vect{s'}$ in the intersection of $\N_{\Omega'}$ with some translate of $\R_+^r$. 
\end{enumerate} 
\end{cor}

\begin{proof}
{(1) $\implies$ (2).} This follows from Corollary~\ref{cor:31}.

{(2) $\implies$ (3).} Obvious.

{(3) $\implies$ (1).} By Corollary~\ref{cor:31}, we may find $\vect{s'}\in \R^r$ and  a $(\delta,R)$-lattice $(\zeta_{j,k},z_{j,k})_{j\in J,k \in K}$ for some $\delta>0$ and some $R>1$ such that, defining $h_k\coloneqq \Im z_{j,k}-\Phi(\zeta_{j,k})$ for every $k\in K$ and for some (hence every $j\in J$), the mapping 
\[
\lambda \mapsto \sum_{j,k} \lambda_{j,k} B^{\vect{s'}}_{(\zeta_{j,k},z_{j,k})} \Delta^{(\vect b+\vect d)/p-\vect s-\vect{s'}}_\Omega(h_k)
\]
induces a continuous linear mapping of $\ell^{p,q}_0(J,K)$ into $A^{p,q}_{\vect s,0}(D)$ (resp.\ of $\ell^{p,q}(J,K)$ into $A^{p,q}_{\vect s}(D)$) and a strict morphism of $\ell^{p,q}_0(J,K)$ onto $\widetilde A^{p,q}_{\vect s,0}(D)$ (resp.\ of $\ell^{p,q}(J,K)$ onto $\widetilde A^{p,q}_{\vect s}(D)$). The assertion follows from Proposition~\ref{prop:23}.

{(1) $\iff$ (4).} This follows Corollary~\ref{cor:29} and Proposition~\ref{prop:30}.
\end{proof}

\begin{deff}\label{def:6}
Take $p,q\in ]0,\infty]$, $\vect s\in \R^r$, and  $\vect{s'}\in \N_{\Omega'}$ such that $\vect{s}+\vect{s'}\in \frac 1 p (\vect b+\vect d)+ \frac{1}{2 q'} \vect{m'}+(\R_+^*)^r$. Then, we define $\widehat A^{p,q}_{\vect s,\vect{s'}}(D)$ as  the Hausdorff locally convex space associated with the space of $f\in \Hol(D)$ such that $f*I^{-\vect{s'}}_{\Omega}\in \widetilde A^{p,q}_{\vect s+\vect{s'}}(D)$, endowed with the corresponding topology.  We define $\widehat A^{p,q}_{\vect s,\vect{s'},0}(D)$ analogously
\end{deff}

Thus, $\widehat A^{p,q}_{\vect s,\vect{s'}}(D)$ can be identified with $\widetilde A^{p,q}_{\vect s}(D)=\widehat A^{p,q}_{\vect s,\vect{0}}(D)$ when $\vect{s}\in \frac 1 p (\vect b+\vect d)+ \frac{1}{2 q'} \vect{m'}+(\R_+^*)^r$. In addition, the new definition of $\widehat A^{2,2}_{\vect s,\vect{s'}}(D)$  coincides with the preceding one, thanks to Proposition~\ref{prop:30}.
Furthermore, arguing as in the proof of Proposition~\ref{prop:54}, one may prove that convolution by $I^{\vect {-s'}}_\Omega$ induces an isomorphism of $\widehat A^{p,q}_{\vect s,\vect{s'}}(D)$ onto $\widetilde A^{p,q}_{\vect s+\vect{s'}}(D)$.

\smallskip

We conclude this section with some necessary conditions for property $(D)^{-\vect s,0}_{p,q}$, which are then also necessary for  the equality $A^{p,q}_{\vect s}(D)=\widetilde A^{p,q}_{\vect s}(D)$ by Propositions~\ref{prop:5} and~\ref{prop:45}.

\CR
\begin{prop}\label{prop:81}
Take $p,q\in ]0,\infty]$ and $\vect s\in\R^r$, and assume that property $(D)^{-\vect s,0}_{p,q}$ holds. Then, the following hold:
\begin{enumerate}
\item[\em(1)] if $\vect s\in \frac 1 p (\vect b+\vect d)+\frac{1}{2 q'}\vect{m'}+(\R_+^*)^r$, then $\vect s\in\left( \frac{1}{2 p}-\frac{1}{2 q}\right)\vect{m'}+\R_+^r$;

\item[\em(2)] if $q> 2$, $p<\infty$, and $n=0$, then $\vect s\in \left(\frac{1}{4}-\frac{1}{2 q}  \right)\vect{m'}+(\R_+^*)^r$.
\end{enumerate}
\end{prop}\CB

Observe that property $(D)^{-\vect s}_{p,q}$  holds if $A^{p,q}_{\vect s}(D)=\widetilde  A^{p,q}_{\vect s}(D)$. In addition, property $(D)^{-\vect s}_{p,q}$ implies property $(D)^{-\vect s,0}_{p,q}$.

The proof is based on~\cite[Proposition 4.34]{BekolleBonamiGarrigosRicci}, which deals with the case in which $p,q\in [1,\infty[$, $\vect s\in \R\vect 1_r$, and $D$ is a symmetric tube domain. Notice that the corresponding assertion of~\cite[Proposition 4.34]{BekolleBonamiGarrigosRicci} is sharper than~{(1)}, since a more refined procedure allows to show that $\vect s\in \left( \frac{1}{2 p}-\frac{1}{2 q}\right)\vect{m'}+(\R_+^*)^r$ (under the assumptions of~\cite[Proposition 4.34]{BekolleBonamiGarrigosRicci}).

\begin{proof}
{(1)} Take $\vect{s'}\in \frac 1 p (\vect b+\vect d)-\frac{1}{2 q}\vect{m'}-\vect s-(\R_+^*)^r $, so that Lemma~\ref{lem:80} implies that $B^{\vect{s'}}_{(0, i e_\Omega)}\in \widetilde A^{p,q}_{\vect s,0}(D)$. 
Since, by Proposition~\ref{prop:45}, $\widetilde A^{p,q}_{\vect s,0}(D)\subseteq A^{p,\infty}_{\vect s}(D)$, Proposition~\ref{prop:60} implies that $\vect {s'}\in \frac 1 p (\vect b+\vect d)-\frac{1}{2 p}\vect{m'}-(\R_+^*)^r$ if $p<\infty$ and $\vect {s'}\in \frac 1 p (\vect b+\vect d)-\frac{1}{2 p}\vect{m'}-\R_+^r$ if $p=\infty$. By the arbitrariness of $\vect{s'}$, this implies that
\[
\vect s\in\left( \frac{1}{2 p}-\frac{1}{2 q}\right)\vect{m'}+\R_+^r.
\]

{(2)} Assume that $q>2$, $p<\infty$, and  $n=0$.
Take $(\lambda_k)_{k\in K}$, $(t_k)$, $(\varphi_k)$, $(\psi_k)$, $C$, and $c$ as in Definition~\ref{def:2}.
In addition, fix a finite subset $K'$ of $K$ such that $B_{\Omega'}(\lambda_k,R\delta)\cap B_{F'}(0,1)\neq \emptyset$ for every $k\in K'$, take $\varphi\in C_c^\infty(F')$ so that
\[
\chi_{B_{\Omega'}(\lambda_k,\delta/2)}\meg \varphi(\,\cdot\,-\lambda_k)\meg \chi_{B_{\Omega'}(\lambda_k,\delta)}
\]
for every $k\in K'$, choose $T\in \Sc_\Omega(\Nc)$ such that $\Fc_\Nc(T)= \varphi$, and define $T_k\coloneqq e^{i\langle\lambda_k,\,\cdot\,\rangle}T$, so that $T_k\in \Sc_\Omega(\Nc)$ and
\[
\Fc_\Nc(T_k)=\varphi(\,\cdot\, -\lambda_k)
\]
for every $k\in K'$ (here we use the assumption that $n=0$). By Lemma~\ref{lem:64}, we may assume that $\varphi_k(\lambda \cdot t_k^{-1})=1 $ for every $\lambda\in B_{\Omega'}(\lambda_k,\delta)$, so that $T_k*\psi_k=T_k$ for every $k\in K'$. 

Now, take a probability space $(X,\mi)$ and a finite family $(r_k)_{k\in K'}$ of Rademacher functions on $X$, that is, $\mi$-measurable functions on $X$ such that
\[
\left(\bigotimes_{k\in K'} r_k\right)_*(\mi)=\frac{1}{2^{\card(K')}}\sum_{\eps\in \Set{-1,1}^{K'}} \delta_\eps
\] 
(cf.~\cite[C.1]{Grafakos}). Then, by Khintchine's inequality there is a constant $C_1>0$ (independent of $K'$) such that
\[
\frac{1}{C_1} \left(\sum_{k\in K'} \abs{a_k}^2  \right)^{1/2} \meg \norm*{\sum_{k\in K'} a_k r_k }_{L^p(\mi)}\meg C_1  \left(\sum_{k\in K'} \abs{a_k}^2  \right)^{1/2} 
\]
for every $(a_k)\in \C^{K'}$ (cf.~\cite[C.2]{Grafakos}).

Now, setting $C_2\coloneqq C \max_{H} \ee^{c\langle \,\cdot\,,e_\Omega \rangle}$ where $H$ is a compact subset of $F'$ such that $B_{\Omega'}(\lambda,R\delta)\subseteq H$ for every $\lambda\in \Omega'$ such that $B_{\Omega'}(\lambda,R\delta)\cap B_{F'}(0,1)\neq \emptyset$ (cf.~Corollary~\ref{cor:18}),
\[
\begin{split}
\norm*{ \sum_{k\in K'} r_k(x) a_k T_k }_{L^p(\Nc)}&\meg C_2 \norm*{\Delta_{\Omega'}^{-\vect s}(\lambda_k) \norm*{ r_k(x) a_k T_k }_{L^p(\Nc)} }_{\ell^q(K')}\\
&= C_2 \norm{T}_{L^p(\Nc)}\norm*{\Delta_{\Omega'}^{-\vect s}(\lambda_k) a_k }_{\ell^q(K')}
\end{split}
\]
for $\mi$-almost every $x\in X$, so that
\[
\norm{a_k}_{\ell^2(K')}\meg C_1 C_2\norm*{\Delta_{\Omega'}^{-\vect s}(\lambda_k) a_k }_{\ell^q(K')}.
\]
Choosing $a_k=\Delta_{\Omega'}^{(q/(q-2))\vect s}(\lambda_k)$ and recalling that $q>2$, we then find
\[
\sum_{k\in K'} \Delta_{\Omega'}^{(2q/(q-2))\vect s}(\lambda_k)\meg (C_1 C_2)^{2 q/(q-2)}.
\]
By Corollary~\ref{cor:34}, the arbitrariness of $K'$ implies that
\[
\int_{\Omega'\cap B_{F'}(0,1)}\Delta_{\Omega'}^{(2q/(q-2))\vect s}(\lambda)\,\dd\nu_{\Omega'}(\lambda)<\infty,
\]
so that $\vect s\in \frac{q-2}{4 q}\vect{m'}+(\R_+^*)^r $ by Lemma~\ref{lem:10} and Proposition~\ref{prop:58}.
\end{proof}

\section{Bergman Projectors}\label{sec:6:7}

In this section we deal with the boundedness properties of the Bergman projectors $P_{\vect s}$. As for atomic decomposition, it is somewhat simpler to deal with integral operators with positive kernels, so that we shall introduce an auxiliary operator $P_{\vect s,+}$.
It then turns out that, when $p,q\Meg 1$, the boundedness of $P_{\vect {s'},+}$ on $L^{p,q}_{\vect s}(D)$ is equivalent to property $\atomic^{p,q}_{\vect s,\vect{s'},+}$ (cf.~Proposition~\ref{prop:31}). 
This will, in particular, imply that property $\atomic^{p,q}_{\vect s,\vect{s'},+}$ implies property $\atomic^{p',q'}_{\vect b+\vect d-\vect s-\vect{s'},\vect{s'},+}$.

For the operators $P_{\vect{s'}}$ a weaker result holds. On the one hand, in full generality one may prove that the boundedness of $P_{\vect{s'}}$ on $L^{p,q}_{\vect s}(D)$ implies property $\atomic^{p,q}_{\vect s,\vect{s'}}$ (cf.~Proposition~\ref{prop:61}).
On the other hand, if property $\atomic^{p,q}_{\vect s,\vect{s'}}$ holds for every $\vect{s'}$ in a translate of $-\R_+^r$, then $P_{\vect{s'}}$ is bounded on $L^{p,q}_{\vect s}(D)$ for every $\vect{s'}$ in a translate of $-\R_+^r$ (cf.~Corollaries~\ref{cor:32} and~\ref{cor:37}).

In addition to that, if $A^{p',q'}_{\vect b+\vect d-\vect s-\vect{s'},0}(D)=\widetilde A^{p',q'}_{\vect b+\vect d-\vect s-\vect{s'},0}(D)$  and some simple necessary conditions are satisfied, then $P_{\vect{s'}}$ induces a continuous linear mapping of $L^{p,q}_{\vect s,0}(D)$ into $\widetilde A^{p,q}_{\vect s}(D)$ (cf.~Theorem~\ref{teo:7}).

\begin{deff}\label{40}
For every \CR $\vect s\in \vect b+\vect d-\frac 1 2 \vect m-(\R_+^*)^r$, \CB define the Bergman projector $P_{\vect s}\colon C_c(D)\to \Hol(D)$ by\index{Bergman projector}
\[
(P_{\vect s} f)(\zeta,z)\coloneqq c_{\vect s}\int_{D} f(\zeta',z') B^{\vect s}_{(\zeta',z')}(\zeta,z)\Delta_\Omega^{-\vect s}(\Im z'-\Phi(\zeta'))\,\dd\nu_D(\zeta',z')
\]
for every $f\in C_c(D)$ and for every $(\zeta,z)\in D$, where $c_{\vect s}\coloneqq \frac{\abs{\Pfaff(e_{\Omega'})}\Gamma_{\Omega'}(-\vect{s}) }{4^m \pi^{n+m} \Gamma_\Omega(\vect b+\vect d-\vect{s})}$. In addition, for every $\vect s\in \R^r$ we define $P_{\vect s,+}\colon C_c(D)\to C(D)$ by
\[
(P_{\vect s,+} f)(\zeta,z)\coloneqq \int_{D} f(\zeta',z') \abs*{B^{\vect s}_{(\zeta',z')}(\zeta,z)}\Delta_\Omega^{-\vect s}(\Im z'-\Phi(\zeta'))\,\dd\nu_D(\zeta',z')
\]
for every $f\in C_c(D)$ and for every $(\zeta,z)\in D$.
\end{deff}

Notice that we dropped the constant $c_{\vect s}$ from the definition of $P_{\vect s,+}$ to simplify the computations. 

In the following result we shall prove some necessary conditions for the boundedness of $P_{\vect{s'}}$ on $L^{p,q}_{\vect s}(D)$. Notice that the condition $p,q\Meg 1$ is a consequence of the fact that $L^{p,q}_{\vect s}(D)$ does \emph{not} embed in $L^1_\loc(\nu_D)$ if $\min(p,q)<1$.

\begin{prop}\label{prop:15}
Take \CR $\vect s\in \R^r$, $\vect {s'}\in \vect b+\vect d-\frac 1 2 \vect m-(\R_+^*)^r$, \CB and $p,q\in ]0,\infty]$, and assume that $P_{\vect {s'}}$ induces an endomorphism of $L^{p,q}_{\vect s,0}(D)$ (resp.\ a continuous linear mapping of $L^{p,q}_{\vect s,0}(D)$ into $L^{p,q}_{\vect s}(D)$). Then, the following hold:
\begin{itemize}
\item $p,q\Meg 1$;

\item $ \vect s\in \sup\left( \frac{1}{2 q} \vect{m}, \frac{1}{p}(\vect b+\vect d)+\frac{1}{2 q'}\vect{m'}  \right)+(\R_+^*)^r$ (resp.\ $\vect s\in  \R_+^r$ and $\vect s\in\Big(\frac{1}{p}(\vect b+\vect d)+\frac{1}{2 }\vect{m'}+(\R_+^*)^r\Big)$ if $q=\infty$);

\item $\vect{s'}\in \frac{1}{\min(p,p')} (\vect b+\vect d)-\frac{1}{2 \min(p,p')}\vect{m'}-(\R_+^*)^r$;

\item $\vect s+\vect{s'}\in \frac 1 p (\vect b+\vect d)-\frac{1}{2q}\vect{m'} -(\R_+^*)^r$;

\item $\vect s+\vect{s'}\in \vect b+\vect d-\frac{1}{2 q'}\vect{m}-(\R_+^*)^r$ if $q'<\infty$ and  $\vect s+\vect{s'}\in \vect b+\vect d-\R_+^r$ if $q'=\infty$.
\end{itemize}
\end{prop}

\CR Similar results hold for $P_{\vect{s'},+}$, without any initial restrictions on $\vect{s'}$.\CB

\begin{proof}
Since the linear mapping $f\mapsto f(0,i e_\Omega) $ is continuous on $A^{p,q}_{\vect s}(D)$, the assumptions imply that the mapping
\[
C_c(\Omega)\ni f\mapsto \int_D f(\zeta,z) B^{\vect{s'}}_{(\zeta,z)}(0,i e_\Omega) \Delta_\Omega^{-\vect{s'}}(\Im z-\Phi(\zeta))\,\dd \nu_D(\zeta,z)\in \C
\]
induces a continuous linear functional on $L^{p,q}_{\vect s,0}(D)$, so that $p,q\Meg 1$ and 
\[
B^{\vect{s'}}_{(\zeta,z)}(0,i e_\Omega) \in L^{p',q'}_{\vect b+\vect d-\vect s-\vect{s'}}(D)
\]
thanks to Proposition~\ref{prop:51}.\footnote{Argue directly to exclude the case in which $\min(p,q)<1$ and $\max(p,q)=\infty$.} Then, Proposition~\ref{prop:60} implies that the following conditions hold:
\begin{itemize}
\item $\vect s\in \frac 1 p(\vect b+\vect d)+\frac{1}{2 q'}\vect{m'}+(\R_+^*)^r$ if $q'<\infty$ and $\vect s\in \frac 1 p(\vect b+\vect d)+\R_+^r$ if $q'=\infty$;

\item $\vect{s'}\in \frac{1}{p'}(\vect b+\vect d)-\frac{1}{2 p'} \vect{m'}-(\R_+^*)^r$ if $p'<\infty$ and  $\vect{s'}\in -\R_+^r$ if $p'=\infty$;

\item $\vect s+\vect{s'}\in \vect b+\vect d-\frac{1}{2 q'}\vect m-(\R_+^*)^r$ if $q'<\infty$ and   $\vect s+\vect{s'}\in \vect b+\vect d-\R_+^r$ if $q'=\infty$.
\end{itemize}

Now, fix $r>0$ such that $\overline B_{E\times F_\C}((0,i e_\Omega),r)\subseteq D$ and choose $\tau\in C([0,r])$ so that $\tau(r)=0$ and $\int_{E\times F_\C} \tau(\abs{(\zeta,z)-(0,i e_\Omega)})\,\dd (\zeta,z)=1$. Then, define
\[
f\colon D\ni (\zeta,z)\mapsto\tau(\abs{(\zeta,z)-(0,i e_\Omega)})\Delta_\Omega^{\vect b+2\vect d+\vect{s'}}(\Im z-\Phi(\zeta))\in \C,
\]
so that clearly $f\in L^{p,q}_{\vect s,0}(D)$.
Now, the mapping $(\zeta',z')\mapsto \overline{B^{\vect{s'}}_{(\zeta',z')}(\zeta,z)}$ is holomorphic for every $(\zeta,z)\in D$, so that
\[
P_{\vect{s'}}(f)=\frac{\abs{\Pfaff(e_{\Omega'})}\Gamma_{\Omega'}(-\vect{s'}) }{4^m \pi^{n+m} \Gamma_\Omega(\vect b+\vect d-\vect{s'})} B^{\vect{s'}}_{(0,i e_\Omega)}
\] 
by Cauchy's integral formula, suitably applied.
Therefore, $B^{\vect{s'}}_{(0,i e_\Omega)}\in L^{p,q}_{\vect s,0}(D)$ (resp.\ $B^{\vect{s'}}_{(0,i e_\Omega)}\in L^{p,q}_{\vect s}(D)$), whence the other assertions thanks to Proposition~\ref{prop:60} and the preceding remarks. 
\end{proof}

In the following result we show that, even though $C_c(D)$ need not be dense in $L^{p,q}_{\vect s}(D)$, when $P_{\vect {s'}}$ induces a continuous linear mapping of $L^{p,q}_{\vect s,0}(D)$ into $L^{p,q}_{\vect s}(D)$, it is possible to define a canonical extension of $P_{\vect {s'}}$ to $L^{p,q}_{\vect s}(D)$ in a rather constructive way.
We then show that $P_{\vect {s'}}$ is self-adjoint in some sense.

\begin{prop}\label{prop:57}
Take \CR $p,q\in [1,\infty]$, $\vect s\in \R^r$, and $\vect{s'}\in \vect b+\vect d-\frac 1 2 \vect m-(\R_+^*)^r$, \CB and assume that $P_{\vect {s'}}$ induces a continuous linear mapping of $L^{p,q}_{\vect s,0}(D)$ into $L^{p,q}_{\vect s}(D)$. 
Define
\[
\langle f\vert g\rangle_{\vect{s'}}\coloneqq \int_D f(\zeta',z') \overline{g(\zeta',z')} \Delta_\Omega^{-\vect{s'}}(\Im z'-\Phi(\zeta'))\,\dd \nu_D(\zeta',z')
\]
for every $f\in L^{p,q}_{\vect s}(D)$ and for every $g\in L^{p',q'}_{\vect b+\vect d-\vect s-\vect{s'}}(D)$.
Then, the following hold:
\begin{enumerate}
\item[\em(1)] define, for every $f\in L^{p,q}_{\vect s}(D)$ and for every $(\zeta,z)\in D$,
\[
(P^{p,q}_{\vect s,\vect{s'}} f)(\zeta,z)\coloneqq c_{\vect{s'}} \big\langle f\big\vert B^{\vect{s'}}_{(\zeta,z)}\big\rangle_{\vect{s'}}
\]
where $c_{\vect{s'}}\coloneqq \frac{\abs{\Pfaff(e_{\Omega'})} \Gamma_{\Omega'}(-\vect{s'})}{4^m \pi^{n+m} \Gamma_\Omega(\vect b+\vect d-\vect s)}$. Then, $P^{p,q}_{\vect s,\vect{s'}}$ is a well-defined continuous linear projector of $L^{p,q}_{\vect s}(D)$ onto $A^{p,q}_{\vect s}(D)$;

\item[\em(2)] define, for every $g\in L^{p',q'}_{\vect b+\vect d-\vect s-\vect{s'}}(D)$ and for every $(\zeta,z)\in D$,
\[
(P^{p',q'}_{\vect b+\vect d-\vect s-\vect{s'},\vect{s'}} g)(\zeta,z)\coloneqq c_{\vect{s'}}\overline{\big\langle  B^{\vect{s'}}_{(\zeta,z)}\big\vert g\big\rangle_{\vect{s'}}}
\]
Then, $P^{p,q}_{\vect s,\vect{s'}}$ is a well-defined continuous linear projector of $L^{p',q'}_{\vect b+\vect d-\vect s-\vect{s'}}(D)$ onto $A^{p',q'}_{\vect b+\vect d-\vect s-\vect{s'}}(D)$;

\item[\em(3)] for every $f\in L^{p,q}_{\vect s}(D)$ and for every $g\in L^{p',q'}_{\vect b+\vect d-\vect s-\vect{s'}}(D)$,
\[
\begin{split}
\langle P^{p,q}_{\vect s,\vect{s'}} f\vert g\rangle_{\vect{s'}}=\big \langle f\big\vert P^{p',q'}_{\vect b+\vect d-\vect s-\vect{s'},\vect{s'}}g\big\rangle_{\vect{s'}}.
\end{split}
\]
\end{enumerate}
\end{prop}

In the sequel we shall simply write $P_{\vect {s'}}$ instead of $P^{p,q}_{\vect s,\vect{s'}}$ or $P^{p',q'}_{\vect b+\vect s-\vect s-\vect{s'},\vect{s'}}$.

\begin{proof}
\textsc{Step I.} Observe first that, arguing as in the proof of Lemma~\ref{lem:35} and using Proposition~\ref{prop:15},  we see that $P^{p,q}_{\vect s,\vect{s'}}$ and $P^{p',q'}_{\vect b+\vect d-\vect s-\vect{s'},\vect{s'}} $ are well defined and induce continuous linear mappings of $L^{p,q}_{\vect s}(D)$ and  $L^{p',q'}_{\vect b+\vect d-\vect s-\vect{s'}}(D)$ into $\Hol(D)$, respectively.

\textsc{Step II.} Let us prove that $P^{p,q}_{\vect s,\vect{s'}}$ is a projector of $L^{p,q}_{\vect s}(D)$ onto $A^{p,q}_{\vect s}(D)$. By~{step I} and Propositions~\ref{prop:9} and~\ref{prop:15}, it will suffice to show that $P^{p,q}_{\vect s, \vect{s'}}$ maps $L^{p,q}_{\vect s}(D)$ into itself.
Observe first that, by assumption, there is a constant $C>0$ such that
\[
\norm{P^{p,q}_{\vect s,\vect{s'}} f}_{L^{p,q}_{\vect s}(D)} =\norm{P_{\vect{s'}}f}_{L^{p,q}_{\vect s}(D)}\meg C \norm{f}_{L^{p,q}_{\vect s}(D)}
\]
for every $f\in C_c(D)$. Now, take $f\in L^{p,q}_{\vect s}(D)$ and choose two positive functions $\varphi\in C_c(\Nc)$ and $\psi\in C_c(\Omega)$ so that $\norm{\varphi}_{L^1(\Nc)}=\norm{\psi}_{L^1(F)}=1$. Define 
\[
\varphi_\rho(\zeta,x)\coloneqq r^{-(2 n+2 m)}\varphi(\rho^{-1} \zeta,\rho^{-2} x) \qquad \text{and}\qquad \psi_\rho(h)\coloneqq \rho^{-m} \psi_h(\rho^{-1}h)
\] 
for every $(\zeta,x)\in \Nc$, for every $h\in \Omega$, and for every $\rho>0$. In addition, define
\[
f_{\rho_1,\rho_2,\rho_3}(\zeta,z) \coloneqq \Delta_\Omega^{-\vect s}(h)\int_\Omega  [ \Delta_\Omega^{\vect s}(h') ( (\chi_{B((0, i e_\Omega),\rho_1)} f)_{h'}  )*\varphi_{\rho_2}](\zeta,x) \psi_{\rho_3}(h-h')\,\dd h'
\]
for every $(\zeta,z)\in D$ and for every $\rho_1,\rho_2,\rho_3>0$, where $x\coloneqq \Re z$ and $h\coloneqq \Im z-\Phi(\zeta)$.

Then, $f_{\rho_1,\rho_2,\rho_3}\in C_c(D)$ and 
\[
\norm{f_{\rho_1,\rho_2,\rho_3}}_{L^{p,q}_{\vect s}(D)}\meg \norm{f}_{L^{p,q}_{\vect s}(D)}
\]
for every $\rho_1,\rho_2,\rho_3>0$. In addition,
\[
\lim_{\rho_1\to +\infty}\lim_{\rho_2\to 0^+} \lim_{\rho_3\to 0^+} f_{\rho_1,\rho_2,\rho_3}=f
\]
in the weak topology $\sigma(L^{p,q}_{\vect s}(D), L^{p',q'}_{\vect b+\vect d-\vect s-\vect{s'}}(D))$, so that
\[
\lim_{\rho_1\to +\infty}\lim_{\rho_2\to 0^+} \lim_{\rho_3\to 0^+} P^{p,q}_{\vect s,\vect{s'} }f_{\rho_1,\rho_2,\rho_3}=P^{p,q}_{\vect s,\vect{s'}}f
\]
pointwise on $D$. Therefore,
\[
\norm{P^{p,q}_{\vect s,\vect{s'}}f}_{L^{p,q}_{\vect s}(D)}\meg C \norm{f}_{L^{p,q}_{\vect s}(D)}
\]
by lower semi-continuity.

\textsc{Step III.} Take $f,g\in C_c(D)$ and observe that
\[
\langle P^{p,q}_{\vect s,\vect{s'}} f\vert g\rangle_{\vect{s'}}= \big\langle f\big\vert P^{p',q'}_{\vect b+\vect d-\vect s-\vect{s'},\vect{s'}}g\big\rangle_{\vect{s'}}.
\]
By the arguments of~{step II}, the same holds for every $f\in L^{p,q}_{\vect s}(D)$. Now, observe that~\cite[Theorem 1 of \S 2]{BenedekPanzone} shows that the sesquilinear form $\langle \,\cdot\, \vert \,\cdot\,\rangle$ induces an antilinear isometry of $  L^{p',q'}_{\vect b+\vect d-\vect s-\vect{s'}}(D)$ into $L^{p,q}_{\vect s}(D)'$. Therefore,
\[
\norm{ P^{p',q'}_{\vect b+\vect d-\vect s-\vect{s'}}g  }_{L^{p',q'}_{\vect b+\vect d-\vect s-\vect{s'}}(D)}\meg C \norm{g}_{L^{p',q'}_{\vect b+\vect d-\vect s-\vect{s'}}(D)}
\]
for every $g\in C_c(D)$. Hence, the arguments of~{step II} show that $P^{p,q}_{\vect b+\vect d-\vect s-\vect{s'},\vect{s'}}$ is a continuous linear projector of $L^{p',q'}_{\vect b+\vect d-\vect s-\vect{s'}}(D)$ onto $A^{p',q'}_{\vect b+\vect d-\vect s-\vect{s'}}(D)$, whence~{(2)}. Assertion~{(3)} follows.
\end{proof}

\begin{prop}\label{prop:31}
Take $p,q\in [1,\infty]$ and $\vect s,\vect{s'}\in \R^r$. Then, the following conditions are equivalent:
\begin{enumerate}
\item[\em(1)] property $\atomic^{p,q}_{\vect s,\vect{s'},0,+}$ (resp.\  $\atomic^{p,q}_{\vect s,\vect{s'},+}$) holds;

\item[\em(2)] $P_{\vect{s'},+}$ induces an endomorphism of $L^{p,q}_{\vect s,0}(D)$ (resp.\ $L^{p,q}_{\vect s}(D)$).
\end{enumerate}

In particular, property $\atomic^{p,q}_{\vect s,\vect{s'},0,+}$ (resp.\ property $\atomic^{p,q}_{\vect s,\vect{s'},+}$)  implies property $\atomic^{p',q'}_{\vect b+\vect d-\vect s-\vect{s'},\vect{s'},+}$.
\end{prop}

\begin{proof}
Take a $(\delta,R)$-lattice $(\zeta_{j,k},z_{j,k})_{j\in J,k\in K}$ on $D$ for some $\delta>0$ and some $R>1$ (cf.~\ref{lem:32}), and define $h_k\coloneqq \Im z_{j,k}-\Phi(\zeta_{j,k})$ for every $k\in K$ and for some (hence every $j\in J$).

{(1) $\implies $ (2).} \CR Let us prove that the linear mapping
\[
T\colon L^{p,q}_{\vect s}(D)\ni f\mapsto \left(\int_{B((\zeta_{j,k},z_{j,k}),R\delta)} f(\zeta,z) \Delta^{\vect s-(\vect b+\vect d)/p}_\Omega(\Im z-\Phi(\zeta))\,\dd \nu_D(\zeta,z)  \right)_{j,k}\in \ell^{p,q}(J,K)
\]
is well-defined and continuous, and induces a continuous linear mapping of $L^{p,q}_{\vect s,0}(D)$ into $\ell^{p,q}_0(J,K)$.
Indeed, define $B_{j,k}\coloneqq B((\zeta_{j,k},z_{j,k}),R\delta)$ for every $j\in J$ and for every $k\in K$. In addition, observe that there are $N\in\N$ and a constant $C_1>0$ such that
\[
\sum_{j,k} \chi_{B_{j,k}}\meg N\chi_D, \qquad \sum_k \chi_{B_\Omega(h_k,R\delta)}\meg N\chi_\Omega,
\]
and such that
\[
\norm{(\chi_{B_{j,k}})_h}_{L^{p'}(\Nc)}\meg C_1 \Delta_\Omega^{-(\vect b+\vect d)/p'}(h) \qquad \text{and} \qquad \norm{\chi_{B_\Omega(h_k,R\delta)}}_{L^{q'}(\nu_\Omega)}\meg C_1
\]
for every $j\in J$, for every $k\in K$, and for every $h\in \Omega$ (cf.~Proposition~\ref{prop:56} and the proof of Theorem~\ref{teo:8}).	
Then,
\[
\begin{split}
	\norm{(T f)_{j,k}}_{\ell^p(J)}&\meg \int_{B_\Omega(h_k,R\delta)} \norm*{\int_\Nc \abs{f_h(\zeta,x)} (\chi_{B_{j,k}})_h(\zeta,x)\,\dd (\zeta,x)}_{\ell^p(J)} \Delta_\Omega^{\vect s+(\vect b+\vect d)/p'}(h)\,\dd \nu_\Omega(h)\\
	&\meg C_1\int_{B_\Omega(h_k,R\delta)} \norm*{\norm*{f_h (\chi_{B_{j,k}})_h}_{L^p(\Nc)}}_{\ell^p(J)}  \Delta_\Omega^{\vect s}(h)\,\dd \nu_\Omega(h)\\
	&=C_1\int_{B_\Omega(h_k,R\delta)} \norm*{f_h \norm*{(\chi_{B_{j,k}})_h}_{\ell^p(J)}}_{L^p(\Nc)}  \Delta_\Omega^{\vect s}(h)\,\dd \nu_\Omega(h)\\
	&\meg C_1 N^{1/p} \int_{B_\Omega(h_k,R\delta)} \norm*{f_h}_{L^p(\Nc)}  \Delta_\Omega^{\vect s}(h)\,\dd \nu_\Omega(h),
\end{split}
\]
where the first inequality follows from the convexity of the norm of $\ell^p(J)$, the second inequality follows from H\"older's inequality, and the last inequality follows from the choice of $N$.
Therefore,
\[
\begin{split}
	\norm{T f}_{\ell^{p,q}(J,K)}&\meg C_1 N^{1/p} \norm*{\int_{B_\Omega(h_k,R\delta)} \norm*{f_h}_{L^p(\Nc)}  \Delta_\Omega^{\vect s}(h)\,\dd \nu_\Omega(h)  }_{\ell^q(K)}\\
	&\meg C_1^2 N^{1/p} \norm*{\norm*{h\mapsto \norm*{f_h}_{L^p(\Nc)}  \Delta_\Omega^{\vect s}(h) \chi_{B_\Omega(h_k,R\delta)}(h)}_{L^q(\nu_\Omega)}  }_{\ell^q(K)}\\
	&=C_1^2 N^{1/p} \norm*{h\mapsto\norm*{f_h}_{L^p(\Nc)}  \Delta_\Omega^{\vect s}(h)\norm*{ \chi_{B_\Omega(h_k,R\delta)}(h) }_{\ell^q(K)} }_{L^q(\nu_\Omega)} \\
	&\meg C_1^2 N^{1/p+1/q} \norm{f}_{L^{p,q}_{\vect s}(D)},
\end{split}
\]
where the second inequality follows from H\"older's inequality, while the last inequality follows from the choice of $N$. To conclude, observe that $T(C_c(D))\subseteq \C^{(J\times K)}$, so that $T(L^{p,q}_{\vect s,0}(D))\subseteq \ell^{p,q}_0(J,K)$ by continuity.

Now, by Theorem~\ref{teo:2} and Corollary~\ref{cor:34} there is a constant $C_2>0$ such that
\[
\begin{split}
&\int_D \abs{f(\zeta',z') B^{\vect{s'}}_{(\zeta',z')}(\zeta,z)} \Delta^{-\vect {s'}}_\Omega(\Im z'-\Phi(\zeta'))\,\dd \nu_D(\zeta',z')\\
&\qquad\qquad\qquad\meg C_2\sum_{j,k} (T\abs{f})_{j,k}\abs{ B^{\vect{s'}}_{(\zeta_{j,k},z_{j,k})}}\Delta_\Omega^{(\vect b+\vect d)/p-\vect s-\vect{s'}}(h_k)
\end{split}
\]
so that by means of property $\atomic^{p,q}_{\vect s,\vect{s'},0,+}$ (resp.\ $\atomic^{p,q}_{\vect s,\vect{s'},+}$) we infer that $P_{\vect{s'},+}$ induces an endomorphism of $ L^{p,q}_{\vect s,0}(D)$ (resp.\ $ L^{p,q}_{\vect s}(D)$).

{(2) $\implies$ (1).} For every $(j,k)\in J\times K$, choose $\tau_{j,k}\in C_c(\Omega)$ so that 
\[
\chi_{B((\zeta_{j,k},z_{j,k}),\delta/2)}\meg \tau_{j,k}\meg \chi_{B((\zeta_{j,k},z_{j,k}),\delta)}.
\]
Then, the proof of Theorem~\ref{teo:5} shows that the mapping
\[
T'\colon \ell^{p,q}(J,K)\ni \lambda \mapsto \sum_{j,k} \lambda_{j,k} \tau_{j,k} \Delta_\Omega^{(\vect b+\vect d)/p-\vect s}\in L^{p,q}_{\vect s}(D)
\]
is well-defined and continuous, and that $T'(\ell^{p,q}_0(J,K))\subseteq L^{p,q}_{\vect s,0}(D)$.
Now, observe that Theorem~\ref{teo:2} and Corollary~\ref{cor:34} imply that there is a constant $C_3>0$ such that
\[
\sum_{j,k} \abs{\lambda_{j,k} B^{\vect{s'}}_{(\zeta_{j,k},z_{j,k})}}\Delta_\Omega^{(\vect b+\vect d)/p-\vect s-\vect{s'}}(h_k) \meg \frac{C_3}{\nu_D(B((0, i e_\Omega),\delta/2))}  P_{\vect {s'},+}(T'\abs{\lambda}).
\]
Since $P_{\vect{s'},+}$ induces an endomorphism of $L^{p,q}_{\vect s,0}(D)$ (resp.\ $L^{p,q}_{\vect s}(D)$), property $\atomic^{p,q}_{\vect s,\vect{s'},0,+}$ (resp.\  $\atomic^{p,q}_{\vect s,\vect{s'},+}$) follows.\CB
\end{proof}

\begin{cor}
Take $\vect{s},\vect{s'}\in \R^r$ and $p,q\in [1,\infty]$. Assume that the following conditions hold:
\begin{itemize}
\item $\vect s\in \frac{1}{2 q}\vect m+\frac{1}{2 q'}\vect{m'}+(\R_+^*)^r $;

\item  $\vect s+\vect{s'}\in \vect b+\vect d-\frac{1}{2 q'}\vect m-\frac{1}{2 q} \vect{m'}-(\R_+^*)^r$.
\end{itemize}
Then $P_{\vect{s'},+}$ induces endomorphisms of $L^{p,q}_{\vect s,0}(D)$ and $L^{p',q'}_{\vect b+\vect d-\vect s-\vect{s'},0}(D)$,  and of $L^{p,q}_{\vect s}(D)$ and $L^{p',q'}_{\vect b+\vect d-\vect s-\vect{s'}}(D)$ by transposition.
\end{cor}

This result covers~\cite[Theorem 2.1]{Nana}, which deals with the case in which $\vect{s'}=-q \vect s+\vect b+\vect d$.

\begin{proof}
This follows from Theorem~\ref{teo:5} and Proposition~\ref{prop:31}.
\end{proof}

\begin{prop}\label{prop:61}
Take $p,q\in [1,\infty]$ and $\vect s,\vect{s'}\in \R^r$, and assume that $P_{\vect{s'}}$ induces an endomorphism of $L^{p,q}_{\vect s,0}(D)$ (resp.\ $L^{p,q}_{\vect s}(D)$). Then properties $\atomic^{p,q}_{\vect s,\vect{s'},0}$, $\atomic^{p,q}_{\vect s,\vect{s'}}$, and $\atomics^{p',q'}_{\vect b+\vect d-\vect s-\vect{s'}}$ (resp.\ property $\atomic^{p,q}_{\vect s,\vect{s'}}$ and $\atomic^{p',q'}_{\vect b+\vect d-\vect s-\vect{s'}}$) hold.
\end{prop}

\begin{proof}
Let $(\zeta_{j,k},z_{j,k})_{j\in J, k \in K}$ be a $(\delta,R)$-lattice for some $\delta>0$ and some $R>1$. Define $h_k\coloneqq \Im z_{j,k}-\Phi(\zeta_{j,k})$ for every $k\in K$ and for some (hence every) $j\in J$.
Choose, for every $(j,k)\in J\times K$, an affine automorphism $\varphi_{j,k}$ of $D$ so that 
\[
\varphi_{j,k}(0,i e_\Omega)=(\zeta_{j,k},z_{j,k}),
\] 
so that
\[
\varphi_{j,k}(B((0, i e _\Omega), r))=B((\zeta_{j,k},z_{j,k}),r)
\] 
for every $r>0$. Fix $\delta'>0$ so that $B_{E\times F_\C}((0, i e_\Omega), \delta') \subseteq B((0, i e_\Omega), \delta)$, and define 
\[
B_{j,k}\coloneqq \varphi_{j,k}(B_{E\times F_\C}((0, i e_\Omega), \delta') )
\]
for every $(j,k)\in J\times K$. In addition, define
\[
\Szego(\lambda)(\zeta,z)\coloneqq \sum_{j,k} \lambda_{j,k} \chi_{B_{j,k}}(\zeta,z) \Delta_\Omega^{\vect{s'}-(\vect b+2 \vect d)}(\Im z-\Phi(\zeta)) \Delta_\Omega^{\vect b+2 \vect d+(\vect b+\vect d)/p-\vect s-\vect{s'}}(h_k)
\]
for every $(\zeta,z)\in D$ and for every $\lambda\in \C^{J\times K}$. Then, arguing as in the proof of Theorem~\ref{teo:5}, we see that $\Szego$ induces a continuous linear mapping $\ell^{p,q}(J,K)\to L^{p,q}_{\vect s}(D)$.
In addition, by Proposition~\ref{prop:57} and holomorphy we see that, for every $\lambda\in \ell^{p,q}(J,K)$,
\[
\begin{split}
P_{\vect{s'}}(\Szego(\lambda))&= c_{\vect{s'}}\sum_{j,k} \lambda_{j,k} \overline{\int_{B_{j,k}} B_{(\zeta,z)}^{\vect{s'}}(\zeta',z')\,\dd (\zeta',z')  }  \Delta_\Omega^{\vect b+2 \vect d+(\vect b+\vect d)/p-\vect s-\vect{s'}}(h_k)\\
&= c_{\vect{s'}} C \sum_{j,k} \lambda_{j,k}  B_{(\zeta_{j,k},z_{j,k})}^{\vect{s'}}(\zeta,z) \Delta_\Omega^{(\vect b+\vect d)/p-\vect s-\vect{s'}}(h_k)
\end{split}
\]
where $C\coloneqq \Hc^{2 n+2 m}(B_{E\times F_\C}((0,i e_\Omega),\delta'))$. Property $\atomic^{p,q}_{\vect s,\vect{s'}}$ follows. 
The proof is then completed by means of Proposition~\ref{prop:57}.
\end{proof}

We now show how the Bergman projectors $P_{\vect{s'}}$ interact with the spaces $\widetilde A^{p,q}_{\vect s}(D)$. 

\begin{teo}\label{teo:7}
Take $p,q\in [1,\infty]$ and $\vect s,\vect{s'}\in \R^r$ such that the following conditions hold:
\begin{itemize}
\item $\vect s\in \sup\left( \frac{1}{2 q}\vect m, \frac 1 p (\vect b+\vect d)+\frac{1}{2 q'}\vect{m'}  \right)+(\R_+^*)^r$;

\item $\vect s+\vect{s'}\in \inf\left(\vect b+\vect d-\frac{1}{2 q'} \vect m, \frac 1 p (\vect b+\vect d)- \frac{1}{2 q}\vect{m'}\right)-(\R_+^*)^r$;

\item $A^{p',q'}_{\vect b+\vect d-\vect s-\vect{s'},0}(D)=\widetilde A^{p',q'}_{\vect b+\vect d-\vect s-\vect{s'},0}(D)$.
\end{itemize}
Then, $P_{\vect{s'}}$ induces a continuous linear mapping of $L^{p,q}_{\vect s,0}(D)$ into $\widetilde A^{p,q}_{\vect s}(D)$.
\end{teo}

The proof is based on~\cite[Proposition 4.28]{BekolleBonamiGarrigosRicci}, which deals with the case in which $\vect s\in \R\vect 1_r$ and $D$ is an irreducible symmetric tube domain.

In particular, under the stated assumptions (and assuming further that $p,q<\infty$), $A^{p,q}_{\vect s}(D)=\widetilde A^{p,q}_{\vect s}(D)$ if and only if $P_{\vect{s'}}$ induces an endomorphism of $L^{p,q}_{\vect s}(D)$.

\begin{proof}
Take $f\in L^{p,q}_{\vect s,0}(D)\cap L^{2,2}_{(\vect b+\vect d-\vect{s'})/2}(D)$, and observe that $ (\vect b+\vect d-\vect{s'})/2\in \frac{1}{4 }\vect m+(\R_+^*)$ thanks to our assumptions, so that $P_{\vect{s'}}$ is the self-adjoint projector of $L^{2,2}_{(\vect b+\vect d-\vect{s'})/2}(D)$ onto $A^{2,2}_{(\vect b+\vect d-\vect{s'})/2}(D)$ by Proposition~\ref{prop:6}. 
Therefore, 
\[
P_{\vect{s'}}f\in A^{2,2}_{(\vect b+\vect d-\vect{s'})/2}(D),
\]
so that there is a unique $T\in B^{(\vect{s'}-\vect b-\vect d)/2}_{2,2}(\Nc,\Omega)$ such that 
\[
P_{\vect{s'}}f=\Ec T,
\]
thanks to Corollary~\ref{cor:29}. 
In addition, Proposition~\ref{prop:6} shows that there is a unique $\tau\in \Lc^2_{(\vect b+\vect d-\vect{s'})/2}(\Omega')$ (cf.~Definition~\ref{def:10}) such that 
\[
\pi_\lambda((\Ec T)_h)= e^{-\langle \lambda,h\rangle} \tau(\lambda)
\]
for almost every $\lambda\in \Omega'$ and for every $h\in \Omega$. 
Further,  clearly $T*\eta^*\in L^2(\Nc)$ for every $\eta\in \Sc_\Omega(\Nc)$.  Arguing as in the proof of Theorem~\ref{teo:11}, we then see that 
\[
\Ec(T*\eta^*)_h=(\Ec T)_h*\eta^*
\]
for every $h\in \Omega$, so that
\[
\pi_\lambda(T*\eta^*)=\tau(\lambda) \pi_\lambda(\eta)^*
\]
for almost every $\lambda\in \Omega'$. It is then easily seen that, for every $\eta\in \Sc_{\Omega,L}(\Nc)$ one has 
\[
T*\eta^*\in L^2(\Nc)\qquad \text{and} \qquad \pi_\lambda(T*\eta^*)=\tau(\lambda) \pi_\lambda(\eta)^*
\]
for almost every $\lambda\in \Omega'$.
By means of Corollary~\ref{cor:7} and Proposition~\ref{prop:55}, we also see that the mapping $\lambda \mapsto \pi_\lambda(\eta)$ belongs to $\Lc^{2}_{-(\vect b+\vect d-\vect{s'})/2}(\Omega')$, so that 
\[
\langle T\vert \eta\rangle=(T*\eta^*)(e)= \frac{2^{n-m}\abs{\Pfaff(e_{\Omega'})}}{\pi^{n+m}}\int_{\Omega'} \tr(\pi_\lambda(T*\eta^*)) \Delta^{-\vect b}_{\Omega'}(\lambda)\,\dd \lambda.
\]
Therefore, Proposition~\ref{prop:58} shows that
\[
\begin{split}
\langle T\vert \eta\rangle&=\frac{ 2^{-\vect{s'}}\abs{\Pfaff(e_{\Omega'})}}{4^m\pi^{n+m} \Gamma_\Omega(\vect b+\vect d-\vect{s'}  )}\int_\Omega \int_{\Omega'} \tr(\pi_\lambda(T*\eta^*)) \ee^{-2\langle \lambda,h\rangle} \Delta^{\vect d-\vect{s'}}_{\Omega'}(\lambda)\,\dd \lambda\\
&\qquad \times\Delta^{\vect b+\vect d-\vect{s'}}_\Omega(h)\,\dd \nu_\Omega(h).
\end{split}
\]
Now, the preceding remarks and  Proposition~\ref{prop:55} show that
\[
\begin{split}
\ee^{-2\langle \lambda,h\rangle}i^{\vect b+\vect d-\vect{s'}}\Delta^{\vect b+\vect d-\vect{s'}}_{\Omega'}(\lambda)\pi_\lambda(T*\eta^*)&=\pi_\lambda((\Ec T)_h) \ee^{-\langle \lambda,h\rangle} \pi_\lambda( \eta*I^{\vect{s'}-(\vect b+\vect d)}_\Omega )^*\\
&=\pi_\lambda((P_{\vect{s'}}f)_h) \pi_\lambda( \Ec(\eta*I^{\vect{s'}-(\vect b+\vect d)}_\Omega )_h )^* 
\end{split}
\]
for almost every $\lambda\in \Omega'$ and for every $h\in \Omega$. Therefore, Corollary~\ref{cor:7} implies that
\[
\langle T\vert \eta\rangle=c'_{\vect{s'}}\int_D (P_{\vect{s'}} f)(\zeta,z) \overline{\Ec(\eta*I^{\vect{s'}-(\vect b+\vect d)}_\Omega)(\zeta,z)} \Delta^{-\vect{s'}}_\Omega(\Im z-\Phi(\zeta))\,\dd \nu_D(\zeta,z)
\]
for a suitable constant $c'_{\vect{s'}}\neq 0$.
The Proposition~\ref{prop:57} then implies that 
\[
\langle T\vert \eta\rangle=c_{\vect{s'}}'\int_D f (\zeta,z) \overline{\Ec(\eta*I^{\vect{s'}-(\vect b+\vect d)}_\Omega)(\zeta,z)} \Delta^{-\vect{s'}}_\Omega(\Im z-\Phi(\zeta))\,\dd \nu_D(\zeta,z),
\]
so that
\[
\abs{\langle T\vert \eta\rangle}\meg \abs{c'_{\vect{s'}}} \norm{f}_{L^{p,q}_{\vect s}(D)} \norm*{ \Ec(\eta*I^{\vect{s'}-(\vect b+\vect d)}_\Omega) }_{L^{p',q'}_{\vect b+\vect d-\vect s-\vect{s'}}(D)}.
\]
Now, by assumption there is a constant $C_1>0$ such that
\[
\norm*{ \Ec(\eta*I^{\vect{s'}-(\vect b+\vect d)}_\Omega) }_{L^{p',q'}_{\vect b+\vect d-\vect s-\vect{s'},}(D)}\meg C_1 \norm*{\eta*I^{\vect{s'}-(\vect b+\vect d)}_\Omega}_{ B^{\vect{s}+\vect{s'}-(\vect b+\vect d)}_{p',q'}(\Nc,\Omega)},
\]
while Theorem~\ref{teo:9} shows that there is a constant $C_2>0$ such that
\[
\norm*{\eta*I^{\vect{s'}-(\vect b+\vect d)}_\Omega}_{B^{\vect{s}+\vect{s'}-(\vect b+\vect d)}_{p',q'}(\Nc,\Omega)}\meg C_2 \norm{\eta}_{B^{\vect{s}}_{p',q'}(\Nc,\Omega)}.
\]
By the arbitrariness of $\eta$, Theorem~\ref{teo:10} shows that there is a constant $C_3>0$ such that
\[
\norm{T}_{B^{-\vect s}_{p,q}(\Nc,\Omega)}\meg C_3   \norm{f}_{L^{p,q}_{\vect s}(D)},
\]
that is,
\[
\norm{P_{\vect{s'}} f}_{\widetilde A^{\vect s}_{p,q}(D)}=\norm{\Ec T}_{\widetilde A^{\vect s}_{p,q}(D)}\meg C_3   \norm{f}_{L^{p,q}_{\vect s}(D)}.
\]
The assertion follows by means of Proposition~\ref{prop:8}.
\end{proof}

We now draw some consequences of Theorem~\ref{teo:7}. In the next result we prove that, if $A^{p,q}_{\vect s,0}(D)=\widetilde A^{p,q}_{\vect s,0}(D)$, then one may improve Theorem~\ref{teo:7} and show that $P_{\vect{s'}}$ induces an endomorphism of $L^{p,q}_{\vect s,0}(D)$.

Observe that the assumption on $\vect{s'}$, which does not appear in Theorem~\ref{teo:7}, is necessary only when $p=\infty$. Otherwise, we would only be able to show that $P_{\vect{s'}}$ induces a continuous linear mapping of $L^{p,q}_{\vect s,0}(D)$ into $L^{p,q}_{\vect s}(D)$.

\CR
\begin{cor}\label{cor:38}
Take $p,q\in [1,\infty]$ and $\vect s,\vect{s'}\in \R^r$ such that the following conditions hold:
\begin{itemize}
\item $\vect s\in \sup\left( \frac{1}{2 q}\vect m, \frac 1 p (\vect b+\vect d)+\frac{1}{2 q'}\vect{m'}  \right)+(\R_+^*)^r$;

%\item $\vect{s'}\in \frac 1 p (\vect b+\vect d)-\frac{1}{2 p} \vect{m'}-(\R_+^*)^r$;

\item $\vect s+\vect{s'}\in \inf\left(\vect b+\vect d-\frac{1}{2 q'} \vect m, \frac 1 p (\vect b+\vect d)- \frac{1}{2 q}\vect{m'}\right)-(\R_+^*)^r$;

\item $A^{p,q}_{\vect s,0}(D)=\widetilde A^{p,q}_{\vect s,0}(D)$;

\item $A^{p',q'}_{\vect b+\vect d-\vect s-\vect{s'},0}(D)=\widetilde A^{p',q'}_{\vect b+\vect d-\vect s-\vect{s'},0}(D)$.
\end{itemize}
Then, $P_{\vect{s'}}$ induces an endomorphism of $L^{p,q}_{\vect s,0}(D)$.
\end{cor}

\begin{proof}
By Theorem~\ref{teo:7}, we know that $P_{\vect{s'}}$ induces a continuous linear mapping of $L^{p,q}_{\vect s,0}(D)$ into $L^{p,q}_{\vect s}(D)$. In addition, Proposition~\ref{prop:81} implies that $\vect s\in \left( \frac{1}{2 p}-\frac{1}{2 q} \right)\vect {m'}+\R_+^r$. Since $\vect s+\vect{s'}\in  \frac 1 p (\vect b+\vect d)- \frac{1}{2 q}\vect{m'}-(\R_+^*)^r$ by assumption, this implies that $\vect{s'}\in \frac 1 p (\vect b+\vect d)-\frac{1}{2 p} \vect{m'}-(\R_+^*)^r$. Therefore, it is readily seen that $P_{\vect{s'}}(C_c(D))\subseteq L^{p,q}_{\vect s,0}(D)$ (cf.~Proposition~\ref{prop:60}), whence the result.
\end{proof}
\CB

We now provide some explicit conditions for the boundedness of  $P_{\vect{s'}}$ on $L^{p,q}_{\vect s}(D)$. This result covers~\cite[Theorems 2.2 and 2.3]{Nana}, which deal with the case $\vect{s'}=-q\vect s+\vect b+\vect d$ with a different proof.

\begin{cor}\label{cor:35}
Take $p,q\in [1,\infty]$ and $\vect s,\vect{s'}\in \R^r$ such that the following conditions hold:
\begin{itemize}
\item $\vect s\in \sup\left( \frac{1}{2 q}\vect m  + \left(\frac{1}{2 \min(p,p')}-\frac{1}{2 q} \right)_+\vect{m'}  , \frac 1 p (\vect b+\vect d)+\frac{1}{2 q'}\vect{m'}  \right)+(\R_+^*)^r$;

\item $\vect s+\vect{s'}\in \inf\left(\vect b+\vect d-\frac{1}{2 q'} \vect m-\left(\frac{1}{2 \min(p,p')}-\frac{1}{2 q'} \right)_+\vect{m'}, \frac 1 p (\vect b+\vect d)- \frac{1}{2 q}\vect{m'}\right)-(\R_+^*)^r$.
\end{itemize}
Then, $P_{\vect{s'}}$ induces a endomorphisms of $L^{p,q}_{\vect s,0}(D)$ and $L^{p,q}_{\vect s}(D)$.
\end{cor}

\begin{proof}
This follows from Corollaries~\ref{cor:29} and~\ref{cor:38}, Proposition~\ref{prop:57}, and Theorem~\ref{teo:7}.
\end{proof}

The following result completes the equivalences of Corollary~\ref{cor:32} in the case $p,q\in[1,\infty]$.

\begin{cor}\label{cor:37}
Take $p,q\in [1,\infty]$ and $\vect s\in \frac 1 p(\vect b+\vect d)+\frac{1}{2 q'}\vect{m'}+(\R_+^*)^r$ such that $\vect s\in \frac{1}{2 q}\vect m+(\R_+^*)^r$ (resp.\ $\vect s\in \R_+^r$ if $q=\infty$). Then, the following conditions are equivalent:
\begin{enumerate}
\item[\em(1)] $A^{p,q}_{\vect s,0}=\widetilde A^{p,q}_{\vect s,0}$ (resp.\ $A^{p,q}_{\vect s}=\widetilde A^{p,q}_{\vect s}$);

\item[\em(2)] $P_{\vect{s'}}$ induces an endomorphism of $L^{p,q}_{\vect s,0}(D)$ (resp.\ $L^{p,q}_{\vect s}(D)$) for every $\vect{s'}$ in a translate of $-\R_+^r$.
\end{enumerate}
\end{cor}

\begin{proof}
{(1) $\implies$ (2).} By Corollary~\ref{cor:29}, there is $\vect{s'_0}\in \frac 1 p (\vect b+\vect d)+\frac{1}{2 p}\vect{m'}-(\R_+^*)^r$ such that
\[
\vect s+\vect{s'_0}\in \inf\left(\vect b+\vect d-\frac{1}{2 q'} \vect m, \frac 1 p (\vect b+\vect d)- \frac{1}{2 q}\vect{m'}\right)-(\R_+^*)^r
\] 
and such that
\[
A^{p,q}_{\vect b+\vect d-\vect s-\vect{s'},0}=\widetilde A^{p,q}_{\vect b+\vect d-\vect s-\vect{s'},0}
\]
for every $\vect{s'}\in \vect{s'_0}-\R_+^r$. Then, Corollary~\ref{cor:38} (resp.\ Proposition~\ref{prop:57} and Theorem~\ref{teo:7}),  implies that  $P_{\vect{s'}}$ induces an endomorphism of $L^{p,q}_{\vect s,0}(D)$ (resp.\ $L^{p,q}_{\vect s}(D)$) for every $\vect{s'}\in \vect{s'_0}-\R_+^r$.

{(2) $\implies$ (1).} This follows from Corollary~\ref{cor:32} and Proposition~\ref{prop:61}.
\end{proof}

\section{Notes and Further Results}

\paragraph{5.3.1} With the notation of Theorem~\ref{teo:3}, the Bergman space $A^{p,q}_{\vect s}(D)$ is said to satisfy the interpolation property if the mapping $S$ is onto (hence an isomorphism for sufficiently fine lattices). Cf.~\cite{Rochberg,BekolleKagou2,BekolleIshiNana} for some results in this direction.

\paragraph{5.3.2} \label{sec:5:3:2} Recall that property $\atomic^{p,q}_{\vect s,\vect{s'},+}$ implies property $\atomic^{p,q}_{\vect{s_0}, \vect{s_0'},+}$ for every $\vect{s_0}\in \vect s+\R_+^r$ and for every $\vect{s'_0}\in \vect{s'}+\vect s-\vect{s_0}-\R_+^r$, thanks to Proposition~\ref{prop:16}. It would be interesting to investigate if the same happens for properties $\atomic^{p,q}_{\vect s,\vect{s'}}$ and $\atomics^{p,q}_{\vect s,\vect{s'}}$. This fact would improve considerably the statement of Corollary~\ref{cor:32}.

Analogous considerations hold for the boundedness of the Bergman projectors $P_{\vect{s'}}$, in view of Proposition~\ref{prop:31}.

\paragraph{5.3.3}   In connection with the discussion of~5.3.2, it would be interesting to understand whether the equality $A^{p,q}_{\vect s}(D)=\widetilde A^{p,q}_{\vect s}(D)$ implies the equality $A^{p,q}_{\vect{s'}}(D)=\widetilde A^{p,q}_{\vect{s'}}(D)$ for every $\vect{s'}\in \vect s+\R_+^r$. This problem is closely related to the boundedness of Riemann--Liouville operators between the spaces $A^{p,q}_{\vect s}(D)$. 
In Corollary~\ref{cor:40}, we proved that convolution with $I^{-\vect{s'}}_{\Omega}$ maps $A^{p,q}_{\vect s}(D)$ into $A^{p,q}_{\vect s+\vect{s'}}(D)$ continuously when $\vect{s'}\in \N_{\Omega'}$. Nonetheless, in the case of general $\vect{s'}\in \R_+^r$, the assertion is less clear.

As observed in~\cite{Bekolleetal}, the fact that $I^{-\vect{s'}}_{\Omega}$ induces an isomorphism of $A^{p,q}_{\vect s}(D)$ onto $A^{p,q}_{\vect s+\vect{s'}}(D)$ can be considered as a generalized Hardy's inequality.

\paragraph{5.3.4} \label{sec:5:3:4} \CR In Proposition~\ref{cor:30} we showed that one may define a continuous sesquilinear \CB form on $\widetilde A^{p,q}_{\vect s}(D)\times \widetilde A^{p',q'}_{(\vect b+\vect d)/\min(1,p)-\vect s-\vect{s'}}(D)$ which extends the mapping
\[
(f,g)\mapsto \int_D f(\zeta,z) \overline{g(\zeta,z)} \Delta_\Omega^{-\vect{s'}}(\Im z-\Phi(\zeta))\,\dd \nu_D(\zeta,z)
\]
on $A^{p,q}_{\vect s}(D)\times  A^{p',q'}_{(\vect b+\vect d)/\min(1,p)-\vect s-\vect{s'}}(D)$, and which induces an antilinear isomorphism
\[
\widetilde A^{p',q'}_{(\vect b+\vect d)/\min(1,p)-\vect s-\vect{s'}}(D)\to \widetilde A^{p,q}_{\vect s,0}(D)'.
\]
Thus, we have canonical injective continuous linear mappings
\[
A^{p',q'}_{(\vect b+\vect d)/\min(1,p)-\vect s-\vect{s'}}(D)\to \widetilde A^{p',q'}_{(\vect b+\vect d)/\min(1,p)-\vect s-\vect{s'}}(D)\to A^{p,q}_{\vect s,0}(D)',
\]
none of which is an isomorphism, in general. 

In a similar way, for every $\vect{s''},\vect{s'''}\in \N_{\Omega'}$ such that $\vect s+\vect{s''},(\vect b+\vect d)/\min(1,p)-\vect s-\vect{s'}+\vect{s'''} \in \frac 1 p (\vect b+\vect d)+\frac{1}{2 q'}\vect{m'}+(\R_+^*)^r$ one may also define a sesquilinear form on $\widehat A^{p,q}_{\vect s,\vect{s''}}(D)\times \widehat A^{p',q'}_{(\vect b+\vect d)/\min(1,p)-\vect s-\vect{s'},\vect{s'''}}(D)$ which induces an antilinear isomorphism of the space $\widehat A^{p',q'}_{(\vect b+\vect d)/\min(1,p)-\vect s-\vect{s'},\vect{s''}}(D)$ onto the space $\widehat A^{p,q}_{\vect s,\vect{s'''},0}(D)'$.

\paragraph{5.3.5}  When $D=\C_+$, the dual of $A^{1,1}_{s}(\C_+)$ has been identified with the classical Bloch space (modulo constants)  $\widehat A^{\infty,\infty}_{0,1}(\C_+)$ by Coifman and Rochberg in~\cite{CoifmanRochberg}. In a series of papers, Békollé, Temgoua Kagou, and Zhu extended this characterization to Bergman spaces of type $L^p$, $p\in ]0,1]$ on more general domains (cf.~\cite{Bekolle2,Kagou,Zhu1,Zhu2,Bekolle} and the references therein). With the previous notation, the dual of $A^{p,p}_{\vect s}(D)$ for $p\in ]0,1]$ has been identified with $\widehat A^{\infty,\infty}_{\vect 0,\vect{s'}}(D)$ by means of the sesquilinear form indicated in~5.3.4 above (for every $\vect{s'}\in \N_{\Omega'}\cap\left(\frac{1}{2}\vect{m'}+(\R_+^*)^r \right) $).
In addition to that, observe that Corollary~\ref{cor:35} easily implies that $P_{\vect{s''}}$ induces a (unique) continuous linear mapping $C_0(D)\to \widehat A^{\infty,\infty}_{\vect 0,\vect{s'}}(D)$ if $\vect{s''}\in \vect b+\vect d-\frac 1 2 \vect m-(\R_+^*)^r$. This mapping can be, in turn, suitably extended to a continuous linear mapping of $L^\infty(D)$ onto $\widehat A^{\infty,\infty}_{\vect 0,\vect{s'}}(D)$.\footnote{This is due to the fact that $(P_{\vect{s''}} f)* I^{-\vect{s'}}_\Omega = c_{\vect{s'},\vect{s''}} P_{\vect{s''}-\vect{s'}} (f (\Delta_\Omega^{-\vect{s'}}\circ \rho) )$ for every $f\in C_c(D)$, where $\rho \colon D\ni (\zeta,z)\mapsto \Im z-\Phi(\zeta)\in \Omega$ and $c_{\vect{s'},\vect{s''}}$ is a suitable constant.} This latter fact extends to general homogeneous Siegel domains of type II the results of the papers cited above.

Notice that the generalized Bloch space defined above is different from the generalized Bloch space $\Bc$ initially introduced by Timoney in~\cite{Timoney}, when $r>1$. If $D'$ is a \emph{bounded} homogeneous domain, the space $\Bc(D')$ is defined in~\cite{Timoney} as the space of holomorphic functions on $D'$ which are Lipschitz for the distance induced by the Bergman metric (modulo constants). Since this space is then invariant under composition with holomorphic automorphisms of $D'$, a similar space may be defined also on the homogeneous Siegel domain $D$. 
When $D=\C_+\times \C_+$, then the function 
\[
f\colon (z_1,z_2)\mapsto \log(z_1)\log(z_2), 
\]
where $\log$ denotes the unique holomorphic function on $\C\setminus \R_-$ which coincides with the  usual logarithm on $\R_+^*$, clearly belongs to $\widehat A^{\infty,\infty}_{\vect 0, \vect 1_2}(D)$ (more precisely, $\partial_1\partial_2 f= f * I^{-\vect 1_2}_{(\R_+^*)^2}\in A^{\infty,\infty}_{\vect 1_2}(D)$, so that $f$ induces an elements of $\widehat A^{\infty,\infty}_{\vect 0, \vect 1_2}(D)$).
Nonetheless, it is not hard to prove that there are no holomorphic functions $g$ on $D$ such that $\partial_1\partial_2 g=0$ and $f+g\in \Bc(D)$, that is (cf.~\cite[the proof of the equivalence (1) $\iff$ (3) of Theorem 3.4]{Timoney}),
\[
\sup\limits_{(z_1,z_2)\in D} \left( (\Im z_1)^2 \abs{\partial_1 (f+g)(z_1,z_2)}^2+ (\Im z_2)^2 \abs{\partial_2 (f+g)(z_1,z_2)}^2\right) <\infty.
\]
Consequently, the spaces $\widehat A^{\infty,\infty}_{\vect 0, \vect 1_2}(D)$ and $\Bc(D)$ are not canonically isomorphic.

\paragraph{5.3.6}  If one is only interested in the Bergman spaces $A^{p,p}_{\vect s}(D)$, then the problems considered in Chapter~\ref{sec:6} and in this one  are completely solved when $p\in ]0,1]$ in view of Theorem~\ref{teo:5} and Corollary~\ref{cor:29}, and partially solved for $p\in ]0,2]$ in view of Corollaries~\ref{cor:29},~\ref{cor:32}, and~\ref{cor:37}.
For $p\Meg 2$, Theorem~\ref{teo:7} shows that $P_{\vect{s'}}$ induces a continuous linear mapping of $L^{p,p}_{\vect s,0}(D)$ into $\widetilde A^{p,p}_{\vect s}(D)$, provided that some \emph{necessary} conditions on $\vect s$ and $\vect{s'}$ are satisfied. 
In full generality, we do not know whether $A^{p,p}_{\vect s}(D)=\widetilde A^{p,p}_{\vect s}(D)$, hence we do not know whether $P_{\vect{s'}}$ is bounded on $L^{p,p}_{\vect s}(D)$. Cf.~\cite{BekolleBonamiGarrigosRicci} for a sharper discussion when $D$ is a tube domain over a light cone.

\appendix

\chapter{Mixed Norm Spaces}\label{sec:app}

In this chapter we collect some results on mixed norm spaces $L^{p,q}(\mi,\nu)$, $p,q\in ]0,\infty]$. Our main focus is in the correspondences between $L^{p,q}(\mi,\nu)$ and $L^q(\nu;L^p(\mi))$ and between $L^{p,q}_0(\mi,\nu)$ and $L^q_0(\nu;L^p_0(\mi))$ (Propositions~\ref{prop:49} and~\ref{prop:50}) and in the characterization of the dual of $L^{p,q}_0(\mi,\nu)$ (Propositions~\ref{prop:51}, and~\ref{prop:52}).

We begin with recalling some notions on locally bounded $F$-spaces.

\begin{deff}\label{def:7}
An $F$-space is a complete metrizable topological vector space.\index{F-space@$F$-space} An $F$-space is called locally bounded if it admits a bounded neighbourhood of $0$.\index{F-space@$F$-space!locally bounded}
\end{deff}

Recall that an $F$-space $Z$ is locally bounded if and only if there are $p\in [1,\infty[$ and an absolutely homogeneous function $\norm{\,\cdot\,}\colon Z\to \R_+$ such that the mapping
\[
(z,z')\mapsto \norm{z-z'}^p
\]
is a distance compatible with the topology of $Z$ (cf.~\cite[Theorem 3.2.1]{Rolewicz}). In particular, $\norm{\,\cdot\,}$ is a continuous quasi-norm on $Z$. In addition, if $q\in [1,\infty[$ and $\norm{\,\cdot\,}'\colon Z\to \R_+$ is another absolutely homogeneous function such that the mapping $(z,z')\mapsto \norm{z-z'}'^q$ is a distance compatible with the topology of $Z$, then there is a constant $C>0$ such that
\[
\frac{1}{C}\norm{z}\meg \norm{z}'\meg C\norm{z}
\]
for every $z\in Z$ (cf.~\cite[Theorem 3.2.13]{Rolewicz}). For this reason, we shall sometimes denote by $\norm{\,\cdot\,}_Z$ one such quasi-norm, even though it is not uniquely defined.

\begin{deff}\label{41}
	If $\mi$ is a measure on a locally compact space $X$, $Z$ is a locally bounded $F$-space, and $f\colon X\to Z$ is a (not necessarily $\nu$-measurable) mapping, by an abuse of notation we shall define\footnote{Recall that, for every positive function $g$ on $X$, the symbol $\int_X^* g\,\dd \mi$ denotes the upper integral of $g$, that is, the greatest lower bound of the upper integrals $\int^*_X h\,\dd \mi$, where $h$ is a lower semi-continuous function and $h\Meg g$. The symbol $\int_X^* h\,\dd \mi$, in turn, denotes the smallest upper bound of the integrals $\int_X \varphi\,\dd \mi$, where $\varphi\in C_c(X)$ and $\varphi\meg h$.}
	\[
	\norm{f}_{L^p(\mi;Z)}\coloneqq \left(\int_X^* \norm{f}_Z^p\,\dd \mi  \right)^{1/p}
	\]
	for every $p\in ]0,\infty[$, and we shall denote by $\norm{f}_{L^\infty(\nu;Z)}$ the essential supremum of $\norm{f}_Z$. We define 
	\[
	\Lc^p(\mi;Z)\coloneqq \Set{f\colon X\to \C\colon \text{$f$ is $\mi$-measurable},   \norm{f}_{L^p(\mi;Z)}<\infty}.
	\]
	Denote by $L^p(\mi;Z)$ the Hausdorff space associated with $\Lc^p(\mi;Z)$. Analogously, define 
	\[
	\Lc^p_\loc(\mi;Z)\coloneqq \Set{f \colon \chi_K f \in \Lc^p(\mi;Z) \text{ for every compact subset $K$ of $X$}}.
	\]
	Denote by $L^p_\loc(\mi;Z)$ the Hausdorff space associated with $\Lc^p_\loc(\mi;Z)$.
	
	We define $L^p_0(\mi;Z)$ and $L^p_{0,\loc}(\mi;Z)$ as the closure of $C_c( X;Z)$ in $L^p(\mi;Z)$ and $L^p_\loc(\mi;Z)$, respectively.\label{42} 
\end{deff}

Notice that, if $Z$ is not a locally convex space, then the mapping $f \mapsto \int_X f\,\dd \mi$, defined on the set of $\mi$-measurable step function in $L^1(\mi;Z)$, is not continuous for the norm $\norm{\,\cdot\,}_{L^1(\mi;Z)}$, in general. In particular, it need not extend to a continuous linear mapping on $L^1(\mi;Z)$.

With standard techniques one then proves the following results.

\begin{prop}\label{prop:46}
	Let $\mi$ be a measure on a locally compact space $X$, $Z$ a locally bounded $F$-space, and take $p\in ]0,\infty]$. Then, $L^p(\mi;Z)$ is a locally bounded $F$-space.
\end{prop}

\begin{prop}\label{prop:47}
	Let $\mi$ be a measure on a locally compact space $X$, $Z$ a locally bounded $F$-space, and take $p\in ]0,\infty[$. Then, 
	\[
	L^p_0(\mi;Z)=L^p(\mi;Z) \qquad \text{and} \qquad L^p_{0,\loc}(\mi;Z)=L^p_\loc(\mi;Z),
	\]
	while 
	\[
	L^\infty_0(\mi;Z)=C_0(\Supp{\mi};Z)\qquad \text{and}\qquad L^\infty_{0,\loc}(\mi;Z)=C(\Supp{\mi};Z).
	\]
\end{prop}

We now pass to mixed norm spaces.

\begin{deff}\label{def:3}
	If $\mi$ and $\nu$ are two Radon measures on two locally compact spaces $X$ and $Y$, respectively, then for every $p,q\in ]0,\infty]$ we define the mixed norm space\index{Mixed norm space} 
	\[
	\Lc^{p,q}(\mi,\nu)\coloneqq \Set{f\colon \text{$f$ is $(\mi\otimes \nu)$-measurable}, \\\norm*{y\mapsto \norm{f(\,\cdot\,,y)}_{L^p(\mi)} }_{L^q(\nu)}<\infty},
	\]
	endowed with the corresponding topology. We denote by $L^{p,q}(\mi,\nu)$ the Hausdorff space associated with $\Lc^{p,q}(\mi,\nu)$, and by $L^{p,q}_0(\mi,\nu)$ the closure of $C_c(X\times Y)$ in $L^{p,q}(\mi,\nu)$.
\end{deff}

\begin{prop}\label{prop:62}
	Let $X$ and $Y$ be two locally compact spaces, and $\mi,\nu$ two Radon measures on $X$ and $Y$, respectively. Take $p,q\in ]0,\infty]$. Then, $L^{p,q}(\mi,\nu)$ is a locally bounded $F$-space.
\end{prop}

We now discuss the relationship between $L^{p,q}(\mi,\nu)$ and $L^q(\nu;L^p(\mi))$. 
Notice that, even though $L^{p,q}_0(\mi,\nu)$ is canonically isomorphic to $L^{q}_0(\nu;L^p_0(\mi))$, the space  $L^{p,q}(\mi,\nu)$ fails to be canonically isomorphic to $L^{q}(\nu;L^p(\mi))$, in general, for a lack of measurability. Indeed, if $f\in L^{\infty,q}(\mi,\nu)$, then the mapping $y\mapsto f(\,\cdot\,,y)\in L^\infty(\mi)$ is not $\nu$-measurable, in general. Roughly speaking, the best one can say is that the mapping $y\mapsto \int_X f(x,y) g(x)\,\dd \mi(x)$ is $\nu$-measurable for every $g\in L^1(\mi)$ (compare this fact with the Dunford--Pettis theorem).
These results can be proved with standard techniques. The proofs are omitted.

\begin{prop}\label{prop:49}
	Let $X$ and $Y$ be two locally compact spaces, and $\mi,\nu$ two Radon measures on $X$ and $Y$, respectively. Take $p,q\in ]0,\infty]$.  Then, there is an isometry 
	\[
	T\colon L^q(\nu;L^p(\mi))\to L^{p,q}(\mi,\nu)
	\]
	such that 
	\[
	(T f)(\,\cdot\,,y)= f(y)
	\]
	for $\nu$-almost every $y\in Y$. If $p<\infty$ or $\nu$ is atomic, then $T$ is onto.
\end{prop}

Recall that a Radon measure $\mi$ is discrete or atomic if $\mi(K)= \sum_{x\in K} \mi(\Set{x})$ for every compact set $K$, and that $\mi$ is diffuse if $\mi(\Set{x})=0$ for every $x$. \index{Measure!discrete}\index{Measure!atomic} \index{Measure!diffuse}
Then, every Radon measure can be written in a unique way as the sum of an atomic and a diffuse measure, cf.~\cite[Proposition 15 of Chapter V, \S 5, No.\ 10]{BourbakiInt1}.

\begin{prop}\label{prop:50}
	Keep the hypotheses and the notation of Proposition~\ref{prop:49}. Then, $T$ induces an isometry of $L^q_0(\nu;L^p_0(\mi))$ onto $L^{p,q}_0(\mi,\nu)$.
	
	If, in addition, $q<\infty$ and $X$ has a countable base, then 
	\[
	L^{\infty,q}_0(\nu,\mi)=\Set{f\in L^{\infty,q}(\nu,\mi)\colon f(\,\cdot\,,y)\in C_0(\Supp{\mi}) \text{ for $\nu$-almost every $y\in Y$}}.
	\]
\end{prop}

We now characterize the dual of $L^{p,q}_0(\mi,\nu)$. For the sake of simplicity, we shall not present a general result, but we shall content ourselves with the cases we are concerned with (namely, the cases in which both $\mi$ and $\nu$ are diffuse measures or counting measures).

\begin{deff}\label{43}
	Define $p'\coloneqq (\max(1,p))'$ for $p\in ]0,\infty]$, so that $p'=\infty$ when $p\in ]0,1]$ while $\frac{1}{p}+\frac{1}{p'}=1$ when $p\in [1,\infty]$.
\end{deff}

\begin{prop}\label{prop:51}
	Let $X$ and $Y$ be two locally compact spaces and $\mi,\nu$ two diffuse Radon measures on $X$ and $Y$, respectively. Take $p,q\in ]0,\infty[$.
	Then, the following hold:
	\begin{enumerate}
		\item[\em(1)] if $\min(p,q)<1$, then $ L^{p,q}(\mi,\nu)'=\Set{0} $;
		
		\item[\em(2)] if $p,q\Meg 1$, then the bilinear mapping
		\[
		L^{p,q}(\mi,\nu)\times L^{p',q'}(\mi,\nu)\ni (f,g)\mapsto \int_{X\times Y} f g\,\dd (\mi\otimes \nu)
		\]
		induces an isometry of $L^{p',q'}(\mi,\nu)$ onto $L^{p,q}(\mi,\nu)'$.
	\end{enumerate}
\end{prop}

The first assertion is proved as in the case $p=q$ (cf.~\cite{Day}). The second assertion is a particular case of~\cite[Theorem 1 of \S 3]{BenedekPanzone}.

\begin{prop}\label{prop:52}
	Let $X$ and $Y$ be two discrete spaces, and $\mi,\nu$ the counting measures on $X$ and $Y$, respectively. Take $p,q\in ]0,\infty]$. Then, the bilinear mapping
	\[
	L^{p,q}_0(\mi,\nu)\times L^{p',q'}(\mi,\nu)\ni (f,g) \mapsto \sum_{(x,y)\in X\times Y} f(x,y) g(x,y)
	\]
	induces an isometry of $L^{p',q'}(\mi,\nu)$ onto $L^{p,q}_0(\mi,\nu)'$.
\end{prop}

When $p,q\Meg 1$, this is a particular case of~\cite[Theorem 1 of \S 3]{BenedekPanzone}. The general case is treated with a similar proof.

\backmatter

\chapter*{Index of Notation}

\begin{center}
  \begin{longtable}{llr}
\textbf{Symbol} & \textbf{Description} & \textbf{Page} \\
\endfirsthead
\textbf{Symbol} & \textbf{Description} & \textbf{Page}\\
\endhead
\endlastfoot
$\N$& natural numbers, starting from $0$ & \\
$\N^*$& natural numbers, starting from $1$ & \\
$\Z$& rational integers & \\
$\R$& real numbers & \\
$\C$& complex numbers & \\
$\R_+$& positive ($\Meg 0$) real numbers & \\
$\R_+^*$& strictly positive ($>0$) real numbers & \\
$A^{(B)}$& families in $A^B$ with finite support & \\
$T_x(M)$& tangent space of $M$ at $x$ & \\
$\varphi_*(X)$ & push-forward of $X$ under  $\varphi$ & \\
$\Lin(X;Y)$&   continuous linear mappings $X\to Y$ & \\
$X'$& (continuous) dual of $X$, $\Lin(X;\C)$ &\\
$\Lin(X)$&  endomorphisms of $X$ & \\
$\Lin^2(X)$ & Hilbert-Schmidt endomorphisms of $X$ & \\
$\Lin^1(X)$ & trace-class endomorphisms of $X$ & \\
$\trasp T$& transpose of $T$ & \\ 
$ \langle \,\cdot\,,\,\cdot\,\rangle $ & bilinear pairing & \\
$ \langle \,\cdot\,\vert\,\cdot\,\rangle $ & sesquilinear pairing & \\
$\Sc$&  Schwartz functions & \\
$\Hol$& holomorphic functions & \\
$C_0$&  continuous functions vanishing at infinity & \\
$\chi_A$& characteristic function of the set $A$ & \\
$I_X$& identity mapping of $X$ & \\
$\Mcal^1$& bounded  Radon measures & \\
$\Hc^n$& suitably normalized Hausdorff measure& \\
$f\cdot \mi$& measure with density $f$ with respect to $\mi$ & \\
$\abs{\mi}$& absolute value of the measure $\mi$ & \\
$\Supp{T}$& support of $T$ & \\
$\dd \pi$& differential of the representation of $\pi$ & \\
$\,\check{\,}$& the inversion $x\mapsto x^{-1}$& \\
$T^*$& $\check{\overline{T}}$& \\
$L_g, R_g$& left and right translations by $g$ & \\
$\vect{1}_k$ & $(1,\dots, 1)\in \N^k$ & \\
$X_\C$& complexification of $X$, $X\otimes_\R \C$ & \\
$\lambda_\C$ & complexification of $\lambda\in X'$, $\lambda\otimes_\R I_\C$ &\\
$\C_+$ & $\R+i \R_+^*$ & \pageref{0}\\
$E$ & vector space over $\C$ of dimension $n$ & \\
$F$ & vector space over $\R$ of dimension $m$ & \pageref{16} \\
$B_E(\zeta,r)$, $B_{F}(x,r)$ & Euclidean balls on $E$ and $F$ & \\
$B_{E\times F_\C}((\zeta,z),r)$  & Euclidean ball on  $E\times F_\C$ & \\
$\Fc_F$, $\Fc_{F'}$ & Fourier transforms on $F$, $F'$ & \\
$\Omega$ & proper open convex cone in $F$ & \pageref{1}\\
$\Omega'$ & dual of $\Omega$ & \pageref{1}\\
$\Phi$ & & \pageref{1}, \pageref{2}\\
$\Phi(\zeta)$ & $\Phi(\zeta,\zeta)$ & \pageref{1}\\
$D$ & Siegel domain of type II & \pageref{3}\\
$b D$ & \v Silov boundary of $D$ & \pageref{4} \\
$\Nc$ & & \pageref{5}\\
$f_h$ & $ (\zeta,x)\mapsto f(\zeta,x+i \Phi(\zeta)+i h) $ & \pageref{12}\\
$W$ & & \pageref{6}\\
$\abs{\Pfaff}$ & Pfaffian & \pageref{7}\\
$\Lambda_+$ & & \pageref{11}\\
$H_\lambda$ & & \pageref{10}\\
$\pi_\lambda$ & Bargmann representation & \pageref{8}\\
$P_{\lambda,0} $ & & \pageref{9}\\
$\Lc$ & Laplace transform & \pageref{13}\\
$e_\Omega$, $e_{\Omega'}$ & base points of $\Omega$ and $\Omega'$ & \pageref{16}\\
$r$ & rank of $\Omega$ & \pageref{16}\\
$T_+$ & & \pageref{14}\\
$\vect m,\vect{m'}$ & & \pageref{15}\\
$\vect b$ & & \pageref{17}\\
$ \Delta_\Omega^{\vect s}$, $\Delta^{\vect s}_{\Omega'}$ & generalized power functions & \pageref{19}\\
$\vect d$ & & \pageref{20}\\
$\nu_\Omega$, $\nu_{\Omega'}$ & invariant measures on $\Omega$ and $\Omega'$ & \pageref{20}\\
$\Gamma_\Omega$, $\Gamma_{\Omega'}$ & Gamma functions on $\Omega$ and $\Omega'$ & \pageref{21}\\
$\N_\Omega$, $\N_{\Omega'}$ & & \pageref{22}\\
$a^{\vect s}$, $i^{\vect s}$ & $a^{\vect s}=\ee^{(\sum_j s_j)\log a}$, $i^{\vect s}= e^{(\sum_j s_j) \pi i/2}$ & \pageref{22}\\
$I^{\vect s}_\Omega$, $I^{\vect s}_{\Omega'}$ & Riemann--Liouville operators & \pageref{23}\\
$\vect s^{\vect{s'}}$ & $\prod_{j=1}^r s_j^{s'_j}$ & \pageref{24}\\
$B^{\vect s}_{(\zeta,z)}$ & & \pageref{25}\\
$L^{p,q}_{\vect s}(D)$, $L^{p,q}_{\vect s,0}(D)$ & & \pageref{26}\\
$d$ & invariant distance on $D$ & \pageref{27}\\
$d_\Omega$, $d_{\Omega'}$ & invariant distances on $\Omega$ and $\Omega'$ &  \pageref{29}\\
$B((\zeta,z),r)$ & invariant ball on $D$ &  \pageref{27}\\
$B_\Omega(h,r)$, $B_{\Omega'}(\lambda,r)$  & invariant balls on  $\Omega$, and $\Omega'$ &   \pageref{29}\\
$\nu_D$ & invariant measure on $D$ & \pageref{27}\\
$A^{p,q}_{\vect s}(D)$, $A^{p,q}_{\vect s,0}(D)$ & weighted Bergman space & \pageref{34}\\
$\Lc^2_{\vect s}(\Omega')$ & & \pageref{35}\\
$\widehat A^{2,2}_{\vect s,\vect{s'}}(D)$ & &\pageref{def:8}\\
$\ell^{p,q}(J,K)$, $\ell^{p,q}_0(J,K)$ & & \pageref{36}\\
$\atomic^{p,q}_{\vect s,\vect{s'}}$, $\atomic^{p,q}_{\vect s,\vect{s'},0}$ & & \pageref{38}\\
$\atomics^{p,q}_{\vect s,\vect{s'}}$, $\atomics^{p,q}_{\vect s,\vect{s'},0}$ & & \pageref{38}\\
$\atomic^{p,q}_{\vect s,\vect{s'},+}$, $\atomic^{p,q}_{\vect s,\vect{s'},0,+}$ & & \pageref{37}\\
$\Sc_\Omega(\Nc)$ & & \pageref{28}\\
$(g\times t)_*\varphi$ & $\Delta_\Omega^{\vect b+\vect d}(t) \varphi\circ (g\times t)^{-1}$ & \pageref{45}\\
$(g\times t)^*\varphi$ & $\Delta_\Omega^{-(\vect b+\vect d)}(t) \varphi\circ (g\times t)$  & \pageref{45}\\
$\Fc_\Nc$ & & \pageref{30}\\
$\Sc_{\Omega,L}(\Nc)$ & & \pageref{31}\\
$\widetilde \Sc_\Omega(\Nc)$ & & \pageref{32}\\
$\Sc_{\Omega,L}'(\Nc)$ & dual of $\Sc_{\Omega,L}(\Nc)$ & \pageref{48}\\
$B^{\vect s}_{p,q}(\Nc,\Omega)$, $\mathring B^{\vect s}_{p,q}(\Nc,\Omega)$ & Besov spaces  & \pageref{33}\\
$\sigma^{\vect s}_{p,q}$ & $\sigma(B^{\vect s}_{p,q}(\Nc,\Omega),
\mathring B^{-\vect s-(1/p-1)_+(\vect b+\vect d)}_{p',q'}(\Nc,\Omega) )$& \pageref{44}\\
$\Ec$ & extension operator & \pageref{39}\\
$\widetilde A^{p,q}_{\vect s}(D)$, $\widetilde A^{p,q}_{\vect s,0}(D)$ & $\Ec(B^{-\vect s}_{p,q}(\Nc,\Omega))$, $ \Ec(\mathring B^{-\vect s}_{p,q}(\Nc,\Omega))$& \pageref{39}\\
$(D)^{\vect s}_{p,q}$, $(D)^{\vect s,0}_{p,q}$ & & \pageref{def:2}\\
$\widehat A^{p,q}_{\vect s,\vect{s'}}(D)$, $\widehat A^{p,q}_{\vect s,\vect{s'},0}(D)$ & & \pageref{def:6}\\
$P_{\vect s}$, $P_{\vect s,+}$ & Bergman projector & \pageref{40}\\
$L^{p}(\mi;Z)$, $L^p_0(\mi;Z)$ & Lebesgue space & \pageref{41}\\
$L^{p}_\loc(\mi;Z)$, $L^p_{0,\loc}(\mi;Z)$ & local Lebesgue space & \pageref{42}\\
$L^{p,q}(\mi,\nu)$ & & \pageref{def:3}\\
$p'$ & $\max(1,p)'$ & \pageref{43}
\end{longtable}
\end{center}

\printindex

\end{document}